\documentclass[reqno,11pt]{amsart}

\usepackage[a4paper,left=23mm,right=23mm,top=30mm,bottom=30mm,marginpar=25mm]{geometry}
\usepackage[utf8]{inputenc}
\usepackage[T1]{fontenc}
\usepackage[english]{babel}
\usepackage{lmodern}
\usepackage{amsmath}
\usepackage{amsthm}
\usepackage{amssymb}
\usepackage{mathrsfs}
\usepackage{stmaryrd}
\usepackage{dsfont}
\usepackage{hyperref}	
\usepackage{graphicx}
\usepackage{color}
\usepackage{longtable}
\usepackage{mathtools}
\usepackage{eurosym}
\usepackage{float}
\usepackage[dvipsnames]{xcolor}
\usepackage{enumitem}
\usepackage{tikz}

\newcommand{\Dcal}{\mathcal{D}}

\newcommand{\Fcal}{\mathcal{F}}

\newcommand{\Rcal}{\mathcal{R}}
\newcommand{\Scal}{\mathcal{S}}
\newcommand{\Tcal}{\mathcal{T}}

\newcommand{\N}{\mathbb{N}}

\newcommand{\R}{\mathbb{R}}

\newcommand{\dl}{\, \mathrm{d}}

\newcommand{\sph}{\mathbb{S}}

\newcommand{\p}{\partial}

\newcommand{\norm}[2]{{\|#1\|}_{#2}}%Norm
\DeclarePairedDelimiter{\abs}{\lvert}{\rvert}

\newcommand{\alp}{\alpha}
\newcommand{\veps}{\varepsilon}
\newcommand{\Om}{\Omega}

\newcommand{\redconvker}{\Theta_\rho^\delta}
\newcommand{\redconvkerrs}{\Tilde{\Theta}_\rho^\delta}
\newcommand{\bla}[1]{\llbracket #1\rrbracket}

\makeatletter
\newcommand{\clippedrightharpoonup}{%
	\mathrel{\mathpalette\clipped@rightharpoonup\relax}%
}
\newcommand{\clipped@rightharpoonup}[2]{%
	\sbox\z@{$\m@th#1\mspace{6mu}$}%
	\clipbox{0pt 0pt {\dimexpr\width-.5\wd\z@} 0pt}{$\m@th#1\rightharpoonup$}%
	\smash{\clipbox{{\wd\z@} 0pt 0pt {-\width}}{$\m@th#1\rightharpoonup$}}%
}
\makeatother

\makeatletter
\newcommand{\clippedto}{%
	\mathrel{\mathpalette\clipped@to\relax}%
}
\newcommand{\clipped@to}[2]{%
	\sbox\z@{$\m@th#1\mspace{6mu}$}%
	\clipbox{0pt 0pt {\dimexpr\width-.5\wd\z@} 0pt}{$\m@th#1\to$}%
	\smash{\clipbox{{\wd\z@} 0pt 0pt {-\width}}{$\m@th#1\to$}}%
}
\makeatother

\newcommand{\ol}[1]{\overline{#1}}

\DeclareMathOperator{\supp}{supp}

\DeclareMathOperator{\diver}{div}

\newcommand{\comment}[1]{{\textcolor{red}{$\langle\hspace{-0.3ex}\langle$Leon: #1$\rangle\hspace{-0.3ex}\rangle$}}}

\newtheorem{Theorem}{Theorem}[section]
\newtheorem{Definition}[Theorem]{Definition}
\newtheorem{Notation}[Theorem]{Notation}

\newtheorem{Corollary}[Theorem]{Corollary}
\newtheorem{Proposition}[Theorem]{Proposition}
\newtheorem{Lemma}[Theorem]{Lemma}

\newtheorem*{Conjecture}{Conjecture}

\theoremstyle{definition}
\newtheorem{Remark}[Theorem]{Remark}
\newtheorem{Setting}[Theorem]{Convention}
\newtheorem{Assumption}[Theorem]{Assumption}

\newenvironment{theorem}[1][]
{ \begin{Theorem}[#1]}
	{\end{Theorem}}

\newenvironment{definition}[1][]
{ \begin{Definition}[#1]}
	{\end{Definition} }

    \newenvironment{notation}[1][]
{ \begin{Notation}[#1]}
	{\end{Notation} }

\newenvironment{corollary}[1][]
{ \begin{Corollary}[#1]}
	{\end{Corollary} }

\newenvironment{proposition}[1][]
{ \begin{Proposition}[#1]}
	{\end{Proposition}}

\newenvironment{lemma}[1][]
{ \begin{Lemma}[#1]}
	{\end{Lemma}}

\newenvironment{remark}[1][]
{ \begin{Remark}[#1]}
	{\end{Remark}}

\numberwithin{equation}{section}

\title[Nonlocal Sobolev spaces as generalized Bessel potential spaces]{Characterizing nonlocal Sobolev spaces as \\generalized Bessel potential spaces }

\author{Leon~Happ}
\address{{\textsc{Leon~Happ.}} Departamento de Matem\'aticas, 
	Universidad Aut\'onoma de Madrid, 
	C/ Francisco Tom\'as y Valiente 7,
	28049 Madrid, Spain}
\email{leon.happ@uam.es}

\author{Carolin~Kreisbeck}
\address{{\textsc{Carolin-Kreisbeck.}} Mathematische-Geographische Fakult\"at, 
	Katholische Universit\"at Eichst\"att-Ingolstadt, 
	Ostenstraße 28,
	85072 Eichst\"att, Germany}
\email{carolin.kreisbeck@ku.de}

\author{Carlos~Mora-Corral}
\address{{\textsc{Carlos~Mora-Corral.}} Departamento de Matem\'aticas, 
	Universidad Aut\'onoma de Madrid, 
	C/ Francisco Tom\'as y Valiente 7,
	28049 Madrid, Spain
    and 
    Instituto de Ciencias Matem\'aticas, CSIC-UAM-UC3M-UCM, 
    28049 Madrid,
    Spain}
\email{carlos.mora@uam.de}

\date{\today}

\begin{document}

\begin{abstract}  
\smallskip

%\commentCK{VERSION BASIC:} In this paper, we study nonlocal Sobolev spaces $H^{\rho,p}(\R^n)$ with $p \in (1,\infty)$, defined via a general class of nonlocal gradients $D^\rho$ associated with a kernel function $\rho$. Our main contribution is to establish three new equivalent characterizations of these spaces: first, in terms of Fourier multipliers; second, through the square root $(-\Delta^\rho)^{1/2}$ of the $\rho$-Laplacian; and third, as potential spaces generated by a suitable convolution kernel. These results generalize the well-known fact that, in the special case of the Riesz kernel, the resulting fractional Sobolev spaces coincide with the classical Bessel potential spaces. In this spirit, our work can also be viewed as introducing a natural notion of generalized Bessel potential spaces.\smallskip

%\commentCK{I tried to flip the story for the abstract, mentioning our characterization early on, and keeping it a bit more matter-of-fact. The introduction can still start with the state-of-the-art and historic review and be a bit more narrative...}

%\commentCK{VERSION 1:}
In this paper, we study nonlocal Sobolev spaces $H^{\rho,p}(\R^n)$ with $p \in (1,\infty)$, defined via a general class of nonlocal gradients $D^\rho$ associated with a kernel function $\rho$. As our main contribution, we establish three new equivalent characterizations of these spaces. First, we describe $H^{\rho,p}(\R^n)$ in terms of Fourier multipliers involving the symbol of the nonlocal potential associated with $D^\rho$. Second, we give a thorough functional analytic treatment of the nonlocal $\rho$-Laplacian $\Delta^\rho$, which allows us to represent general nonlocal Sobolev spaces through its square root $(-\Delta^\rho)^{1/2}$. Third, we show that $H^{\rho,p}(\R^n)$ coincides with the image $g_\rho * L^p(\R^n)$ for a suitable convolution kernel $g_\rho \in L^1(\R^n)$, which confirms that these spaces are genuine potential spaces. These results generalize the well-known fact that, in the special case of the Riesz kernel, the resulting fractional Sobolev spaces coincide with the classical Bessel potential spaces. In this spirit, the present work can also be viewed as introducing a natural notion of generalized Bessel potential spaces.\bigskip

\vspace{8pt}

%\commentCK{Check if ok?}

%\comment{just a collection for now} 46E35, 26A33, 47G20, 31B10, 42B35,47G40, 42B30, 74A70}
% Which ones?

%46E35: Sobolev spaces and other spaces of “smooth” functions, embedding theorems, trace theorems  - Suggestion Carolin: primary
%26A33: Fractional derivatives and integrals 
%31B10: Integral representations, integral operators, integral equations methods in higher dimensions
%47G20: Integro-differential operators 
%42B35: Function spaces arising in harmonic analysis
%47G40: Potential operators 
%42B15: multipliers for harmonic analysis in several variables 

%___________________________________
%74A70: Peridynamics - Suggestion Carolin: Secondary, or remove, maybe remove, it only plays a minor role...
%42B30: H^p spaces - would not include this
%46E30: Spaces of measurable functions (Lp-spaces, Orlicz spaces, Köthe function spaces, Lorentz spaces, rearrangement invariant spaces, ideal spaces, etc.) - Sobolev functions fit better

%%31B10: Integral representations, integral operators, integral equations methods in higher dimensions - others seems to fit better

\vspace{8pt}

\vspace{8pt}

%\noindent\textsc{Date:} \today.
\end{abstract}

	\maketitle
	{\parskip=0em \tableofcontents}

    %%%%%%%%%%%%%%%%%%%%%%%%%%%%%%%%%%%%%%%%%%%%%%%%%%%%%%%%%%%%%%%%%
	%%%%%%%%%%%%%%%%%%%%%%%%%%%%%%%%%%%%%%%%%%%%%%%%%%%%%%%%%%%%%%%%%

%Abstract: Leon Neu - Carolin check
%Einleitung: Carolin kann lesen (C3)
%2. Notation and preliminaries
%2.1 Notation: To be done
%2.2 Recap of nonlocal Sobolev spaces: Done
%2.3 Fourier transform: Done
%2.4 Fourier multiplier theory: Fully done
%2.5 Further hypothesis on the kernel function: Done 
%2.6 Basic Fourier multiplier: Fully done
%3. Fourier characterization: Fully done
%4. Characterization in terms of the $\rho$-Laplacian 
%4.1 Definition and properties: Done, draft (L2)
%4.2 Hilbert space case: Two checks regarding the connection with section 4.3 remain, rest done
%4.3 general p-case: Carolin lesen (C1), Leon entsprechende Ergänzungen oben (L1)
%4.4 Characterization in terms of: Leon nochmal ueberarbeiten (done, shortened substantially)
%5.  und Appendix: Carolin (2)

	\section{Introduction}
    The Riesz fractional gradient was probably first studied in its own right in \cite{Ho59} by Horváth and then very prominently in the more recent works \cite{ShSp15} and \cite{ShSp18} by Shieh and Spector (see their papers also for further historical comments and mentions of earlier appearances of the same object in the literature), even though it is closely related and owes its name to the Riesz potentials introduced by Riesz already in \cite{Ri38} (see also \cite{Ri49}). As established in \cite{ShSp15}, the function spaces associated to the Riesz fractional gradient coincide with the fractional Sobolev spaces known as Bessel potential spaces, going back to the work \cite{ArSm61} by Aronszajn and Smith. These spaces have been the object of intensive study since the 1960s and play a key role in the modern fields of potential theory, harmonic analysis and interpolation theory (see, e.g., the classical references \cite{ArSm61,St70, BeLo76, Tr83, AH96}, as well as the recent works \cite{BGS25, BCGS25}), not mentioning many others.
    However, Riesz fractional gradients are only a model example of a wider class of nonlocal gradients, where the common local gradient is replaced by an integral over the first order difference quotient weighted with a general kernel (\cite{MeSp15, BCMC23, CuKrSch23, BMCS25}). 
   Nonlocal gradients have gained popularity in applied mathematics over the course of the last decades as an alternative approach to models involving conventional differential operators, and important applications comprise peridynamics and nonlocal hyperelasticity (see, e.g., \cite{Si00, MeDu15,BCMC24}).  
    %(their application in the calculus of variations through the connection to the Riesz fractional gradient, and subsequently the investigation into more general nonlocal gradients, is a relatively recent development following the above-mentioned contributions, lies in their potential use in certain models in materials science. In particular, the works of \comment{Silling} in the field of peridynamics sparked .... 
    Mathematically speaking, nonlocal gradients provide a means to enforce relaxed regularity requirements, in contrast with classical Sobolev spaces, upon the functions in question, while a lot of elements in the standard functional analytic toolbox, as a generalized fundamental theorem of calculus, Poincaré-type inequalities, compactness properties, and embedding relations, are usually still available in this broader setting (\cite{CoSt19,CuKrSch23,BMCS25}).
    %\color{red} However, it should be pointed out that, as different as these competing approaches might seem, the canonical limit passages in the calculus of variations (localization, homogenization, discrete-to-continuum analysis etc.) suggest that the right or wrong of local or nonlocal models is often a matter of scale (\comment{But if I am not mistaken, in none of the results, they actually recover a nonlocal gradient-type model in the limit}). \color{black}

    In this work, we address the question whether the nonlocal Sobolev spaces associated to a certain class of generalized nonlocal gradients can be viewed as generalized Bessel potential spaces. In particular, we are interested in whether the well-known equivalent characterizations of Bessel potential spaces carry over to the setting of generalized nonlocal gradients studied in \cite{BMCS25}. Our main findings are summarized in the following subsections.
    
    Throughout the paper, we put an emphasis on a thorough derivation of the results and keep the exposition as self-contained as possible, striving, in particular, to make the paper accessible also for the interested reader not already equipped with a deep understanding of the Fourier theoretic and/or functional theoretic aspects of the matter.

    %\commentCK{I feel we still need to decide on a clear story for the paper. Of course, there are several aspects we cover, but some messages seem mixed... }

    \subsection{A brief overview over the main results}
    We are concerned subsequently with \(\rho\)-nonlocal Sobolev spaces (just nonlocal Sobolev spaces for short) \(H^{\rho,p} (\R^n)\) according to the following framework: for a radial kernel \(\rho:\R^n\to [0,\infty]\) subject to suitable hypothesis, \(H^{\rho,p} (\R^n)\) is given as  
    \begin{equation}\label{eq:spaces_intro}
        H^{\rho,p} (\R^n)
        \coloneqq \{u\in L^p(\R^n): D^\rho u\in L^p(\R^n; \R^n)\},
    \end{equation}
    where the nonlocal derivative is defined for \(\varphi\in C_c^\infty(\R^n)\) by 
    \begin{equation}\label{eq:gen_nlg_intro}
        D^\rho \varphi(x)
        = \int_{\R^n}\frac{\varphi(x)-\varphi(y)}{|x-y|} \frac{x-y}{|x-y|}\rho(x-y)\dl y, \qquad x\in \R^n 
    \end{equation}
    and is extended to \(L^p(\R^n)\) by duality (see Section~\ref{sec:recap_nl_sob} for more details).
    %and Section~\ref{sec:discussion} for more background
    %and an account of the existing literature).
    %\commentCK{Make sure to mention also previous works on special kernels, e.g. Riesz fractional~\cite{ShSp15, ??, Mengesha, ComiStefani, Shilhavy}, finite-horizon fractional \cite{BellidoCuetoMoraCorral, CuKrSch23,...}, etc. and motivation applications in peridynamics and nonlocal hyperelasticity which also inspired the general setting... Also ok in introduction.}
    This setting agrees with the general approach developed in \cite{BMCS25}, which comprises as special cases the classical Riesz fractional gradient, as studied, for example, in \cite{ShSp15, CoSt19, Si20}, and finite horizon fractional gradients appearing in \cite{BCMC23, BCMC24, CuKrSch23, KrSch24}. In particular, in the notation of \eqref{eq:spaces_intro} and \eqref{eq:gen_nlg_intro}, the Riesz fractional gradients correspond to the kernel function \(\rho_s(x)=c(n,s)|x|^{-(n+s-1)}\), for \(s\in (0,1)\) and suitable constants \(c(n,s)\), and the respective nonlocal Sobolev spaces are \(H^{\rho_s,p}(\R^n)\) (see also \eqref{eq:rieszgrad}).
    For further works on related general non-local gradients we refer to \cite{DGLZ13, MeSp15, DEGOK21, DEGMS22}.
    %(see also Section~\ref{sec:discussion}). 
    %We also refer the interested reader to \cite{BMCS25} for a more comprehensive account on the concept of general nonlocal gradients.
    
     We establish in this paper several characterizations of the spaces \(H^{\rho,p} (\R^n)\) and their topological duals, generalizing the corresponding identities for Bessel potential spaces (in Section~\ref{sec:discussion} we recall classical Bessel potential spaces and elaborate on the precise connections to our results).
     %, which, in turn, correspond to the nonlocal Sobolev spaces \(H^{\rho_s,p}(\R^n)\) for the choice of the kernel \(\rho_s(x)=c(n,s)|x|^{-(n+s-1)}\), \(s\in (0,1)\), as first shown in \cite{ShSp15}. 
     As a disclaimer, we note upfront that the characterizations below may differ in their assumptions on the kernel function \(\rho\), in particular, in the requirements regarding the behavior of \(\rho\) in the origin manifesting in the conditions (H3) and (H3') (the precise hypotheses can be found in Section~\ref{sec:hypotheses}).
    
	First, we demonstrate that, under suitable premises already devised in \cite{BMCS25}, it holds 
	\begin{equation}\label{eq:char_pot2}
	H^{\rho,p} (\R^n)
	= \{u\in \mathcal{S}'(\R^n):\, \mathcal{F}^{-1}((1+|\lambda_\rho|^2)^{{1}/{2}}\hat{u})\in L^p(\R^n)\}
	\end{equation}
    for \(\lambda_\rho\) the Fourier symbol of the nonlocal gradient. This result can be found in Theorem~\ref{thm:alt_char}. A similar identity is true for sufficiently regular bounded domains, as demonstrated in Corollary~\ref{cor:alt_charOm}. 
We also derive a related representation of the topological dual space 
	\begin{equation*}
		H^{-,\rho,p'}(\R^n)\coloneqq (H^{\rho,p}(\R^n))'.
	\end{equation*}
	In particular, one has that 
	\begin{equation}\label{eq:char_dual}
	H^{-,\rho,p'}(\R^n)
	= \{u\in \mathcal{S}'(\R^n):\, \mathcal{F}^{-1}((1+|\lambda_\rho|^2)^{-{1}/{2}}\hat{u})\in L^{p'}(\R^n)\}.
	\end{equation}
It also holds a characterization similar to 
%the classical result \cite[Proposition 9.30]{Bre11}, 
the classical result stating that, for any open \(\Om\subset \R^n\), any element in $H^{-,\rho,p'}(\Om)$ can be identified with some $u_0 + \diver^{\rho} U$ for $u_0 \in L^p (\Om)$ and $U \in L^p (\Om; \R^n)$. These two findings are presented in Theorem~\ref{thm:alt_dual_char} and Proposition~\ref{prop:charac_bijection}, respectively.

Moreover, we deduce that 
    \begin{align}
        H^{\rho,p}(\R^n)
        &= \{u\in L^p(\R^n):\, (-\Delta^{\rho})^{1/2} u \in L^p(\R^n)\}
        \label{eq:char_lap_pre}\\
         &= \{u\in \mathcal{S}'(\R^n):\, (1-\Delta^{\rho})^{1/2} u \in L^p(\R^n)\}, \label{eq:char_bessel_pre}
    \end{align}
    where \(\Delta^{\rho}\) is the distributional nonlocal generalized \(\rho\)-Laplacian (see Definition~\ref{de:Laplacianp}) and \((-\Delta^{\rho})^{1/2}\), \((1-\Delta^{\rho})^{1/2}\) are square roots of the respective squared operators. The identity \eqref{eq:char_lap_pre} is proved in Proposition~\ref{prop:charac_lapace2} and Corollary~\ref{prop:charac_lapacep}, while \eqref{eq:char_bessel_pre} is the subject of Corollary~\ref{prop:charac_bessel_op}. 
    We note that while \((1-\Delta^{\rho})^{1/2}\) above is understood as an element of \(\mathcal{S}'(\R^n)\) by duality, more work is required to properly introduce and study \((-\Delta^{\rho})^{1/2}\) as an unbounded operator in \(L^p(\R^n)\). 
    %Here, we identify it with the closure in \(L^p(\R^n)\) of the unbounded operator %(see Definition~\ref{def:ub_mult_op_pl})
    %\[
    % F_{\abs{\lambda_\rho}}u \coloneqq \mathcal{F}^{-1} (\abs{\lambda_{\rho}} \widehat{u}) , \qquad u \in \Scal(\R^n).
    %\]
%In addition, in the Hilbert space case \(p=2\), the operator \((- \Delta^{\rho})^{1/2}\) is uniquely defined as the square root of the positive semidefinite and self-adjoint operator \(-\Delta^\rho\). 

	Finally, we prove that similar to the case of the Riesz fractional gradient and the Bessel potential spaces, we can characterize \(H^{\rho,p}(\R^n)\) as the image of \(L^p(\R^n)\) under a suitable convolution operator. Precisely, we establish that there exists \(g_\rho\in L^1(\R^n)\) such that 
	\begin{equation}\label{eq:char_pot}
	H^{\rho,p}(\R^n)
	=\{g_\rho*f:\, f\in L^p(\R^n)\}.
	\end{equation}
%	where the convolution for vectorial functions is understood component-wise. 
%The convolution kernel \(g_\rho\) is the Fourier inverse of \((1+|\lambda_\rho|^2)^{-{1}/{2}}\) and its particular form can be found in the proof of Proposition~\ref{prop:invers_fourier}. 
%We add that some of the abovementioned results naturally extend to general open domains \(\Om\subset \R^n\).
%Namely, if \(|\p \Om|=0\) we have that (see Corollary~\ref{cor:alt_charOm})
%	\begin{equation}\label{eq:char_pot3}
%	H^{\rho,p}_0(\Om)
%	= \{u\in \mathcal{S}'(\R^n):\, \supp u\subset \ol{\Om} \text{ and } \mathcal{F}^{-1}((1+|\lambda_\rho|^2)^{{1}/{2}}\hat{u})\in L^p(\R^n)\}.
%	\end{equation}
For this result we refer to Theorem~\ref{thm:pot_char}.

Let us point out that throughout the paper we confine our analysis to the range \(p\in (1,\infty)\), since it is the natural setting for the Fourier methods employed in our approach.

%%%%%%%%%%%%%%%%%%   
    
We shall conclude this section by sharing a short observation about the naming of the nonlocal Sobolev spaces treated in this work.
When introduced in \cite{BMCS25}, the authors also discuss possible alternative names for the objects $H^{\rho,p} (\R^n)$ and conclude to give them the name of \emph{$\rho$-nonlocal Sobolev spaces}.
Although they acknowledge that \emph{$\rho$-Bessel spaces} would be another sensible name, this option was ruled out at the time because there was no immediate connection with the Bessel potential.
However, the main result of this work (Theorem~\ref{thm:pot_char}) bridges this exact gap and shows that they indeed define a natural notion of \emph{generalized Bessel potential spaces}.

%For this reason, in this work we have decided to adopt the name
%\emph{generalized Bessel potential spaces}.
%\commentC{The use of the name Bessel was retaken, in \cite{BGS25} and \cite{BCGS25}.} \comment{how is that remark intended exactly? Because I think that they actually only treat the classical Bessel potential spaces.}

Finally, we also stress that we are not the first ones coining the term \emph{generalized Bessel potential spaces}. For example, in \cite{ShEK23} the same name has been used in a different context for other objects, namely, the spaces of functions that are obtained by replacing the Fourier with the Hankel transform and adding a certain weight in the characterization of Bessel spaces of the type \eqref{eq:char_pot2}. 

%%%%%%%%%%%%%%%%%%%%%%%%%%%%%%%%%%%%%%%%%
\subsection{Discussion of results and methodology}\label{sec:discussion}
We start our discussion with a few comments about the Bessel potential spaces. For the following facts, we refer to the classical treatise \cite[Chapter~V.3]{St70} by Stein, the book \cite[Section~7]{Mi96}, and also the recent work \cite{BGS25}, as well as the references therein. %\commentCK{$\rho$-Laplacian and corresponding representation are not mentioned in this dicsussion}\comment{They are, but towards the end. We could add a it there. }
    
    To set the stage, we define for \(s\in\R\) the operators 
    \begin{equation*}
        \Lambda_s:\mathcal{S}'(\R^n)\to \mathcal{S}'(\R^n), 
        \quad f\mapsto \mathcal{F}^{-1}((1+4\pi^2|\cdot|^2)^{-s/2}\hat{f})
    \end{equation*}
    and the $L^1$-function 
    \begin{equation}\label{eq:gs}
		g_s(x)
		\coloneqq \frac{1}{(4\pi)^{s/2}\Gamma(\frac{s}{2})}
        \int_0^\infty e^{-\pi|x|^2/\delta}e^{-\delta/4\pi}\delta^{(s-n)/2}\frac{\dl \delta}{\delta},
        \qquad x\in \R^n.
	\end{equation}
    For \(s>0\), it can be shown that  
    \begin{equation*}
        \widehat{g_s}
        = (1+4\pi^2|\cdot|^2)^{-s/2},
    \end{equation*}
    which implies that 
    \begin{equation*}
        \Lambda_s f
        = g_s* f
    \end{equation*}
    and proves, together with Young's convolution inequality, that \(\Lambda_s\) is a well defined operator from \(L^p(\R^n)\) into itself 
    %for any \(s>0\) and 
    for any \(p\in [1,\infty]\). Moreover, since \(g_s>0\) and \(\mathcal{F}(g_s)(0)=1\), it follows that \(\norm{g_s}{L^1(\R^n)}=\mathcal{F}(g_s)(0)=1\). In particular, \(g_s\) defines a probability measure on \(\R^n\) and, in the terminology of probability theory, the symbol \((1+4\pi^2|\cdot|^2)^{-s/2}\) defines its \textit{characteristic function}.
    
    We are now ready to formally introduce the Bessel potential spaces, which can be done in two equivalent ways, either following Lions and defining them as 
    \begin{equation}\label{eq:bessel_lions}
        \{u\in \mathcal{S}'(\R^n):\, \Lambda_{-s}f\in L^p(\R^n)\}
        = \{u\in \mathcal{S}'(\R^n):\, \mathcal{F}^{-1}((1+4\pi^2|\cdot|^2)^{{s}/{2}}\hat{u})\in L^p(\R^n)\},
    \end{equation}
    or in the spirit of Calderón as potential spaces in the sense
    \begin{equation}\label{eq:bessel_calderon}
        \Lambda_s(L^p(\R^n))
        = \{\Lambda_s f:\, f\in L^p(\R^n)\}
        = \{g_s*f:\, f\in L^p(\R^n)\}.
    \end{equation}
    Indeed, it is easily verified that both definitions agree, using the semigroup property 
    \begin{equation}\label{eq:sg_property}
        \Lambda_{s_0}\Lambda_{s_1}
        = \Lambda_{s_0+s_1} ,
        \qquad 
        s_0,s_1\in \R.
    \end{equation}
We will denote by \(H^{s,p}(\R^n)\) the Bessel potential space of order \(s\in \R\) and exponent \(p\in [1,\infty]\), according to \eqref{eq:bessel_lions} and \eqref{eq:bessel_calderon}. 

%***4/8/2026

    The Bessel potential spaces owe their name to the Bessel potentials associated to the operators \((1-\Delta)^{-s/2}\), whose Fourier symbol is \((1+4\pi^2|\cdot|^2)^{-s/2}\). In particular, \(H^{s,p}(\R^n)\) is commonly also expressed in the form
    \begin{equation}\label{eq:bessel_bessel}
        H^{s,p}(\R^n)
        = \{u\in \mathcal{S}'(\R^n):\, (1-\Delta)^{s/2} u \in L^p(\R^n)\}.
    \end{equation}
    The Bessel potentials are themselves closely related to the Riesz potentials based on the operators \((-\Delta)^{-s/2}\). In particular, it can be shown that for  \(s>0\) the formal quotients
    \begin{equation*}
        \frac{(-\Delta)^{s/2}}{(1-\Delta)^{s/2}}
        \quad\text{and}\quad
        \frac{(1-\Delta)^{s/2}}{1+(-\Delta)^{s/2}},
    \end{equation*}
or in their Fourier representation
    \begin{equation*}
        \frac{(2\pi|\cdot|)^s}{(1+4\pi^2|\cdot|^2)^{s/2}}
        \quad\text{and}\quad
        \frac{(1+4\pi^2|\cdot|^2)^{s/2}}{1+ (2\pi|\cdot|)^s},
    \end{equation*}
    admit Fourier inverses, i.e., there exist finite signed measures \(\mu_1^s, \mu_2^s\) on \(\R^n\) such that
    \begin{equation*}
        \mu_1^s*f=\frac{(-\Delta)^{s/2}}{(1-\Delta)^{s/2}}f
        \quad\text{and}\quad
        \mu_2^s *f =\frac{(1-\Delta)^{s/2}}{1+(-\Delta)^{s/2}}f
    \end{equation*}
    for all \(f\in L^p(\R^n)\), and
    \begin{equation*}
        \widehat{\mu_1^s}=\frac{(2\pi|\cdot|)^s}{(1+4\pi^2|\cdot|^2)^{s/2}}
        \quad\text{and}\quad
        \widehat{\mu_2^s}=\frac{(1+4\pi^2|\cdot|^2)^{s/2}}{1+ (2\pi|\cdot|)^s}.
    \end{equation*}
    Leveraging these relations, the Bessel potential spaces can be represented in terms of the Riesz potential as
     \begin{align}\label{eq:bessel_riesz}
        H^{s,p}(\R^n)
        &= \{u\in L^p(\R^n):\, (-\Delta)^{s/2} u \in L^p(\R^n)\}.
    \end{align}

    This finishes our recapitulation of the different characterizations of the Bessel potential spaces. We transition now to the question why they gained such an importance over the last years in the context of nonlocal gradients. The answer is contained in the observation that there exists yet another way to describe Bessel potential spaces: as established in \cite{ShSp15}, they coincide with the nonlocal Sobolev spaces (see Section~\ref{sec:recap_nl_sob}) based on the kernels 
    \begin{equation}\label{eq:rieszgrad}
		\rho_s(x)
		\coloneqq \frac{1}{\gamma (n,s)|x|^{n-1+s}}
		\quad \text{for}\quad
		\gamma (n,s)\coloneqq \frac{\pi^{\frac{n}{2}}\Gamma\left(\frac{1-s}{2}\right)}{2^s\Gamma\left(\frac{n+1+s}{2}\right)},
        \qquad x\in \R^n\setminus \{0\},
	\end{equation}
    this is, 
    \begin{equation*}
        H^{s,p}(\R^n)=H^{\rho_s,p}(\R^n).
    \end{equation*}
    The path taken in this paper is now opposite to the historical development of the theory of Bessel potential spaces: starting from generalized nonlocal Sobolev spaces, we strive to prove characterizations similar to \eqref{eq:bessel_lions}, \eqref{eq:bessel_calderon}, \eqref{eq:bessel_bessel} and \eqref{eq:bessel_riesz} for very general kernels. In our general theory, the quantity \(2\pi|\cdot|\) is replaced by the expression \(|\lambda_\rho|=2\pi|\cdot|\widehat{Q_\rho}\) (see Section~\ref{sec:fourier}). We demonstrate that the analogues of \eqref{eq:bessel_lions} and \eqref{eq:bessel_bessel}, i.e., 
    \begin{equation*}
	H^{\rho,p} (\R^n)
	= \{u\in \mathcal{S}'(\R^n):\, \mathcal{F}^{-1}((1+|\lambda_\rho|^2)^{{1}/{2}}\hat{u})\in L^p(\R^n)\}
	\end{equation*}
    and
    \begin{align*}
        \begin{split}
        H^{\rho,p}(\R^n)
        = \{u\in \mathcal{S}'(\R^n):\, (1-\Delta^{\rho})^{1/2} u \in L^p(\R^n)\},
        \end{split}
    \end{align*}
    where \((1-\Delta^{\rho})^{1/2}\) is extended by duality to the space of tempered distributions,
    merely rest on the assumption that the generalized fractions 
    	 \begin{equation}\label{eq:general_fractions}
    	 	 \frac{|\lambda_\rho|}{(1+|\lambda_\rho|^2)^{{1}/{2}}}
    	 	\quad\text{ and }\quad 
    	 	\frac{(1+|\lambda_\rho|^2)^{{1}/{2}}}{1+|\lambda_\rho|}
    	 \end{equation}
         constitute Fourier multipliers; in particular, it is not needed that they admit Fourier inverses at this point. As it turns out, under stronger assumptions (the additional hypotheses being the lower bound of \(\rho\) by a Riesz-type fractional gradient in the sense of (H3) in Section~\ref{sec:hypotheses}) they actually do (see Lemma~\ref{lem:gen_pot_frac_inverse}), and the same ideas can be used to show the potential-type characterization of \(H^{\rho,p}(\R^n)\). More concretely, we construct a Fourier inverse in \(L^1(\R^n)\) of the expression \((1+|\lambda_\rho|^2)^{-{1}/{2}}\), %\commentCK{Comment on what happened to $4\pi^2$}, 
         given by  
       \begin{equation*}
            g_{\rho}(x)
    		\coloneqq \frac{1}{(4\pi^2)^{1/2}\Gamma(\frac{1}{2})}
            \int_0^\infty f^\delta_\rho(x)e^{-\delta/(4\pi^2)}\delta^{(1-n)/2}\frac{\dl \delta}{\delta},
            \qquad x\in \R^n,
       \end{equation*}
       for
        \begin{align*}
        f^\delta_\rho(x) \coloneqq \delta^{n/2}\int_{\R^n} e^{-\abs{\lambda_{\rho} (\xi)}^2\delta / (4 \pi^2)}e^{2\pi i x\cdot \xi}\dl \xi,
        \qquad x\in \R^n.
		\end{align*}		
       This is enough to deduce that also Calderon's definition in \eqref{eq:bessel_calderon} of the Bessel spaces finds an analogue in the general sense of nonlocal Sobolev spaces. In particular, it holds
    \begin{equation*}
        H^{\rho,p}(\R^n)
	   =\{g_\rho*f:\, f\in L^p(\R^n)\}
    \end{equation*}
    and the spaces \(H^{\rho,p}(\R^n)\) can indeed be viewed as generalized potential spaces. 
    %In addition, the function \(g_\rho\) can also be used to show that the fractions in \eqref{eq:general_fractions} admit Fourier(-Stieltjes) inverses, again in the sense of finite signed measures on \(\R^n\) (see Lemma~\ref{lem:gen_pot_frac_inverse}).
    
    However, for the moment, we were not able to prove that \(g_\rho\) is positive and \(\norm{g_\rho}{L^1(\R^n)}=1\), which would mirror the situation of the Bessel potential spaces. Consequently, it is not known whether \((1+|\lambda_\rho|^2)^{-{1}/{2}}\) is again the characteristic function of a probability measure. 

    Starting from the same stronger assumptions on \(\rho\) mentioned above, we further show that \(\Delta^\rho\) is well defined as an unbounded operator in \(L^p(\R^n)\), and that \(\ol{F_{\abs{\lambda_\rho}}}\), the closure in \(L^p(\R^n)\) of the unbounded operator %(see Definition~\ref{def:ub_mult_op_pl})
    \[
     F_{\abs{\lambda_\rho}}u \coloneqq \mathcal{F}^{-1} (\abs{\lambda_{\rho}} \widehat{u}) , \qquad u \in \Scal(\R^n),
    \]
    %,    which is the closure of the multiplier operator with Fourier symbol \(\abs{\lambda_\rho}\) and domain \(\Scal(\R^n)\), 
    is a square root of \(-\Delta^{\rho}\).
    %(which is unique for \(p=2\)).  
    In addition, in the Hilbert space case \(p=2\), the operator \((- \Delta^{\rho})^{1/2}\) is uniquely defined as the square root of the positive semidefinite and self-adjoint operator \(-\Delta^\rho\). 
    As a consequence of these results, we find that a nonlocal version of \eqref{eq:bessel_riesz} holds true, this is
    \begin{equation*}
       H^{\rho,p}(\R^n)
        = \{u\in L^p(\R^n):\, (-\Delta^\rho)^{1/2} u \in L^p(\R^n)\} 
        = D(\ol{F_{\abs{\lambda_\rho}}}),
    \end{equation*}
    with \(D(\ol{F_{\abs{\lambda_\rho}}})\) the domain in \(L^p(\R^n)\) of \(\ol{F_{\abs{\lambda_\rho}}}\).

    In summary, we derive analogous expressions to \eqref{eq:bessel_lions}, \eqref{eq:bessel_calderon}, \eqref{eq:bessel_bessel} and \eqref{eq:bessel_riesz} for general nonlocal Sobolev spaces. As a consequence, in light of the above remarks, the operator \((1-\Delta^{\rho})^{-1/2}\) might be interpreted as a generalized Bessel potential operator. 
    Nevertheless, our findings deviate in their generality in some ways from the specif case of the Bessel potential kernels, and we summarize in the following the main differences: 
    \begin{enumerate}[label=(\roman*)]
        \item 
            the operators on \(L^p(\R^n)\) defined through \(f\mapsto g_\rho * f\) do not form a semigroup in a similar sense to the identity (following from \eqref{eq:sg_property})
            \begin{equation*}
                g_{\rho_{s_0}} * (g_{\rho_{s_1}} * f) 
                =g_{\rho_{s_0+s_1}} * f,
                \qquad f\in L^p(\R^n)
            \end{equation*}
            holding for the Bessel potential spaces;
%\commentC{ADD: in the sense that $g_{\rho_1} * (g_{\rho_2} * f) \neq g_{\rho_1 \rho_2} * f$ (or is it $g_{\rho_1} * (g_{\rho_2} * f) \neq g_{\rho_1 \circ \rho_2} * f$?) for kernels $\rho_1, \rho_2$}
        \item 
            for the moment, the conditions on $\rho$ to prove the characterization of \(H^{\rho,p}(\R^n)\) similar to \eqref{eq:bessel_lions} is milder than those required to show that the generalized Bessel potential \((1-\Delta^{\rho})^{-1/2}\) (and also the fractions in \eqref{eq:general_fractions}) admits a Fourier inverse \(g_\rho\in L^1(\R^n)\), and that \(\Delta^\rho\) is a closable operator over \(L^p(\R^n)\) with well defined square root;
        \item 
            we were not able to prove the positivity of \(g_\rho\), which would immediately imply that \(g_\rho\) defines a probability distribution with the corresponding characteristic function given by \((1+|\lambda_\rho|^2)^{-{1}/{2}}\);
        \item 
we suppose that the kernel \(\rho\) has compact support (which is not the case for the Riesz fractional gradient), but this can be assumed without loss of generality, as explained in Section~\ref{sec:hypotheses}). 
\end{enumerate}

%%%%%%%%%%%%%%%%%%%%%%%%%%%%%%%%%%%%%%%%%
\subsection{Outline} 
After fixing the notation used throughout the paper in Section~\ref{sec:notation}, we recall the definition and main properties of nonlocal Sobolev spaces in Section~\ref{sec:recap_nl_sob}. Section~\ref{sec:fourier} then collects the Fourier-analytic aspects of nonlocal gradients needed later on. 
In Section~\ref{sec:fm_theory}, we set up the framework of unbounded operators defined via multiplication in Fourier space, before stating in Section~\ref{sec:hypotheses} the precise hypotheses on the kernel functions that will be assumed in different parts of the paper. The preliminaries conclude with Section~\ref{sec:fourier_mult}, where we derive key estimates and establish the central Fourier multiplier results on which the subsequent analysis relies.

%In Section~\ref{sec:notation} we fix some notation that is employed throughout the paper. Afterwards, we recall in Section~\ref{sec:recap_nl_sob} the definition and important properties of the nonlocal Sobolev spaces introduced in \cite{BMCS25}, devote Section~\ref{sec:fourier} to the prerequisite Fourier aspects of nonlocal gradients, and prepare in Section~\ref{sec:fm_theory} the use of unbounded operators defined through multiplication in Fourier space, before turning in Section~\ref{sec:hypotheses} to the precise hypotheses on the kernel functions of the nonlocal gradients that will be put in place in different parts of the paper. The preliminaries are concluded in Section~\ref{sec:fourier_mult} with the derivation of important estimates and the central Fourier multiplier results that underlie the successive findings.

The presentation of the different characterizations of the nonlocal Sobolev spaces is split into three main parts. First, Section~\ref{sec:charac:fourier} covers what we call the Fourier characterization, where the space is given as the subspace of tempered distributions in terms of a suitable multiplier. Such a description holds for functions defined on the whole space, as well as on bounded domains, which we show in Section~\ref{sec:fourier_full_space} and Section~\ref{sec:fourier_bounded}, respectively. Afterwards, we also derive the Fourier characterization of the dual spaces in Section~\ref{sec:fourier_dual}, alongside their representation through derivatives of \(L^p\)-functions.  

Second, in Section~\ref{sec:char_laplacian} we establish the characterization of the nonlocal Sobolev spaces through the \(\rho\)-Laplacian as subspaces of the \(L^p\)-functions. More precisely, after thoroughly introducing this object in Section~\ref{sec:laplace_prelim} and spending some time on deriving first essential properties of it, we focus in Section~\ref{sec:laplace_hilbert} on the Hilbert space case \(p=2\) and then in Section~\ref{sec:laplace_general} on the general case \(p\in (1,\infty)\). Later on, in Section~\ref{sec:laplace_bp_inverse}, we also directly recast the Fourier characterization from Section~\ref{sec:charac:fourier} in terms of the operator \((1-\Delta^\rho)^{1/2}\).

Third, we conclude the cascade of equivalent characterizations in Section~\ref{sec:charac_bessel}, and justify the name \textit{generalized Bessel potential spaces}, by constructively proving the existence of a suitable potential convolution kernel.

%%%%%%%%%%%%%%%%%%%%%%%%%%%%%%%%%%%%%%%%%
%\subsection{Methodology}
%\comment{Maybe obsolete with previous sectiond and subsequent one.}

%%%%%%%%%%%%%%%%%%%%%%%%%%%%%%%%%%%%%%%%%
\section{Preliminaries}
\subsection{Notation and terminology}\label{sec:notation}
    We denote by \(\N\) the natural numbers and write \(\N_0\) if we include zero; \(\R^n\) is the Euclidean space with dimension \(n\in \N\). 
	We write \(B_r(x)\) for the ball centered at \(x\in \R^n\) and with radius \(r>0\). We set \(B_r\coloneqq B_r(0)\) and write \(\sph^{n-1}\coloneqq \p B_1\) and \(\sph^{n-1}_+\coloneqq\{z\in \sph^{n-1}: z_1>0\}\). The surface area integral over \(\p B_r\) will be denoted by \(\dl \Scal\).
	By \(\Gamma:(0,\infty)\to \R\), \(\Gamma(x)\coloneqq \int_0^\infty t^{x-1}e^{-t}\dl t\) we refer to the respective restriction of the common \(\Gamma\)-function. 
	
	If not mentioned otherwise, we assume throughout the paper that \(p\in (1,\infty)\). Moreover, \(p'\in (1,\infty)\) denotes the dual exponent of $p$ defined through the relation \(1/p+1/p'=1\).
	%Throughout, we take \(p\in (1,\infty)\), and \(p'\in (1,\infty)\) denotes the dual exponent of $p$ defined through the relation \(1/p+1/p'=1\). 
	With \(L^p(\R^n)\) we refer to the usual Lebesgue spaces. 
	We denote the space of compactly supported and smooth functions by \(C_c^\infty(\R^n)\), and its topological dual by \(\Dcal'(\R^n)\). The space of Schwartz functions and its topological dual, the tempered distributions, are abbreviated as \(\Scal(\R^n)\) and \(\Scal'(\R^n)\), respectively. We denote by \(\langle\cdot,\cdot\rangle\) the duality pairing between $L^p(\R^n)$ and $L^{p'}(\R^n)$. Moreover, we write \(\langle\cdot,\cdot\rangle_{\Dcal'\times\Dcal}\), where, for consistency, we use the notation \(\Dcal(\R^n)=C_c^\infty(\R^n)\), and \(\langle\cdot,\cdot\rangle_{\Scal'\times\Scal}\) for the respective duality pairings. 	 We say that an element \(u\) of \(\Dcal'(\R^n)\) or \(\Scal'(\R^n)\) has a representative \(v\) in \(L^p(\R^n)\) if \(\langle u,\varphi \rangle_{\Dcal'\times\Dcal}=\langle v, \varphi\rangle\) for every \(\varphi\in C_c^\infty(\R^n)\) or  \(\langle u,\varphi \rangle_{\Scal'\times\Scal}=\langle v, \varphi\rangle\) for every \(\varphi\in \Scal(\R^n)\), respectively.

	To avoid confusion with the duality pairing, we use the notation \(\bla{\cdot}\coloneqq (1+|\cdot|^2)^{1/2}\) 
	instead of the maybe more common expression \(\langle\cdot\rangle\) for this quantity.
	
	With $\alpha\in \N_0^n$ we always refer to multi-indices, with the convention that $|\alpha|=\alpha_1+\alpha_2+\ldots+\alpha_n$. 
	For another multi-index \(\beta\in \N_0^n\), the inequality \(\beta\le \alp\) is understood componentwise, i.e., \(\beta_i\le \alp_i\) for all \(i\in\{1,\ldots n\}\); in particular, \(|\alp-\beta|=|\alp|-|\beta|\) for \(\beta\le \alp\). Finally, addition of multi-indices is understood componentwise. 
		If we say that all derivatives of a function \(u\in C^\infty(\R^n)\) have (at most) polynomial growth, we mean, as this is usually understood (see, e.g., \cite{Schwartz66}), that for each \(\alp\in \N_0^n\) there exist \(C_\alp>0\) and \(N_\alp\ge 0\) such that 
	\begin{equation*}
	|\p^\alp u(x)|
	\le C_\alp (1+|x|)^{N_\alp},
	\qquad x\in \R^n.
	\end{equation*}
	Accordingly, we also say that a function has polynomial growth as soon as the estimate holds for \(\alp=0\). Both notions naturally extend to vector-valued functions.
	
	We use \(D(A)\) to denote the domain of definition of an unbounded linear operator \(A:D(A)\subset X\to Y\) between two Banach spaces \(X\) and \(Y\). 
	   In the case \(p=2\) and for some densely defined linear operator \(T:D(T)\subset L^2(\R^n)\to L^2(\R^n)\), we write \(T^*\) for the adjoint operator \(T^*:D(T^*)\subset L^2(\R^n)\to L^2(\R^n)\). 
       %which is uniquely defined through the identity \(\langle A u, v\rangle=\langle u,A^* v\rangle\) for all \((u,v)\in D(A)\times D(A^*)\). At several places we write \(1\) for the identity map in \(L^p(\R^n)\).
	Standard references for unbounded operators in Banach spaces are \cite[Chapter 3, Sections 2.1 and 5.1]{Kato76} and \cite[Section 2.1]{CaHa98}.
	We will also use the terminology of closable operators (see, e.g., \cite[Chapter 3, Section 5.3]{Kato76}). For our purposes, it will be enough to know that an unbounded linear operator $T$ in the Banach space $X$ is closable if and only if the following condition holds: if a sequence $(u_j)_j$ in $D(T)$ tends to $0$ in $X$ and there exists $v \in X$ such that $Tu_j \to v$ in $X$, then $v=0$.
	In that case, the closure $\ol{T}$ of $T$ is the unbounded operator in $X$ with domain the set of $u \in X$ for which there exist a sequence $(u_j)_j \subset D(T)$ and $v \in X$ such that $u_j \to u$ and $T u_j \to v$ in $X$.
	The element $v$ is uniquely determined and one sets $\ol{T} u \coloneqq v$. %\commentCK{Make explicit what is the domain of $\ol{F}$.}
	Recall that if $T_1, T_2$ are unbounded operators, the domain of its composition $T_2 T_1$ is the set 
	\begin{align}\label{domain_composition}
	D(T_2T_1) =\{u \in D(T_1): T_1 u\in D (T_2)\},
	\end{align}
	see, e.g., \cite[Chapter 3, Section 5.1]{Kato76}.

    %%%%%%%%%%%%%%%%%%%%%%%%%%%%%%%%%%%%%%%5

    \subsection{Recap on nonlocal Sobolev spaces}\label{sec:recap_nl_sob}
    Throughout this paper \(\rho:\R^n\to [0,\infty]\) will always denote an even (radially symmetric) kernel, whose radial profile is given by \(\ol{\rho}:\R\to [0,\infty]\). As minimum assumptions on the kernel, we adopt from \cite{BMCS25} the hypotheses 
    \begin{equation}\tag{H0}
        \begin{cases}
            \rho:\R^n\to [0,\infty] \text{ is radial with } 
            \rho(x)\in \R \text{ for all } x\in \R^n\setminus\{0\},\\
            \rho\in L^1_{\mathrm{loc}}(\R^n) \text{ with } 
            \displaystyle \int_{\R^n}\min\{1,|x|^{-1}\}\rho(x)\dl x <\infty,\\
            \inf_{B_{\delta_0}}\rho >0 \text{ for some } {\delta_0}>0.
            \end{cases}
    \end{equation}
In particular, (H0) ensures that for any \(u\in C_c^\infty(\R^n)\) and \(U\in C_c^\infty(\R^n;\R^n)\) the nonlocal gradient \(D^\rho_0 u\) defined as
    \begin{equation}\label{eq:def:drho}
        D^\rho_0 u(x)
        \coloneqq \int_{\R^n}\frac{u(x)-u(y)}{|x-y|} \frac{x-y}{|x-y|}\rho(x-y)\dl y , \qquad x \in \R^n
    \end{equation}
lies in \(L^1(\R^n;\R^n)\cap L^\infty (\R^n;\R^n)\),     and the nonlocal divergence 
    \begin{equation}\label{eq:def:divrho}
        \diver^\rho_0 U(x)
        \coloneqq \int_{\R^n}\frac{U(x)-U(y)}{|x-y|} \cdot\frac{x-y}{|x-y|}\rho(x-y)\dl y , \qquad x \in \R^n
    \end{equation}
lies in \(L^1(\R^n)\cap L^\infty (\R^n)\) (see \cite[Lemma~2.3]{BMCS25}). Then, for \(u\in L^1(\R^n)+L^\infty(\R^n)\) the distributional, or weak, nonlocal gradient is the distribution defined through the relation %\commentCK{Specify what an object $D^\rho u$ is?}\comment{like this?}
    \begin{equation}
       % \langle D^\rho u,\Phi\rangle 
       \langle D^\rho u,\Phi\rangle_{\Dcal'\times \Dcal} 
        =-\int_{\R^n} u(x)\diver^\rho_0 \Phi(x)\dl x ,
        \qquad
        \Phi\in C_c^\infty (\R^n;\R^n).
    \end{equation}
    In particular, thanks to \cite[Proposition~3.2]{BMCS25}, for \(u\in C_c^\infty(\R^n)\) its pointwise and distributional nonlocal gradient agree, that is,
    \begin{equation}\label{eq:pw_equals_dual}
        D^\rho_0 u=D^\rho u 
        \quad \text{a.e. in } \Om,
        \qquad  u\in C_c^\infty(\R^n).
    \end{equation}
    Moreover, as the proof of \cite[Lemma~2.3]{BMCS25} shows, \(D^\rho_0 u\) and \(\diver^\rho_0 u\) are actually well defined for any bounded Lipschitz function \(u\in C^{0,1}(\R^n)\cap L^\infty(\R^n)\). In particular, one can quickly check that equality \eqref{eq:pw_equals_dual} extends to such functions (see also the comments after \cite[Definition~3.3]{BMCS25}). If \(\rho\) has compact support, then the same remains true without the $L^{\infty}$-requirement, i.e., \(D^\rho_0 u\) and \(\diver^\rho_0 u\) are well defined and identity \eqref{eq:pw_equals_dual} is valid for all \(u\in C^{0,1}(\R^n)\).
    
    We recall now from \cite[Definition~3.4]{BMCS25} the definitions of the \(\rho\)-nonlocal Sobolev 
    %\commentC{Decide the name for these spaces. This is OK, but somehow contradicts the view taken in the title} 
    spaces \(H^{\rho,p}(\R^n)\) and \(H^{\rho,p}_0(\Om)\) for   \(p\in (1,\infty)\) and $\Omega \subset \R^n$ open:
    \begin{align}
        \begin{split}\label{eq:def_spaces}
            H^{\rho,p}(\R^n)
        &\coloneqq \{u\in L^p(\R^n):\, D^\rho u\in L^p(\R^n;\R^n)\},\\
        H^{\rho,p}_0(\Om)
        &\coloneqq \{u\in H^{\rho,p}(\R^n):\, u=0 \text{ a.e. in } \Om^c\};
        \end{split}
    \end{align}
    the latter complementary-value space is the nonlocal analogue of classical Sobolev spaces with vanishing trace.
    Equipped with the norm 
    \begin{equation}\label{eq:norm_def}
        \norm{u}{H^{\rho,p}(\R^n)}
        \coloneqq \left(\norm{u}{L^p(\R^n)}^p + \norm{D^\rho u}{L^p(\R^n;\R^n)}^p\right)^{1/p} ,
    \end{equation}
    they form separable reflexive Banach spaces, and Hilbert spaces for \(p=2\). Moreover, as shown in \cite[Theorem~3.9]{BMCS25}, the space \(C_c^\infty(\R^n)\) lies dense in \(H^{\rho,p}(\R^n)\) and also  \(C_c^\infty(\Omega)\) in \(H_0^{\rho,p}(\Om)\) for open and bounded \(\Om\) with Lipschitz boundary. 

    An important insight related to the theory above is that for functions in \(C_c^\infty (\R^n)\) the nonlocal gradient can be written as the convolution of the classical gradient with the function \(Q_\rho\) defined as
    \begin{equation*}
        Q_\rho(x)
        \coloneqq \int_{|x|}^\infty \frac{\ol{\rho}(t)}{t}\dl t,
        \qquad x\in \R^n\setminus\{0\}.
    \end{equation*}
    Precisely, for \(u\in C_c^\infty (\R^n)\) it holds  (see \cite[Proposition~2.6~\textit{(i)}]{BMCS25})
    \begin{equation}\label{eq:translation}
        D^\rho_0 u
        = Q_\rho*\nabla u
        = \nabla (Q_\rho*u).
    \end{equation}
    If \(\rho\in L^1(\R^n)\), then also \(Q_\rho\in L^1(\R^n)\).
For further properties of \(Q_\rho\), we refer to \cite[Lemma~2.5]{BMCS25}.

%We introduce the space $\Tcal (\R^n)$ as the set of $u \in C^{\infty} (\R^n)$  such that \comment{\(\delta\)? maybe \(\p\)?} $\delta^{\alpha} u \in L^1 (\R^n) \cap L^{\infty} (\R^n)$ for all $\alp \in \N^n$.
%For $\R^n$-valued functions we use the notation $\Tcal (\R^n, \R^n)$.

%\begin{lemma}\label{le:Tcal}
%Let $\rho \in L^1 (\R^n)$.
%If $u \in \Tcal (\R^n)$ then $D^{\rho} u \in \Tcal (\R^n, %\R^n)$.
%If $u \in \Tcal (\R^n, \R^n)$ then $\diver^{\rho} u \in %\Tcal (\R^n)$.
%\end{lemma}
%\begin{proof}
%Let $u \in \Tcal (\R^n)$ and let $\{ u_j \}_{j \in \N}$ be a sequence in $C^{\infty}_c (\R^n)$ tending to $u$ ***

%The second part is analogous.
%\end{proof}
     \subsection{Fourier transform}\label{sec:fourier}
    We write \(\mathcal{S}(\R^n)\) for the space of Schwartz functions over \(\R^n\) and denote for \(f\in \mathcal{S}(\R^n)\) by \(\mathcal{F}f\) (if sensible, also expressed as \(\hat{f}\) for brevity) the Fourier transform defined as
    \begin{equation*}
        \mathcal{F}f(\xi)
        \coloneqq \int_{\R^n}f(x)e^{-2\pi ix\cdot \xi}\dl x.
    \end{equation*}
    The inverse Fourier transform will be denoted by \(\mathcal{F}^{-1}\) and is defined through the identity \(\mathcal{F}^{-1} f(\xi)=\mathcal{F}f(-\xi)\) for $\xi\in \R^n$. The Fourier transform and its inverse are extended to \(L^1(\R^n)+L^2(\R^n)\) and to the tempered distributions \(\mathcal{S}'(\R^n)\) as usual (see, e.g., \cite[Chapter~2]{Gr14}).    
   Observe that if we denote by \(\widetilde{f}\coloneqq f(-\,\cdot)\) the reflection of a function \(f\), then %by definition of the inverse Fourier transform,
   \begin{equation}\label{eq:FourierReflection}
 \Fcal^{-1}f = \Fcal \widetilde{f} = \color{black}\widetilde{\Fcal f}
\end{equation} 
for any function $f$ for which $\Fcal^{-1} f$ is well defined. 

    A central role in the subsequent discussion will be played by the function 
     \begin{equation*}\label{eq:lambdarho}
        \lambda_\rho(\xi)
        \coloneqq \int_{\R^d}\frac{\rho(x)x}{|x|^2}\left(e^{2\pi i\xi\cdot x}-1\right)\dl x,
         \qquad \xi \in \R^n,
    \end{equation*}
    which also admits the representation (see \cite[Lemma~2.4]{BMCS25})
    \begin{equation*}
        \lambda_\rho(\xi)
        =
        \frac{i\xi}{|\xi|}\int_{\R^n} \frac{\rho(x)x_1}{|x|^2}\sin(2\pi|\xi|x_1)\dl x,
         \qquad \xi \in \R^n\setminus \{0\}.
    \end{equation*}
    We point out that the \(x_1\) in the integrand above can be replaced by any \(x_i\) for \(i=1,\ldots,n\) due to the radial symmetry of $\rho$.
    The importance of \(\lambda_\rho\) roots in the fact that it is the Fourier symbol associated to the nonlocal gradient and divergence, i.e., for \(u\in C_c^\infty(\R^n)\) and \(\Phi\in C_c^\infty(\R^n;\R^n)\) it holds (see \cite[Lemma~2.4]{BMCS25})
    \begin{equation}\label{eq:nlgrad_symbol}
        \mathcal{F}(D^\rho_0 u)
        = \lambda_\rho\hat{u},
        \qquad 
        \mathcal{F}(\diver^\rho_0 \Phi)
        = \lambda_\rho\cdot \hat{\Phi}.
    \end{equation}
    This extends readily to \(u\in \mathcal{S}(\R^n)\) and \(\Phi\in \mathcal{S}(\R^n;\R^n)\), respectively; we note that the Schwartz space is included in \(H^{\rho,p}(\R^n)\) for all \(p\in [1,\infty]\) (see  \cite[Proposition~3.5]{BMCS25}) and thus the Fourier transforms above are well defined. 
    For the Fourier inverse of the nonlocal divergence of \(\Phi\in \mathcal{S}(\R^n;\R^n)\), one has
%\begin{lemma}\label{lem:fourier_identities_as}
           \begin{equation}\label{eq:fourier_identities_as}
    		\mathcal{F}^{-1}(-\diver^\rho_0 \Phi)
            = \lambda_\rho\cdot\mathcal{F}^{-1}{\Phi};
            \end{equation}        
%\end{lemma}
%\begin{proof}
indeed, the fact that \(\lambda_\rho\) is odd, \eqref{eq:FourierReflection} and \eqref{eq:nlgrad_symbol} imply %for all $\xi \in \R^n$, 
   % \begin{equation*}
  %  		\mathcal{F}^{-1}(-\diver^\rho_0 \Phi)(\xi)
   %         =\mathcal{F}(-\diver^\rho_0 \Phi)(-\xi)
   %         = -\lambda_\rho(-\xi) \cdot \mathcal{F}{\Phi}(-\xi)
   %         = \lambda_\rho(\xi)\cdot\mathcal{F}^{-1}{\Phi}(\xi).
   %         \end{equation*}     
      %  \comment{Alternatively: 
        \begin{equation*}
    		\mathcal{F}^{-1}(-\diver^\rho_0 \Phi)
            =\mathcal{F}(-\widetilde{\diver^\rho_0 \Phi})
            = -\widetilde{\lambda_\rho} \cdot\widetilde{ \mathcal{F}{\Phi}}
            = \lambda_\rho\cdot\mathcal{F}^{-1}{\Phi}.
            \end{equation*}   
            %}
%\end{proof}
 %In fact, as we show later on, \(\lambda_\rho\) actually represents the Fourier symbol in the sense of tempered distributions of the weak nonlocal gradient \(D^\rho\) over the whole space \(H^{\rho,p}(\R^n)\) (see Theorem~\ref{thm:alt_char}). 

    %\comment{check reference in the end}. \commentC{Exactly: Theorem~\ref{thm:alt_char} is not the right reference}\comment{I actually added the result to Theorem~\ref{thm:alt_char} (but this is of course up for discussion)}
    
    Furthermore, if \(\rho\in L^1(\R^n)\), then \(Q_\rho\in L^1(\R^n)\) and its Fourier transform is well defined.
    %and it holds (see \cite[Proposition 2.6 (ii)]{BMCS25})
   % \begin{equation}\label{eq:fourier_grad_qp}
    %    \mathcal{F}(D^\rho_0 u)
    %    = 2\pi i\xi\widehat{{Q}_\rho}\hat{u},
    %    \quad u\in C_c^\infty(\R^n).
    %\end{equation}
    In particular, combining \eqref{eq:translation} and \eqref{eq:nlgrad_symbol}, one finds that (see \cite[Proposition~2.6~\textit{(ii)}]{BMCS25})
    \begin{equation*}
        \lambda_\rho(\xi)
        = 2\pi i\xi\widehat{{Q}_\rho}(\xi),
         \qquad \xi \in \R^n
    \end{equation*}
   and 
    \begin{equation}\label{eq:qrho_hat}
        \widehat{{Q}_\rho}(\xi)
        = \frac{1}{2\pi|\xi|}\int_{\R^n} \frac{\rho(x)x_1}{|x|^2}\sin(2\pi|\xi|x_1)\dl x,
         \qquad \xi \in \R^n\setminus \{0\}.
    \end{equation}
    Starting from the identity above, the authors of \cite{BMCS25} studied the properties of \(\widehat{{Q}_\rho}\) and devised sufficient conditions ensuring, in particular, that \(\widehat{{Q}_\rho}\) is positive and satisfies certain lower bounds. We will recall the corresponding results in Section \ref{sec:hypotheses}.
    Later on in Lemma~\ref{lem:lambda_FintoF_pre} it is derived that, as a consequence of these conditions, \(\lambda_\rho\), understood as a multiplier in Fourier space, maps \(\mathcal{S}(\R^n)\) into \(\mathcal{S}(\R^n;\R^n)\). 

     It is also worth mentioning that the Fourier transform preserves radiality (see \cite[Corollary~2.2.12]{Gr14}); thus, since \(\rho\) is radial by hypothesis (H0), \(\widehat{{Q}_\rho}\) is again radial, while \(\lambda_\rho\) is vector radial.

\subsection{Fourier multiplier theory}\label{sec:fm_theory} %\commentCK{Do we have a comment on our default choice of $p$ somewhere?}\comment{we kind of fix it above \eqref{eq:def_spaces}, but we could also make that more explicite}
Let us now turn to $L^p$-Fourier multipliers and, more generally, to unbounded linear operators on $L^p(\R^n)$ defined via multiplication in Fourier space.
%beginning with some preliminaries.
 %   We now turn to the matter of \(L^p\)-Fourier multipliers and, more general, linear unbounded operators on \(L^p(\R^n)\) defined through multiplication in Fourier space, starting with some preliminary considerations. 
   % \commentCK{Slightly streamlined and  rephrased, please check if ok?}\comment{yes, sounds good}
   % We prepare our exposition by some preliminary considerations. The general setup is as follows.
     %  If $w\in \Scal'(\R^n)$ is a tempered distribution, then we say that $\Fcal^{-1}w\in L^p(\R^n)$ \commentCK{Parentheses removed, ok?} if there exists a (necessarily unique) $v\in L^p(\R^n)$ such that
    For a tempered distribution $w \in \Scal'(\R^n)$, we write $\Fcal^{-1}w \in L^p(\R^n)$ if there
exists a (necessarily unique) $v \in L^p(\R^n)$ such that
\begin{equation}\label{eq:def_fourier_charac}
    \langle v, \varphi \rangle = 
   \langle  w, \mathcal{F}^{-1}\varphi\rangle_{\mathcal{S}'\times \mathcal{S}},
    \qquad \varphi\in \Scal(\R^n),
\end{equation}
and we identify $\Fcal^{-1}w$ with $v$. %\commentCK{Parentheses removed, ok?}\comment{yes}
Since this convention appears in two contexts below (see Sections~\ref{sec:charac:fourier}
and~\ref{sec:char_laplacian}), we briefly discuss it here. 
%\commentCK{Is it better to introduce the terminology for  distributions here? Check again for consistency with the definition of the $\rho$-Laplacian.}\commentCK{Strictly speaking, this general observation about (tempered) distributions lying in $L^p(\R^n)$ could as well  become part of the notations (and terminology) section. Let's discuss}
%This convention will appear in two contexts in this paper (see Sections~\ref{sec:charac:fourier} and \ref{sec:char_laplacian}), so we make its precise meaning explicit.
%This notion appears in two contexts below, so its precise meaning is made explicit.
%We use this notion in two contexts and briefly make its meaning precise.
%This will be used in two different contexts, and, even though the resulting notation is standard, we consider it worthwhile to explain the exact meaning of the respective expressions which will appear repeatedly throughout the paper. 

% First, let
 If \(m\in C^\infty(\R^n)\) is such that all its derivatives have polynomial growth and \(u\in \mathcal{S}'(\R^n)\), then $m\hat u$ and $\Fcal^{-1}(m\hat u)$ are again
tempered distributions, defined as usual by duality, i.e.,
\begin{equation*}\label{eq:fourier_first_case}
    %\langle \mathcal{F}^{-1}(m\hat{u}), \eta\rangle \coloneqq 
    \langle \mathcal{F}^{-1}(m\hat{u}), \varphi\rangle_{\mathcal{S}'\times \mathcal{S}}
    = \langle  u, \mathcal{F}(m\mathcal{F}^{-1}\varphi)\rangle_{\mathcal{S}'\times \mathcal{S}},
    \qquad \varphi\in \Scal(\R^n),
\end{equation*}
where we note that since \(m \mathcal{F}^{-1}\varphi\in \mathcal{S}(\R^n)\), its Fourier transform is guaranteed to be well defined. 
Choosing $w = m\hat u$, we understand $\Fcal^{-1}(m\hat u) \in L^p(\R^n)$ in the sense of~\eqref{eq:def_fourier_charac}.
%\commentCK{ I suggest to remove this here: These ideas underlie the Fourier characterization in Section~\ref{sec:charac:fourier}, already incorporated above. }\comment{Yes, let's remove it here then.}
%\commentCK{Check this place as option for symmetry result.}

%Another case will appear. If $m : \R^n \to \R$ is merely a measurable function with polynomial growth but \(u\in \mathcal{S}(\R^n)\), then the product $m \hat{u}$ defines a tempered distribution as

If, on the other hand, $m \colon \R^n \to \R$ is merely measurable and of
polynomial growth, but $u \in \Scal(\R^n)$, then $m\hat{u}$ is integrable and $m\hat{u}$
defines a tempered distribution via 
    \begin{equation*}\label{eq:fourier_second_case}
     \langle m \hat{u} , \varphi \rangle_{\Scal' \times \Scal} 
     = \langle m \hat{u} , \varphi \rangle
     = \int_{\R^n} m \hat{u} \varphi \dl x, \qquad \varphi \in \Scal(\R^n).
    \end{equation*}
By taking $w = m\hat{u}$ in~\eqref{eq:def_fourier_charac}, we can therefore
again give meaning to the expression $\Fcal^{-1}(m\hat{u})$. Moreover, since
$m\hat{u} \in L^1(\R^n)\cap L^2(\R^n)$, this distributional inverse Fourier transform coincides with the pointwise one and, in particular,
$\Fcal^{-1}(m\hat{u}) \in L^2(\R^n)\cap L^\infty(\R^n)$.

For later use, we state an auxiliary result on symmetry relations for multiplication in Fourier space by even or odd functions.% \commentCK{Is this transition sentence ok?}\comment{Yes}

  \begin{lemma}\label{lem:fourier_identities}
Assume that $m_1,m_2: \R^n \to \R$ are measurable functions with polynomial growth such that \(m_1\) is even and \(m_2\) is odd. Then, for any \(u\in \mathcal{S}(\R^n)\), 
%\commentCK{How do we know that the domain of $F_{m_1}$ and $F_{m_2}$ is $\Scal(\R^n)$? cf. Lemma~\ref{le:Fmclosable}, but there, stronger assumptions on $m$ are required. Make consistent. Do we want to avoid writing $F_{m_1}$ etc. here and move further up?}\comment{We shouldn't use $F_{m_1}$ and $F_{m_2}$ here, we introduced it at some point to make the notation more consistent, but it was a mistake as you point out. While the expressions \(F_{m_i}\varphi\) are well defined as functions, they are not as operators in general; but the latter is needed neither here nor in the applications later.}
        \begin{enumerate}[label=(\roman*)]
            \item\label{item:fi_symm}
            $\Fcal^{-1}(m_1\hat{u})
    		= \mathcal{F}(m_1\mathcal{F}^{-1}u)$ and 
            $  \Fcal^{-1}(\Fcal^{-1}(m_1\hat{u}))
    			=m_1\mathcal{F}^{-1}u$;
%If, in addition, \(m_1\) is an \(L^p\)-Fourier multiplier, then the equalities above extend to all \(\varphi\in L^p(\R^n)\). \commentC{Do we need this only for (i) or also for (ii)?}
        \item\label{item:fi_antisymm}
$\Fcal^{-1}(m_2\hat{u}) = -\mathcal{F}(m_2\mathcal{F}^{-1}u)$ and $\Fcal^{-1}(\Fcal^{-1}(m_2\hat{u})) = -m_2\mathcal{F}^{-1}u$.
        \end{enumerate}
   %     \commentCK{ In particular, if $u\in D(F_m)$, then, $F_{m_1}u = \mathcal{F}(m_1\mathcal{F}^{-1}u)$ and $F_{m_2}u = -\mathcal{F}(m_2\mathcal{F}^{-1}u)$. Check! Best if Leon decides. }
	\end{lemma}
	\begin{proof}
Let $m \in \{m_1, m_2\}$ and denote by $\widetilde{m}$ its reflection. 
Note that, as \(m\) has polynomial growth, for \(u\in \mathcal{S}(\R^n)\), the expression \(\mathcal{F}^{-1}(m \hat{u})\) is always well defined, since \(m \hat{u}\in L^1(\R^n)\).
Using~\eqref{eq:FourierReflection},
%that
%\begin{equation}\label{eq:FourierReflection}
% \Fcal^{-1} \psi = \Fcal \widetilde{\psi} = %\tilde{\hat{\psi}}
%\end{equation}
%for any function $\psi$ for which $\Fcal^{-1} \psi$ is defined, \color{black} 
we have for $x \in \R^n$, 
		\begin{align*}
			\begin{split}\label{eq:fi}
			    \mathcal{F}^{-1}(m\hat{u})(x)
%			&= \mathcal{F}(m\hat{\varphi})(-x)\\
			&= \int_{\R^n}m(\xi)\hat{u}(\xi)e^{2\pi i x\cdot \xi}\dl\xi\\
			&=\int_{\R^n}m(\xi)\hat{u}(\xi)e^{-2\pi i x\cdot (-\xi)}\dl\xi\\
			&=\int_{\R^n} \widetilde{m}(\xi) \widetilde{\Fcal u}(\xi)e^{-2\pi i x\cdot \xi}\dl\xi\\
			&=\mathcal{F}(\widetilde{m}\mathcal{F}^{-1} u)(x).
			\end{split}
		\end{align*}
       % \commentCK{Consistent notation $\tilde{m}$ vs. $\widetilde{m}$}\comment{I would choose $\widetilde{m}$}
     %  \commentCK{Fixed inconsistency $\varphi$ vs. $u$, check if ok now? }
With this, we obtain the first equalities of \textit{\ref{item:fi_symm}} and \textit{\ref{item:fi_antisymm}}.
The second ones follow by taking the inverse Fourier transform of both sides of the first one.
		%If \(m\) defines an \(L^p\)-Fourier multiplier, the result translates to the general case that \(\varphi\in L^p(\R^n)\) by density. 
        %\commentC{We have to conciliate this with the language of $\bar{F}$. The LHS is clear: this is $\ol{F_{m_1}} \varphi$. For the RHS one we haven't give a name. Maybe $\ol{G_{m_1}}$ with $G_{m_1} \varphi = \mathcal{F}(m_1\mathcal{F}^{-1}\varphi)$. The second LHS is even worse. Cannot we simply use that $\Fcal^{-1} \Fcal^{-1} \psi= \tilde{\psi}$? Let's see where and how this is used and check whether we need it for $\varphi \in L^p$, or we just need the final formula. Finally, an option that does not require extra notation is $\ol{\Fcal M_m \Fcal^{-1}}$, with a previous definition of the multiplication operator $M_m$}
        %\comment{I checked that we only need the result for functions in \(\Scal\). I also changed the left-hand sides. For the rest I actually prefer to leave them  Thus, I think that we can just leave the statement as it is to keep the notation as easy as possible.}
The lemma is proved.
	\end{proof}

   We next define the following unbounded multiplier operators, which will play a prominent role throughout Section~\ref{sec:char_laplacian}. 
   %\commentCK{Some kind of reference? Where does this come from, or "home-made". Ask Carlos.}\comment{more homemmade. Of course, people have treated this kind of unbounded operators in the literature, but I don't think that exactly this definition is fixed anywhere.}
    
%    \begin{definition}\label{def:ub_mult_op_pl}
%Let $m : \R^n \to \R$ be a measurable function with polynomial growth, and let \(u\in \mathcal{S}(\R^n)\).
%
%We write \(\Fcal^{-1} (m \hat{u}) \in L^p (\R^n)\) if there exists \(v\in L^p(\R^m)\) (necessarily unique) such that 
% If there exists \(v\in L^p(\R^m)\) (necessarily unique) such that  \begin{equation*}
%        \langle v,\varphi\rangle 
%        =\langle m \hat{u} , \Fcal^{-1} \varphi \rangle,
%        \qquad \varphi\in C_c^\infty(\R^n),
%    \end{equation*}
%    we write $\Fcal^{-1} (m \hat{u}) = v$ and say that \(\Fcal^{-1} (m \hat{u}) \in L^p (\R^n)\).
%In this case, we identify $\Fcal^{-1} (m \hat{u}) = v$.
%\commentCK{Can we improve the wording for clarity?}
%\end{definition}
\begin{definition}\label{def:ub_mult_op_pl} 
Let \(p\in (1,\infty)\). Let $m : \R^n \to \R$ be a measurable function with polynomial growth, and let \(u\in \mathcal{S}(\R^n)\).
The unbounded operator $F_m$ in $L^p (\R^n)$ with domain
    \[
     D(F_m) \coloneqq\{ u \in \Scal(\R^n) : \Fcal^{-1} (m \hat{u}) \in L^p (\R^n) \}
    \]
    is defined as
    \[
     F_m u \coloneqq \Fcal^{-1} (m \hat{u}) , \qquad u \in D(F_m).
    \]
    \end{definition}
 Note that in the case of $p\in [2, \infty]$, one automatically has $D(F_m)=\Scal(\R^n)$ as $\Fcal^{-1}(m\hat{u})\in L^p(\R^n)$ for all $u\in \Scal(\R^n)$. 
%    \commentCK{With the above, we would have that $D(F_m)$ is always $\Scal(\R^n)$ for $p\geq 2$. It seems that we should make a comment in this direction... What about closability?}
    %\comment{Yes, let's either include this into the sentence above, where you also mention this, or here, or maybe even at both places. We could point out that the question of the domain therefore is mainly interesting in the case \(p\in (1,2)\). However, I think that the issue of closability is independent of this.}
%    \begin{remark}\label{rem:fi_conincidence}
%    \comment{added}
%        If for \(u\in \Scal(\R^n)\) the expression \(\Fcal^{-1} (m \hat{u})\) exists in \(L^p(\R^n)\) as the pointwise Fourier inverse, then it naturally follows by uniqueness of the \(L^p\)-representative that \(F_m u=\Fcal^{-1} (m \hat{u})\). \commentCK{I suggest to remove this remark. That $m\hat{u}\in L^1(\R^n)$ in this set-up,  and hence, $\F^{-1}(m)$}
%    \end{remark}
%
%\commentCK{Maybe have a paragraph after all. Also, this holds more generally (and is used later) not just for operators $F_m$. Can we find a way to reflect this better?}

    \begin{notation}\label{def:op_extension}
    If the operator $F_m$ in Definition~\ref{def:ub_mult_op_pl} is closable (see Section~\ref{sec:notation}), its closure in \(L^p(\R^n)\) is denoted by $\ol{F_{m}}$. %\commentCK{It seems liek the same notation already introduced above? Am I overlooking something?}
    \end{notation}
    
    %\color{blue}
%\begin{remark}
%   choice of  $\Scal(\R^n)$ or $C_c^\infty(\R^n)$ in the definition of the domain of $F_m$ does not make a difference for the closures $\overline{F_m}$.
%    \end{remark} \color{black}

    %\comment{\(\ol{F_m}\) or \(\ol{F}_m\)? I like the first one more, but the overline %becomes very long for long expressions of the multiplier. I guess we should change to the first one but I leave it for now.} \commentC{Even though I prefer typographically the second, the comment in next paragraph shows that the first one is the only consistent}
 If the domain of $\ol{F_{m}}$ is \(L^p(\R^n)\), then $m$ is called an $L^p$-Fourier multiplier, while \(\ol{F_{m}}\) is called an \(L^p\)-Fourier multiplier operator.
    It is well known that for \(p=2\) the condition \(m\in L^\infty(\R^n)\) is necessary and sufficient for $m$ to be an  $L^2$-Fourier multiplier
%, or, equivalently, that $F_m$ can be uniquely extended to an \(L^2\) Fourier multiplier operator
%or, equivalently, that $\overline{F_m}$ is an \(L^2\) Fourier multiplier operator.
(see \cite[Theorem~2.5.10]{Gr14}). For general \(p\in (1,\infty)\), one needs to require stronger properties for $m$ to extend \(F_m\) to a bounded operator in all of \(L^p(\R^n)\). Sufficient conditions are, for example, devised by the Mihlin-H\"ormander multiplier theorem, which will be used repeatedly throughout the paper (see, e.g.,~\cite[Theorem~6.2.7]{Gr14}).

  The next lemma shows that for smooth multipliers $m$ with derivatives of polynomial growth, the operator $F_m$ has the
following basic properties, which, in particular, allow us to define its closure~$\overline{F_m}$.
  %  For smooth multipliers with derivatives of polynomial growth, one obtains these basic properties $F_m$ has thefollowing basic properties, which in particular allow to define the closure $\overline{F_m}$. 
\begin{lemma}\label{le:Fmclosable}
Let $m : \R^n \to \R$ be a smooth function such that all its derivatives have polynomial growth.
Then the unbounded operator $F_m$ in $L^p(\R^n)$ maps from \(D(F_m)=\mathcal{S}(\R^n)\) into $\Scal(\R^n)$ and is closable. 
Moreover, for all \(u\in D(\ol{F_m})\) the expression \(\mathcal{F}^{-1}(m\hat{u})\) exists as tempered distribution in \(L^p(\R^n)\), and it holds
\begin{equation*}
  \Fcal^{-1}(m\hat{u}) 
  =  \ol{F_m} u.
    %\qquad u\in D(\ol{F_m}).
\end{equation*}
%\commentCK{Do we have a default setting for $p$?} % and it maps into $\Scal(\R^n)$. 
%Moreover, it holds the identity
%\begin{equation*}
%    F_m u = \Fcal^{-1} (m\hat{u}),
%    \qquad u\in \Scal(\R^n),
%\end{equation*}
%where the right-hand side denotes the pointwise Fourier inverse, which lies again in \(\Scal(\R^n)\).
\end{lemma}
\begin{proof}
Let $u \in \mathcal{S}(\R^n)$.
Then $\hat{u} \in \mathcal{S}(\R^n)$, and since all derivatives of $m$ have polynomial growth,  $m \hat{u} \in \mathcal{S}(\R^n)$.
Consequently, \(\Fcal^{-1} (m\hat{u})\) is well defined and  lies in \(\Scal(\R^n)\subset L^p (\R^n)\).
%\footnote{Removed the following, ok? and by Remark~\ref{rem:fi_conincidence} we have $F_mu= \Fcal^{-1} (m\hat{u}) \in \mathcal{S}(\R^n)\subset L^p (\R^n)$.}

Now we prove that $F_m$ is closable. 
Let $(u_j)_j$ be a sequence in $D(F_m)$ such that $u_j \to 0$ in $L^p (\R^n)$ and $F_m u_j \to v$ in $L^p (\R^n)$ for some $v \in L^p (\R^n)$.
%Then $u_j \to 0$ in $\Scal'(\R^n)$, so $\widehat{u_j} \to 0$ in $\Scal'(\R^n)$.
Using that $F_m u_j \to v$ in $\Scal'(\R^n)$, we obtain that for all $\varphi \in \Scal(\R^n)$, 
\[
 \lim_{j\to \infty}\langle F_m u_j , \varphi \rangle
 = \langle v ,\varphi \rangle.
\]
But, as $\Fcal (m \mathcal{F}^{-1}\varphi) \in \Scal(\R^n)$ by the assumption on \(m\),
\[
 \lim_{j\to \infty}\langle F_m u_j , \varphi \rangle
 = \lim_{j\to \infty}\langle m \widehat{u_j} , \mathcal{F}^{-1}\varphi \rangle
 = \lim_{j\to \infty}\langle u_j , \Fcal(m \mathcal{F}^{-1}\varphi)\rangle= 0 .
\]
Thus, $\langle {v} , \varphi \rangle = 0$ for all $\varphi \in \Scal(\R^n)$, yielding $v=0$.

Finally, in order to see the last claim, let $u\in D(\ol{F_m})$ and  $(w_j)_j\subset D(F_m)$ with \(w_j\to u\) in \(L^p(\R^n)\) and such that \(F_m w_j \to \ol{F_m}u\) in \(L^p(\R^n)\). Then, one finds
\begin{align*}
    \langle \ol{F_{m}} u, \varphi\rangle
    =\lim_{j\to \infty}\langle F_m w_j, \varphi \rangle
    =  \lim_{j\to \infty}\langle  w_j, \mathcal{F}(m\mathcal{F}^{-1}\varphi)\rangle
    = \langle  u, \mathcal{F}(m\mathcal{F}^{-1}\varphi)\rangle,
    \qquad \varphi\in \Scal(\R^n),
\end{align*}
so that \(\mathcal{F}^{-1}(m\hat{u})\) as a tempered distribution is equal to \(\ol{F_{m}} u\) in \(L^p(\R^n)\), concluding the proof.
\end{proof}

%%%%%%%%%%%%%%%%%%%%%%%%%%%%%%%%%%%%%%%%%%%%%%%%%%%%%%%%%%
Naturally, and for our later purposes in Section~\ref{sec:char_laplacian} very important, the question arises, how concatenations of unbounded operators as introduced in Definition~\ref{def:ub_mult_op_pl} behave. 
In that respect, one always has at least that \(F_{m_2m_1}\) extends \(\ol{F_{m_2}}F_{m_1}\)
as %(potentially)
unbounded operators. For the reader's convenience, we give below a detailed proof of this property.
We add that this behavior is well known in the case of \(L^p\)-Fourier multipliers (see \cite[Proposition~2.5.13]{Gr14}), where even the equality \(\ol{F_{m_2}}\, \ol{F_{m_1}}=\ol{F_{m_2m_1}}\) holds. 
%\commentCK{Can we show the analogous result for the extensions (under suitable assumptions?), that is, \(\ol{F_{m_2}}\,\ol{F_{m_1}}\) is extended by $\ol{F_{m_2m_1}}$ as unbounded operators? This would automatically yield $\overline{F_{|\lambda_\rho|}}^2$ being extended by $\ol{F_{|\lambda_\rho|^2}}$, and similarly that $\overline{F_{\bla{\lambda_\rho}^2}}$ extends $\overline{F_{\bla{\lambda_\rho}}}^2$, etc...}
%\commentC{With this formula I agree that the only sensible way to write the closure is $\ol{F_{m}}$}

%\commentCK{Introductory sentence missing! }
%***25/7/2026
\begin{lemma}\label{lem:mults_algebra}
        Let $m_1,m_2 : \R^n \to \R$ be measurable functions with polynomial growth, consider the operators $F_{m_1}, {F_{m_2}}$ in $L^p(\R^n)$, and further assume that $F_{m_2}$ is closable. Then
\begin{equation*}
 D(\ol{F_{m_2}}F_{m_1})\subset D(F_{m_2m_1})
\end{equation*}
and 
        \begin{equation}\label{aux:algebra}
            \ol{F_{m_2}}F_{m_1}u= F_{m_2m_1}u,
            \qquad u\in D(\ol{F_{m_2}}F_{m_1}).
        \end{equation}
    \end{lemma}
    \begin{proof}
       % \textit{Step~1: Restriction to functions in \(\mathcal{S}(\R^n)\).}
Let \(u\in D(\ol{F_{m_2}}F_{m_1})\), recalling~\eqref{domain_composition}.
As $F_{m_1}u \in D(\ol{F_{m_2}})$, there exists a sequence \((v_j)_j\subset D (F_{m_2})\) such that \(v_j\to F_{m_1} u\) in \(L^p(\R^n)\) and 
    \begin{equation*}
        F_{m_2} v_j\to \ol{F_{m_2}}F_{m_1} u
        \quad \text{in } L^p(\R^n).
    \end{equation*}
    At the same time, by Definition \ref{def:ub_mult_op_pl}, for all $\varphi\in C_c^\infty(\R^n)$,
    \begin{equation*}
        \langle F_{m_2} v_j, \varphi\rangle
        =\langle m_2 \widehat{v_j}, \mathcal{F}^{-1} \varphi \rangle
        = \langle \widehat{v_j}, m_2 \mathcal{F}^{-1} \varphi \rangle.
    \end{equation*}
Note that \(m_2 \mathcal{F}^{-1} \varphi\in L^{p'}(\R^n)\), since $\mathcal{F}^{-1} \varphi\in \mathcal{S} (\R^n)$ and \(m_2\) has polynomial growth. Hence, we can choose a second sequence \((w_k)_k\subset \mathcal{S}(\R^n)\) with \(w_k\to m_2 \mathcal{F}^{-1}\varphi\) in \(L^{p'}(\R^n)\). On the one hand, this entails that 
    \begin{equation}\label{eq:doublelim_kj}
        \lim_{j\to \infty}\lim_{k\to \infty} \langle \widehat{v_j}, w_k\rangle
        = \lim_{j\to \infty} \langle \widehat{v_j}, m_2 \mathcal{F}^{-1}\varphi\rangle
         = \lim_{j\to \infty} \langle F_{m_2} v_j, \varphi\rangle
         = \langle \ol{F_{m_2}}F_{m_1} u, \varphi\rangle,
    \end{equation}
    while, on the other hand, 
    \begin{equation}\label{eq:weak_conv}
        \lim_{j\to \infty} \langle \widehat{v_j}, w_k\rangle
        = \lim_{j\to \infty} \langle v_j, \widehat{w_k}\rangle
        = \langle F_{m_1}u, \widehat{w_k}\rangle
        = \langle m_1 \hat{u}, w_k\rangle, \qquad k\in \N
    \end{equation}
    and 
    \begin{equation*}
        \lim_{k\to\infty}\langle m_1 \hat{u}, w_k\rangle
        = \langle m_1 \hat{u}, m_2\mathcal{F}^{-1}\varphi\rangle
        = \langle m_2m_1 \hat{u}, \mathcal{F}^{-1}\varphi\rangle,
        %= \langle F_{m_2m_1} u, \varphi\rangle,
    \end{equation*}
    and together
    \begin{equation}\label{eq:doublelim_jk}
        \lim_{k\to \infty}\lim_{j\to \infty} \langle \widehat{v_j}, w_k\rangle
        = \langle m_2m_1 \hat{u}, \mathcal{F}^{-1}\varphi\rangle.
    \end{equation}
As \(m_1\) has polynomial growth and $u \in \mathcal{S} (\R^n)$, it holds \(m_1\hat{u}\in L^p(\R^n)\).
Moreover, similarly to \eqref{eq:weak_conv} one infers that $ \langle \widehat{v_j}, \psi \rangle \to \langle m_1 \hat{u}, \psi \rangle$ for all $\psi \in \mathcal{S} (\R^n)$, and hence, $(\widehat{v_j})_j$ converges weakly to \(m_1\hat{u}\) in \(L^p(\R^n)\). In particular, \((\widehat{v_j})_j\) is bounded in \(L^p(\R^n)\). We conclude that the limit \(\lim_{k\to \infty} \langle \widehat{v_j}, w_k\rangle\) is attained uniformly in \(j\in \N\), which follows from 
    \begin{align*}
        \lim_{k\to \infty} \sup_{j\in \N}
        |\langle \widehat{v_j}, w_k\rangle
        -\langle \widehat{v_j}, m_2\mathcal{F}^{-1}\varphi\rangle|
       % &=\lim_{k\to \infty} \sup_{j\in \N}
        %|\langle \widehat{v_j},w_k- m_2\mathcal{F}^{-1}\varphi\rangle|\\
        %&
        \le \lim_{k\to \infty} ( \sup_{j\in \N}\norm{\widehat{v_j}}{L^p(\R^n)} ) \norm{w_k-m_2\mathcal{F}^{-1}\varphi}{L^{p'}(\R^n)}
        = 0.
    \end{align*}
     By a classical argument,the Moore–Osgood theorem on the interchange of limits, this implies that the two double-limits in \eqref{eq:doublelim_kj} and \eqref{eq:doublelim_jk} agree. By Definition \ref{def:ub_mult_op_pl}, \(u\in D(F_{m_2m_1})\) with 
    \begin{equation*}
        F_{m_2m_1} u 
        = \ol{F_{m_2}}F_{m_1} u,
    \end{equation*}
    which yields the statement.
    \end{proof}
   % \commentCK{Intro sentence missing here: The next lemma shows that  is an operator from Schwartz functions into itself}

    Moreover, one has the inclusion of domains \(D(\overline{F_{m_1}})\subset D(\overline{F_{m_2m_1}})\) of the respective closed operators, if at least \(m_2\) is an \(L^p\)-Fourier multiplier.
\begin{lemma}\label{lem:domain_mulitplier}
  Let $m_1,m_2 : \R^n \to \R$ be measurable functions with polynomial growth, consider the operators $F_{m_1}, F_{m_2}$ in $L^p(\R^n)$, and further assume that $F_{m_1}$ is closable.  If $m_2$ is an $L^p$-Fourier multiplier and also $F_{m_1m_2}$ is closable, then
  \begin{align*}
  D(\overline{F_{m_1}})\subset D(\overline{F_{m_2m_1}}),
  \end{align*}
  and for any \(u\in D(\overline{F_{m_1}})\) and \((u_j)_j\subset D(F_{m_1})\) with \(u_j\to u\) and \(F_{m_1}u_j\to \ol{F_{m_1}}u\) in \(L^p(\R^n)\), it holds
  \begin{equation*}
      F_{m_2 m_1} u_j\to \overline{F_{m_2m_1}} u
      \quad \text{in } L^p(\R^n). 
  \end{equation*}
\end{lemma}
%\begin{lemma}
%  Let $m_1,m_2 : \R^n \to \R$ be measurable functions with polynomial growth, consider the operators $F_{m_1}, F_{m_2}$ in $L^p(\R^n)$, and further assume that $m_2$ is an $L^p$-Fourier multiplier and that $F_{m_1m_2}$ is closable. Then, \(F_{m_1}\) is closable with
%  \begin{align*}
%  D(\overline{F_{m_1}})\subset D(\overline{F_{m_2m_1}}),
%  \end{align*}
%  and for any \(u\in D(\overline{F_{m_1}})\) and \((u_j)_j\subset D(F_{m_1})\) with \(u_j\to u\) and \(F_{m_1}u_j\to \ol{F_{m_1}}u\) in \(L^p(\R^n)\), it holds
%  \begin{equation*}
%      F_{m_2 m_1} u_j\to \overline{F_{m_2m_1}} u
%      \quad \text{in } L^p(\R^n). 
%  \end{equation*}
%\end{lemma}
\begin{proof} 
First observe that since $\ol{F_{m_2}}$ is an $L^p$-multiplier operator, and thus, $D(\ol{F_{m_2}})=L^p(\R^n)$, it holds in view of Definition~\ref{def:ub_mult_op_pl} that 
$D(F_{m_1})= \{u\in \Scal(\R^n): F_{m_1}\in L^p(\R^n)\} = D(\ol{F_{m_2}}F_{m_1})$. Consequently,
Lemma~\ref{lem:mults_algebra} implies in particular that $D(F_{m_1})\subset D(F_{m_2m_1})$. 
%Indeed, by definition (see~\eqref{domain_composition}), $D(\overline{F_{m_2}}F_{m_1})=\{u\in D(F_{m_1}): F_{m_1}\in L^p(\R^n)\} = D(F_{m_1})$ . 
%\footnote{\commentCK{Strictly speaking, it suffices also to observe that $D(F_{m_1})\subset D(\ol{F_{m_2}}F_{m_1})$.}}

Let $u\in D(\overline{F_{m_1}})$ and choose a sequence $(u_j)_j\subset D(F_{m_1})$ with $u_j\to u$ and $F_{m_1}u_j\to \ol{F_{m_1}}u$ in $L^p(\R^n)$. Since $\overline{F_{m_2}}F_{m_1}u_j=F_{m_1m_2}u_j$ for all $j\in \N$ by~\eqref{aux:algebra}, and because $\overline{F_{m_2}}$ is an $L^p$-Fourier multiplier operator, one deduces that 
\begin{align*}
F_{m_2m_1}u_j=\overline{F_{m_2}}F_{m_1}u_j\to \overline{F_{m_2m_1}}u \quad\text{in $L^p(\R^n)$,}
\end{align*} 
which yields $u\in D(\ol{F_{m_2m_1}})$.
\end{proof}

%One could be tempted to think that the above equality carries over to the operators \(\ol{F_{m_2}}\, \ol{F_{m_1}}\) and \(\ol{F_{m_2m_1}}\) if \(F_{m_1}\) is closable, but this is not true in general, this is, Lemma~\ref{lem:mults_algebra} does not automatically imply that \(D(\ol{F_{m_2}}\, \ol{F_{m_1}})\subset D(\ol{F_{m_2m_1}})\). 
We conclude this section with a result that, in the cases we will be interested in, reduces proving an identity of the type \(\ol{F_{m_2}}\, \ol{F_{m_1}}=\ol{F_{m_2m_1}}\) of operators in \(L^p(\R^n)\) to showing the equality \(D(\ol{F_{m_2}}\, \ol{F_{m_1}})= D(\ol{F_{m_2m_1}})\) of their domains. 

\begin{lemma}\label{lem:closed_ops_equality}
    Let $m_1,m_2 : \R^n \to \R$ be measurable functions with polynomial growth, consider the operators $F_{m_1}, {F_{m_2}}$ in $L^p(\R^n)$, and further assume that \(F_{m_1},F_{m_2}\) are closable and that \(\mathcal{F}(m_2\mathcal{F}^{-1}\varphi)\), \(\mathcal{F}(m_2m_1\mathcal{F}^{-1}\varphi)\) lie in \(L^{p'}(\R^n)\) for any \(\varphi\in \Scal(\R^n)\). Then, if
\begin{equation*}
 D(\ol{F_{m_2}}\, \ol{F_{m_1}})= D(\ol{F_{m_2m_1}}),
\end{equation*}
it follows that
\begin{equation*}
    \ol{F_{m_2}}\, \ol{F_{m_1}}
    = \ol{F_{m_2m_1}}
\end{equation*}
as operators in \(L^p(\R^n)\).
\end{lemma}
\begin{proof} 
    Let \(u\in  D(\ol{F_{m_2}}\, \ol{F_{m_1}})= D(\ol{F_{m_2m_1}})\), and take \((u_j)_j, (v_j)_j\subset \mathcal{S}(\R^n)\) such that 
\begin{align*}
    u_j\to u, 
    \quad F_{m_1} u_j\to \ol{F_{m_1}} u , 
    \quad v_j\to\ol{F_{m_1}} u, 
    \quad F_{m_2} v_j\to \ol{F_{m_2}}\,\ol{F_{m_1}}  u 
    \quad \text{in } L^p(\R^n).
\end{align*} 
    We compute
    \begin{align*}
    \langle \ol{F_{m_2}}\, \ol{F_{m_1}} u,\varphi\rangle
    &= \lim_{j\to \infty} \langle F_{m_2}v_j,\varphi\rangle
    =
    \lim_{j\to \infty} 
        \langle {v_j},\mathcal{F}(m_2\mathcal{F}^{-1}\varphi)\rangle\\
       & = 
        \lim_{j\to \infty}
        \langle F_{m_1} u_j,\mathcal{F}(m_2\mathcal{F}^{-1}\varphi)\rangle
        = \lim_{j\to \infty}
        \langle  u_j,\mathcal{F}(m_2m_1\mathcal{F}^{-1}\varphi)\rangle\\
        & = \langle u, \mathcal{F}(m_2m_1\mathcal{F}^{-1}\varphi)\rangle 
        = \langle F_{m_2m_1}u, \varphi \rangle,
        \quad \varphi\in \Scal(\R^n),
    \end{align*}
    where we have used in particular that \(\mathcal{F}(m_2\mathcal{F}^{-1}\varphi), \mathcal{F}(m_2m_1\mathcal{F}^{-1}\varphi)\in L^{p'}(\R^n)\) by assumption. 
    %\commentCK{Do we also use Plancherel-type arguments above or should these identities be considered in the sense of distributions, which is not clear under the current assumptions on $m_1, m_2$...? Discuss with Leon. }\comment{We need that \(\mathcal{F}(m_2\mathcal{F}^{-1}\varphi)\) and \(\mathcal{F}(m_2m_1\mathcal{F}^{-1}\varphi)\) lie in \(L^{p'}(\R^n)\). We need to find the correct additional assumptions to add to the statement. This should also resolve the issue with the Plancherel-type argument that you mentioned.} due to the polynomial growth of \(m_2\).
    This shows that \(\ol{F_{m_2}}\, \ol{F_{m_1}}  = \ol{F_{m_2m_1}}\) as operators in \(L^p(\R^n)\).
\end{proof}

    \subsection{Further hypotheses on the kernel function}\label{sec:hypotheses}
    The properties of the function \(\widehat{Q_\rho}\) depend strongly on those of the kernel \(\rho\). We assume in the following that \(\rho\in L^1(\R^n)\). 

    It is straightforward, using the coarea formula, to rewrite the expression of \(\widehat{Q_\rho}\) in \eqref{eq:qrho_hat} as (see \cite[Lemma~4.3]{BMCS25}) 
    \begin{equation*}
        \widehat{{Q}_\rho}(\xi)
        = \frac{1}{\pi|\xi|}\int_{\sph_+^{n-1}}
        \int_{0}^\infty r^{n-2}\ol{\rho}(r)\sin(2\pi|\xi|r z_1)\dl r\dl \mathcal{H}^{n-1}(z),
         \qquad \xi \in \R^n\setminus \{0\}.
    \end{equation*}
    Basically by a variant of Pólya's criterion for the non-negativity of oscillatory integrals of \(\sin\)-type (see \cite{Po18}), it was shown in \cite[Lemma~4.3]{BMCS25} that \(\widehat{Q_\rho}>0\) if \(r\mapsto r^{n-2}\ol{\rho}(r)\) is strictly decreasing (actually, in \cite[Lemma~4.3]{BMCS25} it is only assumed that the function \(r\mapsto r^{n-2}\ol{\rho}(r)\) is decreasing, not \textit{strictly} decreasing; but this is compensated there by the additional doubling property found below in assumption (H1)). In \cite{BMCS25}, the positivity of \(\widehat{Q_\rho}\) is then used to formally invert the nonlocal derivative in Fourier space by introducing the symbol \(\widehat{Q_\rho}^{-1}\). 
    Additional conditions on \(\rho\) are then imposed to derive Poincaré inequalities, compact embeddings and a nonlocal version of the fundamental theorem of calculus. 
    Most of these conditions prove equally relevant to our purposes. We therefore recall those that we shall use, and comment on their implications afterwards.
    %Most~\footnote{What about (H4)...maybe weaker formulation is ok, maybe even these...} of these conditions prove to be pivotal for our purposes as well, we therefore recollect them here and comment on their implications afterwards.\commentCK{Find good phrasing to remove footnote.}
   % \comment{Just change to a semicolon in the middle? Here are also two other options:}
   % \comment{Most of these conditions prove equally relevant to our purposes. We therefore recall those that we shall use, and comment on their implications afterwards.}
  %  \comment{Since most of these conditions will play a pivotal role in our analysis as well, we recall them here - restricting ourselves to those relevant for our purposes - and comment on their implications afterwards.}
Let \({\delta_0}>0\) such that \(\inf_{B_{\delta_0}}\rho>0\) as in (H0).
    \begin{enumerate}
        \item[(H1)]
            The function \(f_\rho:(0,\infty)\to \R, r\mapsto r^{n-2}\ol{\rho}(r)\) is decreasing, and there exists \(\mu\in (0,1)\) such that \(\mu f_\rho(r/2)\ge f_\rho(r)\) for \(r\in (0,{\delta_0})\).
        \item[(H2)]
            The function \(f_\rho\) is smooth in \((0,\infty)\), and for \(r\in (0,{\delta_0})\) it holds
            \begin{equation*}
                -C\frac{d}{dr} f_\rho(r)\ge \frac{f_\rho(r)}{r}
                \quad\text{and}\quad
                \Bigg|\frac{d^k}{dr^k} f_\rho(r)\Bigg| \le C_k \frac{f_\rho(r)}{r^k},
                \qquad  k\in \N.
            \end{equation*}
        \item[(H3)]
            There exists \(\sigma\in (0,1)\) such that the function \(g_\rho:(0,\infty)\to \R, r\mapsto r^{n+\sigma-1}\ol{\rho}(r)\) is almost decreasing on \((0,{\delta_0})\).
    \end{enumerate}
    The expression \textit{almost decreasing} in (H3) means that \(g_\rho(r_1)\ge Cg_\rho(r_2)\) for all \(r_1\le r_2\in (0,{\delta_0})\) and some fixed constant \(C>0\). Qualitatively speaking, (H3) encodes that \(\rho\) lies around the origin above the kernel \(|\cdot|^{-(n+\sigma-1)}\) (which gives rise to the Riesz fractional gradient of order \(\sigma\)), and thus, has at least a singularity of the corresponding order in the origin. 
    We also introduce as a weaker replacement of condition (H3)
    \begin{equation}\tag{H3'}
			\liminf_{t\downarrow 0} t^{n-1}\overline{\rho}(t)>0.
	\end{equation}
    For more details on the different hypotheses we refer the reader to \cite{BMCS25}.
%\commentC{Here is where we should cite Felix' paper}

    Hypothesis~(H1) ensures not only that \(\widehat{{Q}_\rho}\) is positive, as remarked above, but it also guarantees the validity of the estimate (see \cite[Lemma~4.3]{BMCS25})
    \begin{equation*}\label{eq:first_lower_est}
			\widehat{Q_\rho}(\xi)
			\ge C\frac{\ol{\rho}(1/\xi)}{|\xi|^n},
			\qquad \xi\in B_{2/{\delta_0}}^c.
	\end{equation*}
    Building on this estimate, as outlined at the start of the proof of \cite[Lemma~4.11]{BMCS25}, under the additional hypothesis~(H3') it holds 
    \begin{equation}\label{eq:lower_prelim}
			\widehat{Q_\rho}(\xi)
			\ge C|\xi|^{-1},
			\qquad  \xi\in B_{2/{\delta_0}}^c.
	\end{equation}
    Moreover, with hypothesis~(H3) one can strengthen this lower bound. In fact, under conditions~(H0), (H1) and (H3) one has
    (see in \cite[Eq.\ (5.1)]{BMCS25})
    \begin{equation}\label{eq:lest_qrho_pre}
		\widehat{Q_\rho}(\xi)
		\ge C|\xi|^{\sigma-1}
		\ge C|\xi|^{-1},
		\qquad  \xi\in B_{2/{\delta_0}}^c.
	\end{equation}
    Finally, assuming additionally that \(\rho\) has compact support, so that \(\widehat{Q_\rho}\) is  smooth by the Paley--Wiener theorem, and imposing hypothesis~(H2), allows to deduce upper bounds for the derivatives \(\p^\alp \widehat{Q_\rho}\) (see \cite[Lemma~4.10]{BMCS25}). In particular, if (H0), (H1), (H2) and (H3') are satisfied, then
	\begin{equation}\label{eq:uest_qrho_alp_pre}
		|\p^\alp \widehat{Q_\rho}(\xi)|
		\le C_\alp|\xi|^{-|\alp|}\widehat{Q_\rho}(\xi),
		\qquad  \xi\in B_{2/{\delta_0}}^c.
	\end{equation}
    If \(\rho\) has compact support, then, by the smoothness of \(\widehat{Q_\rho}\), all estimates above that are valid on \(B_{2/{\delta_0}}^c\) also hold true on \(B_1^c\) (up to adapting the constants in the inequalities). For this reason, we will drop the specific \({\delta_0}>0\) in the estimates from here on and work throughout the paper exclusively on \(B_1\) and its complement \(B_1^c\).

    %%%%%%%%%%%%%%%%%%%%%%%%%%%%%%%%%%%%%%%

%***30/7/2026

%\commentCK{Consider removing the following headline}
%    \subsection{On the condition of compact support}\label{sec:compact_supp}
    As we have seen above, assuming that \(\rho\) has compact support simplifies the analysis in some aspects. Not only does this imply that \(\rho\in L^1(\R^n)\) and \(Q_\rho\in L^1(\R^n)\), since \(\rho\in L^1_{\textrm{loc}}(\R^n)\) according to hypothesis~(H0), but more importantly, it also guarantees that \(Q_\rho\) again has compact support and \(\widehat{Q_\rho}\) is smooth. 
    %As it turns out, the assumption of compact support is not too restrictive, since for any two kernels \(\rho_1,\rho_2\) satisfying (H0) and \((\rho_1-\rho_2)/|\cdot|\in L^1(\R^n)\) it holds (see \cite[Proposition~3.10]{BMCS25})
    %\begin{equation*}
    %    H_0^{\rho_1,p}(\Om)=H_0^{\rho_2,p}(\Om)
    %\end{equation*}
    %for all open \(\Om\subset \R^n\) and all \(p\in [1,\infty]\). In particular, if \(\rho_1,\rho_2\) agree locally around the origin, then the associated nonlocal Sobolev spaces coincide thanks to (H0). Thus, we can always replace \(\rho\) subject to (H0) by \(\chi\rho\) with a radial \(\chi\in L^\infty(\R^n)\) such that \(\chi|_{B_r}\equiv 1\) for some \(r>0\); in this case, \(H_0^{\rho,p}(\Om)=H_0^{\chi\rho,p}(\Om)\) (see also \cite[Example~3.11]{BMCS25}). In addition, by their local nature, also hypothesis (H1)--(H3) and (H3') are left unaltered under the multiplication with \(\chi\) if, in addition, in view of (H1), \(\chi\) is decreasing. For completeness, we also mention that \cite[Proposition~3.10]{BMCS25} proves that if one of the spaces \(H_0^{\rho_1,p}(\Om)\), \(i=1,2\), admits a Poincaré inequality, then so does the other.
As it turns out, the assumption of compact support is not too restrictive, since for any kernel \(\rho\) 
satisfying (H0) and any radially decreasing \(\chi\in L^\infty(\R^n)\) such that \(\chi|_{B_r}\equiv 1\) for some \(r>0\), we have \(H_0^{\rho,p}(\Om)=H_0^{\chi\rho,p}(\Om)\) (see \cite[Example~3.11]{BMCS25}). In addition, if $\rho$ satisfies any of the conditions (H1)--(H3) and (H3'), then so does $\chi \rho$. For completeness, we also mention that \cite[Proposition~3.10]{BMCS25} proves that if \(H_0^{\rho,p}(\Om)\), admits a Poincaré inequality, then so does \(H_0^{\chi\rho,p}(\Om)\).

    %%%%%%%%%%%%%%%%%%%%%%%%%%%%%%%%%%%%%%%%%
	\subsection{Basic Fourier multiplier results involving the symbol $\lambda_\rho$}\label{sec:fourier_mult}
    %Some estimates and a Fourier multiplier result
 We will start this section by presenting important estimates involving the Fourier symbol \(\lambda_\rho\) (see Section~\ref{sec:fourier}) of the nonlocal derivative \(D^\rho_0\). These estimates build upon the inequalities in Section~\ref{sec:hypotheses}, and their proof, only consisting of careful computations, is deferred  to Appendix~\ref{sec:appendix_ests}.
  %  Building on the estimates in Section~\ref{sec:hypotheses}, we derive the following key inequalities, \commentCK{maybe something like...estimates on the Fourier symbols of the nonlocal gradient and the $\rho$ Laplacian...or...estimates involving $\lambda_\rho$... or something else for just a bit more context.} whose proof is deferred to Appendix~\ref{sec:appendix_ests}. 
	
	\begin{lemma}\label{lem:estimates}
	 Assume that \(\rho\) has compact support and satisfies assumptions (H0)--(H2) and (H3'). Then, for all \(\xi\in B_1^c\) and all multi-indices \(\alp\in \N^n_0\),
		\begin{align}
			&|\p^\alp \lambda_\rho(\xi)|,|\p^\alp |\lambda_\rho(\xi)||
			\le C_\alp |\lambda_\rho(\xi)||\xi|^{-|\alp|},\label{eq:second_estimate}\\
			&|\p^\alp (|\lambda_\rho(\xi)|)^{\gamma}|, |\p^\alp (1+|\lambda_\rho(\xi)|)^{\gamma}|
			\le  C_\alp |\lambda_\rho(\xi)|^{\gamma}|\xi|^{-|\alp|},\label{eq:third_estimate}\\
			&|\p^\alp (1+|\lambda_\rho(\xi)|^2)|
			\le C_\alp |\lambda_\rho(\xi)|^2|\xi|^{-|\alp|},\label{eq:fourth_estimate}\\
			&|\p^\alp (1+|\lambda_\rho(\xi)|^2)^{\gamma}|
			\le C_\alp  |\lambda_\rho(\xi)|^{2\gamma}|\xi|^{-|\alp|},
			\qquad  \gamma \in \R,\label{eq:fifth_estimate}
			%\le C_\alp |\xi|^{-\sigma-|\alp|}
		\end{align}
		 with constants \(C_\alp>0\).
	\end{lemma}
	%\begin{remark}
	%	To ensure the inequality in \eqref{eq:lest_qrho}, and subsequently \eqref{eq:uest_qrho_alp}, one can also assume that 
	%	\begin{equation*}
	%		\liminf_{t\downarrow 0} t^{n-1}\overline{\rho}(t)>0
	%	\end{equation*}
	%	instead of \comment{hypothesis (H3) in \cite{BMCS25}}. Indeed, as outlined at the start of the proof of \cite[Lemma~4.11]{BMCS25}, as a consequence of \cite[Lemma~4.3 and Eq.~(4.7)]{BMCS25}, it holds
	%	\begin{equation*}
	%			\widehat{Q_\rho}(\xi)
	%		\ge C|\xi|^{-1}
	%		\qquad \text{ for } |\xi|\ge 1.
	%	\end{equation*}
	%\end{remark}

     % \comment{
      It is important for our purposes to note that the operators $D^{\rho}_0$  and $\diver^{\rho}_0$ map Schwartz functions again to Schwartz functions. 
      %$\Scal(\R^n)$ into $\Scal (\R^n;\R^n)$ and $\Scal(\R^n; \R^n)$ into $\Scal (\R^n)$, respectively. 
      Indeed, in view of~\eqref{eq:nlgrad_symbol}, it suffices to prove that \(\lambda_\rho\) is smooth and all its derivatives have polynomial growth, as  Lemma~\ref{lem:lambda_FintoF_pre} below shows. 
      %To be precise, as pointed out after \eqref{eq:pw_equals_dual} and \eqref{eq:nlgrad_symbol}, for \(\varphi\in \mathcal{S}(\R^n)\) and \(\Phi\in \mathcal{S}(\R^n;\R^n)\), the expressions \(D^\rho_0 \varphi\) and \(\diver^\rho_0 \Phi\) are well defined and their Fourier symbol is given in both cases by \(\lambda_\rho\). Hence, their Fourier inverses then lie indeed in the respective spaces claimed above.
      %}
    \begin{lemma}\label{lem:lambda_FintoF_pre}
    %\comment{Alternative for later use in Section~\ref{sec:char_laplacian}: 
    Assume that \(\rho\) has compact support. If \(\rho\) satisfies assumptions (H0), then $\lambda_\rho\in C^\infty(\R^n)$ has linear growth, and if it additionally fulfills (H1)--(H2) and (H3'), then 
    \begin{equation*}
        |\p^\alp \lambda_\rho(\xi)|
        \le C_{\alp} |\xi|^{1-|\alp|}
    \end{equation*}
     for all \(\xi\in \R^n\) and multi-indices \(\alp\in \N^n_0\), for some suitable constants \(C_\alp>0\). % which depend also on \(\norm{\rho}{L^1(\R^n)}\).
   % Assume that \(\rho\) has compact support and satisfies assumptions (H0)--(H2) and (H3').  Then $\lambda_\rho\in C^\infty(\R^n)$ has linear growth and   \begin{equation*}
    %    |\p^\alp \lambda_\rho(\xi)|
    %    \le C_{\alp} |\xi|^{1-|\alp|}
   % \end{equation*}
  %   for all \(\xi\in \R^n\) and multi-indices \(\alp\in \N^n_0\), for some suitable constants \(C_\alp>0\), which depend also on \(\norm{\rho}{L^1(\R^n)}\).
\end{lemma}
\begin{proof}
    The property that \(\lambda_\rho\in C^\infty(\R^n)\) is a direct consequence of the fact that \(\rho\) has compact support. 
    Let $\xi\in B_1^c$. According to \cite[Lemma~2.4]{BMCS25} it holds
    \begin{equation*}
        |\lambda_\rho(\xi)|
        \le \int_{\R^n} \rho(x)\min\{|x|^{-1}, 2\pi|\xi|\}\dl x,
    \end{equation*}
    which, using that \(\rho\in L^1(\R^n)\), can be further estimated by
    \begin{equation*}
        |\lambda_\rho(\xi)|
        \le 2\pi \norm{\rho}{L^1(\R^n)}|\xi|,
        \qquad \xi\in B_1^c.
    \end{equation*}
    Exploiting also  (H1)--(H2) and (H3'), it follows then immediately from Lemma~\ref{lem:estimates} that \(|\p^\alp \lambda_\rho(\xi)|\le C|\xi|^{1-|\alp|}\) for all multi-indices \(\alp\in \N^n_0\), and for some constants \(C_\alp>0\) depending on \(\norm{\rho}{L^1(\R^n)}\). 
    We conclude by noting that the inequality extends to all \(\xi\in \R^n\), since \(\lambda_\rho\in C^\infty(\R^n)\).
%
%  \commentCK{Make very condensed version and include intro into sentence. Standard!}  As pointed out after \eqref{eq:pw_equals_dual} and \eqref{eq:nlgrad_symbol}, for \(\varphi\in \mathcal{S}(\R^n)\) and \(\Phi\in \mathcal{S}(\R^n;\R^n)\) the expressions \(D^\rho_0 \varphi\) and \(\diver^\rho_0 \Phi\) are well defined and their Fourier symbol is given in both cases by \(\lambda_\rho\), i.e.,
 %   \begin{equation*}
%        \mathcal{F}(D^\rho_0 \varphi)
%        = \lambda_\rho\hat{\varphi}
%        \in \mathcal{S}(\R^n;\R^n),
 %       \quad 
%        \mathcal{F}(\diver^\rho_0 \Phi)
 %       = \lambda_\rho\cdot \hat{\Phi}
%        \in \mathcal{S}(\R^n).
 %   \end{equation*}
% We note that $\lambda_\rho\hat{\varphi}\in \mathcal{S}(\R^n;\R^n)$ and $\lambda_\rho\cdot \hat{\Phi} \in \mathcal{S}(\R^n)$ since \(\lambda_\rho\) is smooth and all its derivatives have polynomial growth. 
%    Accordingly, the Fourier inverses of the expressions above are again elements of \(\mathcal{S}(\R^n;\R^n)\) and \(\mathcal{S}(\R^n)\), respectively. 
\end{proof}

	Lemma~\ref{lem:estimates} allows us to deduce the following Fourier multiplier result by means of the Mihlin-H\"ormander multiplier theorem. 
    %Recall that $\lambda_\rho(\xi) = 2\pi i \xi \widehat{Q}_{\rho}(\xi)$.
    We set from now on 
    \begin{equation*}
    \bla{\lambda_\rho}\coloneqq (1+|\lambda_\rho|^2)^{1/2}.
    %quad m_\rho\coloneqq \frac{\abs{\lambda_\rho}}{\bla{\lambda_\rho}}.
    \end{equation*}
	\begin{proposition}[Central Fourier multipliers]\label{prop:fourier_multipliers}
		Let \(p\in (1,\infty)\). Assume that \(\rho\) has compact support  and satisfies assumptions (H0)--(H2) and (H3'). Then, the symbols 
		\begin{equation*}\label{eq:fm_2}
		\textit{(i)}%\quad m_\rho,
        \quad \frac{\abs{\lambda_\rho}}{\bla{\lambda_\rho}},
		\qquad \textit{(ii)}\quad \frac{\bla{\lambda_\rho}}{1+\abs{\lambda_\rho}},
		\qquad \textit{(iii)}\quad \frac{1}{\bla{\lambda_\rho}},
      \qquad \textit{(iv)}\quad \frac{1}{\bla{\lambda_\rho}^2}
		\end{equation*}
		are \(L^p\)-Fourier multipliers with corresponding $L^p$-Fourier multiplier operators 
        \begin{equation*}
            \ol{F_{\abs{\lambda_\rho}/\bla{\lambda_\rho}}}, \,
            \ol{F_{\bla{\lambda_\rho}/(1+\abs{\lambda_\rho})}}, \,
            \ol{F_{1/\bla{\lambda_\rho}}} \text{ and }
            \ol{F_{1/\bla{\lambda_\rho}^2}}.
        \end{equation*}
	\end{proposition}
    %\commentCK{Can we exploit that squares of Fourier multipliers via H\"ormander's theorem are Fourier multipliers? Using the Leibniz rule... Then $|\lambda_\rho|^2/\bla{\lambda_\rho}^2$ would be multipliers and analogously for the square of \textit{(ii)}. Can we conclude from that that $D(\ol{F_{|\lambda_\rho|^2}}) = D(\ol{F_{\bla{\lambda_\rho}^2}}))$? This is basically also the idea behind the new addition $(iv)$.}
	\begin{proof}
		\textit{(i).}
			We define \( m_\rho:\R^n\to \R\) by 
		\begin{equation*}
		%f_1\coloneqq \abs{\lambda_\rho} , 
        %\qquad f_2\coloneqq \bla{\lambda_\rho} , \qquad 
        m_\rho\coloneqq \frac{\abs{\lambda_\rho}}{\bla{\lambda_\rho}}.
		\end{equation*}
        In order to apply the Mihlin-H\"ormander multiplier theorem, we verify subsequently that 
		\begin{equation}\label{eq:MH_requirements}
		|\p^\alp m_\rho(\xi)|\le C_\alp|\xi|^{-|\alp|}, \qquad \xi\in \R^n\setminus\{0\}
		\end{equation} 
		for all multi-indices \(\alp\in \N^n_0\). We note that by their smoothness, the functions \(\p^\alp f_i\), \(i=1,2\), are well defined on \(B_1\setminus\{0\}\), and it is readily seen that \(m_\rho\) and all its derivatives are bounded there. So, it will be sufficient to establish inequality \eqref{eq:MH_requirements} for \(\xi\in B_1^c\). To that end we first note that, actually for all \(\xi\in \R^n\),
		\begin{equation}\label{eq:final_est0}
		|m_\rho(\xi)|
		\le 1.
		\end{equation}
		Choose now some \(\alp\in \N^n_0\setminus \{0\}\). Then, with Leibniz' rule, it holds
		%\begin{align}
		%\begin{split}\label{eq:initial_est}
		%|\p^\alp m_\rho(\xi)|
		%&\le C_\alp\Big[\sum_{\beta\le \alp, 1\le |\beta|\le |\alp|-1}
		%|\p^\beta f_1(\xi)||\p^{\alp-\beta} f_2^{-1}(\xi)|\\
		%&\qquad\quad + |\p^\alp f_1(\xi)|| f_2^{-1}(\xi)|
		%+|f_1(\xi)||\p^{\alp} f_2^{-1}(\xi)|\Big].
		%\end{split}
		%\end{align}
        \begin{align}
		\begin{split}\label{eq:initial_est}
		|\p^\alp m_\rho(\xi)|
		&\le C_\alp\sum_{\beta\le \alp}
		|\p^\beta \abs{\lambda_\rho(\xi)}||\p^{\alp-\beta} \bla{\lambda_\rho(\xi)}|, \qquad \xi \in B_1^c.
		\end{split}
		\end{align}
        Resorting to \eqref{eq:second_estimate} and \eqref{eq:fifth_estimate} in Lemma~\ref{lem:estimates}, we conclude that the right-hand side is bounded by
\begin{equation}
\label{eq:final_est2}
		%\sum_{\beta\le \alp}
		%|\p^\beta f_1(\xi)||\p^{\alp-\beta} f_2^{-1}(\xi)| 
        %\le 
        C_\alp\sum_{\beta\le \alp} |\lambda_\rho(\xi)||\xi|^{-|\beta|}
		|\lambda_\rho(\xi)|^{-1}|\xi|^{-|\alp|+|\beta|} = C_\alp |\xi|^{-|\alp|}, 
        \qquad \xi \in B_1^c.
\end{equation}
%for $\xi\in B_1^c$.
	Combining estimate \eqref{eq:final_est0} with \eqref{eq:initial_est} and \eqref{eq:final_est2}, proves \eqref{eq:MH_requirements}. Consequently, thanks to the Mihlin-H\"ormander multiplier theorem, %(see \cite[Theorem~6.2.7]{Gr14}),
\(m\) is an \(L^p\)-Fourier multiplier, for any \(p\in (1,\infty)\).\smallskip
%In particular, by the density of the Schwartz functions \(\mathcal{S}(\R^n)\) in \(L^p(\R^n)\), \(m(\xi)\) can be uniquely extended to a bounded operator in \(L^p(\R^n)\). 
		
		\textit{(ii).}
		The proof is obtained similarly to case \textit{(i)}. In fact, once more with Lemma~\ref{lem:estimates}, now estimates \eqref{eq:third_estimate} and \eqref{eq:fourth_estimate}, we can compute again
        \begin{align*}
		\begin{split}
		&\sum_{\beta\le \alp}
		|\p^\beta \bla{\lambda_\rho}||\p^{\alp-\beta} (1+\abs{\lambda_\rho(\xi)})^{-1}|\\
		&\quad \le C_\alp\sum_{\beta\le \alp} |\lambda_\rho(\xi)||\xi|^{-|\beta|}
		|\lambda_\rho(\xi)|^{-1}|\xi|^{-|\alp|+|\beta|}\\
		&\quad = C_\alp |\xi|^{-|\alp|},
        \qquad  \xi\in B_1^c.
		\end{split}
		\end{align*}
        %, with the roles of \(f_1\) and \(f_2\) exchanged. To that end, note that due to \eqref{eq:lest_qrho} it holds \(|f_2|\le C_0 \) for all \(|\xi|\ge 1\). 
		
		\textit{(iii) and (iv).} 
		It follows from \eqref{eq:lower_prelim} that \(|\lambda_\rho|^{-1}\) is bounded, and we get from Lemma~\ref{lem:estimates} immediately that 
		\begin{align*}
			&|\p^\alp \bla{\lambda_\rho}^{-1}|
			\le C_\alp |\lambda_\rho(\xi)|^{-1}|\xi|^{-|\alp|}
			\le C_\alp |\xi|^{-|\alp|},\\
            &|\p^\alp \bla{\lambda_\rho}^{-2}|
			\le C_\alp |\lambda_\rho(\xi)|^{-2}|\xi|^{-|\alp|}
			\le C_\alp |\xi|^{-|\alp|},
            \qquad  \xi\in B_1^c
		\end{align*}
    for all \(\alp\in \N^n_0\). Thus, one can again apply the Mihlin-H\"ormander multiplier theorem in both cases. 
	\end{proof}

    Our final result in this section collects four specific multiplier operators involving the
symbol $\lambda_\rho$. Their closures will turn out to be central for the
later analysis, since they correspond to the $\rho$-Laplacian and its square
roots. 
\begin{lemma}\label{lem:closability_sr}
Let \(p\in (1,\infty)\). Assume that \(\rho\) has compact support and satisfies (H0)--(H2) and (H3'). Then the following holds:
\begin{enumerate}[label=(\roman*)]
    \item 
   %  Then the unbounded operator Fm in Lp(Rn) maps from D(Fm) = S(Rn) into S(Rn) and is closable.
        The unbounded operators $F_{\bla{\lambda_\rho}}$, $F_{\bla{\lambda_\rho}^2} $, $F_{1/\bla{\lambda_\rho}^2} $ and $F_{|\lambda_\rho|^2}$ in \(L^p(\R^n)\) map their domain $\Scal(\R^n)$ into itself and are closable. 
        %and the identities
      %  \begin{equation}\label{eq:im_scal}
      %      F_{\bla{\lambda_\rho}} u=\Fcal^{-1}(\bla{\lambda_\rho}\hat{u}), \quad
       %     F_{|\lambda_\rho|^2}u=\Fcal^{-1}(\abs{\lambda_\rho}^2\hat{u}) , \qquad u\in \Scal(\R^n)
     %   \end{equation}
      %  hold a.e., where the the pointwise Fourier inverses \(\Fcal^{-1}(\bla{\lambda_\rho}\hat{u})\) and \(\Fcal^{-1}(\abs{\lambda_\rho}^2\hat{u})\) are again elements of \(\Scal(\R^n)\).
    \item 
        The unbounded operators $F_{\abs{\lambda_{\rho}}}$ and $F_{1+\abs{\lambda_{\rho}}}$ in $L^p(\R^n)$ are closable with domains \(\Scal(\R^n)\). %\commentCK{Homework Carolin: Check whether the first part of the statement is really needed here, cf. Remark~\ref{rem:fi_conincidence}. }
     \end{enumerate}
        %Moreover, the identity
      %  \begin{equation*}
      %      F_{\abs{\lambda_{\rho}}}u
       %     = \Fcal^{-1}(\abs{\lambda_\rho} \hat{u}),
       %     \qquad u\in \Scal(\R^n),
      %  \end{equation*}
      %   holds a.e., 
      %  where the right-hand side denotes the pointwise Fourier inverse, and 
    Moreover, it holds
    \begin{alignat}{3}
            &\norm{\ol{F_{\bla{\lambda_{\rho}}}}u}{L^p(\R^n)}
             &&\le \frac{1}{2} \norm{u}{L^p(\R^n)} + \frac{1}{2}\norm{\ol{F_{\bla{\lambda_\rho}^2}} u}{L^p(\R^n)},
             \qquad &&u\in D(\ol{F_{\bla{\lambda_\rho}}})\cap D(\ol{F_{\bla{\lambda_\rho}^2}}), \label{eq:aux_lp_1}\\
           &\norm{F_{\abs{\lambda_{\rho}}}u}{L^p(\R^n)}
            && \le \frac{1}{2} \norm{u}{L^p(\R^n)} + \frac{1}{2}\norm{F_{|\lambda_\rho|^2} u}{L^p(\R^n)},
             \qquad &&u\in \Scal(\R^n). \label{eq:aux_lp_2}
    \end{alignat}

%\commentCK{Mention $L^p$! moreover, with same domain. Also, maybe collect in this lemma also the statements on $F_{\bla{\lambda_\rho}}$ and $F_{|\lambda_\rho|^2}$, which are both closable with domain $\Scal(\R^n)$, the latter even mapping into $\Scal(\R^n)$ etc. to make these easy to reference.} \commentCK{Maybe prove here that $D(F_{|\lambda_\rho|}) =\Scal(\R^n)$ instead of later... Can we conclude this from Lemma~\ref{lem:mults_algebra} }
%\comment{I think we could use Lemma~\ref{lem:mults_algebra} as you suggest but then we have to invoke Proposition~\ref{prop:fourier_multipliers} anyway to guarantee the closability of the operator \(F_{m_2}\) in the statement of the Lemma.}
\end{lemma}
\begin{proof}
%Let $T_1$ be the unbounded operator in $L^p (\R^n)$ with domain $D(T_1) := \Scal$ and expression
%\[
% T_1 u \coloneqq  \Fcal^{-1} ( \bla{\lambda_\rho} \hat{u}) .
%\]
%\textit{Proof of (i).}
\textit{Step~1: proof of (i) and \eqref{eq:aux_lp_1}.} Part \textit{(i)} is an immediate consequence of Lemma~\ref{le:Fmclosable} in combination with Lemma~\ref{lem:lambda_FintoF_pre}. 
%the unbounded operators $F_{\bla{\lambda_\rho}}$ and $F_{|\lambda_\rho|^2}$ have domain equal to \(\mathcal{S}(\R^n)\) and are closable. The remaining claims follow also directly from the same results. 
Furthermore, with the aim of proving \eqref{eq:aux_lp_1}, let \(u\in  D(\ol{F_{\bla{\lambda_\rho}}})\cap D(\ol{F_{\bla{\lambda_\rho}^2}})\) and choose a sequence, existing by definition of the closure, \((u_j)_j\subset \Scal(\R^n)\) such that \(u_j\to u\) and \(F_{\bla{\lambda_\rho}^2} u_j \to \ol{F_{\bla{\lambda_\rho}^2}} u\). Since \(\ol{F_{1/\bla{\lambda_\rho}}}\) is an \(L^p\)-Fourier multiplier operator according to Proposition~\ref{prop:fourier_multipliers}, we infer with Lemma~\ref{lem:domain_mulitplier} that \(F_{\bla{\lambda_\rho}} u_j \to \ol{F_{\bla{\lambda_\rho}}} u\). Then, as a consequence of Young's inequality, one has
\begin{align*}
    \begin{split}
        &|\langle\ol{ F_{\bla{\lambda_{\rho}}} u}, \varphi\rangle|
        = \lim_{j\to \infty} |\langle F_{\bla{\lambda_{\rho}}} u_j, \varphi\rangle|
        =
      \lim_{j\to \infty}|\langle \bla{\lambda_{\rho}} \widehat{u_j}, \mathcal{F}^{-1}\varphi\rangle|\\
    &\le \lim_{j\to \infty}\left(\frac{1}{2}|  \langle \widehat{u_j}, \mathcal{F}^{-1}\varphi\rangle|
    + \frac{1}{2} | \langle \bla{\lambda_{\rho}}^2 \widehat{u_j}, \mathcal{F}^{-1}\varphi\rangle|\right)
    = \frac{1}{2}|  \langle u, \varphi\rangle|
    + \frac{1}{2} | \langle \ol{F_{\bla{\lambda_{\rho}}^2}} {u}, \varphi\rangle|,
    \qquad \varphi \in C^{\infty}_c (\R^n).
    \end{split}
\end{align*}
Taking the supremum over all \(\varphi\in C_c^\infty(\R^n)\) such that \(\norm{\varphi}{L^{p'}(\R^n)}\le 1\), we obtain estimate \eqref{eq:aux_lp_1}.
\smallskip

\textit{Step~2: proof of (ii) and and \eqref{eq:aux_lp_2}.}
%\textit{Proof of (ii).}
We observe upfront that 
%\comment{maybe we include here also the second property; this should clarity the question about $F_{\abs{\lambda_{\rho}}} u_j \to v$ and \(F_{1+\abs{\lambda_{\rho}}} u_j \to v\) below}
\begin{equation}\label{eq:dom_coincidence}
    D(F_{|\lambda_\rho|}) = D(F_{1+|\lambda_\rho|})
    \quad\text{and}\quad
    F_{1+\abs{\lambda_{\rho}}} u=u +F_{\abs{\lambda_{\rho}}} u,
    \qquad u\in D(F_{|\lambda_\rho|}),
\end{equation}
following the linearity of the Fourier transform.
To establish the claim that $F_{\abs{\lambda_{\rho}}}$ and $F_{1+\abs{\lambda_{\rho}}}$ have domain equal to \(\mathcal{S}(\R^n)\), it thus suffices to show that \(D ({F_{\abs{\lambda_\rho}}})=\Scal(\R^n)\), which reduces to proving that \(\Scal(\R^n)\subset D({F_{\abs{\lambda_\rho}}})\), since the inverse inclusion holds by definition.

To that end, we pick \(u\in \mathcal{S}(\R^n)\) and point out that the product \(\abs{\lambda_\rho} \hat{u}\) lies in \(L^2(\R^n)\), since \(\lambda_\rho\) has polynomial growth thanks to Lemma~\ref{lem:lambda_FintoF_pre}, and, hence, the Fourier inverse is well defined  and \(\mathcal{F}^{-1}(\abs{\lambda_\rho} \hat{u})\in L^2(\R^n)\). 
%Without mentioning this again, similar arguments assure that the Fourier inverse of related expressions below will be well defined. 
At the same time, it holds, similar to \textit{Step~1}, with Young's inequality and part~\textit{(i)} that
%, for any $\varphi \in C^{\infty}_c (\R^n)$, 
\begin{align*}
    \begin{split}
        |\langle \mathcal{F}^{-1}(\abs{\lambda_{\rho}} \hat{u}), \varphi\rangle|
    =
     |\langle \abs{\lambda_{\rho}} \hat{u}, \mathcal{F}^{-1}\varphi\rangle|
    &\le 
    %\frac{1}{2}|  \langle \hat{u}, \mathcal{F}^{-1}\varphi\rangle|
    %+ \frac{1}{2} | \langle \abs{\lambda_{\rho}}^2 \hat{u}, \mathcal{F}^{-1}\varphi\rangle|\\
    %&= 
    \frac{1}{2}|  \langle u, \varphi\rangle|
    + \frac{1}{2} | \langle F_{\abs{\lambda_{\rho}}^2} {u}, \varphi\rangle|,
    \qquad \varphi \in C^{\infty}_c (\R^n),
    \end{split}
\end{align*}
and taking again the supremum over all \(\varphi\in C_c^\infty(\R^n)\) such that \(\norm{\varphi}{L^{p'}(\R^n)}\le 1\), yields
\begin{equation}\label{eq:aux_lp_22}
    \norm{\mathcal{F}^{-1}(\abs{\lambda_{\rho}}\hat{u})}{L^p(\R^n)}
    \le \frac{1}{2} \norm{u}{L^p(\R^n)} + \frac{1}{2}\norm{F_{\abs{\lambda_{\rho}}^2} {u}}{L^p(\R^n)}.
\end{equation}
%where $\mathcal{F}^{-1}(\abs{\lambda_{\rho}}^2 \widehat{u}) =F_{|\lambda_\rho|^2}u\in L^p(\R^n)$ according to~\textit{(i)}.
This shows that \(\mathcal{F}^{-1}(\abs{\lambda_\rho} \hat{u})\in L^p(\R^n)\), and hence, $u\in D(F_{|\lambda_\rho|})$. Consequently, \(D (F_{\abs{\lambda_\rho}})= D(F_{1+|\lambda_\rho|})=\mathcal{S}(\R^n)\), and \eqref{eq:aux_lp_2} is a direct consequence of \eqref{eq:aux_lp_22}.

It remains to prove that $F_{\abs{\lambda_\rho}}$ and $F_{1+\abs{\lambda_\rho}}$
are both closable. 
%We give here the proof for $F_{\abs{\lambda_\rho}}$, noting that the argument for
%$F_{1+\abs{\lambda_\rho}}$ is analogous.
%It remains to prove that $F_{\abs{\lambda_{\rho}}}$ and $F_{1+\abs{\lambda_{\rho}}}$ are closable, we focus here on the proof for $F_{\abs{\lambda_{\rho}}}$, as the other proof is similar.
With \eqref{eq:dom_coincidence} in mind, choosing now a sequence $(u_j)_{j}\subset D(F_{\abs{\lambda_{\rho}}}) =\Scal(\R^n)$ such that $u_j \to 0$ in $L^p (\R^n)$, there exists $v \in L^p (\R^n)$ such that $F_{\abs{\lambda_{\rho}}} u_j \to v$ and \(F_{1+\abs{\lambda_{\rho}}} u_j \to v\) in $L^p (\R^n)$. %\commentCK{Should the second $v$ be $w$ or $\tilde{v}, etc.?$ in the first step and then argue that $w=v$! Insert argument for clarity!}
Proposition~\ref{prop:fourier_multipliers}~\textit{(ii)} implies the continuity of the \(L^p\)-Fourier multiplier operator \(\overline{F_{\bla{\lambda_\rho}/(1+|\lambda_\rho|)}}\), and so we find with Lemma~\ref{lem:mults_algebra} that 
\[
 F_{\bla{\lambda_\rho}}  u_j = \overline{F_{\bla{\lambda_\rho}/(1+|\lambda_\rho|)}} F_{1+|\lambda_\rho|} u_j \to \overline{F_{\bla{\lambda_\rho}/(1+|\lambda_\rho|)}} v \quad \text{in } L^p (\R^n).
\]
 As the operator $F_{\bla{\lambda_\rho}} $ is closable in \(L^p(\R^n)\) % \commentCK{Terminology: Closable in $L^p$ vs. operator $F_m$ in $L^p$ is closable... be consistent}
 according to part~\textit{(i)}, we conclude that the limit above is zero, thus, \(F_{\bla{\lambda_\rho}}  u_j\to 0\) in $L^p (\R^n)$.
Finally, applying Proposition~\ref{prop:fourier_multipliers}\,\textit{(i)} and again Lemma~\ref{lem:mults_algebra} yields
\begin{equation*}
    F_{\abs{\lambda_{\rho}}} u_j =  
    \overline{F_{|\lambda_\rho|/\bla{\lambda_\rho}}} F_{\bla{\lambda_\rho}} u_j \to 0 \quad \text{in } L^p (\R^n).
\end{equation*}
Hence, $v=0$ and, consequently, $F_{\abs{\lambda_{\rho}}}$ is closable. 
 Similarly, because $F_{1+\abs{\lambda_{\rho}}} u_j \to 0$ in $L^p (\R^n)$ as well, we have also proved that $F_{1+\abs{\lambda_{\rho}}}$ is closable. 
\end{proof}

\section{Fourier characterization}\label{sec:charac:fourier}

In this section, we establish the characterization of \(H^{\rho,p}(\R^n)\) via Fourier transform, and 
then conclude a similar result for \(H^{\rho,p}_0 (\Om)\) with $\Om \subset \R^n$ open.
Finally, the dual of \(H^{\rho,p}_0 (\Om)\) is also characterized. 

	\subsection{Fourier characterization of \(H^{\rho,p}(\R^n)\)}\label{sec:fourier_full_space}
In the following, we will devise a representation of \(H^{\rho,p}(\R^n)\) in terms of a Fourier characterization over the space \(\mathcal{S}'(\R^n)\). In order for this characterization to make sense, we rely, as outlined in Section~\ref{sec:fourier}, on the fact that \(\bla{\lambda_\rho}\) lies in \(C^\infty(\R^n)\) and %\color{red} itself and \color{black}\comment{I took out "itself" everywhere}
all its derivatives have polynomial growth thanks to Lemma~\ref{lem:lambda_FintoF_pre}.
    
	%We also recall that \(\lambda_\rho = 2\pi i\xi\widehat{Q_\rho}\).
	\begin{theorem}[Fourier characterization]\label{thm:alt_char}
		Let \(p\in (1,\infty)\) and assume that \(\rho\) has compact support and fulfills hypotheses (H0)--(H2) and (H3'). It holds
		\begin{align}
		\begin{split}\label{eq:fourier_identities}
		    H^{\rho,p}(\R^n)
		&= \{u\in \mathcal{S}'(\R^n):\, \mathcal{F}^{-1}(\bla{\lambda_\rho}\hat{u})\in L^p(\R^n)\}\\
        &=\{u\in L^p(\R^n):\, \mathcal{F}^{-1}(\bla{\lambda_\rho}\hat{u})\in L^p(\R^n)\}
		\end{split}
		\end{align}
		with 
		\begin{align}\label{eq:alt_char}
		c_p\norm{\mathcal{F}^{-1}(\bla{\lambda_\rho}\hat{u})}{L^{p}(\R^n)}
		\le \norm{u}{H^{\rho,p}(\R^n)}
		\le C_p\norm{\mathcal{F}^{-1}(\bla{\lambda_\rho}\hat{u})}{L^{p}(\R^n)}
		\end{align}
		for all \(u\in H^{\rho,p}(\R^n)\) and some constants \(c_p,C_p>0\). 
        
        Moreover, 
        \begin{equation}\label{eq:gen_fourier_symbol}
            D^\rho u
            =\mathcal{F}^{-1}(\lambda_\rho \hat u),
            \qquad  u\in H^{\rho,p}(\R^n),
        \end{equation}
        where the right-hand side is understood in the sense of tempered distributions.
	\end{theorem}

    %\commentCK{Should we mention here already that the representation in Theorem~\ref{thm:alt_char} can equivalently be rephrased as $
	%	H^{\rho,p}(\R^n)
	%	= \{u\in L^p(\R^n):\, \mathcal{F}^{-1}(\bla{\lambda_\rho}\hat{u})\in L^p(\R^n)\}$...It is done later, which is ok, I guess.}

	%\begin{remark}
	%	Due to the introductory considerations of this section, even though the statement above is only formulated for kernels \(\rho\) with compact support, we can always reduce to that case by means of \cite[Proposition~3.10]{BMCS25} \comment{Only if (H0)--(H3) are also preserved}. \commentC{The fact that we can assume without loss of generality that $\rho$ has compact support is explained in detail in the manuscript with Felix}
	%\end{remark}
%***6/8/1977
	\begin{proof} %\commentCK{Maybe start with brief idea?}\comment{I am not sure whether this is needed, since we also provide some guidelines at the beginnings of both inclusions.}\commentCK{OK!}
		We define 
		\begin{equation*}
		X^{\rho,p}(\R^n)
		\coloneqq \{u\in \mathcal{S}'(\R^n):\, \mathcal{F}^{-1}(\bla{\lambda_\rho}\hat{u})\in L^p(\R^n)\}.
		\end{equation*}
        We point out that the inclusions \(L^p(\R^n)\subset \mathcal{S}'(\R^n)\) and \(H^{\rho,p}(\R^n)\subset \mathcal{S}'(\R^n)\) are fulfilled automatically.
        Moreover, as soon as we know that \(H^{\rho,p}(\R^n)= X^{\rho,p}(\R^n)\), then the second equality in \eqref{eq:fourier_identities} follows directly from \(H^{\rho,p}(\R^n)\subset L^p(\R^n)\). Henceforth, we will use \(C>0\) to denote a generic positive constant that may change throughout the proof. 
        
	First, let us show that \(H^{\rho,p}(\R^n)\subset X^{\rho,p}(\R^n)\) continuously, i.e., that the symbol \(\bla{\lambda_\rho}\) defines a bounded operator from \(H^{\rho,p}(\R^n)\subset L^p(\R^n)\) into \(L^p(\R^n)\). We pick some arbitrary \(u\in C_c^\infty(\R^n) \subset H^{\rho,p}(\R^n)\) and observe that \(\mathcal{F}^{-1}((1+|\lambda_\rho|)\hat{u})\) is well defined as an element in \(\Scal'(\R^n)\) due to the polynomial growth of \(\lambda_\rho\) (see Lemma~\ref{lem:lambda_FintoF_pre}). Our goal is to demonstrate that \(\mathcal{F}^{-1}((1+|\lambda_\rho|)\hat{u})\in L^p(\R^n)\) because Proposition~\ref{prop:fourier_multipliers}\,\text{(ii)} then entails 
		\begin{equation}\label{est98}
		\norm{\mathcal{F}^{-1}(\bla{\lambda_\rho}\hat{u})}{L^p(\R^n)}
		\le C \norm{\mathcal{F}^{-1}((1+|\lambda_\rho|)\hat{u})}{L^p(\R^n)}.
		\end{equation}	
		%		Let \(\varphi\in C_c^\infty(\R^n)\). We compute
		%		\begin{equation}\label{eq:start}
		%			|\langle \mat,,hcal{F}^{-1}((1+|\lambda_\rho|^2)^{\frac{1}{2}}\hat{u}),\varphi\rangle|
		%			= |\langle \hat{u},(1+|\lambda_\rho|^2)^{\frac{1}{2}}\hat\varphi\rangle|
		%			\le \langle \hat{u},(1+|\lambda_\rho|^2)^{\frac{1}{2}}\psi\rangle
		%			\le \langle \hat{u},(1+|\lambda_\rho|)\psi\rangle,
		%		\end{equation}
		%		where we define \(\psi\coloneqq \text{sign} (\hat{u})|\hat\varphi|\) and first used Parceval's relation for the \(L^2\) inner product.  Note that since \(u\in C_c^\infty(\R^n)\subset \mathcal{S}(\R^n)\), it holds \(\hat{u}\in \mathcal{S}(\R^n)\) and then \(\mathcal{F}^{-1}((1+|\lambda_\rho|^2)^{\frac{1}{2}}\hat u )\in \mathcal{S}(\R^n)\). For the first term it follows directly from Hölder's inequality and another application of Parceval's identity that \(\langle \hat{u},\psi\rangle\le \norm{u}{L^p}\norm{\varphi}{L^{p'}}\). 	
        
        By definition, we obtain for any \(\varphi\in C_c^\infty (\R^n)\) 
        \begin{equation*}
            \langle \mathcal{F}^{-1}((1+|\lambda_\rho|)\hat{u}),\varphi\rangle_{\Scal'\times \Scal}
            =\langle u,\Fcal((1+|\lambda_\rho|)\mathcal{F}^{-1}\varphi)
		      =	\langle \hat{u},\Fcal^{-1}\varphi\rangle +\langle \hat{u},|\lambda_\rho|\Fcal^{-1}\varphi\rangle.
        \end{equation*}
       The first contribution on the right-hand side is bounded with the help of H\" older's inequality. %by \(\norm{ u}{L^p(\R^n)}\norm{\varphi}{L^{p'}(\R^n)}\).
		In order to estimate the remaining summand on the right-hand side, 
        %we turn to the Fourier version of the nonlocal fundamental theorem of calculus in \cite[Proposition~4.4]{BMCS25}. In particular, by the version of Lemma~\ref{lem:aux_ftoc} we have
        we observe that, using the identity \(\mathcal{F}(D^\rho_0 u)=\lambda_\rho \hat{u}=2\pi i \widehat{Q_\rho} \hat{u}\) (see~\eqref{eq:nlgrad_symbol}), it can be rewritten as
        \begin{align*}
        \left\langle \frac{-i\xi }{2\pi|\xi|^2\widehat{Q_\rho}(\xi)} \cdot \mathcal{F}(D^\rho_0 u),\abs{\lambda_\rho}\Fcal^{-1}\varphi  \right\rangle
        =
		\left\langle\sum_{j=1}^n\frac{-i\xi_j }{|\xi|}(\mathcal{F}(D^\rho_0 u))_j,\Fcal^{-1}\varphi\right\rangle
		= \langle \mathcal{F}(\mathcal{R}D^\rho_0 u),\Fcal^{-1}\varphi\rangle,
		\end{align*}
		%\begin{align*}
		      %\langle \mathcal{F}^{-1}(|\lambda_\rho|\hat{u}),\varphi\rangle
		%\langle \hat{u},|\lambda_\rho|\Fcal^{-1}\varphi\rangle
		      %&= \lim_{r\downarrow 0}\int_{B_r^c} \frac{-i\xi }{2\pi|\xi|^2\hat Q_\rho}\cdot\widehat{D^\rho_0 u}|\lambda_\rho|\Fca^l^{-1}\varphi\dl \xi\\
	       	%&= \langle\sum_{j=1}^n\frac{-i\xi_j }{|\xi|}(\widehat{D^\rho_0 u})_j,\Fca^l^{-1}\varphi\rangle\\
		%= \langle \Fcal(\Rcal D^\rho_0 u),\Fcal^{-1}\varphi\rangle,
		%\end{align*}
		where 
        \begin{equation}\label{eq:riesz_vec}
            \mathcal{R}U\coloneqq \sum_{j=1}^n  \mathcal R_j U_j,
            \qquad U\in L^p(\R^n;\R^n)
        \end{equation}
       % \(\Rcal D^\rho_0 u = \sum_{j=1}^n\mathcal{F}({R_j}(D^\rho_0 u)_j)\) 
       for \(\Rcal_j\) the \(j\)-th Riesz transforms with the Fourier symbol \(\xi\mapsto -i\xi_j/|\xi|\). Importantly, the Riesz transforms are bounded operators on \(L^p(\R^n)\) for \(p\in(1,\infty)\) and \(j\in \{1, \ldots, n\}\) (cf. \cite[Corollary~5.2.8]{Gr14}). 
		Since \(\rho\) has compact support, we know that \(D^\rho_0 u\in C_c^\infty (\R^n)\) and, in particular, \(\Rcal D^\rho_0 u\in L^2(\R^n)\). Thus, H\" older's inequality yields
		\begin{align*}
		|\langle \Fcal(\Rcal D^\rho_0 u),\Fcal^{-1}\varphi\rangle|
		= |\langle\Rcal D^\rho_0 u,\varphi\rangle|
		\le \norm{\Rcal D^\rho_0 u}{L^p(\R^n)}\norm{\varphi}{L^{p'}(\R^n)}
		\le C\norm{D^\rho_0 u}{L^p(\R^n;\R^n)}\norm{\varphi}{L^{p'}(\R^n)}.
		\end{align*}
		Taking now the supremum of the left-hand side over all \(\varphi\in C_c^\infty (\R^n)\) with \(\norm{\varphi}{L^{p'}(\R^n)}\le 1\), we arrive at
		\begin{equation*}
		\norm{\mathcal{F}^{-1}((1+|\lambda_\rho|)\hat{u})}{L^p(\R^n)}
		\le C(\norm{ u}{L^p(\R^n)}+\norm{D^\rho_0 u}{L^p(\R^n;\R^n)}),
		%<+\infty,
		\end{equation*}
		%where the finiteness follows from \(u\in C_c^\infty(\R^n) \subset H^{\rho,p}(\R^n)\).
        where \(\mathcal{F}^{-1}((1+|\lambda_\rho|)\hat{u})\) has a representative in \(L^p(\R^n)\) due to Riesz' representation theorem (see, e.g., \cite[Theorem~4.11]{Bre11} ).
		In view of~\eqref{est98}, this proves the claim for functions in \(C_c^\infty(\R^n)\). For arbitrary \(u\in H^{\rho,p}(\R^n)\), the statement can be deduced by the density of \(C_c^\infty(\R^n)\) in \(H^{\rho,p}(\R^n)\) (see \cite[Theorem~3.9]{BMCS25}).
		
		%		For the other inclusion, we use that \(C_c^\infty(\R^n)\) is actually dense in \(\mathcal{S}'(\R^n)\) in the weak\(^*\)-topology (see \comment{Grafakos, 2.3.23}).

		For the other inclusion, \( X^{\rho,p}(\R^n)\subset H^{\rho,p}(\R^n)\), let now \(u\in X^{\rho,p}(\R^n)\). 
        %New:
        %Lemma~\ref{lem:lp_implication} in combination with Proposition~\ref{prop:fourier_multipliers}~\textit{(iii)} for \(m=\bla{\lambda_\rho}\) shows that \(u\in L^p(\R^n)\).
        %end new
       We demonstrate first that \(u\) can be regarded as an element in \(L^p(\R^n)\), arguing by means of Proposition~\ref{prop:fourier_multipliers}~\textit{(iii)}. 
		  	For any \(\varphi\in \Scal(\R^n)\), we define the auxiliary function
	\begin{equation*}
	f_1
	\coloneqq \mathcal{F}^{-1}(\bla{\lambda_\rho}^{-1}\hat{\varphi}),
	\end{equation*}
    and observe that, by virtue of the second identity of Lemma~\ref{lem:fourier_identities}~\textit{\ref{item:fi_symm}} for the even multiplier \(\bla{\lambda_\rho}^{-1}\), it holds \(\bla{\lambda_\rho}^{-1}\mathcal{F}^{-1}\varphi= \Fcal^{-1} f_1\). Then
	we compute
	\begin{align*}
	\langle u,\varphi\rangle_{\Scal'\times \Scal}
	&= \langle u, \mathcal{F}(\mathcal{F}^{-1}\varphi)\rangle_{\Scal'\times \Scal}\\
	&= \langle u, \mathcal{F}(\bla{\lambda_\rho}(\bla{\lambda_\rho}^{-1}\mathcal{F}^{-1}\varphi))\rangle_{\Scal'\times \Scal}\\
	&=\langle u, \mathcal{F}(\bla{\lambda_\rho}\mathcal{F}^{-1}f_1)\rangle_{\Scal'\times \Scal}\\
	&= \langle \mathcal{F}^{-1}(\bla{\lambda_\rho}\hat{u}), f_1\rangle.
	\end{align*}
	Combining the assumption \(\mathcal{F}^{-1}(\bla{\lambda_\rho}\hat{u})\in L^p(\R^n)\) with Hölder's inequality and Proposition~\ref{prop:fourier_multipliers}~\textit{(iii)}, we conclude that
	\begin{align}
	\begin{split}\label{eq:bounds}
	|\langle u,\varphi\rangle_{\Scal'\times \Scal}|
	\le \norm{\mathcal{F}^{-1}(\bla{\lambda_\rho}\hat{u})}{L^p(\R^n)}
	\norm{f_1}{L^{p'}(\R^n)}
	\le C\norm{\mathcal{F}^{-1}(\bla{\lambda_\rho}\hat{u})}{L^p(\R^n)}\norm{\varphi}{L^{p'}(\R^n)}.
	\end{split}
	\end{align}
    %for some constant \(C>0\). 
    
	Choose now some \(\varphi\in L^{p'}(\R^n)\) and \((\varphi_k)_k\subset C_c^\infty(\R^n)\) such that \(\varphi_k\to \varphi\) in \(L^{p'}(\R^n)\). It follows by linearity from \eqref{eq:bounds} that \((\langle u,\varphi_k\rangle_{\Scal'\times \Scal})_k\) constitutes a Cauchy sequence in \(\R\), which attains a well defined limit. In particular, setting \(\langle u,\varphi\rangle\coloneqq \lim_{k\to \infty}\langle u,\varphi_k\rangle_{\Scal'\times \Scal}\), the expression \(\langle u,\varphi\rangle\) is independent of the choice of the approximating sequence. At the same time, we obtain along with \eqref{eq:bounds} that 
	\begin{equation*}
	|\langle u,\varphi\rangle_{\Scal'\times\Scal}|
	\le C \norm{\mathcal{F}^{-1}(\bla{\lambda_\rho}\hat{u})}{L^p(\R^n)}
	\norm{\varphi}{L^{p'}(\R^n)}.
	\end{equation*}
	The Riesz representation theorem yields that \(u\) indeed coincides with an element in \(L^p(\R^n)\). 
        
		It remains to prove that \(D^\rho u\in L^p(\R^n; \R^n)\).
		Precisely, we claim that \(\mathcal{F}^{-1}(\lambda_\rho \hat u)\) is well defined as an element in \(L^p(\R^n)\) and that \(D^\rho u=\mathcal{F}^{-1}(\lambda_\rho \hat u)\). Indeed, similar to the steps 
        %taken to prove Lemma~\ref{lem:lp_implication},
        taken above to establish that \(u\) coincides with an \(L^p\)-function, 
        one can also show that \(\mathcal{F}^{-1}(\lambda_\rho \hat u)\in L^p(\R^n)\).
        To see this, we note that, for any \(\Phi\in C_c^\infty(\R^n;\R^n)\) and recalling \eqref{eq:riesz_vec}, the function  
       % \commentCK{The function $m_\rho$ appears here again.}\comment{Still, I would refrain from using \(m_\rho\), since then we would probably need to introduce another symbol for the inverse fraction, and maybe this becomes too confusing.}
		\begin{equation*}
		f_2\coloneqq -\mathcal{F}^{-1}\left(\frac{|\lambda_\rho|}{\bla{\lambda_\rho}} \Fcal(\Rcal \Phi)\right)
		\end{equation*}
		is an element of \(L^p(\R^n)\) with 
        \begin{equation*}
             \norm{f_2}{L^p(\R^n)}
             \le C \norm{\Phi}{L^p(\R^n;\R^n)},
         \end{equation*}
         %for some \(C>0\), 
         by virtue of Proposition~\ref{prop:fourier_multipliers}~\textit{(i)} and the boundedness of the Riesz transforms.
        Then it holds, with the second identity of Lemma~\ref{lem:fourier_identities}~\textit{\ref{item:fi_antisymm}} applied componentwise,
         \begin{align*}
		 \langle \mathcal{F}^{-1}(\lambda_\rho \hat u),\Phi\rangle_{\Scal'\times \Scal}
		 &= \langle u, \mathcal{F}(\lambda_\rho\cdot \mathcal{F}^{-1}\Phi)\rangle_{\Scal'\times \Scal}\\
		 &= \left\langle u, \mathcal{F}\left(\bla{\lambda_\rho}\Big[\frac{i\xi}{|\xi|}\frac{|\lambda_\rho|}{\bla{\lambda_\rho}}\cdot \mathcal{F}^{-1}\Phi\Big]\right)\right\rangle_{\Scal'\times \Scal}\\
		 &=\langle u, \mathcal{F}(\bla{\lambda_\rho}\mathcal{F}^{-1}(f_2))\rangle_{\Scal'\times \Scal}\\
        &= \left\langle \mathcal{F}^{-1}(\bla{\lambda_\rho}\hat{u}), f_2\right\rangle.
		 \end{align*}
        Following, mutatis mutandis, 
        the computations in 
        \eqref{eq:bounds} and the steps thereafter, 
        %the remaining steps in the proof of Lemma~\ref{lem:lp_implication} 
        and applying Riesz' representation theorem, it follows as claimed that \(\mathcal{F}^{-1}(\lambda_\rho \hat u)\) has a representative in \(L^p(\R^n)\).
        
        Finally, we verify that \(D^\rho u=\mathcal{F}^{-1}(\lambda_\rho \hat u)\). We recall that since \(\rho\) has compact support, it holds \(\diver^\rho_0 \Phi \in C_c^\infty (\R^n)\) for any \(\Phi\in C_c^\infty (\R^n;\R^n)\) and the integral 
		\begin{equation*}
		\int_{\R^n} u(x)\diver^\rho_0 \Phi(x)\dl x
		\end{equation*}
		is well defined. With 
        %Lemma~\ref{lem:fourier_identities_as} 
        \eqref{eq:fourier_identities_as} we find
		\begin{equation*}
		-\langle u,\diver^\rho_0 \Phi\rangle 
        %\int_{\R^n} u(x)\diver^\rho_0 \Phi(x)\dl x
		= \langle \hat{u}, \mathcal{F}^{-1}(-\diver^\rho_0 \Phi)\rangle_{\Scal'\times \Scal}
		=	\langle \hat{u}, \lambda_\rho\cdot\mathcal
        F^{-1}{\Phi}\rangle_{\Scal'\times \Scal}
		%= \int_{\R^n}\mathcal{F}^{-1}(\lambda_\rho\hat{u})\cdot  \Phi\dl x
        =\langle \mathcal{F}^{-1}(\lambda_\rho\hat{u}),  \Phi\rangle
		\end{equation*}
		for all \(\Phi\in C_c^\infty (\R^n;\R^n)\). By definition, this establishes the sought identity \(D^\rho u=\mathcal{F}^{-1}(\lambda_\rho \hat u)\), which is exactly \eqref{eq:gen_fourier_symbol}.
	\end{proof}
    %\begin{remark}\label{rmk:gen_fourier_symbol}
    %\comment{maybe something like this? And then we could mention this also in the preliminaries}
    %    As can be seen from the last part of the proof, the function \(\lambda_\rho\) represents not only the Fourier symbol of the nonlocal gradient \(D^\rho_0\) acting on functions in \(C_c^\infty(\R^n)\), but also in the general case for \(D^\rho\) defined over \(H^{\rho,p}(\R^n)\). More precisely, for any \(u\in H^{\rho,p}(\R^n)\) it holds the identity \(D^\rho u=\mathcal{F}^{-1}(\lambda_\rho \hat u)\), where the right-hand side is understood in the sense of tempered distributions.
%\commentC{Yes. In fact, I would put it in the statement. But I doubt whether it should be in this result or in a separate one (to be placed much earlier than here)}\comment{I added it to the statement. If you agree, we can delete the remark.}
%    \end{remark}
    \begin{remark}[Hilbert space case]\label{rk:Hp=2}
        In the case \(p=2\), by Plancherel's identity, the alternative characterization 
        \begin{equation*}
			H^{\rho}(\R^n)
		= \{u\in \mathcal{S}'(\R^n):\, \bla{\lambda_\rho}\hat{u}\in L^2(\R^n)\}
	   \end{equation*}
       holds true (see also \cite[Proposition 4.7]{BGR26} for a recent related result).
%Under a different point of view, this is done in the recent work \cite[Proposition 4.7]{BGR26}. \commentCK{different point of view reads a little vague}
Moreover, the statement of Proposition~\ref{prop:fourier_multipliers} can be easily proved in this case, since the space of \(L^2\)-Fourier multipliers coincides with \(L^\infty\) (see \cite[Theorem~2.5.10]{Gr14}), noting that the assumptions on $\rho$ can be weakened, as done, for example, in \cite[Section~4.2 and Theorem~7.2~\textit{(ii)}]{BMCS25}.
    \end{remark}
    \begin{remark}[Kernels without compact support] We recall from Section~\ref{sec:hypotheses} that for any kernel \(\rho\) satisfying conditions (H0)--(H2) and (H3'), the space \(H^{\rho,p}(\R^n)\) is equally described as \(H^{\tilde{\rho},p}(\R^n)\) for some suitable kernel \(\tilde{\rho}\) of compact support still fulfilling conditions (H0)--(H2) and (H3').
   %With this in mind, we will assume that \(\rho\) has compact support and, in particular, that \(\rho\in L^1(\R^n)\).\color{black}
        %Let \(\rho\) be any kernel satisfying conditions (H0)--(H2) and (H3'). 
        %Then, by \cite[Proposition~3.10~(i)]{BMCS25}, it holds \(H_0^{\rho,p}(\Om)=H_0^{\tilde{\rho},p}(\R^n)\) with equivalent norms for 
        Then, Theorem~\ref{thm:alt_char} entails that 
        \begin{align*}
		c\norm{\mathcal{F}^{-1}(\bla{\lambda_{\tilde{\rho}}}\hat{u})}{L^{p}(\R^n)}
		\le \norm{u}{H^{\rho,p}(\R^n)}
		\le C\norm{\mathcal{F}^{-1}(\bla{\lambda_{\tilde{\rho}}}\hat{u})}{L^{p}(\R^n)}
		\end{align*}
		for all \(u\in H^{\rho,p}(\R^n)\) and some constants \(c,C>0\).
    \end{remark}

	 %%%%%%%%%%%%%%%%%%%%%%%%%%%%%%%%%

	\subsection{Fourier characterization of \(H^{\rho,p}_0(\Om)\)}\label{sec:fourier_bounded}

Let for the rest of the section $\Om \subset \R^n$ be open.
%As in \cite[Definition 3.4]{BMCS25}, we define 
Recalling~\eqref{eq:def_spaces}, the space $H^{\rho,p}_0(\Om)$ is defined as the set of $u \in H^{\rho,p} (\R^n)$ such that $u = 0$ a.e.\ in $\Om^c$.
If $\Om$ has a Lipschitz boundary, this set can be described as the closure of $C^{\infty}_c (\Om)$ in $H^{\rho,p} (\R^n)$; see \cite[Theorem 3.9]{BMCS25}.
The Fourier characterization of $H^{\rho,p} (\R^n)$ given by Theorem \ref{thm:alt_char}, translates directly to the space $H^{\rho,p}_0(\Om)$. 
%\commentCK{Since this is just a simple consequence of the characterization in Section 3.1, the question is whether an individual section is needed, but for structuring reasons, it is fine I guess. }

\begin{corollary}[Fourier characterization for complementary-value spaces]\label{cor:alt_charOm}
		Let \(p\in (1,\infty)\), let $\Om \subset\R^n$ be open such that $|\partial \Om| = 0$ and assume that \(\rho\) has compact support and satisfies hypotheses (H0)--(H2) and (H3'). Then
		\begin{equation*}
		H^{\rho,p}_0 (\Om)
		= \{u\in \mathcal{S}'(\R^n): \supp u \subset \overline{\Om} , \, \mathcal{F}^{-1}(\bla{\lambda_\rho}\hat{u})\in L^p(\R^n)\}
		\end{equation*}
		and, for all \(u\in H^{\rho,p}_0(\Omega)\), inequalities \eqref{eq:alt_char} hold.
	\end{corollary}
\begin{proof}
As $|\partial \Om| = 0$, the subspace $H^{\rho,p}_0 (\Om)$ can be equivalently described as the set of $u \in H^{\rho,p} (\R^n)$ such that $\supp u \subset \overline{\Om}$.
Since, for a locally integrable function, its support as a function coincides with its support as a distribution (see, e.g., \cite[Sections I.1 and I.3]{Schwartz66}), we conclude by Theorem \ref{thm:alt_char} the desired characterization.
\end{proof}

\begin{remark}
Instead of \cite[Definition 3.4]{BMCS25}, an equally natural definition of $H^{\rho,p}_0 (\Om)$ would have been the set of $u \in H^{\rho,p} (\R^n)$ such that $u = 0$ a.e.\ in $\ol{\Om}^c$; in other words, the set of $u \in H^{\rho,p} (\R^n)$ such that $\supp u \subset \ol{\Om}$.
In that case, the analogue of the characterization given in Corollary \ref{cor:alt_charOm} holds without the assumption $|\p \Om| = 0$.
\end{remark}

%\begin{remark}
%    The assumption that $|\partial \Om| = 0$ is a consequence of the definition of $H^{\rho,p}_0 (\Om)$ in \comment{\ldots}. If instead we define 
%    \begin{equation*}
%        H^{\rho,p}_0 (\Om)
%        \coloneqq \{u \in H^{\rho,p} (\R^n):\, 
%        u = 0 \text{ a.e. in } \overline{\Om}^c\},
%    \end{equation*}
%    then the condition that $|\partial \Om| = 0$ is suspended.
%\end{remark}

\subsection{Fourier characterization of the dual space}\label{sec:fourier_dual}
		For \(p\in (1,\infty)\) we denote by \(p'\in (1,\infty)\) its dual exponent fulfilling \(1/p+1/{p'}=1\), and define the topological dual space of \(H^{\rho,p}_0(\Om)\) by
		\begin{equation*}
			H^{-,\rho,p'}(\Om)
			\coloneqq (H^{\rho,p}_0(\Om))'.
		\end{equation*}
As $C^{\infty}_c (\Om)$ is dense in $H^{\rho,p}_0 (\Om)$ due to \cite[Theorem~3.9]{BMCS25} (under assumption (H0)), we have $H^{-,\rho,p'} (\Om) \subset \Dcal' (\Om)$.
%\color{red} However,[Rephrase for less focus and insert a little intuitive explanation at the end of the sentence.] \color{black}
%\comment{like this?}
We point out that, if \(\Om\) is bounded, then $H^{-,\rho,p'} (\Om)$ is in general not a subset of $\Scal' (\R^n)$, since $H^{-,\rho,p'} (\Om)$ only restricts the behavior of the functionals acting on functions defined over \(\Om\).  
    
	%	it holds (with \comment{BMCS, Theorem~4.5})
	%	\begin{equation*}
	%		H^{\rho,p}_0(\Om)\hookrightarrow L^2(\Om)\hookrightarrow 	H^{-,\rho,q}_0(\Om)
	%	\end{equation*}
	%	compactly if \(\lim_{t\to 0} t^{n-1}\bar{\rho}(t)=\infty\).
	The embedding \(H^{\rho,p}_0(\Om)\hookrightarrow L^{p}(\Om)\) is dense, hence the dual embedding $L^{p'}(\Om)\hookrightarrow H^{-,\rho,p'}(\Om)$ holds. Moreover, if $\Om$ is bounded, by \cite[Theorem~4.11]{BMCS25} (provided that (H0)--(H2) and (H3') are fulfilled), the embedding \(H^{\rho,p}_0(\Om)\hookrightarrow L^{p}(\Om)\) is compact, and, consequently, so is the embedding $L^{p'}(\Om)\hookrightarrow H^{-,\rho,p'}(\Om)$ (see, e.g., \cite[Theorem 6.4]{Bre11}).

	The Fourier characterization of $H^{\rho,p}(\R^n)$ given in Theorem~\ref{thm:alt_char} has a natural counterpart for its dual. 
	\begin{theorem}[Dual Fourier characterization]\label{thm:alt_dual_char} 
		Let \(p\in (1,\infty)\),  and assume that \(\rho\) has compact support and fulfills hypotheses (H0)-(H2) and (H3'). Then
		\begin{equation*}
		H^{-,\rho,p'}(\R^n)
		= \{u\in \mathcal{S}'(\R^n):\, \mathcal{F}^{-1}(\bla{\lambda_\rho}^{-1}\hat{u})\in L^{p'}(\R^n)\}
		\end{equation*}
		and, for all \(u\in H^{-,\rho,p'}(\R^n)\),
		\begin{align*}
			C_{p}^{-1} \norm{\mathcal{F}^{-1}(\bla{\lambda_\rho}^{-1}\hat{u})}{L^{p'}(\R^n)}
			\le \norm{u}{H^{-,\rho,p'}(\R^n)}
			\le c_{p}^{-1} \norm{\mathcal{F}^{-1}(\bla{\lambda_\rho}^{-1}\hat{u})}{L^{p'}(\R^n)} ,
		\end{align*}
		where \(c_{p},C_{p}>0\) are the constants from Theorem~\ref{thm:alt_char}.
	\end{theorem}

    We will deduce this statement as a special case of the following general result, which characterizes the duals of spaces defined via suitable multiplication in Fourier space.  %generalizes the relation between Bessel potential spaces and their duals. 
    The proof relies on standard techniques; nevertheless, since we were unable to locate this precise statement in the literature, we include it here for the sake of completeness.
    
    Henceforth, let \(m\in C^\infty(\R^n)\) be positive, even, and such that all its derivatives have polynomial growth. Further assume that also \(m^{-1}=1/m\) has polynomial growth.
    %is an \(L^p\)-Fourier multiplier for all \(p\in (1,\infty)\). 
    Set 
    \begin{equation*}
        H_{m}^p(\R^n)
        \coloneqq \{u\in \mathcal{S}'(\R^n):\, \mathcal{F}^{-1}(m\hat{u})\in L^{p}(\R^n)\}
    \end{equation*}
    and assume that \(\norm{\cdot}{H_m^p(\R^n)}\) is a norm equivalent to \(\norm{\mathcal{F}^{-1}(m\mathcal{F}(\cdot))}{L^p(\R^n)}\), i.e.,
    \begin{equation}\label{eq:norm_h}
        c_{H_m^p} \norm{\mathcal{F}^{-1}(m\hat{u})}{L^p(\R^n)}
        \le\norm{u}{H_{m}^p(\R^n)} 
        \le C_{H_m^p} \norm{\mathcal{F}^{-1}(m\hat{u})}{L^p(\R^n)}
    \end{equation}
    for some \(c_{H_m^p},C_{H_m^p}>0\). We point out that \(H_{m}^p(\R^n)\) is well defined following the assumptions on \(m\). 
    Finally, let \((H_{m}^p(\R^n))'\) be the topological dual space of \(H_{m}^p(\R^n)\). 
    \begin{lemma}\label{lem:aux_dual}
    Let \(p\in (1,\infty)\). If \(C_c^\infty(\R^n)\) is dense in \(H_{m}^p(\R^n)\), then  
        \begin{equation*}
            (H_{m}^p(\R^n))'
            = \{u\in \mathcal{S}'(\R^n):\, \mathcal{F}^{-1}(m^{-1}\hat{u})\in L^{p'}(\R^n)\},
        \end{equation*}
        and, for all \(u\in H_{m}^p(\R^n)\),
		\begin{align*}
			C_{H_{m}^p}^{-1} \norm{\mathcal{F}^{-1}(m^{-1}\hat{u})}{L^{p'}(\R^n)}
			\le \norm{u}{(H_{m}^p(\R^n))'}
			\le c_{H_{m}^p}^{-1} \norm{\mathcal{F}^{-1}(m^{-1}\hat{u})}{L^{p'}(\R^n)}.
		\end{align*}
    \end{lemma}
    \begin{proof}
        We define 
        \begin{equation*}
        X'
        \coloneqq \{u\in \mathcal{S}'(\R^n):\, \mathcal{F}^{-1}(m^{-1}\hat{u})\in L^{p'}(\R^n)\}
        \end{equation*}
        and prove first that \( (H_{m}^p(\R^n))'\subset X'\). 
        Applying the definition of the Fourier transform and its inverse for tempered distributions, it follows with the first identity of Lemma~\ref{lem:fourier_identities}~\textit{\ref{item:fi_symm}} that for any \(u\in (H_{m}^p(\R^n))'\) and any \(\varphi \in C_c^\infty(\R^n)\), 
        \begin{align}
        \begin{split}\label{eq:dual1_gen}
        \langle \mathcal{F}^{-1}(m^{-1}\hat{u}), \varphi\rangle_{\Scal'\times\Scal}
        = \langle u, \mathcal{F}(m^{-1}\mathcal{F}^{-1}\varphi)\rangle
        = \langle u, \mathcal{F}^{-1}(m^{-1}\hat{\varphi})\rangle.
        \end{split}
        \end{align}
        Moreover, additionally assuming that \(\norm{\varphi}{L^{p}}\le 1\), we have thanks to \eqref{eq:norm_h}
        \begin{equation*}\label{eq:norm_inequ_first}
        \norm{\mathcal{F}^{-1}(m^{-1}\hat{\varphi})}{H_{m}^p(\R^n)}
        \le C_{H_{m}^p} \norm{\mathcal{F}^{-1}(m(m^{-1}\hat{\varphi}))}{L^{p}(\R^n)}
        \le C_{H_{m}^p}
        \end{equation*}
        and, hence,     
        \begin{equation*}
        |\langle u, \mathcal{F}^{-1}(m^{-1}\hat{\varphi})\rangle|
        \le C_{H_{m}^p}\norm{u}{(H_{m}^p(\R^n))'}.
        \end{equation*}
        Taking the supremum over all \(\varphi\in C_c^\infty(\R^n)\) with \(\norm{\varphi}{L^{p}}\le 1\), it follows from \eqref{eq:dual1_gen}, together with the Riesz representation theorem, that
        \begin{equation*}
        C_{H_m^p}^{-1}\norm{\mathcal{F}^{-1}(m\hat{u})}{L^{p'}(\R^n)}
        \le \norm{u}{(H_{m}^p(\R^n))'},
        \end{equation*}
        and we obtain that \( (H_{m}^p(\R^n))'\subset X'\).
        
        For the reverse inclusion \(X'\subset (H_{m}^p(\R^n))'\), we pick \(u\in X'\) and again \(\varphi\in C_c^\infty(\R^n)\subset H_{m}^p(\R^n)\), where the inclusion holds by assumption. By virtue of \eqref{eq:norm_h}, the function 
        \begin{equation*}
        f
        \coloneqq \mathcal{F}^{-1}(m\hat{\varphi})
        \end{equation*}
        is an element of \(\Scal(\R^n)\) by the assumptions on \(m\) and 
        \begin{equation}\label{eq:potential1_gen}
            \norm{f}{L^{p}(\R^n)}
            \le c_{H_{m}^p}^{-1} \norm{\varphi}{H_{m}^p(\R^n)}.
        \end{equation}
        Now, again with the first identity of Lemma~\ref{lem:fourier_identities}~\textit{\ref{item:fi_symm}},
        \begin{align*}
        \langle u, \varphi\rangle
        = \langle u, \mathcal{F}^{-1}(m^{-1}\mathcal{F}(f))\rangle
        = \langle u, \mathcal{F}(m^{-1}\mathcal{F}^{-1}(f))\rangle
        =\langle \mathcal{F}^{-1}(m^{-1}\hat{u}),f\rangle.
        \end{align*}
        If \(\varphi\in C_c^\infty(\R^n)\) is chosen such that \(\norm{\varphi}{H_{m}^p(\R^n)}\le 1\), it follows from the above computations, Hölder's inequality and the norm estimate in \eqref{eq:potential1_gen} that 
        \begin{equation*}
        |\langle u, \varphi\rangle|
        \le c_{H_{m}^p}^{-1} \norm{\mathcal{F}^{-1}(m^{-1}\hat{u})}{L^{p'}(\R^n)}.
        \end{equation*}
        Taking in the left-hand side the supremum over all \(\varphi\in C_c^\infty(\R^n)\), and using the density of \(C_c^\infty(\R^n)\) in \(H_{m}^p(\R^n)\), yields
        \begin{equation*}
        \norm{u}{(H_{m}^p(\R^n))'}\le c_{H_{m}^p}^{-1} \norm{\mathcal{F}^{-1}(m^{-1}\hat{u})}{L^{p'}(\R^n)}.
        \end{equation*} 
        We conclude that \( (H_{m}^p(\R^n))'=X'\), as claimed.
    \end{proof}
    The proof of Theorem~\ref{thm:alt_dual_char} is now obtained as a simple corollary.
    \begin{proof}[Proof of Theorem~\ref{thm:alt_dual_char}]
        The result follows from Lemma~\ref{lem:aux_dual} with \(m\coloneqq \bla{\lambda_\rho}\) so that \(H_{m}^p(\R^n)=H^{\rho,p}(\R^n)\). We note that the validity of the estimates in \eqref{eq:norm_h} follows from Theorem~\ref{thm:alt_char}, while the fact that \(C_c^\infty(\R^n)\) is dense in \(H^{\rho,p}(\R^n)\) was shown in \cite[Theorem~3.9]{BMCS25} and only relies on hypothesis~(H0). 
    \end{proof}
    
Although not of Fourier-type, we end this section with an alternative characterization of the dual through derivatives of $L^p$-functions.

Let $\Om \subset \R^n$ be open.
Given $u_0 \in L^p (\Om)$ and $U \in L^p (\Om; \R^n)$, then $u_0 + \diver^{\rho} U \in H^{-, \rho, p} (\Om)$, since for all $\varphi \in H^{\rho, p'} (\Om)$,
\[
 \langle u_0 + \diver^{\rho} U, \varphi \rangle 
 = \int_{\Om} \left( u_0 \varphi - U\cdot D^{\rho} \varphi \right) \dl x \leq \norm{u_0}{L^p (\Om)} 
 \norm{\varphi}{L^p (\Om)} + \norm{U}{L^p (\Om; \R^n)} \norm{D^{\rho} \varphi}{L^{p'} (\Om; \R^n)},
\]
showing that 
\begin{equation*}
    \norm{u_0 + \diver^{\rho} U}{H^{-,\rho, p} (\Om)} 
    \le \left( \norm{u_0}{L^p (\Om)}^p + \norm{U}{L^p (\Om; \R^n)}^p \right)^{1/p}.
\end{equation*}
On the contrary, any element in \(H^{-,\rho, p} (\Om)\) can be identified with an element of \(L^p(\Om)+L^p(\Om;\R^n)\), and the proof of this property is the same as the classical one (see, e.g., Proposition~8.14 and Proposition~9.20 in \cite{Bre11}).

\begin{proposition}[Dual characterization through derivatives]\label{prop:charac_bijection}
Let $\Om \subset \R^n$ be open, \(p\in (1,\infty)\), and assume that \(\rho\) satisfies (H0).
For every $F \in H^{-,\rho, p} (\Om)$ there exist $u_0 \in L^p (\Om)$ and $U \in L^p (\Om; \R^n)$ such that $F = u_0 + \diver^{\rho} U$ and
\[
 \tilde{c}_p \norm{F}{H^{-,\rho, p} (\Om)} 
 \le  \max\{\norm{u_0}{L^p (\Om)}, \norm{U}{L^p (\Om; \R^n)}\}
 \le \tilde{C}_p \norm{F}{H^{-,\rho, p} (\Om)}
\]
for some constants \(\tilde{c}_p, \tilde{C}_p>0\).
When $\Om$ is bounded and \(\rho\) satisfies additionally (H1), (H2) and (H3'), one can take $u_0=0$.
\end{proposition}
\begin{proof}
    The proof can be obtained following the steps in \cite[Proposition~8.14]{Bre11}. Note that the norm defined in \eqref{eq:norm_def} is equivalent to the one used in \cite[Proposition~8.14]{Bre11}, so we find the equivalent dual norm characterization claimed in the statement; the constants \(\tilde{c}_p, \tilde{C}_p>0\) only depend on \(p\in (1,\infty)\). To derive the last claim we use that for any \(\rho\) subject to hypotheses (H0)--(H2) and (H3') it holds with \cite[Theorem~4.11~\textit{(i)}]{BMCS25} a Poincaré-Wirtinger-type inequality in \(H_0^{\rho, p'}(\Om)\). With this result one can again adapt the respective steps in the proof of \cite[ Proposition~8.14]{Bre11} accordingly. 
\end{proof}

We conclude with the comment that unlike in the classical case, if $u \in C_c^\infty(\R^n)$, then $D^{\rho} u$ is not, in general, a distribution.
Indeed, if $T$ is a distribution then $D^{\rho} T$ is defined, for all $\varphi \in C^{\infty}_c (\R^n; \R^n)$, as
\[
 \langle D^{\rho} T, \varphi \rangle = - \langle T, D^{\rho} \varphi \rangle,
\]
but $D^{\rho} \varphi$ is not, in general, in $C^{\infty}_c (\R^n)$, since it does not need to be compactly supported if the horizon of \(\rho\) is not finite.
This may be solved by working in the space $\Tcal (\R^n)$ of $C^{\infty} (\R^n)$ functions $u$ such that $\partial^{\alpha} u \in L^1 (\R^n) \cap L^{\infty} (\R^n)$ for all $\alpha \in \N^n_0$.
In this case, $D^{\rho} \varphi \in \Tcal (\R^n)$ when $\varphi \in \Tcal (\R^n)$, so $D^{\rho} T \in \Tcal (\R^n)'$ whenever $T \in \Tcal (\R^n)'$ (see, e.g., \cite{MSP99, Si20}).

\section{Characterization in terms of the $\rho$-Laplacian}\label{sec:char_laplacian}
In this section, we prove a characterization of $H^{\rho,p} (\R^n)$ using the square root of the negative  $\rho$-Laplacian, in the spirit of \eqref{eq:bessel_bessel} for the Bessel space.
%Although a full analysis of the negative $\rho$-Laplacian is beyond the scope of this article,\footnote{Be less apologetic...} a precise definition is required in order to develop the relation with the nonlocal spaces and, in particular, to introduce its square root.\footnote{Comment on compact support commented out.} 
 After introducing the \(\rho\)-Laplacian and deriving its main properties that are important for our purposes, the presentation is split into the cases \(p=2\) and \(p\in (1,\infty)\). For the latter, we shall use the language of unbounded operators in Banach spaces, and for \(p=2\), the special theory of unbounded operators in Hilbert spaces (see, e.g., \cite{Kato76, ReSi80, Rudin91, CaHa98}).

%We assume throughout most of this section that $\rho$ has compact support.This implies that $D^{\rho}_0 u \in C^{\infty}_c (\R^n)$ when $u \in C^{\infty}_c (\R^n)$.

%\comment{I was reconsidering how to present the different cases of \(p\). I am now inclined to think that the best option would be to present the definition and the density result in the general \(p\) case and then to split the sections fromthere onward. This way, we only need to give the definition once in a unified way and also the density result can be provided directly for all \(p\).}

\subsection{Definition and first properties of the $\rho$-Laplacian}\label{sec:laplace_prelim} 
%\commentCK{Suggestion for restructuring: Move (i) before the definition, do not mention the unbounded operator property, and move "important preliminary relations" up, before introducing $\rho$-Laplacian as a unbounded operator.  }
%\begin{definition}[$\rho$-Laplacian]\label{de:Laplacianp}
%Let \(p\in (1,\infty)\) and \(\rho\) satisfy (H0).
%\begin{enumerate}[label=(\roman*)]
%\item
%The $\rho$-Laplacian for smooth functions is denoted by $\Delta^{\rho}_0$ and defined as the unbounded operator in $L^p (\R^n)$ with domain $C^{\infty}_c (\R^n)$ %through\footnote{Consider using $\varphi$ here...}
%\[
% \Delta^{\rho}_0 u \coloneqq \diver^{\rho}_0 D^{\rho}_0 u , \qquad u \in C^{\infty}_c (\R^n) .
%\]
%
%\item
%Given $u \in L^p (\R^n)$, we say that $\Delta^{\rho} u \in L^p (\R^n)$ if there exists $v \in L^p (\R^n)$ such that
%\begin{equation}\label{eq:vp}
% \int_{\R^n} v \varphi = \int_{\R^n} u \Delta^{\rho}_0 \varphi , \qquad \varphi \in C^{\infty}_c (\R^n).
%\end{equation}
%The $\rho$-Laplacian $\Delta^{\rho}$ is the unbounded operator in $L^p (\R^n)$ with domain 
%\[
% D(\Delta^{\rho}) = \{ u \in L^p (\R^n) : \Delta^{\rho} u \in L^p (\R^n) \}
%\]
%defined as $\Delta^{\rho} u = v$, where $v$ is the unique $L^p (\R^n)$ function satisfying \eqref{eq:vp}.
%\end{enumerate}
%\end{definition}
In analogy to the definition of the pointwise and distributional nonlocal gradient and divergence in Section~\ref{sec:recap_nl_sob}, we introduce in the following the pointwise and distributional nonlocal \(\rho\)-Laplacian. 
%In particular/to be precise, ...; and maybe here already: For simplicity we will assume throughout the section that \(\rho\) has compact support.
%}
For \(u\in C_c^\infty(\R^n)\), we define its pointwise $\rho$-Laplacian $\Delta^{\rho}_0$ in the natural way as 
%\commentCK{With our notation convention, is it more consistent to write $\varphi$ or $u$. Since it will be used as test function further down, probably ok.}\comment{Here it is indeed really complicated. Nevertheless, I would actually maybe rather use \(u\). Also, this is what we use in the introduction when introducing \(\diver^{\rho}_0\) and \(D^{\rho}_0 \).}
%the unbounded operator in $L^p (\R^n)$ with domain $C^{\infty}_c (\R^n)$ through\footnote{Consider using $\varphi$ here...}
\begin{align}\label{def:Drho0}
 \Delta^{\rho}_0 u \coloneqq \diver^{\rho}_0 D^{\rho}_0 u,
\end{align}
with $D_0^\rho$ and $\diver_0^\rho$ given by \eqref{eq:def:drho} and \eqref{eq:def:divrho}.
Note that right-hand side in~\eqref{def:Drho0} is well defined for any \(\rho\) subject to (H0), since \(D^\rho_0\) is a bounded Lipschitz function (see comments after \eqref{eq:pw_equals_dual}).
%We recall that if $\rho$ additionally has compact support, then we even have $D^{\rho}_0 u \in C^{\infty}_c (\R^n)$ for any $u \in C^{\infty}_c (\R^n)$.
We will assume from here on that \(\rho\) has compact support, and recall that if this is the case, then we have $D^{\rho}_0 u, \Delta_0^\rho u \in C^{\infty}_c (\R^n)$ for any $u \in C^{\infty}_c (\R^n)$.

We shall now first collect some preliminary relations. 
By \eqref{eq:nlgrad_symbol}, it holds
\begin{equation}\label{eq:FourierLaplacian_pre}
 \Fcal(\Delta_0 ^{\rho} u) 
 =\lambda_{\rho} \cdot \mathcal{F}(D^\rho_0 u) 
 = \lambda_{\rho} \cdot \lambda_{\rho} \hat{u} 
 = - \abs{\lambda_{\rho}}^2 \hat{u} ,
 \qquad u\in C_c^\infty(\R^n),
\end{equation}
since $\lambda_{\rho} $ is purely imaginary on \(\R^n\). 

In addition, as demonstrated in the proof of \cite[Lemma~3.6]{BMCS25}, if \(\rho\) fulfills (H0), then for any \(k\in C^{\infty}_c (\R^n), U\in C^{\infty}_c (\R^n;\R^n)\) it holds $k * \diver^{\rho}_0 U = \diver^{\rho}_0 (k * U)$. This remains true for kernels \(k\in L^1(\R^n)\) and, by a similar calculation, one gets an analogous identity for the gradient operator \(D^\rho_0\). Hence, for \(\rho\) with compact support, it follows from \(\Delta^{\rho}_0  = \diver^{\rho}_0 D^{\rho}_0\) that 
\begin{equation}\label{eq:commutation}
    k * D_0^{\rho} u 
    = D_0^{\rho} (k * u),\quad
    k * \Delta_0^{\rho} u 
    = \Delta_0^{\rho} (k * u),
    \qquad k\in L^1(\R^n), u\in C_c^\infty(\R^n).
\end{equation}

For \(u\in L^1(\R^n)+L^\infty(\R^n)\), we define the distributional, or weak, nonlocal \(\rho\)-Laplacian $\Delta^\rho u$ as the element in \(\Dcal'(\R^n)\) such that
    \begin{equation}\label{def_weakLaplacian}
        \langle \Delta^\rho u,\varphi\rangle_{\Dcal'\times \Dcal} 
        =\langle u , \Delta^{\rho}_0 \varphi \rangle,
        \qquad \varphi\in C_c^\infty(\R^n).
    \end{equation}
 In what follows, we understand the \(\rho\)-Laplacian as an unbounded operator in $L^p(\R^n)$. 

\begin{definition}[$\rho$-Laplacian]\label{de:Laplacianp}
Let \(p\in (1,\infty)\) and let \(\rho\)  have compact support and satisfy (H0).
   % \commentCK{Since we carefully introduce above the meaning of (tempered) distributions lying in $L^p(\R^n)$ and with Leon's planned comment earlier on $\varphi\in \Scal(\R^n)$ vs. $\varphi\in C_c^\infty(\R^n)$, I think we should adjust this first part of the definition for consistency. }\comment{One thing to keep in mind, however, is that the point of the above definition is that \(\Delta^\rho\) is introduced in its "very" weak form (partially integrating twice, not only ones). In some way we need to make that clear.} \commentCK{Agreed. My comment before was intended towards the expression $\Delta^\rho u\in L^p(\R^n)$, which can also be interpreted for distributions...Is it an option to define $\Delta^\rho$ first distributionally?  }\comment{In \(\Dcal'\) then? We could do that, I guess. }
The $\rho$-Laplacian $\Delta^{\rho}$ is the unbounded operator in $L^p (\R^n)$ with domain 
\[
 D(\Delta^{\rho}) = \{ u \in L^p (\R^n) : \Delta^{\rho} u \in L^p (\R^n) \},
\]
where $\Delta^\rho u$ is defined as in~\eqref{def_weakLaplacian} for $u\in L^p(\R^n)$.
%$\Delta^{\rho} u = v$, where $v$ is the $L^p$-function satisfying \eqref{eq:vp}.
\end{definition}
Strictly speaking, the $\rho$-Laplacian in $L^p (\R^n)$ depends on $p$,
%requires a notation indicating its dependence on $p$, 
but, for the sake of simplicity, we will plainly denote it by $\Delta^{\rho}$.  %\commentCK{Comment on the use of notation, which is the same for the weak and "our" $\rho$-Laplacian.}\comment{Actually, I think there is no reason, since there is no danger of confusion. We do nothing else for the nonlocal gradient.}
%\commentCK{Do we have a good terminology to use when we want to indicate a specific underlying $p$. Later I saw something like "in the sense of $L^p(\R^n)$", is there a better option?}% and make the $p$ explicit in the statements.
%\comment{In fact,what you refer to is the only point where me make the dependence specific, I think. Carlos and me were also talking about this issue, but didn't come up with a better alternative. We could of course introduce this dependency into the notation for this case and write \(\Delta^{\rho}_p\) or \(\Delta^{\rho,p}\) or something similar.}

\begin{remark}
a) Observe that for \(\rho\) with compact support a direct computation shows that \(C_c^\infty(\R^n)\subset D(\Delta^\rho)\) and  $\Delta_0^\rho$ and $\Delta^\rho$ coincide on $C_c^\infty(\R^n)$.
First, recall that if $u \in C^{\infty}_c (\R^n)$, then $D^{\rho} u \in C^{\infty}_c (\R^n)$ and $\Delta_0^{\rho} u \in C^{\infty}_c (\R^n)$.
Furthermore,
\begin{equation}\label{eq:Delta0Symmetric}
 \langle\Delta_0^{\rho} u, \varphi \rangle
 =- \langle D_0^{\rho} u, D_0^{\rho} \varphi\rangle 
 = \langle u, \Delta_0^{\rho} \varphi\rangle,
 \qquad \varphi\in C_c^\infty(\R^n).
\end{equation}
Therefore, $u \in D (\Delta^{\rho})$ and $\Delta^{\rho} u = \Delta_0^{\rho} u$.\smallskip
%\comment{Added:}

b) More generally, 
%if \(\rho\) has compact support
%\comment{I realized that we need to change the statement of Lemma~\ref{lem:lambda_FintoF_pre} a bit, since we need to stress that (H0) and compact support of \(\rho\) alone imply that \(\lambda_\rho\) has linear growth; I devised a draft}, 
 it also holds that \(\Scal(\R^n)\subset D(\Delta^\rho)\).
To be precise, combining the linear growth of \(\lambda_\rho\) due to Lemma~\ref{lem:lambda_FintoF_pre} with the first identity of Lemma~\ref{lem:fourier_identities}~\textit{\ref{item:fi_symm}} and \eqref{eq:FourierLaplacian_pre}, one obtains     for any \(u\in \Scal(\R^n)\) that
     \begin{equation}\label{eq:rholap_S}
         \langle \mathcal{F}^{-1}(-\abs{\lambda_\rho}^2\hat{u}) , \varphi \rangle
        = \langle u,\mathcal{F}^{-1}(-\abs{\lambda_\rho}^2\hat{\varphi})\rangle 
        = \langle u, \Delta_0^{\rho} \varphi\rangle,
        \qquad \varphi\in C_c^\infty(\R^n).
     \end{equation}
     In particular, one can conclude from these considerations that it makes no difference for our subsequent discussion whether  \(\Delta^\rho_0\) is introduced in \eqref{def:Drho0} for functions  in \(C_c^\infty(\R^n)\) or \(\Scal(\R^n)\). 
     \end{remark}

The following auxiliary estimate and the resulting approximation result will be used repeatedly below. The proof can be found in Appendix~\ref{sec:appendix_aux}.
\begin{lemma}\label{lem:lap_est}
    Let \(p\in (1,\infty)\) and assume that \(\rho\) satisfies (H0)  %--(H2) and (H3') 
    and has compact support.
    %Let \(g\in C^\infty(\R^n)\) be such that 
    %\begin{equation*}
    %    g, \Delta^\rho g
    %    \in L^p(\R^n).
    %\end{equation*}
    %Then, 
    For any \(g\in C^\infty(\R^n)\cap H^{\rho,p}(\R^n)\)\(\cap D(\Delta^\rho)\) and \(f\in C_c^\infty(\R^n)\) it holds 
		\begin{equation}\label{eq:est_lap}
		\norm{\Delta^\rho_0 (fg)-f\Delta^\rho g}{L^p(\R^n)}
		\le C\norm{\nabla f}{W^{1,\infty}(\R^n;\R^n)}\norm{g}{H^{\rho,p}(\R^n)}
		\end{equation}
		for some constant \(C=C(\rho)>0\).
\end{lemma}
In particular, for \(\sigma \in C_c^\infty(\R^n)\) with \(\sigma\equiv 1\) in \(B_1\), and the smooth sequence of cut-off functions $(\sigma_j)_{j}$ defined as \(\sigma_j(\cdot)\coloneqq \sigma(\cdot/j)\), one gets
        \begin{equation}\label{eq:simul_approx_lem}
            \sigma_jg\to g,\quad
           \Delta^\rho_0(\sigma_jg)\to \Delta^\rho g
            \qquad \text{in } L^p(\R^n).
        \end{equation}
         In order to see this, we observe that \(\sigma_jg\to g\) and \(\sigma_j\Delta^{\rho} g\to \Delta^{\rho} g\) by dominated convergence, and that \(\norm{\nabla \sigma_j}{W^{1,\infty}}\le C/j\) for some uniform constant \(C>0\). The claimed convergence of the \(\rho\)-Laplacian follows then immediately from the above-proved estimate~\eqref{eq:est_lap} for the choice \(f= \sigma_j\), since
\begin{equation*}
		\norm{\Delta^{\rho}_0 (\sigma_j g)-\sigma_j\Delta^{\rho} g}{L^p(\R^n)}
		\le C\norm{\nabla \sigma_j}{W^{1,\infty}(\R^n;\R^n)}\norm{g}{H^{\rho,p}(\R^n)}.
\end{equation*} 

We demonstrate next that \(D(\Delta^{\rho})\subset H^{\rho, p}(\R^n)\).
The proof is inspired by the techniques of \cite[Theorem~5.2]{BMCS25} (see also Remark~\ref{rmk:ftoc} below).

\begin{proposition}\label{prop:poincare2}
    Let \(p\in (1,\infty)\) %\commentCK{Also here, the role of $p$ for $\Delta^\rho$ is somewhat implicit, but maybe ok here...} 
    and assume that \(\rho\) has compact support and satisfies (H0)--(H3). 
    Then there exists a vector radial function \(K_\rho\in L^1(\R^n;\R^n)\) such that 
    \begin{equation}\label{eq:fi_lap_to_grad}
        \widehat{K_\rho}
        = \frac{\lambda_\rho}{1+|\lambda_\rho|^2}.
    \end{equation}
    
    Moreover, %for $u\in D(\Delta^\rho)$ 
    \begin{align}\label{identity19}D^\rho u= K_\rho*(u-\Delta^\rho u)\in L^p(\R^n), \quad u\in D(\Delta^\rho)\end{align}
   % \begin{equation*}
    %    K_\rho*(\Phi-\Delta^\rho_0\Phi)
     %   = \diver^\rho_0 \Phi,
     %   \qquad \Phi\in C_c^\infty(\R^n;\R^n).
    %\end{equation*}
    and, in particular, 
    \begin{equation*}
        D(\Delta^\rho)\subset H^{\rho,p}(\R^n).
    \end{equation*}
\end{proposition}
\begin{proof}
    By \(C>0\) we will denote a constant that may change from line to line and possibly depends on \(\rho\) and, occasionally, on some multi-index \(\alp\in\N_0^n\).\smallskip
    
    \textit{Step~1: Preliminary estimates.}
    We recall initially that due to the assumptions and by \eqref{eq:lest_qrho_pre} it holds \(\lambda_\rho(\xi)\ge C|\xi|^{\sigma} \ge C\) for all \( \xi\in B_1^c\). In particular, since
    \begin{equation*}
        \frac{|\lambda_\rho|}{1+|\lambda_\rho|^2}
        \le C \frac{1}{|\lambda_\rho|}
        \le C,
        \qquad \xi\in B_1^c
    \end{equation*}
    and because \(\lambda_\rho \in C^\infty(\R^n)\) (see Lemma~\ref{lem:lambda_FintoF_pre}), 
    the right-hand side of \eqref{eq:fi_lap_to_grad} is uniformly bounded in \(\R^n\) and defines a tempered distribution, as which its Fourier inverse is well defined.     Moreover, with the help of Lemma~\ref{lem:estimates} one can compute that
    \begin{align}\label{eq:est_aux_lap_to_grad}
        \begin{split}
            \Bigg| \p^\alp \frac{\lambda_\rho(\xi)}{1+|\lambda_\rho(\xi)|^2}\Bigg|
        &\le \sum_{\beta\le \alp}
		|\p^\beta \lambda_\rho(\xi)||\p^{\alp-\beta} (1+|\lambda_\rho(\xi)|^2)^{-1}| \\
        &\le C\sum_{\beta\le \alp} |\lambda_\rho(\xi)||\xi|^{-|\beta|}
		|\lambda_\rho(\xi)|^{-2}|\xi|^{-|\alp|+|\beta|} \\
        &= C |\xi|^{-|\alp|} |\lambda_\rho(\xi)|^{-1}\\
        &\le C |\xi|^{-|\alp|-\sigma},
        \qquad \xi\in B_1^c
        \end{split}
    \end{align}
    for any \(\alp\in \N_0^n\).\smallskip

    \textit{Step~2: Validating that $K_\rho\in L^1(\R^n)$. } We show now that
    \begin{equation*}
        K_\rho(x)\coloneqq \Fcal^{-1}\left(\frac{\lambda_\rho}{1+|\lambda_\rho|^2}\right)(x)
        = \int_{\R^n} \frac{\lambda_\rho(\xi)}{1+|\lambda_\rho(\xi)|^2}e^{2\pi i x\cdot\xi } \dl \xi
        = \int_{\R^n} \frac{\lambda_\rho(\xi)}{1+|\lambda_\rho(\xi)|^2}e^{2\pi i |x|\xi_1 } \dl \xi
    \end{equation*}
    exists pointwise for a.e. \(x\in \R^n\) and lies actually in \(L^1(\R^n)\); note that the last identity above follows from the symmetry of \(\lambda_\rho\) and similar arguments as in \cite[Lemma~2.4]{BMCS25}. To that end, we will derive the estimates
    \begin{align}
            & |K_\rho(x)|\le C\frac{1}{|x|^{n-\sigma}}, \qquad x\in B_1\setminus \{0\},\label{eq:krho_xb1}\\
        &|K_\rho(x)| \le C\frac{1}{|x|^{n+1}}, \qquad x\in B_1^c. \label{eq:krho_xb1c}
    \end{align}
    Together, those yield
    \begin{equation}\label{eq:final_est_krho}
             \int_{\R^n}|K_\rho(x)|\dl x
             \le C\int_{B_1}\frac{1}{|x|^{n-\sigma}}\dl x + C\int_{B_1^c} \frac{1}{|x|^{n+1}}\dl x
             <\infty,
        \end{equation}
    which proves indeed that $K_\rho\in L^1(\R^n)$.

    We assume first that \(x\in B_1^c\) and show \eqref{eq:krho_xb1c}. With repeated integration by parts, we compute
    \begin{equation*}
        \int_{\R^n} \frac{\lambda_\rho(\xi)}{1+|\lambda_\rho(\xi)|^2}e^{2\pi i |x|\xi_1 } \dl \xi
        = \frac{1}{(2\pi i |x|)^{n+1}}\left[ \int_{B_{1}} + \int_{B_{1}^c} \right] \p^{n+1}_{\xi_1}\left(\frac{\lambda_\rho(\xi)}{1+|\lambda_\rho(\xi)|^2}\right)e^{2\pi i |x|\xi_1 } \dl \xi.
    \end{equation*}
    Inserting in the integral above over \(B_1^c\) the estimate~\eqref{eq:est_aux_lap_to_grad} with the choice \(\alp=(n+1,0,\dots,0)\) and using on \(B_1\) the smoothness of \(\lambda_\rho\), yields
    \begin{equation}\label{eq:lg_xoutside}
        \int_{\R^n} \frac{\lambda_\rho(\xi)}{1+|\lambda_\rho(\xi)|^2}e^{2\pi i |x|\xi_1 } \dl \xi
        \le C\frac{1}{|x|^{n+1}}\Bigg[1+\int_{B_{1}^c} |\xi|^{-(n+1)}\dl \xi\Bigg]
		\le C\frac{1}{|x|^{n+1}},
        \qquad  x\in B_1^c.
    \end{equation}
    
    Now we study the remaining case that \(x\in B_1\setminus \{0\}\) and derive \eqref{eq:krho_xb1}. The first contribution in the sum 
    \begin{equation*}
        \int_{\R^n} \frac{\lambda_\rho(\xi)}{1+|\lambda_\rho(\xi)|^2}e^{2\pi i |x|\xi_1 } \dl \xi
        = \left[ \int_{B_{1}} + \int_{B_{1}^c} \right] \frac{\lambda_\rho(\xi)}{1+|\lambda_\rho(\xi)|^2}e^{2\pi i |x|\xi_1 } \dl \xi
    \end{equation*}
    is, again by the smoothness of \(\lambda_\rho\), bounded by
    \begin{equation*}
		\Bigg|\int_{B_{1}}\frac{\lambda_\rho(\xi)}{1+|\lambda_\rho(\xi)|^2}e^{2\pi i |x|\xi_1}\dl \xi\Bigg|
		\le C
		\le C\frac{1}{|x|^{n-\sigma}},
        \qquad x\in B_1\setminus \{0\},
		\end{equation*}
        and we derive subsequently a similar bound also for the second integral. To that end, we partially integrate \((n-1)\)-times to obtain 
        \begin{equation}\label{eq:lq_B1c}
        \begin{split}
			\int_{B_1^c} \frac{\lambda_\rho(\xi)}{1+|\lambda_\rho(\xi)|^2}e^{2\pi i |x|\xi_1}\dl \xi
			&= \frac{(-1)^{n-1}}{(2\pi i|x|)^{n-1}}
			\int_{B^c_1}
			\p^{n-1}_{\xi_1}\left(\frac{\lambda_\rho(\xi)}{1+|\lambda_\rho(\xi)|^2}\right)
			e^{2\pi i |x|\xi_1}\dl \xi
			+  R(x),
		\end{split}
        \end{equation}
		where the boundary integrals resulting from the partial integrations are collected in the remainder term
		\begin{equation*}
		R(x) 
		= \sum_{k=1}^{n-1} \frac{(-1)^{k-1}}{(2\pi i|x|)^{n-1}}
		\int_{\p B_1} \p^{k-1}_{\xi_1}\left(\frac{\lambda_\rho(\xi)}{1+|\lambda_\rho(\xi)|^2}\right)\p^{n-1-k}_{\xi_1}e^{2\pi i |x| \xi_1}\dl \mathcal{S}(\xi).
		\end{equation*} 
        With \eqref{eq:est_aux_lap_to_grad} we get 
        \begin{align}
            \begin{split}\label{eq:lg_gen_fin_est_R}
            |R(x) |
            \le C  \sum_{k=1}^{n-1} 
            \frac{1}{|x|^{k}}
			\int_{\p B_{1}}|\xi|^{-(k-1+\sigma)}\dl \mathcal{S}(\xi)
            \le C \frac{1}{|x|^{n-\sigma}},
            \qquad x\in B_1\setminus \{0\} .
            \end{split}
		\end{align}
        To estimate the right-hand side in~\eqref{eq:lq_B1c} further, we apply in the following the decomposition of \(B_1^c\) into \(B_{1/|x|}\setminus B_1\) and \(B_{1/|x|}^c\).  Partially integrating once more on the latter domain it holds 
    \begin{align}
    \label{aux54}	
    \begin{split}& \frac{(-1)^{n-1}}{(2\pi i|x|)^{n-1}} 
			\int_{B^c_{1/|x|}}
			\p^{n-1}_{\xi_1}\left(\frac{\lambda_\rho(\xi)}{1+|\lambda_\rho(\xi)|^2}\right)
			e^{2\pi i |x|\xi_1}\dl \xi\\
			&= 
			\frac{(-1)^{n}}{(2\pi i|x|)^{n}} 
			\int_{\p B_{1/|x|}}
			\p^{n-1}_{\xi_1}\left(\frac{\lambda_\rho(\xi)}{1+|\lambda_\rho(\xi)|^2}\right)
			e^{2\pi i |x|\xi_1}\dl \mathcal{S}(\xi)\\
			&+\frac{(-1)^{n}}{(2\pi i|x|)^{n}} 
			\int_{B^c_{1/|x|}}
			\p^{n}_{\xi_1}\left(\frac{\lambda_\rho(\xi)}{1+|\lambda_\rho(\xi)|^2}\right)
			e^{2\pi i |x|\xi_1}\dl \xi.
            \end{split}
		\end{align}  
        
        Estimating now the remaining contributions on \(B_1^c\) with the bounds from \eqref{eq:est_aux_lap_to_grad} leads in view of~\eqref{eq:lq_B1c}  to
		\begin{align}
            \begin{split}\label{eq:lg_gen_fin_est_1}
                \frac{1}{|x|^{n-1}}
			\int_{B_{1/|x|} \setminus B_{1}}
			\Big|\p^{n-1}_{\xi_1}\left(\frac{\lambda_\rho(\xi)}{1+|\lambda_\rho(\xi)|^2}\right)\Big|\dl \xi
			\le C\frac{1}{|x|^{n-1}}
			\int_{B_{1/|x|}}
			\frac{1}{|\xi|^{n-1+\sigma}}\dl \xi
			= C
			\frac{1}{|x|^{n-\sigma}},
            \end{split}
		\end{align}
		and, with the help of~\eqref{aux54}, one gets
		\begin{align}
		  \begin{split}\label{eq:lg_gen_fin_est_2}
		      \frac{1}{|x|^{n}} 
			\int_{\p B_{1/|x|}}
			\Big|\p^{n-1}_{\xi_1}\left(\frac{\lambda_\rho(\xi)}{1+|\lambda_\rho(\xi)|^2}\right)\Big|
			\dl \mathcal{S}(\xi)
			\le C\frac{1}{|x|^{n}} 
			\int_{\p B_{1/|x|}}
			\frac{1}{|\xi|^{n-1+\sigma}}
			\dl \mathcal{S}(\xi)
			= C\frac{1}{|x|^{n-\sigma}}, 
		  \end{split}
		\end{align}
as well as,
		\begin{align}
		  \begin{split}\label{eq:lg_gen_fin_est_3}
		      &\frac{1}{|x|^{n}}
		\int_{B^c_{1/|x|}}
		\Big|\p^{n}_{\xi_1}\left(\frac{\lambda_\rho(\xi)}{1+|\lambda_\rho(\xi)|^2}\right)\Big|\dl \xi
		\le C\frac{1}{|x|^{n}}
		\int_{B^c_{1/|x|}}
		\frac{1}{|\xi|^{n+\sigma}}\dl \xi
		= C\frac{1}{|x|^{n-\sigma}},
		  \end{split}
		\end{align}
        where always \(x\in B_1\setminus \{0\}\).

        Overall, from \eqref{eq:lg_xoutside} for \(x\in B_1^c\) and \eqref{eq:lq_B1c} for \(x\in B_1\setminus \{0\}\) together with estimates \eqref{eq:lg_gen_fin_est_R}--\eqref{eq:lg_gen_fin_est_3}, we obtain \eqref{eq:final_est_krho}
        %\begin{equation*}
        %     \int_{\R^n}|K_\rho(x)|\dl x
         %    \le \int_{B_1}\frac{1}{|x|^{n-\sigma}}\dl x + \int_{B_1^c} \frac{1}%{|x|^{n+1}}\dl x
        %     <\infty
        %\end{equation*}
        and \(K_\rho\in L^1(\R^n;\R^n)\). Note, in addition, that \(K_\rho\) is vector radial, since this property is preserved by Fourier inversion. \smallskip

\textit{Step~3: Proof of~\eqref{identity19} and $D(\Delta^\rho)\subset H^{\rho, p}(\R^n)$.}
           Next, the identity \(K_\rho*(\Phi-\Delta^\rho_0\Phi)= \diver^\rho_0 \Phi\) for \(\Phi\in C_c^\infty(\R^n;\R^n)\) (where \(K_\rho*\Phi\) is defined as the trace of the matrix \(((K_\rho)_i*\Phi_j)_{ij}\), and similarly for \(K_\rho*\Delta^\rho_0\Phi\)) follows by definition of \(K_\rho\) and \eqref{eq:nlgrad_symbol}, \eqref{eq:FourierLaplacian_pre}, stating the Fourier symbols of \(D^\rho_0,\diver^\rho_0\) and \(\Delta^\rho_0\). Precisely, we have
        \begin{equation}\label{eq:lap_to_grad}
            K_\rho*(\Phi-\Delta^\rho_0 \Phi)
            = \Fcal^{-1} \left(\frac{\lambda_\rho}{1+|\lambda_\rho|^2}\cdot\Fcal(\Phi-\Delta^\rho_0 \Phi)\right)
            = \Fcal^{-1}(\lambda_\rho \cdot \widehat{\Phi})
            = \diver^\rho_0 \Phi.
        \end{equation}

        By duality, this also implies the property \(D(\Delta^\rho)\subset H^{\rho,p}(\R^n)\), as we show subsequently. In order to see that, we will need additional properties of \(K_\rho*\Phi\), namely that \(K_\rho*\Phi\in H^{\rho,p}(\R^n)\) and \(K_\rho*\Phi\in D(\Delta^\rho)\) for any \(\Phi\in C_c^\infty(\R^n;\R^n)\). Indeed, \(D^\rho_0 (K_\rho*\Phi)\) is pointwise well defined (see \eqref{eq:pw_equals_dual} and the comments thereafter) and an element of \(L^p(\R^n)\) by Young's convolution inequality and \(K_\rho,Q_\rho\in L^1(\R^n)\), since, with \eqref{eq:commutation} and \eqref{eq:translation},
        \begin{equation*}
            D^\rho_0(K_\rho*\Phi)
            = K_\rho*D^\rho_0 \Phi
            = K_\rho*(Q_\rho * \nabla \Phi),
            \qquad j\in \N.
        \end{equation*}
        A similar computation shows that \(\Delta^\rho K_\rho*\Phi\) exists pointwise and that \(K_\rho*\Phi\in D(\Delta^\rho)\).
        
        This allows us now to apply \eqref{eq:simul_approx_lem} and to conclude that 
         %for a smooth sequence of cut-off functions $\{\sigma_j \}_{j>0}$ obtained via the scalings \(\sigma_j(\cdot)\coloneqq \sigma(\cdot/j)\) for some \(\sigma \in C_c^\infty(\R^n)\) with \(\sigma\equiv 1\) in \(B_1\), we can approximate \(K_\rho* \Phi\) by \((\sigma_j(K_\rho* \Phi))_j\subset C_c^\infty(\R^n)\) such that 
         there exists \((\sigma_j)_j\subset C_c^\infty(\R^n)\) such that \((\sigma_j(K_\rho* \Phi))_j\subset C_c^\infty(\R^n)\) and
        \begin{equation}\label{eq:simul_approx}
            \sigma_j(K_\rho* \Phi)\to K_\rho*\Phi,
            \quad \Delta^\rho_0(\sigma_j(K_\rho* \Phi))\to \Delta^\rho (K_\rho* \Phi)
            \qquad \text{in } L^p(\R^n).
        \end{equation}
       %For the second convergence we use that \(\sigma_j\Delta^{\rho} (K_\rho* \Phi)\to \Delta^{\rho} (K_\rho* \Phi)\) by dominated convergence, in combination with the inequality \(\norm{\nabla \sigma_j}{W^{1,\infty}}\le C/j\) and the estimate from Lemma~\ref{lem:lap_est}, which in our case for the choice \(f= \sigma_j\) and \(g=K_\rho* \Phi\) reads as 
%\begin{equation*}
%		\norm{\Delta_0^{\rho} (\sigma_j (K_\rho* \Phi))-\sigma_j\Delta^{\rho} (K_\rho* \Phi)}{L^p(\R^n)}
%		\le C\norm{\nabla \sigma_j}{W^{1,\infty}(\R^n;\R^n)}\norm{K_\rho* \Phi}{H^{\rho,p}(\R^n)}.
%\end{equation*} 

        We have collected now all ingredients for the final argument.
        Using first the fact that \(K_\rho\) is vector radial together with the convolution identity from \cite[Proposition 4.16]{Bre11} and then later the commutation relation in \eqref{eq:commutation}, it holds 
        \begin{align*}
            &\langle K_\rho*(u-\Delta^\rho u),\Phi \rangle 
            = -\langle u-\Delta^\rho u,K_\rho*\Phi \rangle 
            = -\lim_{j\to \infty}\langle u-\Delta^\rho u,\sigma_j (K_\rho* \Phi) \rangle \\
            &= -\lim_{j\to \infty}\langle u,(\text{Id}-\Delta^\rho_0)(\sigma_j (K_\rho* \Phi)) \rangle 
            = -\langle u,K_\rho*(\Phi-\Delta^\rho_0\Phi) \rangle 
            = -\langle u,\diver^\rho_0 \Phi \rangle.
        \end{align*}
        Furthermore, by Young's convolution inequality and \(K_\rho\in L^1(\R^n)\), we deduce that \(K_\rho*(u-\Delta^\rho u)\in L^p(\R^n)\) whenever \(u\in D(\Delta^\rho)\). It follows then by definition that \(D^\rho u= K_\rho*(u-\Delta^\rho u)\in L^p(\R^n)\) exists, implying that \(u\in H^{\rho,p}(\R^n)\). 
        \end{proof}
    \begin{remark}[Hilbert space case]\label{rem:poincare2_p2} 
        For \(p=2\), the inclusion \(D(\Delta^\rho)\subset H^{\rho,2}(\R^n)\) holds true under the milder conditions that \(\rho\) has compact support and merely satisfies (H0), (H1) and (H3'). In fact, using Fourier multiplier theory in \(L^2(\R^n)\) instead, one does not need to rely on the property that \({\lambda_\rho}/({1+|\lambda_\rho|^2})\) admits a Fourier inverse in \(L^1(\R^n)\). We sketch here the arguments. 
        
        Due to \eqref{eq:lower_prelim}, the assumptions on \(\rho\) assure that \({\lambda_\rho}/({1+|\lambda_\rho|^2})\in L^\infty(\R^n)\) and this implies that this fraction is an \(L^2\)-Fourier multiplier (see \cite[Theorem~2.5.10]{Gr14}). 
        We define
         \begin{equation*}
            R_\rho: L^2(\R^n)\to L^2(\R^n;\R^n), 
            \quad u\mapsto \Fcal^{-1} \left(\frac{\lambda_\rho}{1+|\lambda_\rho|^2}\hat u\right)
        \end{equation*}
        and extend it to a map \(L^2(\R^n;\R^n)\to L^2(\R^n)\) by inserting the scalar product in Fourier space. It follows from Parseval's identity that the adjoint operator of \(R_\rho\) is given by \(- R_\rho\). Now, for any \(\Phi\in C_c^\infty(\R^n;\R^n)\) the functions
        \begin{equation}\label{eq:Finverses_aux}
            \Fcal^{-1}\left(\frac{\lambda_\rho\otimes \lambda_\rho}{1+|\lambda_\rho|^2} \cdot\hat{\Phi}\right),
            \quad -\Fcal^{-1}\left(\frac{\lambda_\rho|\lambda_\rho|^2}{1+|\lambda_\rho|^2} \cdot\hat{\Phi}\right)
        \end{equation} 
        lie in \(L^2(\R^n;\R^n)\) and \(L^2(\R^n)\), respectively, because \({\lambda_\rho\otimes \lambda_\rho}/({1+|\lambda_\rho|^2})\), \({\lambda_\rho|\lambda_\rho|^2}/({1+|\lambda_\rho|^2})\) have polynomial growth (which follows again from \eqref{eq:lower_prelim} and the fact that \(\lambda_\rho\) has linear growth as shown in Lemma~\ref{lem:lambda_FintoF_pre}, where the arguments merely use \(\rho\in L^1(\R^n)\)). Additionally, the expressions in \eqref{eq:Finverses_aux} are equal to \(D^\rho R_\rho \Phi\) and \(\Delta^\rho R_\rho \Phi\), respectively, as can be quickly verified using Parseval's identity. We conclude from this that \(R_\rho \Phi\in H^{\rho, 2}(\R^n)\) and, with similar arguments as in \eqref{eq:lap_to_grad}, that \((1-\Delta^\rho) R_\rho \Phi=\diver^\rho_0 \Phi\).  Hence, again by a suitable approximation and \eqref{eq:simul_approx_lem},
        \begin{align*}
            \langle R_\rho(u-\Delta^\rho u),\Phi \rangle 
            =-\langle u-\Delta^\rho u, R_\rho\Phi \rangle 
            = -\langle u,(\text{Id}-\Delta^\rho)R_\rho\Phi \rangle 
            = \langle u,\diver^\rho_0 \Phi \rangle.
        \end{align*}
        Thus, we have shown that \(D^\rho u= R_\rho(u-\Delta^\rho u)\in L^2(\R^n;\R^n)\) and \(u\in H^{\rho,2}(\R^n)\).   
    \end{remark}

    \begin{remark}[Refinement of a nonlocal fundamental theorem of calculus]\label{rmk:ftoc}
        We comment briefly on \cite[Theorem~5.2]{BMCS25} where the authors show a nonlocal analogue of the classical fundamental theorem of calculus. In particular, they establish that for compactly supported \(\rho\) subject to (H0)--(H4) (where (H4) is an additional hypothesis that they impose) there exists a function \(V_\rho\in C^\infty(\R^n\setminus \{0\};\R^n)\cap L^1_{\text{loc}}(\R^n;\R^n)\) such that 
        \begin{equation*}
            u(x)
            = \int_{\R^n} V_\rho(x-y)\cdot D^\rho_0 u(y)\dl y,
            \qquad x\in \R^n, u\in C_c^\infty(\R^n).
        \end{equation*}
        Later in \cite[Corollary~5.3]{BMCS25}, this is extended to \(u\in H_0^{\rho,p}(\Om)\) for open and bounded \(\Om\subset \R^n\). 

        Inspecting the proof of \cite[Theorem~5.2]{BMCS25} one finds that hypothesis~(H4) is only used to derive that an expression similar to 
        \begin{equation*}
            \frac{1}{|x|^{n-1}}
			\int_{B_{1/|x|} \setminus B_{1}}
			\frac{1}{|\xi|^{n}\widehat{Q_\rho}(\xi)}\dl \xi,
            \qquad x\in B_1
        \end{equation*}
        is finite and integrable over \(B_1\). From \eqref{eq:lest_qrho_pre} (which does not require (H4)) and the computations in \eqref{eq:lg_gen_fin_est_1} it follows now that 
        \begin{equation*}
            \frac{1}{|x|^{n-1}}
			\int_{B_{1/|x|} \setminus B_{1}}
			\frac{1}{|\xi|^{n}\widehat{Q_\rho}(\xi)}\dl \xi
            \le C\frac{1}{|x|^{n-1}}
			\int_{B_{1/|x|}}
			\frac{1}{|\xi|^{n-1+\sigma}}\dl \xi
			= C
			\frac{1}{|x|^{n-\sigma}},
            \qquad x\in B_1,
        \end{equation*}
        where, as desired, the right-hand side is integrable. Consequently, the nonlocal fundamental theorem of calculus presented in \cite[Theorem~5.2]{BMCS25} and \cite[Corollary~5.3]{BMCS25} is valid for the larger class of kernels being compactly supported and only constrained by hypotheses (H0)--(H3).
    \end{remark}

   % The result of Proposition~\ref{prop:poincare2} allows to find simultaneous approximations as in \eqref{eq:simul_approx} \commentCK{To avoid repeating the same argument twice one could move this up.} for mollifications of functions in \(D(\Delta^\rho)\), and to prove that $C^{\infty}_c (\R^n)$ is a core for $\Delta^{\rho}$. 
   % \comment{
    The result of Proposition~\ref{prop:poincare2} allows us to prove that $C^{\infty}_c (\R^n)$ is a core for $\Delta^{\rho}$. 
    %}

    \begin{lemma}[Core property]\label{le:corep}
    Let \(p\in (1,\infty)\) and assume that \(\rho\) has compact support and satisfies (H0)--(H3). For each $u \in D (\Delta^{\rho})$ there exists a sequence $(u_j)_j$ in $C^{\infty}_c (\R^n)$ such that 
    \begin{equation*}
        u_j \to u ,\quad \Delta_0^{\rho} u_j \to \Delta^{\rho}  u 
         \qquad \text{in } L^p(\R^n).
    \end{equation*}
\end{lemma}
\begin{proof}
Let $u \in D(\Delta^{\rho})$.
Let $( \eta_{\veps})_{\veps>0}\subset C_c^\infty(\R^n)$ be a standard family of mollifiers.
Then $\eta_{\veps} * u \to u$ in $L^p (\R^n)$ for \(\veps\to 0\). Thanks to \eqref{eq:commutation}, we have $\eta_{\veps} * \Delta_0^{\rho} \varphi = \Delta_0^{\rho} (\eta_{\veps} * \varphi)$ for any \(\varphi\in C_c^\infty(\R^n)\).
%As proved in \cite[Lemma~3.6]{BMCS25}, it holds $\eta_{\veps} * D^{\rho}_0 \varphi = D^{\rho}_0 (\eta_{\veps} * \varphi)$ for all $\varphi \in C^{\infty}_c (\R^n)$, and similarly for the divergence operator \(\diver^\rho_0\). Hence, we find that $\eta_{\veps} * \Delta_0^{\rho} \varphi = \Delta_0^{\rho} (\eta_{\veps} * \varphi)$. 
Together with \cite[Proposition 4.16]{Bre11} and the fact that $\eta_{\veps}$ is radial, we obtain 
\begin{equation*}
     \langle\eta_{\veps} * u, \Delta_0^{\rho} \varphi\rangle 
     = \langle u,\eta_{\veps} * \Delta_0^{\rho} \varphi\rangle 
     = \langle u,\Delta_0^{\rho} (\eta_{\veps} * \varphi)\rangle 
     = \langle \Delta^{\rho} u, \eta_{\veps} * \varphi\rangle
     = \langle\eta_{\veps} * \Delta^{\rho} u, \varphi \rangle
\end{equation*}
for all $\varphi \in C^{\infty}_c (\R^n)$.
Therefore, \(\eta_{\veps} * u\in D(\Delta^\rho)\)   and $\Delta^{\rho} (\eta_{\veps} * u) = \eta_{\veps} * \Delta^{\rho} u$, and consequently, \( \Delta^{\rho} (\eta_{\veps} * u) \to \Delta^{\rho} u\) in \( L^p(\R^n)\) as $\veps\to 0$. Applying Proposition~\ref{prop:poincare2}, further yields \(\eta_{\veps} * u\in H^{\rho,p}(\R^n)\).
%Now let once more $\{ \sigma_j \}_{j>0}$ be the family of smooth functions in \(C_c^\infty(\R^n)\) defined as \(\sigma_j(\cdot)\coloneqq \sigma(\cdot/j)\) for \(\sigma \in C_c^\infty(\R^n)\) with \(\sigma\equiv 1\) in \(B_1\).
%By dominated convergence, for each $\veps>0$, we have $\sigma_j (\eta_{\veps} * u) \to \eta_{\veps} * u$ in $L^p (\R^n)$. Furthermore, with Proposition~\ref{prop:poincare2} we know that \(\eta_{\veps} * u\in H^{\rho,p}(\R^n)\) and from Lemma~\ref{lem:lap_est} we infer with the choice \(f= \sigma_R\) and \(g=\eta_{\veps} * u\) that
%\begin{equation}\label{eq:aux:est_applied}
%		\norm{\Delta_0^{\rho} (\sigma_j (\eta_{\veps} * u))-\sigma_j\Delta^{\rho} (\eta_{\veps} * u)}{L^p(\R^n)}
%		\le C\norm{\nabla \sigma_j}{W^{1,\infty}(\R^n;\R^n)}\norm{\eta_{\veps} * u}{H^{\rho,p}(\R^n)}
%\end{equation}
%for some constant \(C=C(\rho)>0\). Since, by construction, \(\norm{\nabla \sigma_j}%{W^{1,\infty}}\to 0\), the right-hand side vanishes for \(j\to \infty\). Thus, again by dominated convergence, for each \(\veps>0\),
%\begin{equation*}
 %   \lim_{j\to\infty} \Delta_0^{\rho} (\sigma_j (\eta_{\veps} * u))
  %  = \lim_{j\to\infty} \sigma_j\Delta^{\rho} (\eta_{\veps} * u)
   % = \Delta^{\rho} (\eta_{\veps} * u)
   % \qquad \text{in } L^p(\R^n).
%\end{equation*}

Thus, for each fixed $\veps>0$, we infer from \eqref{eq:simul_approx_lem} that there exists \((\sigma_j)_j\subset C_c^\infty(\R^n)\) such that \((\sigma_j(\eta_\veps* u))_j\subset C_c^\infty(\R^n)\) and\begin{equation*}
    \sigma_j (\eta_\veps* u)\to \eta_\veps* u,\quad
    \Delta^\rho_0(\sigma_j (\eta_\veps* u))\to \Delta^\rho (\eta_\veps* u)
    \qquad \text{in } L^p(\R^n).
\end{equation*}

Finally, by choosing a suitable diagonal sequence \((\sigma_{j(\veps)} (\eta_{\veps} * u))_\veps\) depending on \(\veps\), and letting \(\veps\to 0\), we obtain the desired result.
\end{proof}
\begin{remark}\label{rem:core_p2}
 a)   Since in the above proof, hypothesis~(H3) is only used to apply Proposition~\ref{prop:poincare2} to obtain that \(\eta_{\veps} * u\in H^{\rho,p}(\R^n)\), it follows from Remark~\ref{rem:poincare2_p2} that for \(p=2\), Lemma~\ref{le:corep} is true for all kernels \(\rho\) with compact support and fulfilling merely (H0), (H1) and (H3').\smallskip
%\end{remark}
%\begin{remark}\label{rem:sim_approy_smooth}

  b) For \(u\in C^\infty(\R^n)\cap D(\Delta^\rho)\) the proof of Lemma~\ref{le:corep} reduces to applying the cut-off argument from \eqref{eq:simul_approx_lem} which only relies on Lemma~\ref{lem:lap_est}. Thus, the claimed simultaneous approximation exists as soon as hypothesis~(H0) is satisfied and \(\rho\) has compact support.  \smallskip
%\end{remark}
%\begin{remark}\label{rem:simul_approx_pre}

c)
    It was shown in \cite[Theorem~3.9]{BMCS25} that by the very same construction used above in Lemma~\ref{le:corep}, but here only requiring (H0), any function in \(u\in H^{\rho,p}(\R^n)\) can be approximated by a sequence in \(C_c^\infty(\R^n)\). Since additionally \(D(\Delta^{\rho})\subset H^{\rho,p}(\R^n)\) by virtue of Proposition~\ref{prop:poincare2}, this proves that for each $u \in D (\Delta^{\rho})$ there exists a sequence $(u_j)_j\subset C^{\infty}_c (\R^n)$ such that $u_j \to u$, \(D^\rho_0 u_j\to D^\rho u\), and $\Delta_0^{\rho} u_j \to \Delta^{\rho} u$ simultaneously in $L^p$-norm.
\end{remark}

%For the remaining part of this section, we split the exposition into the special instance \(p=2\) and the general case \(p\in (1,\infty)\).

\subsection{The Hilbert space case}\label{sec:laplace_hilbert}
In this section, we will write \(H^{\rho}(\R^n)\) for \(H^{\rho,2}(\R^n)\) to ease the notation. 
We have the following essential properties of the $\rho$-Laplacian $\Delta^{\rho}$ introduced in Definition~\ref{de:Laplacianp} when $p=2$. %\commentCK{Can we make a brief comment relating to the role of closed operators in the Hilbert space case, to make a bridge with the general $p$ case, maybe at the end?... Here or below in 4.3?}

\begin{proposition}[Basic properties of the $\rho$-Laplacian for $p=2$]\label{prop:FourierLaplacian}
Assume that \(\rho\) has compact support and satisfies (H0), (H1) and (H3').
Then \(\Delta^{\rho}\) has the following properties:
\begin{enumerate}[label=(\roman*)]
\item\label{item:LaplacianNegative} $\Delta^{\rho}$ is self-adjoint and dissipative. 

\item\label{item:LaplacianPartial}
 For all $u\in D (\Delta^{\rho})$ and  $v\in H^{\rho}(\R^n)$, it holds that 
\begin{equation}\label{eq:DeltaSymmetric}
 \langle\Delta^{\rho} u, v\rangle =
 -\langle D^\rho u, D^\rho v\rangle.
\end{equation}

\item\label{item:LaplacianFourier} 
    The domain of \(\Delta^{\rho}\) fulfills
    \begin{equation*}
        D(\Delta^{\rho}) = \{u \in L^2 (\R^n) : \abs{\lambda_{\rho}}^2 \widehat{u} \in L^2 (\R^n) \},
    \end{equation*}
   and
\begin{equation}\label{eq:LaplaceFourier}
  \Delta^{\rho} u = \Fcal^{-1}(- \abs{\lambda_{\rho}}^2 \widehat{u}) , \qquad u \in D(\Delta^{\rho}) .
\end{equation}
%\commentCK{Relation with operators $F_{...}$ not mentioned, probably on purpose? Ok.}\comment{Yes, the idea was not to resort to the operator framework in this section, since it is not necessary; the remark below was intended to provide a bridge.}
\end{enumerate}
%\commentCK{Make signs consistent with general $p$-case, $\Delta^\rho$ vs. $-\Delta^\rho$}
\end{proposition}
\begin{proof}
%The proof is given in several steps. \smallskip
\textit{\textit{\ref{item:LaplacianNegative}} and \textit{\ref{item:LaplacianPartial}}.}
%Let $u \in C^{\infty}_c (\R^n)$.
%We have that $D^{\rho} u \in C^{\infty}_c (\R^n)$ and $\Delta_0^{\rho} u \in C^{\infty}_c (\R^n)$.
%In addition, for any $\varphi \in C^{\infty}_c (\R^n)$, 
%\begin{equation}\label{eq:Delta0Symmetric}
% \int_{\R^n} \Delta_0^{\rho} u \varphi = - \int_{\R^n} D_0^{\rho} u \cdot D_0^{\rho} \varphi = \int_{\R^n} u \Delta_0^{\rho} \varphi .
%\end{equation}
%Therefore, $u \in D (\Delta^{\rho})$ and $\Delta^{\rho} u = \Delta_0^{\rho} u$. \commentCK{More Step~1 until here up after definition of $\Delta_\rho$. }
We start by showing that $\Delta^{\rho}$ is self-adjoint.
As shown in \eqref{eq:Delta0Symmetric}, we have \( \langle\Delta_0^{\rho} u, v \rangle= \langle u, \Delta_0^{\rho} v\rangle\) for any \(u,\varphi\in C_c^\infty(\R^n)\). This entails that $\Delta^{\rho}$ is an extension of $\Delta_0^{\rho}$ and that $\Delta_0^{\rho}$ is symmetric. In addition, invoking the core property of Lemma~\ref{le:corep}, we find that \(\Delta^\rho\) is symmetric as well. Let us demonstrate now that $(\Delta^{\rho}_0)^* = \Delta^{\rho}$, where we understand \(\Delta^{\rho}_0\) as an unbounded operator in \(L^2(\R^n)\) with domain \(C_c^\infty(\R^n)\).
%\textit{Step~2: $(\Delta^{\rho}_0)^* = \Delta^{\rho}$.}
The adjoint $(\Delta_0^{\rho})^*$ exists because $C^{\infty}_c (\R^n)$ is dense in $L^2 (\R^n)$.
By \cite[Section 2.4]{CaHa98} and Riesz' representation theorem, a function $u \in L^2 (\R^n)$ belongs to $D((\Delta_0^{\rho})^*)$ if and only if there exists $v \in L^2 (\R^n)$ such that
\[
 \langle u, \Delta_0^{\rho} \varphi\rangle 
 = \langle v, \varphi \rangle , \qquad \varphi \in C^{\infty}_c (\R^n) ,
\]
which, according to Definition \ref{de:Laplacianp}, holds if and only if $u \in D(\Delta^{\rho})$.
Moreover, in this case, $(\Delta_0^{\rho})^* u = v = \Delta^{\rho} u$. Together with the symmetry of \(\Delta^\rho\) it follows \(\Delta^{\rho} \subset (\Delta^{\rho} )^*\subset (\Delta_0^{\rho})^* = \Delta^{\rho}\); in other words, $\Delta^{\rho}$ is self-adjoint.
%\textit{Step~2: proof of \textit{\ref{item:LaplacianNegative}} and \textit{\ref{item:LaplacianPartial}}.}  

We are now in the position to prove \textit{\ref{item:LaplacianPartial}}.
Taking \(u\in D(\Delta^\rho)\), we find with the help of Lemma~\ref{le:corep}, Remark~\ref{rem:core_p2}\,a) and c) a sequence \((u_j)\subset C_c^\infty(\R^n)\) with the property that $u_j \to u$, \(D^\rho_0 u_j\to D^\rho u\), and $\Delta_0^{\rho} u_j \to \Delta^{\rho} u$ in $L^2 (\R^n)$. We conclude that 
\begin{equation*}
    \langle\Delta^\rho u, v \rangle
    = \lim_{j\to \infty} \langle\Delta^\rho_0 u_j, v \rangle
    = \lim_{j\to \infty} \langle\diver^\rho_0 D^\rho_0 u_j, v\rangle
    =-\lim_{j\to \infty} \langle D^\rho_0 u_j, D^\rho v \rangle
    = - \langle D^\rho u, D^\rho v\rangle .
\end{equation*}
This implies \eqref{eq:DeltaSymmetric} and \textit{\ref{item:LaplacianPartial}}.
%\textit{Step~4: }
In order to conclude the proof of \textit{\ref{item:LaplacianNegative}}, it remains to show that \(\Delta^\rho\) is dissipative. 
%Let $u \in D(\Delta^{\rho})$ and let $(u_j)_j$ be as before.
But it follows immediately from \eqref{eq:DeltaSymmetric} that
\begin{equation}\label{eq:negative}
 \langle\Delta_0^{\rho} u_j, u_j \rangle= - \norm{D^{\rho} u_j}{L^2 (\R^n)}^2 \leq 0 .
\end{equation}
%By the convergence properties of $(u_j)_j$ given by Lemma~\ref{le:corep} we obtain that
%\[
% \int_{\R^n} \Delta_0^{\rho} u_j u_j \to \int_{\R^n} \Delta^{\rho} u u ,
%\]
%so
%\[
% \int_{\R^n} \Delta^{\rho} u u \leq 0
%\]
By Lemma~\ref{le:corep} this carries over to general  $u \in D(\Delta^{\rho})$,
and, hence, $\Delta^{\rho}$ is indeed dissipative. \smallskip

%\textit{Step~3: determining \(D(\Delta^{\rho})\) and validity of \eqref{eq:LaplaceFourier}.}
\textit{\textit{\ref{item:LaplacianFourier}}.}  
We define in $L^2 (\R^n)$ the unbounded operator \(A\) with domain \(D(A) = \{u \in L^2 (\R^n) : \abs{\lambda_{\rho}}^2 \widehat{u} \in L^2 (\R^n) \}\) as \( A u \coloneqq\mathcal{F}^{-1} (-\abs{\lambda_{\rho}}^2 \widehat{u})\), \(u \in D(A)\). 
%\[
% D(A) = \{u \in L^2 (\R^n) : \abs{\lambda_{\rho}}^2 \widehat{u} \in L^2 (\R^n) \}
%\]
%as
%\[
% A u = \mathcal{F}^{-1} (-\abs{\lambda_{\rho}}^2 \widehat{u}) , \qquad u \in D(A) .
%\]
%As the Fourier transform is a unitary operator in $L^2 (\R^n)$, the operator $A$ is unitarily equivalent to the multiplication operator $M$ defined as the unbounded operator in $L^2 (\R^n)$ with domain $D (M) = \{ f \in L^2 (\R^n) : \abs{\lambda_{\rho}}^2 f \in L^2 (\R^n)\}$ through
%\[
% (M f) (\xi) = -\abs{\lambda_{\rho}}^2 f (\xi) .
%\]
By properties of the Fourier transform, the operator $A$ is unitarily equivalent to the Fourier multiplication operator with multiplier \(-\abs{\lambda_{\rho}}^2\). 
Hence, as $-\abs{\lambda_{\rho}}^2$ is real, the operator $A$ is self-adjoint (see, e.g., \cite[Section VIII.3, Proposition 1]{ReSi80}).

Now, for each $u \in C^{\infty}_c (\R^n)$, we have that $\hat{u} \in \Scal(\R^n)$, so $\abs{\lambda_{\rho}}^2 \hat{u} \in L^2 (\R^n)$ because \(\lambda_\rho\) has linear growth according to Lemma~\ref{lem:lambda_FintoF_pre}. In addition, \(\Fcal(\Delta_0 ^{\rho} u) = - \abs{\lambda_{\rho}}^2 \widehat{u}\) according to \eqref{eq:FourierLaplacian_pre}.
%using \eqref{eq:nlgrad_symbol}, we find that
%\begin{equation}\label{eq:FourierLaplacian}
% \widehat{\Delta_0 ^{\rho} u} =\lambda_{\rho} \cdot \widehat{D^{\rho} u} = \lambda_{\rho} \cdot \lambda_{\rho} \widehat{u} = - \abs{\lambda_{\rho}}^2 \widehat{u} ,
%\end{equation}
%since $\lambda_{\rho} (\xi)$ is purely imaginary for real $\xi$.
Therefore, $\Delta^{\rho}_0$ and $A$ coincide in $C^{\infty}_c (\R^n)$.
As $A$ is a self-adjoint extension of the symmetric operator $\Delta^{\rho}_0$, we obtain that $(\Delta_0^{\rho})^*$ is an extension of $A$, so, by \textit{Step~1}, $\Delta^{\rho}$ is an extension of $A$.

Finally, we show that $D(\Delta^{\rho}) \subset D (A)$.
%Let $u \in D(\Delta^{\rho})$ and let $(u_j)_j$ be again the sequence of Lemma~\ref{le:corep} and Remark~\ref{rem:core_p2}.
%By \eqref{eq:FourierLaplacian_pre}, $\widehat{\Delta_0 ^{\rho} u_j} = - \abs{\lambda_{\rho}}^2 %\widehat{u_j}$ for each $j\in \N$.
%By Lemma~\ref{le:corep} and the continuity of the Fourier transform in $L^2 (\R^n)$ we have %that, for a subsequence, $\widehat{\Delta_0 ^{\rho} u_j} \to \widehat{\Delta ^{\rho} u}$ in $L^2 (\R^n)$ and a.e., and $\widehat{u_j} \to \widehat{u}$ in $L^2 (\R^n)$ and a.e. Therefore, $\widehat{\Delta ^{\rho} u} = - \abs{\lambda_{\rho}}^2 \widehat{u}$ a.e.\ and $\abs{\lambda_{\rho}}^2 \widehat{u} \in L^2 (\R^n)$.
%Thus, $u \in D(A)$.
%We conclude that $\Delta^{\rho} = A$ as unbounded operators. This proves \ref{item:LaplacianFourier}. 
Since $\Fcal(\Delta_0 ^{\rho} u)  = - \abs{\lambda_{\rho}}^2 \widehat{u}$ for \(u\in C_c^\infty(\R^n)\), it follows for $u \in D(\Delta^{\rho})$ with the approximating sequence provided by Lemma~\ref{le:corep} that $\Fcal(\Delta_0 ^{\rho} u)  = - \abs{\lambda_{\rho}}^2 \widehat{u}$ a.e.\ and, in particular, that $\abs{\lambda_{\rho}}^2 \widehat{u} \in L^2 (\R^n)$.
Thus, $u \in D(A)$.
We conclude that $\Delta^{\rho} = A$ as unbounded operators. This proves \textit{\ref{item:LaplacianFourier}}.  
%\commentC{Unify notation: dissipative / negative semidefinite / negative}
\end{proof}
\begin{remark}
 a)   In the language introduced in Section~\ref{sec:fourier}, the statement of Proposition~\ref{prop:FourierLaplacian}~\textit{\ref{item:LaplacianFourier}} can be equivalently expressed as $\Delta^{\rho} = \ol{F_{-\abs{\lambda_{\rho}}^2}}$. We show below in Proposition~\ref{prop:FourierLaplacianp2} that this identity extends to general \(p\in (1,\infty)\).\smallskip

b) An alternative proof of the self-adjointness of $\Delta^{\rho}$ can also be obtained following the chain of reasoning in \cite[Section 2.6.1]{CaHa98} and arguing via the Lax-Milgram theorem in \(H^{\rho}(\R^n)\).
\end{remark}

In view of~Proposition~\ref{prop:FourierLaplacian}, it follows from a well-known result (see, e.g., \cite[Theorem 13.31]{Rudin91}) that there exists a unique self-adjoint, positive semidefinite, unbounded operator in $L^2(\R^n)$ whose square is $- \Delta^{\rho}$. It will be denoted by $(- \Delta^{\rho})^{1/2}$. Moreover, $(- \Delta^{\rho})^{1/2}$ can be identified with the unbounded Fourier multiplication operator with symbol \(|\lambda_\rho|\). While the latter property can also be derived from the multiplication operator form of the spectral theorem and the associated functional calculus (see, e.g., \cite{ReSi80}), we provide, for the reader's convenience, a brief self-contained proof of this fact. %\comment{Headline added similar to general \(p\) case.}

%By general theory (e.g., \cite[Theorem 13.31]{Rudin91} or \cite[Chapter 5, Theorem 3.35]{Kato}), there exists a unique $m$-accretive square root of $- \Delta^{\rho}$, which is also selfadjoint and positive semidefinite. It is denoted by $(- \Delta^{\rho})^{1/2}$.

\begin{lemma}[Square root of $-\Delta^\rho$]\label{le:sqrtLaplacian}
Assume that \(\rho\) has compact support and satisfies (H0), (H1) and (H3').
Then
\[
  D((- \Delta^{\rho})^{1/2}) = \{u \in L^2 (\R^n) : \abs{\lambda_{\rho}} \widehat{u} \in L^2 (\R^n) \}
\]
and
\[
 (- \Delta^{\rho})^{1/2} u = \mathcal{F}^{-1} (\abs{\lambda_{\rho}} \widehat{u}) , \qquad u \in D((- \Delta^{\rho})^{1/2}) .
\]
\end{lemma}
\begin{proof}
We define the operator $B$ in $L^2 (\R^n)$ with domain $D(B) = \{u \in L^2 (\R^n) : \abs{\lambda_{\rho}} \widehat{u} \in L^2 (\R^n) \}$
by \( B u \coloneqq \mathcal{F}^{-1} (\abs{\lambda_{\rho}} \widehat{u})\), \(u \in L^2 (\R^n)\), which is self-adjoint because \(\abs{\lambda_{\rho}}\) is real (see, e.g., \cite[Section VIII.3, Proposition 1]{ReSi80}).
Since
\[
 \langle B u, u \rangle = \langle \abs{\lambda_{\rho}} \widehat{u}, \ol{\, \widehat{u}\, } \rangle =\int_{\R^n} \abs{\lambda_{\rho} (\xi)} \abs{\widehat{u} (\xi)}^2 d \xi \geq 0,
 \qquad u \in D(B),
\]
$B$ is positive semidefinite.
Its square $B^2$ is defined in
\[
 D (B^2) = \{ u \in D(B) : Bu \in D(B) \} = \{ u \in L^2 (\R^n) : \, \abs{\lambda_{\rho}} \widehat{u} \in L^2 (\R^n) , \, \abs{\lambda_{\rho}}^2\widehat{u} \in L^2 (\R^n) \} .
\]
Thanks to H\"older's inequality, for any $f \in L^2 (\R^n)$ with $\abs{\lambda_{\rho}}^2 f \in L^2 (\R^n)$ we have
\[
 \norm{\abs{\lambda_{\rho}} f}{L^2 (\R^n)} \leq \norm{\abs{f}^{1/2}}{L^4 (\R^n)} \norm{\abs{\lambda_{\rho}} \abs{f}^{1/2}}{L^4 (\R^n)} = \norm{f}{L^2 (\R^n)}^{1/2} \norm{\abs{\lambda_{\rho}}^2 f}{L^2 (\R^n)}^{1/2} .
\]
Therefore,
\[
D (B^2) = \{ u \in L^2 (\R^n) :  \abs{\lambda_{\rho}}^2\widehat{u} \in L^2 (\R^n) \} .
\]
This operator acts, for each $u \in D(B^2)$, as
\[
 B^2 u = B (\mathcal{F}^{-1} (\abs{\lambda_{\rho}} \widehat{u})) = \mathcal{F}^{-1} (\abs{\lambda_{\rho}}^2 \widehat{u} ) .
\]
By Proposition~\ref{prop:FourierLaplacian}~\textit{\ref{item:LaplacianFourier}}, $B^2 = -\Delta^{\rho}$.
Because $B$ is self-adjoint and positive semidefinite according to Proposition~\ref{prop:FourierLaplacian}~\textit{\ref{item:LaplacianNegative}}, by uniqueness (see, e.g., \cite[Theorem 13.31]{Rudin91}), $B = (- \Delta^{\rho})^{1/2}$.
\end{proof}
We have prepared the following characterization of \(H^{\rho}(\R^n)\) in terms of the 
square root of the generalized $\rho$-Laplace operator. Our proof is based on the Fourier characterization of \(H^{\rho,p}(\R^n)\) obtained in Theorem~\ref{thm:alt_char}
\begin{proposition}\label{prop:charac_lapace2}
Assume that \(\rho\) has compact support and satisfies (H0)--(H2) and (H3').
Then 
\[
H^{\rho}(\R^n)
= \{u\in L^2(\R^n):\,  |\lambda_\rho|\hat{u} \in L^2(\R^n)\} 
= D ((-\Delta^{\rho})^{1/2}) .
\]
%\commentCK{Consitent use of sign in $\rho$-Laplacian in $p=2$ and general $p$-case.}
\end{proposition}
\begin{proof}
Let $u \in H^{\rho}(\R^n)$.
Then $u\in L^2(\R^n)$ and, by Theorem \ref{thm:alt_char}, $\mathcal{F}^{-1} (\bla{\lambda_\rho}\hat{u})\in L^2(\R^n)$.
        By Proposition \ref{prop:fourier_multipliers}~\textit{(i)}, $\mathcal{F}^{-1}(|\lambda_\rho|\hat{u})\in L^2(\R^n)$ and, hence, $|\lambda_\rho|\hat{u}\in L^2(\R^n)$.
Conversely, let $u\in L^2(\R^n)$ satisfy $|\lambda_\rho|\hat{u}\in L^2(\R^n)$.
Then, $\mathcal{F}^{-1}(|\lambda_\rho|\hat{u})\in L^2(\R^n)$.
As $u\in L^2(\R^n)$, then $\mathcal{F}^{-1}((1+|\lambda_\rho|)\hat{u})\in L^2(\R^n)$. By Proposition~\ref{prop:fourier_multipliers}~\textit{(ii)}, $\mathcal{F}^{-1}(\bla{\lambda_\rho}\hat{u})\in L^2(\R^n)$.
The first equality of the statement is proved; the second follows from Lemma~\ref{le:sqrtLaplacian}.
\end{proof}

\subsection{The general-$p$ case}\label{sec:laplace_general}

The theory of square roots of unbounded operators in Banach spaces is considerably more involved than in Hilbert spaces (see, e.g., \cite{Ha06}) and lies beyond the scope of this work.
Yet, we can define the $\rho$-Laplacian and its square root in $L^p (\R^n)$ in a rather canonical way.
For the rest of the section, we fix \(p\in (1,\infty)\) and, as for \(p=2\) in the previous section, turn to  the problem of representing  $H^{\rho, p}(\R^n)$ in terms of the domain of \((-\Delta^\rho)^{1/2}\). 
%via its Fourier characterization. \color{black}

We remind the reader of the terminology and basic properties introduced in Section~\ref{sec:fourier_mult} of unbounded operators in \(L^p(\R^n)\) defined through a multiplier in Fourier space and closable unbounded operators (especially, Definition~\ref{def:ub_mult_op_pl} and the results thereafter).
Also recall that, as stated in Proposition~\ref{prop:fourier_multipliers}, for $m$ being any of the four functions \(\abs{\lambda_\rho}/\bla{\lambda_\rho}, \bla{\lambda_\rho}/(1+\abs{\lambda_\rho})\), \(1/\bla{\lambda_\rho}\) or \(1/\bla{\lambda_\rho}^2\), the corresponding closure $\ol{F_m}$ is a bounded linear operator in $L^p (\R^n)$. 

The following proposition is the general counterpart of Proposition~\ref{prop:FourierLaplacian} and collects again the main properties of the \(\rho\)-Laplacian. It varies from the result in the Hilbert space case, since there is no concept of self-adjointness in $L^p (\R^n)$.

%\commentCK{Again $C_c^\infty(\R^n)$ vs. $\Scal(\R^n)$, discuss.}

\begin{proposition}[Basic properties of the $\rho$-Laplacian]\label{prop:FourierLaplacianp2}
Assume that \(\rho\) has compact support and satisfies (H0)--(H3).
The operator \(\Delta^{\rho}\) has the following properties:
\begin{enumerate}[label=(\roman*)]
\item\label{item:Laplacianclosurep}  
The closure of $\Delta_0^{\rho}$ in $L^p (\R^n)$ is $\Delta^{\rho}$.
\item\label{item:Laplacianpartialp} 
For all $u\in D (\Delta^{\rho})$ in the sense of \(L^p(\R^n)\) and $v\in D (\Delta^{\rho})$ in the sense of \(L^{p'}(\R^n)\) it holds
\begin{equation}\label{eq:DeltaSymmetricp}
 \langle \Delta^{\rho} u, v\rangle 
 = \langle u, \Delta^{\rho} v\rangle.
\end{equation}
\item\label{item:LaplacianFourierp} 
    It holds $\Delta^{\rho} = \ol{F_{-\abs{\lambda_{\rho}}^2}}$ in \(L^p(\R^n)\), in particular,
\begin{equation}\label{eq:LaplaceFourierp}
     %\widehat{\Delta^{\rho} u}
     %= - \abs{\lambda_{\rho}}^2 \widehat{u} , \qquad u \in \mathcal{S}(\R^n) .
     \Delta^{\rho} u
     =\mathcal{F}^{-1}(-\abs{\lambda_\rho}^2\hat{u}),
     \qquad u \in \mathcal{S}(\R^n).
\end{equation}
%\commentCK{Mention $p$?? e.g. add in $L^p(\R^n)$}\comment{yes, added.} \commentCK{To make analogous to $p=2$, maybe $\Delta^{\rho} u     =\Fcal^{-1}(- \abs{\lambda_{\rho}}^2 \widehat{u}) $ in last line..}\comment{I agree, added.}
%\commentCK{Make presentation consistent with Hilbert space case.}
\end{enumerate}
\end{proposition}
\begin{proof}
%The proof is again presented in multiple steps.\smallskip
\textit{\ref{item:Laplacianclosurep}.}
We show first that \(\Delta_0^{\rho}\) is closable as an operator in \(L^p(\R^n)\).
Let $(u_j)_j$ be a sequence in $C^{\infty}_c (\R^n)$ such that $u_j \to 0$ in $L^p (\R^n)$ and $\Delta_0^{\rho} u_j \to v$ in $L^p (\R^n)$ for some $v \in L^p (\R^n)$.
Then, for all $\varphi \in C^{\infty}_c (\R^n)$,
\[
 \langle \Delta_0^{\rho} u_j, \varphi\rangle \to\langle v, \varphi\rangle
 \quad\text{ and }\quad 
    \langle\Delta_0^{\rho} u_j, \varphi \rangle
 = \langle u_j, \Delta_0^{\rho} \varphi\rangle\to 0,
\]
so \(\langle v, \varphi \rangle = 0\), which implies that $v=0$.

It remains to verify that the closure of $\Delta_0^{\rho}$ in $L^p (\R^n)$ is $\Delta^{\rho}$.
Let therefore $\overline{\Delta_0^{\rho}}$ denote the closure of $\Delta_0^\rho$ in $L^p(\R^n)$. We start by proving that $ D(\overline{\Delta^{\rho}_0})\subset D(\Delta^{\rho})$ and that \(\Delta^{\rho}\) extends \(\overline{\Delta^{\rho}_0}\). Let $(u_j)_j$ be a sequence in $C^{\infty}_c (\R^n)$ such that $u_j \to u$ in $L^p (\R^n)$ and $\Delta^{\rho}_0 u_j \to w$ in $L^p (\R^n)$ for some $u, w \in L^p (\R^n)$.
Then
\[
 \langle u, \Delta_0^{\rho} \varphi\rangle
 = \lim_{j \to \infty} \langle u_j, \Delta_0^{\rho} \varphi \rangle
 = \lim_{j \to \infty} \langle \Delta_0^{\rho} u_j, \varphi \rangle
 =\langle w, \varphi\rangle ,
 \qquad \varphi \in C^{\infty}_c (\R^n),
\]
so $u \in D (\Delta^{\rho})$ and $\Delta^{\rho} u = w = \overline{\Delta^{\rho}_0} u$.
Thus, $\Delta^{\rho}$ is an extension of $\overline{{\Delta}^{\rho}_0}$.

The inclusion $D(\Delta^{\rho}) \subset D(\overline{\Delta^{\rho}_0})$ and the fact that \(\overline{\Delta^{\rho}_0}\) extends \(\Delta^{\rho}\) follow directly from the statement of Lemma~\ref{le:corep}. 
This proves part~\textit{\ref{item:Laplacianclosurep}}.\smallskip

\textit{\ref{item:Laplacianpartialp}.}
Let $u\in D (\Delta^{\rho})$ in the sense of \(L^p(\R^n)\) and $v\in D (\Delta^{\rho})$ in the sense of \(L^{p'}(\R^n)\). By the core property in Lemma~\ref{le:corep}, there exists a sequence \((v_j)_j\subset C_c^\infty(\R^n)\) such that \(v_j\to v\), \(\Delta^\rho_0 v_j\to \Delta^\rho v\) in \(L^{p'}(\R^n)\). It follows that 
\begin{equation*}
    \langle \Delta^\rho u, v  \rangle
    = \lim_{j \to \infty} \langle \Delta^\rho u, v_j \rangle  
    = \lim_{j \to \infty} \langle  u, \Delta_0^{\rho}v_j \rangle
    =\langle  u, \Delta^\rho v\rangle.
\end{equation*}
This proves \eqref{eq:DeltaSymmetricp} and thus, \textit{\ref{item:Laplacianpartialp}}.\smallskip

\textit{\ref{item:LaplacianFourierp}.}
Note that the closure $\ol{F_{-\abs{\lambda_{\rho}}^2}} = -\ol{F_{\abs{\lambda_{\rho}}^2}} $ is well defined thanks to Lemma~\ref{lem:closability_sr}, which also implies that $D(F_{-\abs{\lambda_{\rho}}^2}) =\Scal(\R^n)$. Comparing this with \eqref{eq:FourierLaplacian_pre}, we observe that $\Delta^{\rho}_0$ and $F_{-\abs{\lambda_{\rho}}^2}$ coincide in $C^{\infty}_c (\R^n)$.
Since, as was shown in \textit{\ref{item:Laplacianclosurep}}, the closure of $\Delta^{\rho}_0$ in \(L^p(\R^n)\) with respect to \(C_c^\infty(\R^n)\) is exactly $\Delta^{\rho}$, it follows with the inclusion \(C_c^\infty(\R^n)\subset \Scal(\R^n)\) that \(D(\Delta^{\rho})\subset D(\ol{F_{-\abs{\lambda_{\rho}}^2}})\) and that $\ol{F_{-\abs{\lambda_{\rho}}^2}}$ extends \(\Delta^{\rho}\).

To see the reverse connection, let \((u_j)_j\subset \Scal(\R^n)\) and \(u\in L^p(\R^n)\) such that \(u_j\to u\) and \(F_{-\abs{\lambda_{\rho}}^2} u_j\to \ol{F_{-\abs{\lambda_{\rho}}^2}} u\) in \(L^p(\R^n)\). By identities \eqref{eq:pw_equals_dual} and \eqref{eq:nlgrad_symbol} and the comments thereafter, it holds \eqref{eq:LaplaceFourierp} and \(\Delta^{\rho} u_j = \diver^\rho _0 D^\rho_0 u_j = F_{-\abs{\lambda_{\rho}}^2} u_j\), and, similarly to the proof of \textit{\ref{item:Laplacianclosurep}},
\begin{equation*}
     \langle u, \Delta_0^{\rho} \varphi \rangle
     = \lim_{j \to \infty} \langle u_j, \Delta_0^{\rho} \varphi \rangle
     = \lim_{j \to \infty} \langle \Delta^{\rho} u_j, \varphi \rangle 
     = \lim_{j \to \infty} \langle F_{-\abs{\lambda_{\rho}}^2} u_j,\varphi \rangle
     =\langle \ol{F_{-\abs{\lambda_{\rho}}^2}} u, \varphi\rangle, \quad \varphi\in C_c^\infty(\R^n).
\end{equation*}
This entails that \(\Delta^\rho u\) exists in \(L^p(\R^n)\) and simultaneously that \(D(\ol{F_{-\abs{\lambda_{\rho}}^2}})\subset D(\Delta^{\rho})\) with $\Delta^{\rho} = \ol{F_{-\abs{\lambda_{\rho}}^2}}$ on \(D(\ol{F_{-\abs{\lambda_{\rho}}^2}})\). We have thus demonstrated \eqref{eq:LaplaceFourierp} and \textit{\ref{item:LaplacianFourierp}}, where we note that \eqref{eq:LaplaceFourierp} alone already follows from \eqref{eq:rholap_S}.
\end{proof}

In what follows, our goal is to show that \(\ol{F_{\abs{\lambda_\rho}}}\) defines a square root of \(-\Delta^{\rho}\), generalizing the result of Lemma~\ref{le:sqrtLaplacian}. The proof is more delicate than in the Hilbert space case, and we will establish several auxiliary results building up to it. With Proposition~\ref{prop:FourierLaplacianp2}~\textit{\ref{item:LaplacianFourierp}} in mind, it will suffice to establish that \(\ol{F_{\abs{\lambda_\rho}}}^2= \ol{F_{\abs{\lambda_\rho}^2}}\) as operators in \(L^p(\R^n)\). 
%\commentCK{Switch the previous and the following sentence?} 

First, we will show the representation of the nonlocal Sobolev space \(H^{\rho,p}(\R^n)\) as the domain of the operators \(\ol{F_{\bla{\lambda_\rho}}}\) and \(\ol{F_{\abs{\lambda_\rho}}}\), respectively. The result is based on the Fourier characterization in Theorem~\ref{thm:alt_char} and the density of smooth compactly supported functions in $ H^{\rho,p}(\R^n)$.
%\comment{Change proof strategy: Keep Step~1, but show in Step~2 directly the identity between the domains using Fourier multiplier property, maybe in a separate lemma.}
\begin{lemma}\label{lem:hrhop_asdomains}
    Assume that \(\rho\) has compact support and fulfills hypotheses (H0)--(H2) and (H3'). Then it holds
    \begin{equation*}
       H^{\rho,p}(\R^n)
       =  D(\ol{F_{\bla{\lambda_\rho}}})
        = D(\ol{F_{\abs{\lambda_\rho}}}). 
        %\color{blue} \quad \text{or maybe?? Which do we like better?}   H^{\rho,p}(\R^n)= D(\ol{F_{\bla{\lambda_\rho}}})
        %= D(\ol{F_{\abs{\lambda_\rho}}}). 
    \end{equation*}
\end{lemma}
\begin{proof}
The proof will be given in two steps, corresponding to the two claimed identities. \smallskip

\textit{Step~1.}
We start by showing that \(H^{\rho,p}(\R^n) =  D(\ol{F_{\bla{\lambda_\rho}}})\).
We use that \(H^{\rho,p}(\R^n)\) in the characterization of Theorem~\ref{thm:alt_char} can be written as 
\begin{align}\label{eq:hrhop_reform}
H^{\rho,p}(\R^n)
= \{u\in L^p(\R^n):\, \mathcal{F}^{-1}(\bla{\lambda_\rho}\hat{u})\in L^p(\R^n)\}, 
\end{align}
where \(\mathcal{F}^{-1}(\bla{\lambda_\rho}\hat{u})\) is understood in the sense of \(\mathcal{S}'(\R^n)\). 
In view of ~\eqref{eq:hrhop_reform}, the inclusion \(H^{\rho,p}(\R^n) \subset  D(\ol{F_{\bla{\lambda_\rho}}})\) follows from 
%That~$H^{\rho, p}(\R^n)$ coincides with $D(\ol{F_{\bla\lambda_\rho}})$ follows from 
the fact that \(C_c^\infty(\R^n)\) is dense in \(H^{\rho,p}(\R^n)\) (see \cite[Theorem~3.9]{BMCS25}). Indeed, the density result and the equivalent norm characterization of Theorem~\ref{thm:alt_char} imply that for every \(u\in H^{\rho,p}(\R^n)\), there exists \((u_j)_j\subset C_c^\infty(\R^n)\subset \Scal(\R^n)\) such that 
\begin{equation*}
     F_{\bla{\lambda_\rho}} u_j= \mathcal{F}^{-1}(\bla{\lambda_\rho}\widehat{u_j})
  \to  \mathcal{F}^{-1}(\bla{\lambda_\rho}\hat{u}) \quad \text{in $L^p(\R^n)$},
\end{equation*}
showing that \(u\in D(\ol{F_{\bla\lambda_\rho}})\). The reverse inclusion is a mere consequence of the property that the tempered distribution \(\mathcal{F}^{-1}(\bla{\lambda_\rho}\hat{u})\) coincides with the function \(\ol{F_{\bla{\lambda_\rho}}} u\), whenever \(u\in D(\ol{F_{\bla{\lambda_\rho}}})\), following from Lemma~\ref{le:Fmclosable} for \(m=\bla{\lambda_\rho}\). \smallskip
%\comment{We could also add something like this to Section~\ref{sec:fm_theory}. I included something there after the Notation.}\smallskip

\textit{Step~2.} We will show in this step that 
\begin{align}\label{aux31}
D(\ol{F_{\bla{\lambda_\rho}}})
        = D(\ol{F_{\abs{\lambda_\rho}}}),
        \end{align}
        which will follow essentially from the properties of the $L^p$-Fourier multipliers in Proposition~\ref{prop:fourier_multipliers}. 
To be precise, by applying Lemma~\ref{lem:domain_mulitplier}  twice, first with $m_1=\bla{\lambda_\rho}$ and $m_2=|\lambda_\rho|/\bla{\lambda_\rho}$ (see Proposition~\ref{prop:fourier_multipliers}\, \textit{(i)}), and then with $m_1=1+|\lambda_\rho|$ and $m_2=\bla{\lambda_\rho}/(1+|\lambda_\rho|)$ (see Proposition~\ref{prop:fourier_multipliers}~\textit{(ii)}),  we conclude that
\begin{align}\label{eq:double_inclusion}
D(\overline{F_{\bla{\lambda_\rho}}})\subset D(\overline{F_{|\lambda_\rho|}}),
\quad D(\overline{F_{1+|\lambda_\rho|}})\subset D(\overline{F_{\bla{\lambda_\rho}}}).
\end{align}
All that remains is to prove the inclusion $D(\ol{F_{|\lambda_\rho|}})\subset D(\ol{F_{1+|\lambda_\rho|}})$, and we will see that this is a simple consequence of the notion of domain. Indeed, let $u\in D(\ol{F_{|\lambda_\rho|}})$, implying, by definition of \(\ol{F_{\abs{\lambda_\rho}}}\), that there exists \((u_j)_j\subset \mathcal{S}(\R^n)\) such that \(u_j\to u\) and \(F_{|\lambda_\rho|}u_j\to \ol{F_{\abs{\lambda_\rho}}} u\) in \(L^p(\R^n)\). Together, it holds
\begin{align}
F_{1+|\lambda_\rho|}u_j = u_j + F_{|\lambda_\rho|}u_j\to u + F_{|\lambda_\rho|}u =\ol{F_{1+|\lambda_\rho|}}u \quad \text{ in $L^p(\R^n)$},
\end{align}
showing $u\in D(\ol{F_{1+|\lambda_\rho|}})$. 
In combination with \eqref{eq:double_inclusion}, this yields~\eqref{aux31}, as desired. 
\end{proof}

We prove below that the unbounded operators \(\ol{F_{\abs{\lambda_\rho}}}\) and \(\ol{F_{\bla{\lambda_\rho}}}\) commute with each other.
\begin{lemma}\label{lem:comm_relation}
Assume that \(\rho\) has compact support and fulfills hypotheses (H0)--(H2) and (H3'). Then
\begin{equation*}
    \ol{F_{\bla{\lambda_\rho}}} \, \ol{F_{\abs{\lambda_\rho}}}
        = \ol{F_{\abs{\lambda_\rho}}}\, \ol{ F_{\bla{\lambda_\rho}}}
\end{equation*}
as operators in \(L^p(\R^n)\).
\end{lemma}
\begin{proof}
\textit{Step~1.}
    We first prove the commutation relation on the level of functions in \(\Scal(\R^n)\). Precisely, we claim that it holds the identity
    \begin{equation}\label{eq:comm_relation}
        \ol{F_{\bla{\lambda_\rho}}} F_{\abs{\lambda_\rho}} u
        = F_{\abs{\lambda_\rho}} F_{\bla{\lambda_\rho}} u,
        \qquad u\in \Scal(\R^n).
    \end{equation}
    and, in particular, that the left-hand side is always well defined as an element in \(L^p(\R^n)\). Note that the right-hand side exists due to Lemma~\ref{lem:closability_sr}, which shows that \(F_{\bla{\lambda_\rho}}\) maps \(\Scal(\R^n)\) into itself and that \(D(F_{\abs{\lambda_\rho}}) =\Scal(\R^n)\). To obtain \eqref{eq:comm_relation}, we compute, for \(u\in \Scal(\R^n)\),
    \begin{equation*}
        \langle F_{\abs{\lambda_\rho}} F_{\bla{\lambda_\rho}} u, \varphi\rangle 
        = \langle \bla{\lambda_\rho} \hat{u}, \abs{\lambda_\rho}\Fcal^{-1}\varphi\rangle%_{\Scal'\times \Scal}
        =  \langle \abs{\lambda_\rho} \hat{u}, \bla{\lambda_\rho}\Fcal^{-1}\varphi\rangle%_{\Scal'\times \Scal},
        ,\qquad \varphi\in \Scal(\R^n),
    \end{equation*}
    and, by definition, the right-hand side equals \(\Fcal^{-1}(\bla{\lambda_\rho} F_{\abs{\lambda_\rho}} u)\) as a tempered distribution, i.e., it lies in \(H^{\rho,p}(\R^n)\) by virtue of Theorem~\ref{thm:alt_char}. With Lemma~\ref{lem:hrhop_asdomains}, we have that \(F_{\abs{\lambda_\rho}} u\in D(\ol{F_{\bla{\lambda_\rho}}})\), and equality \eqref{eq:comm_relation} is then a direct consequence of the fact that, as a result of Lemma~\ref{le:Fmclosable} for \(m=\bla{\lambda_\rho}\), \(\ol{F_{\bla{\lambda_\rho}}} v\) agrees with \(\Fcal^{-1}(\bla{\lambda_\rho} \hat{v})\) for any \(v\in D(\ol{F_{\bla{\lambda_\rho}}})\).\smallskip

    \textit{Step~2.}
    In this step we are concerned with deriving the inclusion 
    \begin{equation}\label{eq:211}
        D(\ol{F_{\bla{\lambda_\rho}}} \, \ol{F_{\abs{\lambda_\rho}}})
        \subset D(\ol{F_{\abs{\lambda_\rho}}}\, \ol{ F_{\bla{\lambda_\rho}}})
    \end{equation}
    and the identity
    \begin{equation}\label{eq:212}
     \ol{F_{\abs{\lambda_\rho}}}\, \ol{ F_{\bla{\lambda_\rho}}}u
        =  \ol{F_{\bla{\lambda_\rho}}} \, \ol{F_{\abs{\lambda_\rho}}}u,
        \qquad u\in D(\ol{F_{\bla{\lambda_\rho}}} \, \ol{F_{\abs{\lambda_\rho}}}).
    \end{equation}
     Let therefore \(u\in D(\ol{F_{\bla{\lambda_\rho}}} \, \ol{F_{\abs{\lambda_\rho}}})\) and choose \((u_j)_j\subset \mathcal{S}(\R^n)\) such that \(u_j\to u\) and \( F_{|\lambda_\rho|} u_j\to \ol{F_{|\lambda_\rho|}} u \) in \(L^p(\R^n)\).
   %\begin{align}\label{many_convergences}
   %     u_j\to u, 
   %     \quad F_{|\lambda_\rho|} u_j\to \ol{F_{|\lambda_\rho|}} u 
        %\quad v_j\to \ol{F_{|\lambda_\rho|}} u, 
        %\quad F_{|\lambda_\rho|} v_j\to \ol{F_{|\lambda_\rho|}}^2  u,
   %     \quad \text{in } L^p(\R^n).
   %\end{align} 
   With Lemma~\ref{lem:domain_mulitplier} applied to $m_1=1+|\lambda_\rho|$ and $m_2=\bla{\lambda_\rho}/(1+|\lambda_\rho|)$, it also holds \(u\in D(\ol{F_{\bla{\lambda_\rho}}})\) with \(F_{\bla{\lambda_\rho}} u_j\to \ol{F_{\bla{\lambda_\rho}}} u\) in \(L^p(\R^n)\). 
   Invoking \textit{Step~1}, yields that
    \begin{equation}\label{eq:213}
        F_{\abs{\lambda_\rho}} F_{\bla{\lambda_\rho}} u_j
        = \ol{F_{\bla{\lambda_\rho}}} F_{\abs{\lambda_\rho}} u_j,
        \qquad j\in \N.
    \end{equation}
    Thus, since \(F_{\abs{\lambda_\rho}} u_j\to \ol{F_{|\lambda_\rho|}} u\in D(\ol{F_{\bla{\lambda_\rho}}})\) and \(F_{\bla{\lambda_\rho}} u_j\to \ol{F_{\bla{\lambda_\rho}}} u\), we find in the limit indeed that \(\ol{F_{\bla{\lambda_\rho}}} u\in D(\ol{F_{\abs{\lambda_\rho}}})\) with \(\ol{F_{\abs{\lambda_\rho}}}\, \ol{ F_{\bla{\lambda_\rho}}}u,
        =  \ol{F_{\bla{\lambda_\rho}}} \, \ol{F_{\abs{\lambda_\rho}}}u\). We point out that the regularity of the approximating sequence \((F_{\bla{\lambda_\rho}} u_j)_j\) is not important for this argument to hold; instead, the convergence of the left-hand side in \eqref{eq:213} rests merely upon the fact that the graph of \(\ol{F_{\abs{\lambda_\rho}}}\) is closed by definition of \(\ol{F_{\abs{\lambda_\rho}}}\) being a closed operator. We have hence shown \eqref{eq:211} and \eqref{eq:212} at the same time. \smallskip

    \textit{Step~3.}
    At this point, the overall claim follows as soon as we prove that
    \begin{equation}\label{eq:221}
        D(\ol{F_{\abs{\lambda_\rho}}}\, \ol{ F_{\bla{\lambda_\rho}}})
        \subset D(\ol{F_{\bla{\lambda_\rho}}} \, \ol{F_{\abs{\lambda_\rho}}})
    \end{equation}
    and 
    \begin{equation}\label{eq:222}
       \ol{F_{\bla{\lambda_\rho}}} \, \ol{F_{\abs{\lambda_\rho}}}u
        = \ol{F_{\abs{\lambda_\rho}}}\, \ol{ F_{\bla{\lambda_\rho}}}u,
        \qquad u\in D(\ol{F_{\abs{\lambda_\rho}}}\, \ol{ F_{\bla{\lambda_\rho}}}).
    \end{equation}
    %Precisely, for \(u\in D(\ol{F_{\abs{\lambda_\rho}}}\, \ol{ F_{\bla{\lambda_\rho}}})\), we need to establish that \(u\in D(\ol{F_{\abs{\lambda_\rho}}})\) and \(\ol{F_{\abs{\lambda_\rho}}} u \in D(\ol{F_{\bla{\lambda_\rho}}})\), and 
    The approach is similar to \textit{Step~2}. More precisely, for \(u\in D(\ol{F_{\abs{\lambda_\rho}}}\, \ol{ F_{\bla{\lambda_\rho}}})\) pick \((u_j)_j\subset \mathcal{S}(\R^n)\) with \(u_j\to u\) and \( F_{\bla{\lambda_\rho}} u_j\to \ol{F_{\bla{\lambda_\rho}}} u \) in \(L^p(\R^n)\), infer from Lemma~\ref{lem:domain_mulitplier}, now with $m_1=\bla{\lambda_\rho}$ and $m_2=|\lambda_\rho|/\bla{\lambda_\rho}$, that \(u\in D(\ol{F_{\abs{\lambda_\rho}}})\) with \( F_{|\lambda_\rho|} u_j\to \ol{F_{|\lambda_\rho|}} u \), and use \textit{Step~1} to conclude that \eqref{eq:213} holds true again. In the limit we thus find exactly \eqref{eq:221} and \eqref{eq:222}.
\end{proof}

%\commentCK{Polish this introductory paragraph...}
Next, we establish the identity 
\(\ol{F_{\bla{\lambda_\rho}}}^2=\ol{F_{\bla{\lambda_\rho}^2}}\) as operators in \(L^p(\R^n)\), which is the analogue of the sought equality for \(\bla{\lambda_\rho}\) instead of \(\abs{\lambda_\rho}\). We deduce the requisite core property of the operator \(\ol{F_{\bla{\lambda_\rho}^2}}\) with the help of the \(L^p\)-Fourier multiplier operator \(\ol{F_{1/\bla{\lambda_\rho}^2}}\). 
%Another important ingredient of the proof is that Fourier multipliers defined on tempered distributions are particularly well-behaved in terms of compositions. 
%\commentCK{Check if a result like Lemma~\ref{lem:mults_algebra} for the extensions could be helpful to reduce/shorten a few steps in the following proofs. } 

\begin{lemma}\label{lem:vspace}
    Assume that \(\rho\) has compact support and fulfills hypotheses (H0)--(H2) and (H3'). %\commentCK{Improve phrasing... moreover, don't we show that the two operators coincide as announced before the lemma?}\comment{Ah sorry, yes.I will add that.}
   % Then it holds in \(L^p(\R^n)\) the equality
		%\begin{equation*}
		  %    D (\ol{F_{\bla{\lambda_\rho}}}^2) 
        %     % =X^{\rho,p}(\R^n)
        %      =D(\ol{F_{\bla{\lambda_\rho}^2}}).
		%\end{equation*}
        Then 
        \begin{equation*}
            \ol{F_{\bla{\lambda_\rho}}}^2=\ol{F_{\bla{\lambda_\rho}^2}}
        \end{equation*}
        as operators in \(L^p(\R^n)\).
\end{lemma}
\begin{proof}   
We define, inspired by the Fourier characterization in Theorem~\ref{thm:alt_char}, the auxiliary space
\begin{equation*}
    X^{\rho,p}(\R^n)
    \coloneqq \{u\in L^p(\R^n):\, \Fcal^{-1}(\bla{\lambda_\rho}^2 \hat{u})\in L^p(\R^n)\},
\end{equation*}
%Note that, in general, except in the case of the Riesz fractional gradient and the Bessel potential spaces, \(X^{\rho,p}(\R^n)\) is different from \(H^{\rho^2,p}(\R^n)\).
and we proof in \textit{Steps~1}~and~\textit{2} the statement \(D (\ol{F_{\bla{\lambda_\rho}}}^2) =D(\ol{F_{\bla{\lambda_\rho}^2}})\) by showing the two identities
\begin{equation}\label{eq:first_aux_id}
    D (\ol{F_{\bla{\lambda_\rho}}}^2) = X^{\rho,p}(\R^n)
\end{equation}
and 
\begin{equation}\label{eq:second_aux_id}
    X^{\rho,p}(\R^n)=D(\ol{F_{\bla{\lambda_\rho}^2}}).
\end{equation}
Appealing to Lemma~\ref{lem:closed_ops_equality}  for \(m_1=\bla{\lambda_\rho}=m_2\), this automatically implies that \(\ol{F_{\bla{\lambda_\rho}}}^2=\ol{F_{\bla{\lambda_\rho}^2}}\) as operators in \(L^p(\R^n)\).
%The proof of these identities are the content of the subsequent two steps. 

\textit{Step~1: proof~of~\eqref{eq:first_aux_id}.} 
%We use that \(H^{\rho,p}(\R^n)\subset L^p(\R^n)\) in the characterization of Theorem~\ref{thm:alt_char} can be written as \comment{maybe as a remark after Theorem~\ref{thm:alt_char}?}
%\begin{align}\label{eq:hrhop_reform}
%H^{\rho,p}(\R^n)
%= \{u\in L^p(\R^n):\, \mathcal{F}^{-1}(\bla{\lambda_\rho}\hat{u})\in L^p(\R^n)\}
%= D(\ol{F_{\bla\lambda_\rho}}),
%\end{align}
%where the last equality follows from the fact that \(C_c^\infty(\R^n)\) is dense in \(H^{\rho,p}(\R^n)\) (see \cite[Theorem~3.9]{BMCS25}). In fact, the density result and the equivalent norm characterization of Theorem~\ref{thm:alt_char} imply that for every \(u\in H^{\rho,p}(\R^n)\), there exists \((u_j)_j\subset \Scal(\R^n)\) such that 
%\begin{equation*}
%    \mathcal{F}^{-1}(\bla{\lambda_\rho}\widehat{u_j})
%    = F_{\bla{\lambda_\rho}} u_j\to  \mathcal{F}^{-1}(\bla{\lambda_\rho}\hat{u}),
%\end{equation*}
%showing that \(u\in D(\ol{F_{\bla\lambda_\rho}})\). The reverse inclusion is a mere consequence of the property that the tempered distribution \(\mathcal{F}^{-1}(\bla{\lambda_\rho}\hat{u})\) coincides with the function \(\ol{F_{\bla\lambda_\rho}} u\), whenever \(u\in D(\ol{F_{\bla\lambda_\rho}})\). \comment{We could also add something like this to Section~\ref{sec:fm_theory}.}
%Furthermore, 
We recall that the domain of \(\ol{F_{\bla{\lambda_\rho}}}^2\) is given by
\begin{align*}
    \begin{split}
         D (\ol{F_{\bla{\lambda_\rho}}}^2) 
         = D (\ol{F_{\bla{\lambda_\rho}}} \, \ol{F_{\bla{\lambda_\rho}}}) 
        &= \{ u \in D(\ol{F_{\bla{\lambda_\rho}}}) : \ol{F_{\bla{\lambda_\rho}}}u \in D(\ol{F_{\bla{\lambda_\rho}}}) \},
    \end{split}
\end{align*}
cf.~\eqref{domain_composition}. With the help of Theorem~\ref{thm:alt_char} 
%\commentCK{Maybe refer to comment on subspace of $L^p(\R^n)$.}
%and \eqref{eq:hrhop_reform}
and Lemma~\ref{lem:hrhop_asdomains}, we can reformulate the right-hand side as 
%\commentCK{Do we really need $H^{\rho, p}(\R^n)$ in this reasoning or can it be dropped?}
\begin{align*}
    &\{ u\in D(\ol{F_{\bla{\lambda_\rho}}}) : \ol{F_{\bla{\lambda_\rho}}}u \in D(\ol{F_{\bla{\lambda_\rho}}})\}\\
    &\qquad = \{ u \in H^{\rho,p}(\R^n): \Fcal^{-1}(\bla{\lambda_\rho}\hat{u}) \in H^{\rho,p}(\R^n) \}\\
    &\qquad = \{ u \in L^p(\R^n): \Fcal^{-1}(\bla{\lambda_\rho}\hat{u}), \Fcal^{-1}(\bla{\lambda_\rho}\Fcal(\Fcal^{-1}(\bla{\lambda_\rho}\hat{u})) \in L^p(\R^n) \}\\
    &\qquad = \{ u \in L^p(\R^n): \Fcal^{-1}(\bla{\lambda_\rho}\hat{u}), \Fcal^{-1}(\bla{\lambda_\rho}^2\hat{u}) \in L^p(\R^n)\}.
\end{align*}
Let \(u\in L^p(\R^n)\) be an element of the right-hand side above. 
Then, by \eqref{eq:aux_lp_1} in Lemma~\ref{lem:closability_sr}, in combination with Lemma~\ref{le:Fmclosable} for \(m=\bla{\lambda_\rho}\) and \(m=\bla{\lambda_\rho}^2\), we have that 
\begin{equation*}
    \norm{\Fcal^{-1}(\bla{\lambda_\rho}\hat{u})}{L^p(\R^n)}
    \le \frac{1}{2}\norm{u}{L^p(\R^n)} +\frac{1}{2}\norm{\Fcal^{-1}(\bla{\lambda_\rho}^2\hat{u})}{L^p(\R^n)}.
\end{equation*}
%\commentCK{Can this final step be dropped with the alternative definition of $X^{\rho, p}(\R^n)$ as subspace of $L^p(\R^n)$?}For \(u\in \Scal'(\R^n)\), the fact that \(\Fcal^{-1}(\bla{\lambda_\rho}^2\hat{u})\in L^p(\R^n)\) implies \(u\in L^p(\R^n)\) by Proposition~\ref{prop:fourier_multipliers}~\textit{(iv)} together with Lemma~\ref{lem:lp_implication} \comment{added} for the choice \(m=\bla{\lambda_\rho}^2\).
Thus, we conclude
\begin{equation*}
D (\ol{F_{\bla{\lambda_\rho}}}^2) 
=\{ u \in L^p(\R^n):\Fcal^{-1}(\bla{\lambda_\rho}^2\hat{u}) \in L^p(\R^n) \}
= X^{\rho,p}(\R^n).
\end{equation*}
\smallskip
%\commentCK{Step~1 seems somewhat natural. Is there an alternative arument of Step~2?}

\textit{Step~2: proof~of~\eqref{eq:second_aux_id}.}  The inclusion \(D(\ol{F_{\bla{\lambda_\rho}^2}})\subset X^{\rho,p}(\R^n)\) holds true since Lemma~\ref{le:Fmclosable} for \(m=\bla{\lambda_\rho}^2\) implies that \(\mathcal{F}^{-1}(\bla{\lambda_\rho}^2\hat{u})=\ol{F_{{\bla\lambda_\rho}^2}} u\) for any \(u\in D(\ol{F_{\bla{\lambda_\rho}^2}})\) . 
We thus turn to the reverse inclusion. Observe that proving \(D(\ol{F_{\bla{\lambda_\rho}^2}})\subset X^{\rho,p}(\R^n)\) amounts to finding for each \(u\in X^{\rho,p}(\R^n)\) an approximating sequence \((u_j)_j\subset \Scal(\R^n)\) such that 
\begin{equation}\label{eq:core_bla_square}
    u_j\to u, \quad F_{\bla{\lambda_\rho}^2}u_j\to \ol{F_{\bla{\lambda_\rho}^2}} u,
    \quad \text{in } L^p(\R^n).
\end{equation}
To this end, since \(\Fcal^{-1}(\bla{\lambda_\rho}^2 \hat{u}) \in L^p(\R^n)\) by definition of \(X^{\rho,p}(\R^n)\), there exists, by the density of \(C_c^\infty(\R^n)\) in \(L^p(\R^n)\), a sequence \((w_j)_j\subset C_c^\infty(\R^n)\subset \Scal(\R^n)\) such that \(w_j\to \Fcal^{-1}(\bla{\lambda_\rho}^2 \hat{u})\) in \(L^p(\R^n)\). We define \(u_j\coloneqq F_{1/\bla{\lambda_\rho}^2} w_j\) for \(j\in \N\), note that \((F_{1/\bla{\lambda_\rho}^2} w_j)_j\subset \Scal(\R^n)\) according to Lemma~\ref{lem:closability_sr}~\textit{(i)}, and claim that this sequence satisfies the convergence properties of \eqref{eq:core_bla_square}. Indeed, by Proposition~\ref{prop:fourier_multipliers}~\textit{(iv)}, \(F_{1/\bla{\lambda_\rho}^2}\) is an \(L^p\)-Fourier multiplier and we obtain that 
\begin{equation*}
    u_j = F_{1/\bla{\lambda_\rho}^2} w_j
    \to \ol{F_{1/\bla{\lambda_\rho}^2}} \Fcal^{-1}(\bla{\lambda_\rho}^2 \hat{u})
    = u
    \quad \text{in } L^p(\R^n),
\end{equation*}
where the final equality is due to the identity \(\ol{F_{1/\bla{\lambda_\rho}^2}}v=\Fcal^{-1}(1/\bla{\lambda_\rho}^2 \hat{v})\) for any \(v\in D(\ol{F_{1/\bla{\lambda_\rho}^2}})= L^p(\R^n)\), resulting from Lemma~\ref{le:Fmclosable} with the choice \(m=1/\bla{\lambda_\rho}^2\). Furthermore, it holds with Lemma~\ref{lem:mults_algebra}, 
\begin{equation*}
    F_{\bla{\lambda_\rho}^2}u_j
    = F_{\bla{\lambda_\rho}^2}F_{1/\bla{\lambda_\rho}^2} w_j
    = w_j 
    \to \Fcal^{-1}(\bla{\lambda_\rho}^2 \hat{u})
    \quad \text{in } L^p(\R^n).
\end{equation*}
Hence, \eqref{eq:core_bla_square} is fulfilled and \eqref{eq:second_aux_id} is proved.
\end{proof}
%\begin{remark}
%    \comment{The trick we are using here is that on the level of tempered distributions we immediately have \(\Fcal^{-1}(\bla{\lambda_\rho}\Fcal(\Fcal^{-1}(\bla{\lambda_\rho}\hat{u}))=\Fcal^{-1}(\bla{\lambda_\rho}^2\hat{u})\). Moreover, in the space \(X^{\rho,p}(\R^n)\) it is also easy to construct the sought approximations. }
%\end{remark}

After these preparatory considerations, we will now show that \(\ol{F_{\abs{\lambda_\rho}}}^2= -\Delta^{\rho}\) as operators in \(L^p(\R^n)\), deriving the statement from the preceding lemmas.
%Lemmas~\ref{lem:hrhop_asdomains} and~\ref{lem:vspace} 
%combined with the repeated application of the Fourier multiplier results in Proposition~\ref{prop:fourier_multipliers}. 
%The proof still requires some care and it relies on the special properties of the operators \(\ol{F_{\bla{\lambda_\rho}}}\) and \(\ol{F_{\abs{\lambda_\rho}/\bla{\lambda_\rho}}}\).

\begin{proposition}\label{le:sqrtLaplacianp}
Assume that \(\rho\) has compact support and satisfies (H0)--(H3). Then 
\begin{equation*}
    \ol{F_{\abs{\lambda_\rho}}}^2= -\Delta^{\rho}
\end{equation*}
as operators in $L^p(\R^n)$.
\end{proposition}
\begin{proof}
Subsequently, \textit{Step~1} proves the equality $D(\ol{F_{\abs{\lambda_\rho}^2}}) = D(\ol{F_{\bla{\lambda_\rho}^2}})$, while \textit{Step~2} establishes that $D(\ol{F_{\abs{\lambda_\rho}}}^2) = D(\ol{F_{\bla{\lambda_\rho}}}^2)$. Thereafter, the proof is concluded in \textit{Step~3}. \smallskip

\smallskip

\textit{Step~1:} Similarly to the \textit{Step~2} in the proof of Lemma~\ref{lem:hrhop_asdomains}, applying Lemma~\ref{lem:domain_mulitplier} twice, namely with $m_1=\bla{\lambda_\rho}^2$ and  $m_2=|\lambda_\rho|^2/\bla{\lambda_\rho}^2$, as well as $m_1=(1+|\lambda_\rho|)^2$ and  $m_2=\bla{\lambda_\rho}^2/(1+|\lambda_\rho|)^2$
(the choices for $m_2$ are indeed $L^p$-Fourier multipliers as squares of the $L^p$-Fourier multipliers in Proposition~\ref{prop:fourier_multipliers}\,\textit{(i)} and~\textit{(ii)}, respectively; see \cite[Proposition~2.5.13]{Gr14}) leads to
\begin{align*}
D(\ol{F_{(1+|\lambda_\rho|)^2}}) \subset  D(\ol{F_{\bla{\lambda_\rho}^2}}) \subset D(\ol{F_{|\lambda_\rho|^2}}) .
\end{align*}
Identity in this line of inclusions follows from the observation that  $D(\ol{F_{|\lambda_\rho|^2}})\subset D(\ol{F_{(1+|\lambda_\rho|)^2}})$. To see the latter,  let $u\in D(F_{|\lambda_\rho|^2})$, meaning that there exists \((u_j)_j\subset \mathcal{S}(\R^n)\) such that \(u_j\to u\) and \(F_{|\lambda_\rho|^2}u_j\to \ol{F_{\abs{\lambda_\rho}^2}} u\) in \(L^p(\R^n)\). Appealing to \eqref{eq:aux_lp_2} in Lemma~\ref{lem:closability_sr}, we know that \((F_{|\lambda_\rho|}u_j)_j\) defines a Cauchy sequence in \(L^p(\R^n)\), and hence converges to \(\ol{F}_{\abs{\lambda_\rho}}\) by the closedness of the operator. 
%\begin{equation*}
%    \norm{F_{|\lambda_\rho|}v}{L^p(\R^n)}
%    \le \frac{1}{2}\norm{v}{L^p(\R^n)} +\frac{1}{2}\norm{F_{|\lambda_\rho|^2}v}{L^p(\R^n)},
%    \qquad v\in D(F_{|\lambda_\rho|^2})
%\end{equation*}
%by virtue of Young's inequality and the dual representation of the norms, we also find that \(F_{|\lambda_\rho|}u_j\to \ol{F_{\abs{\lambda_\rho}}} u\) in \(L^p(\R^n)\). 
Therefore, 
\begin{align}
F_{(1+|\lambda_\rho|)^2}u_j 
= u_j + 2F_{|\lambda_\rho|} u_j +  F_{|\lambda_\rho|^2}u_j
\to u + 2 \ol{F_{|\lambda_\rho|}} u + \ol{ F_{|\lambda_\rho|^2}u}=\ol{F_{(1+|\lambda_\rho|)^2}}u \quad \text{ in $L^p(\R^n)$},
\end{align}
showing $u\in D(\ol{F_{1+|\lambda_\rho|}})$. 
Overall, we have proved that  $D(\ol{F_{\abs{\lambda_\rho}^2}}) = D(\ol{F_{\bla{\lambda_\rho}^2}})$.

\smallskip
\textit{Step~2:} 
To show that $D(\ol{F_{\abs{\lambda_\rho}}}^2) \subset D(\ol{F_{\bla{\lambda_\rho}}}^2)$, let \(u\in D(\overline{F_{\abs{\lambda_\rho}}}^2)\). Then, by definition, \(u\in D(\overline{F_{\abs{\lambda_\rho}}})\) and \(\overline{F_{\abs{\lambda_\rho}}}u\in D(\overline{F_{\abs{\lambda_\rho}}})\). 
  %  We define \(\tilde{v}_j\coloneqq F_{\bla{\lambda_\rho}} v_j\), \(j\in \N\), and note that these function are again elements of \(\Scal(\R^n)\) by Lemma~\ref{lem:closability_sr}~\textit{(i)}. It follows with the help of Lemma~\ref{lem:mults_algebra} and Proposition~\ref{prop:fourier_multipliers}~\textit{(ii)} that 
  %  \begin{align}\label{many_convergences_ud}
  %       F_{\bla{\lambda_\rho}} u_j= \ol{F_{\bla{\lambda_\rho}/(1+\abs{\lambda_\rho})}} F_{1+\abs{\lambda_\rho}} u_j
  %       = \ol{F_{\bla{\lambda_\rho}/(1+\abs{\lambda_\rho})}} (u_j +F_{\abs{\lambda_\rho}} u_j)\to \ol{F_{\bla{\lambda_\rho}}} u 
        %\quad F_{\bla{\lambda_\rho}} \tilde{v}_j
        %= \ol{F_{\bla{\lambda_\rho}/(1+|\lambda_\rho|)}}  v_j\to \ol{F_{|\lambda_\rho|}}^2  u,
  %      \quad \text{in } L^p(\R^n).
  %  \end{align}
  It follows immediately from Lemma~\ref{lem:hrhop_asdomains} that \(u\in D(\overline{F_{\bla{\lambda_\rho}}})\) and \(\overline{F_{\abs{\lambda_\rho}}}u\in D(\overline{F_{\bla{\lambda_\rho}}})\). Furthermore, we obtain from Lemma~\ref{lem:comm_relation} that \(\ol{F_{\bla{\lambda_\rho}}} u\in D(\ol{F_{\abs{\lambda_\rho}}})\), and, again with Lemma~\ref{lem:hrhop_asdomains}, we conclude that \(\ol{F_{\bla{\lambda_\rho}}} u\in D(\ol{F_{\bla{\lambda_\rho}}})\). This verifies that \(u\in D (\ol{F_{\bla{\lambda_\rho}}}^2) \) and proves $D(\ol{F_{\abs{\lambda_\rho}}}^2) \subset D(\ol{F_{\bla{\lambda_\rho}}}^2)$. 

   The reverse inclusion $D(\ol{F_{\bla{\lambda_\rho}}}^2) \subset D(\ol{F_{\abs{\lambda_\rho}}}^2)$ is derived similarly, with the roles of \(|\lambda_\rho|\) and \(\bla{\lambda_\rho}\) exchanged. First, from \(u\in D(\ol{F_{\bla{\lambda_\rho}}}^2)\), one has \(u\in D(\ol{F_{\bla{\lambda_\rho}}})\) and \(\ol{F_{\bla{\lambda_\rho}}}u\in D(\ol{F_{\bla{\lambda_\rho}}})\). Then, Lemma~\ref{lem:hrhop_asdomains} entails that \(u\in D(\ol{F_{\abs{\lambda_\rho}}})\) and \(\ol{F_{\bla{\lambda_\rho}}}u\in D(\ol{F_{\abs{\lambda_\rho}}})\), and Lemma~\ref{lem:comm_relation} combined with Lemma~\ref{lem:hrhop_asdomains} yield \(\ol{F_{\abs{\lambda_\rho}}}u\in D(\ol{F_{\abs{\lambda_\rho}}})\), i.e., \(u\in D(\ol{F_{\abs{\lambda_\rho}}}^2)\). This concludes \textit{Step~2}, establishing that \(D(\ol{F_{\bla{\lambda_\rho}}}^2)= D(\ol{F_{\abs{\lambda_\rho}}}^2)\).
\smallskip

\textit{Step~3.}
Using \textit{Steps~1}~and~\textit{2}, we infer from Lemma~\ref{lem:vspace} that 
\begin{equation}\label{eq:domain_lambda_squared}
    D (\ol{F_{\abs{\lambda_\rho}}}^2) 
   = D(\ol{F_{\abs{\lambda_\rho}^2}}).
    %\quad\text{ and }\quad
    %\ol{F_{\abs{\lambda_\rho}}}^2
    %=\ol{F_{\abs{\lambda_\rho}^2}}.
\end{equation}
%let \(u\in D(\ol{F_{\abs{\lambda_\rho}}}^2)=D(\ol{F_{\abs{\lambda_\rho}^2}})\), and take \((u_j)_j, (v_j)_j\subset \mathcal{S}(\R^n)\) such that 
%\begin{align*}
%    u_j\to u, 
%    \quad F_{|\lambda_\rho|} u_j\to \ol{F_{|\lambda_\rho|}} u , 
%    \quad v_j\to \ol{F_{|\lambda_\rho|}} u, 
%    \quad F_{|\lambda_\rho|} v_j\to \ol{F_{|\lambda_\rho|}}^2  u
%    \quad \text{in } L^p(\R^n).
%\end{align*} 
%    We compute 
%    \begin{align*}
%   &\langle \ol{F_{\abs{\lambda_\rho}}}^2,\varphi\rangle
%    = \lim_{j\to \infty} \langle F_{|\lambda_\rho|}v_j,\varphi\rangle
%    =
%    \lim_{j\to \infty} 
%        \langle {v_j},\mathcal{F}(|\lambda_\rho|\mathcal{F}^{-1}\varphi)\rangle\\
%       & = 
%        \lim_{j\to \infty}
%        \langle F_{|\lambda_\rho|} u_j,\mathcal{F}(|\lambda_\rho|\mathcal{F}^{-1}\varphi)\rangle
%        = \lim_{j\to \infty}
%        \langle  u_j,\mathcal{F}(|\lambda_\rho|^2\mathcal{F}^{-1}\varphi)\rangle
%        = \langle u, \mathcal{F}(|\lambda_\rho|^2\mathcal{F}^{-1}\varphi)\rangle,
%        \qquad \varphi\in \Scal(\R^n)
%    \end{align*}
%    using that \(\mathcal{F}(|\lambda_\rho|\mathcal{F}^{-1}\varphi)\in L^2(\R^n)\) since \%(\abs{\lambda_\rho}\) has polynomial growth.
Lemma~\ref{lem:closed_ops_equality} with \(m_1=\abs{\lambda_\rho}=m_2\) (where the requisite property that \(\mathcal{F}(\abs{\lambda_\rho}\mathcal{F}^{-1}\varphi), \mathcal{F}(\abs{\lambda_\rho}^2\mathcal{F}^{-1}\varphi)\in L^{p'}(\R^n)\) for all \(\varphi\in \Scal(\R^n)\) is a direct consequence of Lemma~\ref{lem:closability_sr}) shows that \(\ol{F_{\abs{\lambda_\rho}}}^2
    =\ol{F_{\abs{\lambda_\rho}^2}}\), and it therefore follows along with Proposition \ref{prop:FourierLaplacianp2}\,\textit{(iii)} the identity $\ol{F_{\abs{\lambda_\rho}}}^2 = -\Delta^{\rho}$ as unbounded operators in \(L^p(\R^n)\).
\end{proof}
With the result of Proposition~\ref{le:sqrtLaplacianp} in mind, we define the square root of $-\Delta^\rho$ as an unbounded operator in $L^p(\R^n)$ as follows. 
\begin{definition}[Square root of $-\Delta^\rho$]\label{def:sqrt_lap}
    For \(\rho\) with compact support and satisfying (H0)--(H3) we define the operator
    \begin{equation*}
        (-\Delta^{\rho})^{1/2}:D((-\Delta^{\rho})^{1/2})\to L^p(\R^n)
    \end{equation*}
    by
    \begin{equation*}
      D((-\Delta^{\rho})^{1/2})\coloneqq D(\ol{F_{\abs{\lambda_\rho}}}), \qquad (-\Delta^{\rho})^{1/2}
        \coloneqq \ol{F_{\abs{\lambda_\rho}}}. 
    \end{equation*}
\end{definition}

%\comment{From before:}Finally, we are ready to prove \color{red} a generalization of  the well-known representation of classical Bessel potential spaces, see~e.g.~\cite{}. As in the $p=2$ analogue of Proposition 4.14, the main idea of the proof is again to leverage the alternative Fourier characterization of \(H^{\rho,p}(\R^n)\) provided in Theorem~\ref{thm:alt_char}.  \color{black}

%\commentCK{Idea was not made explicit before...could be phrased clearer, or more vaguely, which could also help... Maybe also give an idea of previous proofs.}\comment{I added something before Proposition~\ref{prop:charac_lapace2}} 
%However, some care is required to connect for \(u\in H^{\rho,p}(\R^n)\) the expressions \(\mathcal{F}^{-1}(\bla{\lambda_\rho}\hat{u})\) -- understood as a tempered distribution -- and \(\ol{F_{\abs{\lambda_\rho}}}u\) with each other. 
%In particular, a priori, \(F_{\bla{\lambda_\rho}}u\) might exists without \(u\) being an element of \(D(\ol{F_{\abs{\lambda_\rho}}}u)\).

%\comment{Or now maybe just:}
Finally, we are ready to conclude, as also described in Section~\ref{sec:discussion}, a generalization of the well-known representation of classical Bessel potential spaces through the square root of the negative fractional Laplace operator.
%, see~e.g.~\cite{}.\commentCK{Insert reference.}
%\comment{Maybe as corollary now?}
\begin{corollary}[Characterization in terms of \((-\Delta^{\rho})^{1/2}\)]\label{prop:charac_lapacep}
Assume that \(\rho\) has compact support and satisfies (H0)--(H3). Then
\[
H^{\rho,p}(\R^n)
%= \{u\in L^p(\R^n):\, F_{\abs{\lambda_\rho}}u\in L^p(\R^n)\} 
%= \{u\in L^p(\R^n):\, (-\Delta^{\rho})^{1/2} u \in L^p(\R^n)\} 
= \{u\in L^p(\R^n):\,  \ol{F_{|\lambda_\rho|}}u \in L^p(\R^n)\}
= D((-\Delta^{\rho})^{1/2}).
\]
\end{corollary}
\begin{proof}
The statement is a direct consequence of Definition~\ref{def:sqrt_lap} together with Proposition~\ref{lem:hrhop_asdomains}.
\end{proof}

%%%%%%%%%%%%%%%%%%%%%%%%%%%%%%%%%%%%%%%%%%

\subsection{Characterization in terms of \((1-\Delta^\rho)^{1/2}\)}\label{sec:laplace_bp_inverse}
%\comment{Draft for a new and much shortened version of the section. It contains all the information below but with just one new definition which is now justified a posteriori, without building up to it through a sequence of other definitions. The arguments are also cleaner than before and incorporate some of the points you were raising below.}
Motivated by the Fourier characterization in Theorem~\ref{thm:alt_char}, we define for \(\rho\) with compact support and subject to hypotheses~(H0)-(H2) and (H3') the operator 
\begin{equation*}
    (1-\Delta^\rho)^{1/2}:\mathcal{S}'(\R^n)\to\mathcal{S}'(\R^n)
\end{equation*}
as
\begin{equation}\label{eq:sqr_tempered}
        (1-\Delta^\rho)^{1/2}u
        \coloneqq \mathcal{F}^{-1}(\bla{\lambda_\rho}\hat{u}),
        \qquad  u\in \Scal'(\R^n).
    \end{equation}
    In fact, this definition is compatible with the framework developed in the preceding subsections. By virtue of the first identity of Lemma~\ref{lem:fourier_identities}~\textit{\ref{item:fi_symm}}, we have the equality
    \begin{align*}
        \langle (1-\Delta^\rho)^{1/2}u,\varphi\rangle_{\mathcal{S}'\times \mathcal{S}} 
         = \langle u,\mathcal{F}(\bla{\lambda_\rho}\mathcal{F}^{-1}\varphi)\rangle_{\mathcal{S}'\times \mathcal{S}} 
        = \langle u, F_{\bla{\lambda_\rho}}\varphi\rangle_{\mathcal{S}'\times \mathcal{S}},
        \qquad  u\in \mathcal{S}'(\R^n),\varphi\in \mathcal{S}(\R^n).
    \end{align*}
    In addition, the operator \(F_{\bla{\lambda_\rho}}\) studied in Lemma~\ref{lem:closability_sr}~\textit{(i)} has a well defined square in \(\Scal(\R^n)\), which equals \(1+F_{\abs{\lambda_\rho}^2}=1-\Delta^\rho\), thanks to \eqref{eq:rholap_S}. Hence, %Jointly, both properties imply that
    \((1-\Delta^\rho)^{1/2}\) is indeed the extension by duality to \(\Scal'(\R^n)\) of a square root of the operator \(1-\Delta^\rho\) in \(\Scal(\R^n)\). %\comment{If the kernel satisfies additionally hypothesis~(H3), then we even also recover the identity \(\ol{F_{\bla{\lambda_\rho}}}^2 = 1+\ol{F_{\abs{\lambda_\rho}^2}}= 1-\Delta^\rho\) as operators in \(L^p(\R^n)\), where the first identity follows from a similar approach as the proof of Lemma~\ref{le:sqrtLaplacianp} (but a bit simpler) and the second one is a direct consequence of Proposition~\ref{prop:FourierLaplacianp2}~\textit{\ref{item:LaplacianFourierp}}. However, I think that I would omit this point since it is not really essential for the idea of this subsection and it still involves some computations.}
    
    With these considerations, we can reformulate Theorem~\ref{thm:alt_char} as follows. 
    %have the following result, which is just a reformulation of Theorem~\ref{thm:alt_char}.
    %\comment{Maybe better as corollary?}
    \begin{corollary}[Characterization in terms of \((1-\Delta^{\rho})^{1/2}\)]\label{prop:charac_bessel_op}
 Let \(p\in(1,\infty)\) and assume that \(\rho\) has compact support and satisfies (H0)--(H2) and (H3'). We have the identity
    \begin{equation*}
        H^{\rho,p}(\R^n)
        = \{u\in \mathcal{S}'(\R^n):\, (1-\Delta^{\rho})^{1/2} u \in L^p(\R^n)\}.
    \end{equation*} 
\end{corollary}

	\section{Characterization as potential spaces}\label{sec:charac_bessel}

	The subject of this section is the identification of \(H^{\rho,p} (\R^n)\) with a suitable potential space, which is a generalization of the identity between \(H^{\rho_s,p} (\R^n)\) for the Riesz fractional gradient \(\rho_s\) and the Bessel potential spaces (see \cite[Theorem~1.7]{ShSp15}).
	In particular, we prove that, under sufficient  conditions on the kernel function $\rho$, \(H^{\rho,p}(\R^n)\) can likewise be regarded as the image of \(L^p(\R^n)\) under convolution with a suitable generalized potential kernel \(g_\rho\), defined in Proposition~\ref{prop:invers_fourier} below. 

  %\commentCK{In the same spirit as...mention introductory Section~\ref{sec:discussion}..., the challenge therefore is..}
  
    As outlined in the introductory Section~\ref{sec:discussion} and in analogy to the case of the Riesz fractional gradient, the candidate for the function \(g_\rho\) is the Fourier inverse of  %\commentCK{Why? Motivation missing.}
	\begin{equation*}
	\bla{\lambda_\rho}^{-1}=(1+ \abs{\lambda_\rho}^2)^{-{1}/{2}}.
	%=(1+4\pi^2 \abs{\lambda_{\rho} }^2)^{-{1}/{2}}.
	\end{equation*}
    Thus, the challenge consists in the following in proving that this Fourier inverse exists, and, moreover, that it lies in \(L^1(\R^n)\). 
	%admits a Fourier inverse in \(L^1(\R^n)\).
    For reasons that will become apparent in Lemma~\ref{lem:gen_pot_frac_inverse}, we show  a more general result. 
	\begin{proposition}\label{prop:invers_fourier}
		Assume that \(\rho\) have compact support and satisfy assumptions (H0)--(H3). 
		 Then for each \(\ell\in\N\) there exists a function \(g_{\rho,\ell}\in L^1(\R^n)\) such that %\commentCK{$\alpha$ is always used for multi-indices, maybe use something different here, e.g. $\ell$ }
		\begin{equation*}
		\widehat{g_{\rho,\ell}}
		= \bla{\lambda_\rho}^{-\ell},
        %=(1+\abs{\lambda_\rho}^2)^{-{\ell}/{2}},
		\end{equation*}
		and the family \(({g_{\rho,\ell}})_{\ell}\) is uniformly bounded in \(L^1(\R^n)\).
		\end{proposition}
	\begin{proof}
		We proceed in several steps. We first formally introduce the candidates for the functions \({g_{\rho,\ell}}\), before proving that they actually constitute elements of \(L^1(\R^n)\). For the latter purpose, we derive explicit pointwise upper bounds, guaranteeing that the \(L^1\)-norm of \(g_{\rho,\ell}\) is bounded independently of \(\ell\in\N\).
		
		\textit{Step~1: construction of the functions $g_{\rho, \ell}$.} 
		We recall initially that according to hypothesis~(H3), we know that (see \eqref{eq:lest_qrho_pre} together with the comments below~\eqref{eq:uest_qrho_alp_pre})
		\begin{equation}\label{eq:lowerboundlambda}
		\abs{\lambda_\rho (\xi)}^2
		\ge C|\xi|^{2\sigma},
        \qquad \xi\in B_1^c
		\end{equation}
		for \(\sigma\in (0,1)\) from (H3) and a constant \(C>0\).
As a consequence of this and the polynomial growth of the derivatives of $|\lambda_{\rho}|^2$, following Lemma~\ref{lem:lambda_FintoF_pre}, for each $\delta>0$, the function
\begin{equation}\label{eq:aux_h}	h_\rho^\delta\coloneqq e^{- \abs{\lambda_{\rho}}^2\delta/(4 \pi^2)}\delta^{n/2}
		\end{equation}
lies in \(\mathcal{S}(\R^n)\). Consequently, the Fourier inverse \(f^\delta_\rho\coloneqq \mathcal{F}^{-1} h^\delta_\rho\) exists and is again an element of \(\mathcal{S}(\R^n)\). Furthermore, appealing to the identity
%\comment{should we use large or small brackets with the \(\ell\) inserted?}
		\begin{equation}\label{eq:int_id}
		t^{-\ell/2}
		=\frac{1}{\Gamma(\frac{\ell}{2})}\int_0^\infty e^{-t\delta}\delta^{\ell/2}\frac{\dl \delta}{\delta}, \qquad  t>0, \ell>0 
		\end{equation}
		it follows with the choice \(t=(1+\abs{\lambda_{\rho}(\xi)}^2)/4\pi^2\) for $\xi\in \R^n$  that 
		\begin{align*}
		\frac{1}{(4\pi^2)^{\ell/2}\Gamma(\frac{\ell}{2})}\int_0^\infty h^\delta_\rho(\xi)e^{-\delta/4\pi^2}\delta^{(\ell-n)/2}\frac{\dl \delta}{\delta}
		&=\frac{1}{(4\pi^2)^{\ell/2}\Gamma(\frac{\ell}{2})}\int_0^\infty e^{- \abs{\lambda_{\rho}(\xi)}^2\delta/(4 \pi^2)}e^{-\delta/4\pi^2}\delta^{\ell/2}\frac{\dl \delta}{\delta}\\
		&=\frac{1}{(4\pi^2)^{\ell/2}\Gamma(\frac{\ell}{2})}\int_0^\infty e^{-(1+\abs{\lambda_{\rho}(\xi)}^2)\delta/4\pi^2}\delta^{\ell/2}\frac{\dl \delta}{\delta}\\
		&= (1+\abs{\lambda_{\rho} (\xi)}^2)^{-{\ell}/{2}},
        \qquad \xi\in \R^n, \ell\in \N.
		\end{align*}
		Thus, defining 
		\begin{equation}\label{eq:def_conv_kernel}
		g_{\rho,\ell}(x)
		\coloneqq \frac{1}{(4\pi^2)^{\ell/2}\Gamma(\frac{\ell}{2})}\int_0^\infty f^\delta_\rho(x)e^{-\delta/4\pi^2}\delta^{(\ell-n)/2}\frac{\dl \delta}{\delta}, \qquad x\in \R^n,
		\end{equation}
		  Fubini's theorem, together with the fact that \(g_{\rho,\ell}\in L^1(\R^n)\) for all \(\ell\in\N\), as will be established below, yields that 
%		\begin{align*}
%		\mathcal{F}\left(\frac{1}{(4\pi^2)^{\alp/2}\Gamma(\frac{\alp}{2})}\int_0^\infty f^\delta_\rho(x)e^{-\delta/4\pi^2}\delta^{(\alp-n)/2}\frac{\dl \delta}{\delta}\right)
%		= (1+4\pi^2 \abs{\lambda_{\rho} (\xi)}^2)^{-{\alp}/{2}}.
%		\end{align*} 
		\begin{align*}
		\widehat{g_{\rho,\ell}}
		= (1+\abs{\lambda_{\rho}}^2)^{-{\ell}/{2}}
        = \bla{\lambda_\rho}^{-\ell}.
		\end{align*}

		\textit{Step~2: sufficient conditions for \(g_{\rho,\ell}\in L^1(\R^n)\).} From here on, \(C>0\) denotes a generic, possibly changing constant which may depend on \(\rho\) but is, in particular, independent of \(\delta>0\) and \(\ell\in \N\). 
        
		It remains to verify the property that \(({g_{\rho,\ell}})_{\ell}\) is uniformly bounded in \(L^1(\R^n)\). 
		By definition, we have
%\comment{up to some preconstant depending on our definition}
        \begin{align*}
        f^\delta_\rho(x) = \delta^{n/2}\int_{\R^n} e^{-\abs{\lambda_{\rho} (\xi)}^2\delta / (4 \pi^2)}e^{2\pi i x\cdot \xi}\dl \xi, \qquad x\in \R^n.
		\end{align*}		
We will show in \textit{Steps~3}, \textit{4} and \textit{5} below that there exist \(\eta_1, \eta_2, \eta_3 \in [0,1/2)\) and $b \in (0,1)$
%, and some constant \(C>0\) 
such that 
%\commentCK{Is there a reason not to simplify the first estimate to $..\leq C\delta^{n/2-\eta}...$ for some $\eta\in (0,1/2)$ (then given as the maximum of $\eta_1,\eta_2, \eta_3$)? Write $\eta_\delta$?}
		\begin{equation}\label{eq:close_orign}
		 \Big|f_\rho^\delta(x)\Big|
		\le C\sum_{i=1}^3\delta^{n/2-\eta_i}\frac{1}{|x|^{n-b}}, \qquad x \in B_1\setminus \{0\}, \, \delta\in (0,1),
		\end{equation} 
		\begin{equation}\label{eq:away_orign}
		\Big|f_\rho^\delta(x)\Big|
		\le C\delta^{n/2}\frac{1}{|x|^{n+1}} , \qquad x \in B_1^c, \, \delta\in (0,1),
		\end{equation} 
        and
        %\comment{this is the new estimate; since we only need it for \(\delta\in (1,\infty)\), I restricted the range of \(\delta\) above.}
        \begin{equation}\label{eq:delta_large}
            \norm{f_\rho^\delta}{L^1(\R^n)}
            \le C\delta^{n/2},
            \qquad \delta\in (1,\infty).
        \end{equation}
We continue with the proof assuming these estimates. %Note that \(\delta^{n/2}\le \delta^{(n-1)/2}\) for \(\delta\in (0,1)\), while \(\delta^{(n-1)/2}\le \delta^{n/2}\) for \(\delta\in (1,\infty)\).
%Since \((\alp-3)/2> -1\) as long as \(\alp\ge 2\), it follows in these cases that 
%\begin{align*}
%\int_0^1 |f^\delta_\rho(x)|e^{-\delta/4\pi^2}\delta^{(\alp-n)/2}\frac{\dl \delta}{\delta}
%&\le C \left(\chi_{x\in B_1}\frac{1}{|x|^{n-\sigma}}
%            +\chi_{x\in B_1^c}\frac{1}{|x|^{n+1}}\right)
%            \int_0^1 e^{-\delta/4\pi^2}\delta^{(\alp-3)/2}\dl \delta\\
%            &\le C\left(\chi_{x\in B_1}\frac{1}{|x|^{n-\sigma}}
%           +\chi_{x\in B_1^c}\frac{1}{|x|^{n+1}}\right),
%\end{align*}
%independently of \(\alp\in \N_{\ge 2}\),
     %   \begin{align*}
     %       \int_{\R^n} |g_{\rho,\alp}(x)|\dl x
     %       \le \frac{1}{(4\pi^2)^{\alp/2}\Gamma(\frac{\alp}{2})}
     %       &\left(\int_{B_1}\frac{1}{|x|^{n-\sigma}}\dl x
     %       + \int_{B_1^c}\frac{1}{|x|^{n+1}}\dl x\right)\\
     %       &\left(\int_0^1 e^{-\delta/4\pi^2}\delta^{(\alp-3)/2}\dl \delta
     %       + \int_1^\infty \right)
     %   \end{align*}
%\begin{align*}
%		\int_{\R^n} |g_{\rho,\alp}(x)|\dl x
%		&\le C \left(1+\frac{1}{(4\pi^2)^{\alp/2}\Gamma(\frac{\alp}{2})}\int_1^\infty e^{-\delta/4\pi^2}\delta^{\alp/2}\frac{\dl \delta}{\delta}\right)\\
%		&\le C\left(1+\frac{1}{(4\pi^2)^{\alp/2}\Gamma(\frac{\alp}{2})}\int_0^\infty e^{-\delta/4\pi^2}\delta^{\alp/2}\frac{\dl \delta}{\delta}\right),
%		\end{align*} 
%		where the right-hand side is now uniformly bounded due to \eqref{eq:int_id}.
Observing that
\begin{align*}
    \int_{\R^n} &|g_{\rho,\ell}(x)|\dl x \\
    &\leq \frac{1}{(4\pi^2)^{\ell/2}\Gamma(\frac{\ell}{2})} \left[ \int_{B_1} \int_0^1 + \int_{B_1^c} \int_0^1 + \int_{\R^n} \int_1^{\infty}\right] \abs{f^\delta_\rho(x)} e^{-\delta/4\pi^2}\delta^{(\ell-n)/2}\frac{\dl \delta}{\delta} \dl x, 
\end{align*}
we then estimate using \eqref{eq:close_orign}, 
\begin{align*}
     \int_{B_1} \int_0^1 \abs{f^\delta_\rho(x)} e^{-\delta/4\pi^2}\delta^{(\ell-n)/2}\frac{\dl \delta}{\delta} \dl x 
     &\leq C \sum_{i=1}^3\int_{B_1} \frac{1}{|x|^{n-b}} \dl x \int_0^1 e^{-\delta/4\pi^2}\delta^{(\ell-2)/2-\eta_i} \dl \delta \\
     &\leq C \int_0^1 \delta^{-1/2-\max_{i\in\{1,2,3\}}\eta_i} \dl \delta \leq C ,
\end{align*}
using \eqref{eq:away_orign},
\[
 \int_{B_1^c} \int_0^1 \abs{f^\delta_\rho(x)} e^{-\delta/4\pi^2}\delta^{(\ell-n)/2}\frac{\dl \delta}{\delta} \dl x \leq C \int_{B_1^c} \frac{1}{|x|^{n+1}} \dl x \int_0^1 e^{-\delta/4\pi^2}\delta^{\ell/2} \dl \delta \leq C ,
\]
and using \eqref{eq:delta_large} and Fubini's theorem,
\begin{align*}
 \int_{\R^n} \int_1^{\infty} \abs{f^\delta_\rho(x)} e^{-\delta/4\pi^2}\delta^{(\ell-n)/2}\frac{\dl \delta}{\delta} \dl x 
 &\leq C \int_1^{\infty} e^{-\delta/4\pi^2}\delta^{\ell/2} \frac{\dl \delta}{\delta} \\
 &\leq C \int_0^{\infty} e^{-\delta/4\pi^2}\delta^{\ell/2} \frac{\dl \delta}{\delta} = C (4 \pi^2)^{\ell/2} \Gamma \left( \frac{\ell}{2}\right).
\end{align*}
Altogether and along with the property of the $\Gamma$-function that 
%fact that (see~\eqref{})  \commentCK{Insert reference...or just by the properties of the $\Gamma$-function}
		\begin{equation*}
			\sup_{\ell \in (0, \infty)} \frac{1}{(4\pi^2)^{\ell/2}\Gamma(\frac{\ell}{2})} < \infty,
		\end{equation*}
we conclude 
\[
 \int_{\R^n} |g_{\rho,\ell}(x)|\dl x \leq C
\]
with a constant independent of $\ell\in \N$.

We will now turn to the remaining proof of \eqref{eq:close_orign}, \eqref{eq:away_orign} and \eqref{eq:delta_large}. For the first two estimates, we will use a quantified form of the fact that, as an inverse Fourier transform of an exponential function, \(f_\rho^\delta\) decays faster than any polynomial. In particular, the corresponding approach refines the techniques applied in Proposition~\ref{prop:poincare2}. To derive \eqref{eq:delta_large}, we will apply different strategy, relying on a standard Fourier estimate for the \(L^1\)-norm . 
%In fact, the same approach as for the other estimates does not work for \(\delta\in (1,\infty)\) (see Remark~\ref{}).  

To unburden the notation, the following estimates will be carried out for the auxiliary function 
\begin{align}\label{Theta_rhodelta}
\redconvker:=e^{-|\lambda_\rho|^2\delta/(4\pi^2)} = h_\rho^\delta \delta^{-n/2},
\end{align}
whose relation with $f_\rho^\delta$ is given through $f_\rho^\delta= \Fcal^{-1} h_\rho^\delta = \delta^{n/2} \Fcal^{-1} \redconvker$.

	%	\commentCK{Introduce the function $\Theta_\rho^\delta$ in a suitable place!}

        \textit{Step~3: proof of \eqref{eq:close_orign}.}
        We will demonstrate the sought estimate for general \(\delta\in (0,\infty)\).
		Using the radial nature of $\redconvker$, and, hence, exploiting the rotational invariance of the integrand, we write
		\begin{equation*}
			\Fcal^{-1} \redconvker(x) = \left[ \int_{B_{1}} +\int_{B_1^c} \right] \redconvker(\xi)e^{2\pi i |x|\xi_1}\dl \xi, \qquad x\in B_1\setminus\{0\}.
		\end{equation*} 
       It follows directly that  %	for \(x\in B_1\), 
		\begin{equation*}
\Bigg|\int_{B_{1}}\redconvker(\xi)e^{2\pi i |x|\xi_1}\dl \xi\Bigg|
		\le C
		\le C\frac{1}{|x|^{n-b}}, \qquad x\in B_1\setminus\{0\}, 
		\end{equation*}
		for any \(b\in (0,1)\). For the integral over \(B_1^c\), 
		applying now the change of variables \(\xi\mapsto \delta^{-1/(2\sigma)}\xi\) and partially integrating \(n-1\) times, we find the expression
		\begin{equation}\label{eq:B1c}
        \begin{split}
			\int_{B_1^c} \redconvker(\xi)e^{2\pi i |x|\xi_1}\dl \xi
			&=\delta^{-n/(2\sigma)}
			\int_{B^c_{\delta^{1/(2\sigma)}}} \redconvker\left(\frac{\xi}{\delta^{1/(2\sigma)}}\right)
			e^{2\pi i |x|\delta^{-1/(2\sigma)}\xi_1}\dl \xi \\
			&= \delta^{-1/(2\sigma)}\frac{(-1)^{n-1}}{(2\pi i|x|)^{n-1}}
			\int_{B^c_{\delta^{1/(2\sigma)}}}
			\p^{n-1}_{\xi_1}\redconvker\left(\frac{\xi}{\delta^{1/(2\sigma)}}\right)
			e^{2\pi i |x|\delta^{-1/(2\sigma)}\xi_1}\dl \xi
			+  R(x) 
			%&= \delta^{-1/(2\sigma)}\frac{(-1)^{n-1}}{(2\pi i|x|)^{n-1}}
			%\int_{\delta^{1/(2\sigma)}<|\xi|\le \delta^{a}/|x|}
			%\p^{n-1}_{\xi_1}g\left(\frac{\xi}{\delta^{1/(2\sigma)}}\right)
			%e^{2\pi i |x|\delta^{-1/(2\sigma)}\xi_1}\dl \xi\\
			%&+\frac{(-1)^{n}}{(2\pi i|x|)^{n}} 
			%\int_{\p B_{\delta^{a}/|x|}}
			%\p^{n-1}_{\xi_1}g\left(\frac{\xi}{\delta^{1/(2\sigma)}}\right)
			%e^{2\pi i |x|\delta^{-1/(2\sigma)}\xi_1}\dl \mathcal{H}^{n-1}(\xi)\\
			%&+\frac{(-1)^{n}}{(2\pi i|x|)^{n}} 
			%\int_{|\xi|> \delta^{a}/|x|}
			%\p^{n}_{\xi_1}g\left(\frac{\xi}{\delta^{1/(2\sigma)}}\right)
			%e^{2\pi i |x|\delta^{-1/(2\sigma)}\xi_1}\dl \xi\\
			%&+  R(x),
		\end{split}
        \end{equation}
		for $x\in B_1\setminus\{0\}$, with the remainder% \commentCK{Check exponents!}
		\begin{equation*}
		R(x) 
		= \delta^{-n/(2\sigma)}\sum_{k=1}^{n-1} \frac{(-1)^{k-1}\delta^{(n-1)/(2\sigma)}}{(2\pi i|x|)^{n-1}}
		\int_{\p B_{\delta^{1/(2\sigma)}}} \p^{k-1}_{\xi_1}\redconvker\left(\frac{\xi}{\delta^{1/(2\sigma)}}\right)\p^{n-1-k}_{\xi_1}e^{2\pi i |x| \delta^{-1/(2\sigma)}\xi_1}\dl \mathcal{S}(\xi),
		\end{equation*}
        comprising all boundary integrals arising from the \(n-1\) partial integrations. 
        
        Fix
\[b
 \coloneqq \frac{\sigma}{4 - 3 \sigma}, \qquad
\gamma:= \frac{1}{4-3\sigma}, 
\]
and note that $b \in (0, \sigma) \subset (0,1)$ and $\gamma=b/\sigma$.
By virtue of estimate \eqref{eq:new_est}, the remainder can be bounded by
		\begin{align}
            \begin{split}\label{eq:gen_fin_est_R}
                |R(x) |
			&\le C  \delta^{-n/(2\sigma)}\sum_{k=1}^{n-1} 
            \frac{\delta^{k/(2\sigma)}}{|x|^{k}}
			\int_{\p B_{\delta^{1/(2\sigma)}}}|\xi|^{-(k-1+\gamma\sigma)}\dl \mathcal{S}(\xi)\\
            &\le C  \delta^{-n/(2\sigma)}\sum_{k=1}^{n-1} 
            \frac{\delta^{k/(2\sigma)}}{|x|^{k}}\delta
            ^{(n-1)/(2\sigma)}
			\int_{\p B_{1}}\delta^{-(k-1+\gamma\sigma)/(2\sigma)}|\xi|^{-(k-1+\gamma\sigma)}\dl \mathcal{S}(\xi)\\
            &\le C \frac{\delta^{-\gamma/2}}{|x|^{n-b}},
            \qquad x\in B_1\setminus \{0\},
            \end{split}
		\end{align}
where the change of variables \(\xi\mapsto \delta^{1/(2\sigma)}\xi\) was used in the second step. 

Fix
\[
 a \coloneq \frac{4 - 3 \sigma}{4 \sigma}
\]
and note that $a \in (1/4, \infty)$ and $a b = 1/4$.
For simplicity, we assume in the following that \(\delta^{a}/|x|>\delta^{1/(2\sigma)}\), entailing that the integral over $B^c_{\delta^{1/(2\sigma)}}$ in the right-hand side of \eqref{eq:B1c} is subdivided into two domains: $B_{\delta^{a}/|x|} \setminus B_{\delta^{1/(2\sigma)}}$ and $B^c_{\delta^{a}/|x|}$. If instead \(\delta^{a}/|x|\le\delta^{1/(2\sigma)}\), one might just extend the integral to \(B^c_{\delta^{a}/|x|}\) and proceed according to the steps below.
In the first domain, we leave the expression unchanged, while in the second one, we partially integrate once more, revealing the identity
        \begin{align*}
			& \delta^{-1/(2\sigma)}\frac{(-1)^{n-1}}{(2\pi i|x|)^{n-1}} 
			\int_{B^c_{\delta^{a}/|x|}}
			\p^{n-1}_{\xi_1}\redconvker\left(\frac{\xi}{\delta^{1/(2\sigma)}}\right)
			e^{2\pi i |x|\delta^{-1/(2\sigma)}\xi_1}\dl \xi\\
			&= 
			\frac{(-1)^{n}}{(2\pi i|x|)^{n}} 
			\int_{\p B_{\delta^{a}/|x|}}
			\p^{n-1}_{\xi_1}\redconvker\left(\frac{\xi}{\delta^{1/(2\sigma)}}\right)
			e^{2\pi i |x|\delta^{-1/(2\sigma)}\xi_1}\dl \mathcal{S}(\xi)\\
			&+\frac{(-1)^{n}}{(2\pi i|x|)^{n}} 
			\int_{B^c_{\delta^{a}/|x|}}
			\p^{n}_{\xi_1}\redconvker\left(\frac{\xi}{\delta^{1/(2\sigma)}}\right)
			e^{2\pi i |x|\delta^{-1/(2\sigma)}\xi_1}\dl \xi, \qquad x\in B_1\setminus\{0\}.
		\end{align*}
Again with \eqref{eq:new_est}, and recalling that
\[
 \gamma=b/\sigma , \qquad a b = 1/4 , \qquad -\frac{1}{2\sigma} +a(1-b)
			=\frac{1-2\sigma}{2\sigma} ,
\]
we ultimately deduce for $x\in B_1\setminus\{0\}$ the estimates 
		\begin{align}
            \begin{split}\label{eq:gen_fin_est_1}
                \delta^{-1/(2\sigma)}\frac{1}{|x|^{n-1}}
			\int_{B_{\delta^{a}/|x|} \setminus B_{\delta^{1/(2\sigma)}}}
			\Big|\p^{n-1}_{\xi_1}\redconvker\left(\frac{\xi}{\delta^{1/(2\sigma)}}\right)\Big|\dl \xi
			&\le C\delta^{-1/(2\sigma)}\frac{1}{|x|^{n-1}}
			\int_{B_{\delta^{a}/|x|}}
			\frac{1}{|\xi|^{n-1+b}}\dl \xi\\
			&= C\delta^{(1-2\sigma)/(2\sigma)}
			\frac{1}{|x|^{n-b}},
            \end{split}
		\end{align}
		and 
		\begin{align}
		  \begin{split}\label{eq:gen_fin_est_2}
		      \frac{1}{|x|^{n}} 
			\int_{\p B_{\delta^{a}/|x|}}
			\Big|\p^{n-1}_{\xi_1}\redconvker\left(\frac{\xi}{\delta^{1/2}}\right)\Big|
			\dl \mathcal{S}(\xi)
			\le C\frac{1}{|x|^{n}} 
			\int_{\p B_{\delta^{a}/|x|}}
			\frac{1}{|\xi|^{n-1+b}}
			\dl \mathcal{S}(\xi)
			= C\delta^{-1/4}\frac{1}{|x|^{n-b}}, 
		  \end{split}
		\end{align}
and finally,
		\begin{align}
		  \begin{split}\label{eq:gen_fin_est_3}
		      &\frac{1}{|x|^{n}}
		\int_{B^c_{\delta^{a}/|x|}}
		\Big|\p^{n}_{\xi_1}\redconvker\left(\frac{\xi}{\delta^{1/(2\sigma)}}\right)\Big|\dl \xi
		\le C\frac{1}{|x|^{n}}
		\int_{B^c_{\delta^{a}/|x|}}
		\frac{1}{|\xi|^{n+b}}\dl \xi
		= C\delta^{-1/4}
		\frac{1}{|x|^{n-b}}.
		  \end{split}
		\end{align}
        Since \(\gamma \in (0,1)\), all powers of \(\delta\) in \eqref{eq:gen_fin_est_R}, \eqref{eq:gen_fin_est_2} and \eqref{eq:gen_fin_est_3} have an exponent in the range \((-1/2,0)\), and because \(1-2\sigma>-\sigma\), the exponent of \(\delta\) in \eqref{eq:gen_fin_est_1} lies in \((-1/2,\infty)\). Thus, keeping in mind in view of \eqref{eq:gen_fin_est_1} that \(\delta\in (0,1)\), it holds 
        \begin{equation*}
        \Big|\int_{\R^n} \redconvker(\xi)e^{2\pi i x\cdot \xi}\dl \xi\Big|
		\le C\sum_{i=1}^3\delta^{-\eta_i}\frac{1}{|x|^{n-b}} 
        \quad\text{and}\quad
		\Big|f_\rho^\delta(x)\Big|
		\le C\sum_{i=1}^3\delta^{n/2-\eta_i}\frac{1}{|x|^{n-b}},
        \qquad  x\in B_1\setminus \{0\}
		\end{equation*}
        for some \(\eta_1, \eta_2, \eta_3 \in [0,1/2)\). 
        %\commentCK{Consider introducing an $\eta_\delta$ as max or min $\eta_1, \eta_2,\eta_3$.}

        \textit{Step~4: proof of \eqref{eq:away_orign}.}
		%\comment{New:}
        In order to prove \eqref{eq:away_orign}, we fix \(\delta\in (0,1)\). Integrating by parts, we write
        \begin{align*}
		 \Fcal^{-1} \redconvker(x) & =\int_{\R^n} \redconvker (\xi)e^{2\pi i x\cdot \xi}\dl \xi
		=\int_{\R^n} \redconvker(\xi) e^{2\pi i |x|\xi_1}\dl \xi\\
		 & = \frac{1}{(2\pi i|x|)^{n+1}} \left[ \int_{B_{1}} + \int_{B_{1}^c} \right] \p^{n+1}_{\xi_1}\redconvker(\xi)e^{2\pi i |x|\xi_1}\dl \xi, \qquad x\in \R^n.
		\end{align*}
		With the choice \(\alp=(n+1,0,\dots,0)\in \N_0^n\) in \eqref{eq:der_small_xi_new1}, we infer that 
		\begin{equation*}
		|\p^{n+1}_{\xi_1} \redconvker(\xi)|
		\le C,
        \qquad
        \xi\in B_1,
		\end{equation*}
		and from \eqref{eq:first_der:est} we get
		\begin{equation*}
		\Big|\p^{n+1}_{\xi_1}\redconvker\left(\xi\right)\Big|
		\le C|\xi|^{-(n+1)},
        \qquad  \xi\in B_1^c.
		\end{equation*}
		Hence, we can estimate 
		\begin{equation*}
		 |  \Fcal^{-1} \redconvker(x)| %=\Big|\int_{\R^n} \redconvker(\xi)e^{2\pi i x\cdot \xi}\dl \xi\Big|
		\le C\frac{1}{|x|^{n+1}}\Bigg[1+\int_{B_{1}^c} |\xi|^{-(n+1)}\dl \xi\Bigg]
		\le C\frac{1}{|x|^{n+1}},
        \qquad  x\in B_1^c,
		\end{equation*}
		proving \eqref{eq:away_orign}.% in the case \(\delta\in (0,1)\).
        
        %For \(\delta\in (1,\infty)\), we employ the change of variables \(\xi\mapsto\delta^{-1}\xi\) ... \comment{New approach coming}
        %\eqref{eq:der_small_xi_new1}, we infer that 
		%\begin{equation*}
		%|\p^{n+1}_{\xi_1} \redconvker(\xi)|
		%\le C,
        %\qquad
        %\xi\in B_1,
		%\end{equation*}
		%and from 
		%Hence, we can estimate 
		%\begin{equation*}
		 % | \Fcal^{-1} \redconvker(x)| %=\Big|\int_{\R^n} \redconvker(\xi)e^{2\pi i x\cdot \xi}\dl \xi\Big|
		%\le C\frac{1}{|x|^{n+1}}\Bigg[1+\int_{B_{1}^c} |\xi|^{-(n+1)}\dl \xi\Bigg]
		%\le C\frac{1}{|x|^{n+1}},
        %\qquad  x\in B_1^c,
		%\end{equation*}
		%proving \eqref{eq:away_orign}.
        %\comment{end}

        %\comment{New:}
        \textit{Step~5: proof of \eqref{eq:delta_large}.}
        %The proof of \eqref{eq:away_orign} cannot be extended to the range \(\delta\in (1,\infty)\), since estimate \comment{...} breaks down.
       As already pointed out, we resort for the proof of \eqref{eq:delta_large} to a different technique. We use that it follows from the relation $f_\rho^\delta=\delta^{n/2} \Fcal^{-1} \redconvker$ and the composition of the changes of variables \(\xi\mapsto \delta^{-1/2}\) and \(x\mapsto \delta^{1/2} x\) that
        \begin{equation}\label{eq:dl_first}
            \int_{\R^n} |f_\rho^\delta(x)|\dl x
            %= \int_{\R^n}\delta^{n/2} |\Fcal^{-1} \redconvker(x)|\dl x
            = \delta^{n/2} \int_{\R^n} |(\Fcal^{-1} \redconvker (\delta^{-1/2}\cdot))(x)| \dl x.
        \end{equation}
        For simplicity, we define 
        \begin{equation*}
            \redconvkerrs(\xi)
            \coloneqq \redconvker (\delta^{-1/2}\xi)
            = e^{(-|\xi|^2\widehat{Q_\rho}(\delta^{-1/2} \xi))/(4\pi^2)},
            \qquad\xi\in \R^n.
        \end{equation*}
        It then follows from \eqref{eq:dl_first} and the rotational symmetry of \(\widehat{Q_\rho}\) that 
        \begin{equation*}
            \norm{f_\rho^\delta}{L^1(\R^n)}
            = \delta^{n/2} \norm{\Fcal \redconvkerrs}{L^1(\R^n)}.
        \end{equation*}
        Next, we choose some \(k>n/2\) and observe that it holds by Hölder's inequality and Plancherel's identity that
        \begin{align*}
            \norm{\Fcal \redconvkerrs}{L^1(\R^n)}
            &\le \norm{(1+|\cdot|^2)^{-k/2}}{L^2(\R^n)}\norm{ (1+|\cdot|^2)^{k/2}\Fcal \redconvkerrs}{L^2(\R^n)}\\
            &=\norm{(1+|\cdot|^2)^{-k/2}}{L^2(\R^n)}\norm{ \redconvkerrs}{W^{k,2}(\R^n)},
        \end{align*}
        where \(\norm{(1+|\cdot|^2)^{-k/2}}{L^2(\R^n)}\) is finite. In the last step, we also used the characterization of the integer order Sobolev spaces \(W^{k,2}(\R^n)\) as Bessel potential spaces (see, e.g. \cite[Section~1.3.1]{Gr14m}). We conclude that in order to prove \eqref{eq:delta_large}, it suffices to show that \
        \begin{equation}\label{eq:dl_second}
            \norm{ \redconvkerrs}{W^{k,2}(\R^n)}
            <\infty.
        \end{equation}
        Splitting \(\R^n\) into \(B_{\delta^{1/2}}\) and \(B_{\delta^{1/2}}^c\), we use \eqref{eq:large_delta1} and \eqref{eq:large_delta2} to find for all \(\alp\in \N_0^n\) with \(|\alp|\le k\) the estimates
        \begin{equation*}
            \int_{B_{\delta^{1/2}}} |\p^\alpha \redconvkerrs(\xi)|^2\dl \xi
            \le C_\alp^2\int_{\R^n}(1+|\xi|)^{4|\alp|} e^{-2c|\xi|^2}\dl \xi
        \end{equation*}
        and
        \begin{equation*}
            \int_{B_{\delta^{1/2}}} |\p^\alpha \redconvkerrs(\xi)|^2\dl \xi
            \le C_\alp^2 \int_{\R^n}e^{-2c|\xi|^{2\sigma}}\dl \xi.
        \end{equation*}
        Since the right-hand sides are bounded independently of \(\delta\in (1,\infty)\), we indeed verify \eqref{eq:dl_second}.
		\end{proof}

    %\commentCK{Make the message of this remark clear and streamline. Maybe a  introductory sentence would also help to provide context? The following discussed a potential alternative proof strategy in lines the classical approach in the standard Bessel case. }
		\begin{remark}\label{rmk:fi_pos}
        We comment here on the difference in the constriction of the functions \((g_{\rho,\ell})_l\) between the general nonlocal gradients case and the Riesz fractional gradient with kernel \(\rho_s\) (see Section~\ref{sec:discussion}). Precisely, for the Riesz fractional gradient it can be quickly observed that \(g_{\rho_s,\ell}\) is non-negative and that \(\norm{g_{\rho_s,\ell}}{L^1(\R^n)}=1\) for any \(\ell\in \N\) (see~\cite[Chapter~7]{Mi96}). In particular, the integrability already follows from the non-negativity, since then
			\begin{equation*}
				\norm{g_{\rho_s,\ell}}{L^1(\R^n)}
				 = \widehat{g_{\rho_s, \ell}}(0) = (1+4\pi^2|\xi|^2)^{-s/2}|_{\xi=0}  =1,
                 \qquad \ell\in \N.
                %... (1+4\pi^2 \abs{\lambda_{\rho}}^2)^{-{\ell}/{2}}|_{\xi=0}=1
			\end{equation*}
			%for all \(\ell\in \N\).
            Following the same idea for any \(\rho\) as in Proposition~\ref{prop:invers_fourier}, one encounters the problem that the non-negativity of the \(g_{\rho,\ell}\) seems difficult to verify in the general case. In fact, with \eqref{eq:aux_h} and \eqref{eq:def_conv_kernel}, one would like to show that \(f^\delta_\rho= \mathcal{F}^{-1} h^\delta_\rho\) is non-negative, where, as the function $\abs{\lambda_{\rho}}$ is even and with the help of the coarea formula, the expression for \(f^\delta_\rho\) can be written as an oscillating integral of \(\cos\)-type, this is
            %the integral            
			%\begin{equation*}
			%\int_{\R^n} e^{-\abs{\lambda_{\rho} (\xi)}^2\delta / (4 \pi^2)}\sin(2\pi i x\cdot \xi)\dl \xi
			%\end{equation*}
			%vanishes for any \(x\in \R^n\), as it has an odd integrand (the finiteness of the integral is due to \eqref{eq:lest_qrho_pre}).
			%With this, using the evenness of the remaining integrand, and by virtue of the coarea formula, the function \(f^\delta_\rho\) can be rewritten as
			\begin{equation*}
			f^\delta_\rho(x)
			=2\delta^{n/2}\int_{\sph^{n-1}_{+}}\int_{0}^\infty r^{n-1}e^{- r^2q(r)^2\delta}\cos(2\pi r |x|z) \dl r\dl \mathcal{S}(z),
            \qquad x\in \R^n,
			\end{equation*}
			denoting by \(q\) the radial representative of \(\widehat{Q_\rho}\). 
            But this remains an open problem for the moment, and instead, we demonstrate that the functions \(g_{\rho,\ell}\) lie in $L^1(\R^n)$ via explicit decay and norm estimates for $f^{\delta}_{\rho}$ (cf.~\eqref{eq:close_orign}, \eqref{eq:away_orign} and~\eqref{eq:delta_large}).
            
            The positivity of the functions \(g_{\rho,\ell}\) is also interesting in a different regard: if \(g_{\rho,\ell}\) is positive, then it defines a probability distribution and \((1+ \abs{\lambda_{\rho}}^2)^{-{\ell}/{2}}\) can be seen as its characteristic function (see again Section~\ref{sec:discussion} for the analogous statement in the case of the Bessel potential spaces). Connected to the above-mentioned difficulties, the question if this holds in the case of general nonlocal gradients also remains open.
            \end{remark}
	
	We now turn to the main result of this section. It relies on the alternative characterization of \(H^{\rho,p}(\R^n)\) obtained in Theorem~\ref{thm:alt_char} together with the existence of $g_{\rho,1}$ according to Proposition~\ref{prop:invers_fourier}.
	\begin{theorem}[Potential space characterization]\label{thm:pot_char}
	Let \(p\in (1,\infty)\) and assume that \(\rho\) has compact support and satisfies assumptions (H0)--(H3). Denote by \(g_{\rho}\in L^1(\R^n)\) the function \(g_{\rho,1}\) from Proposition~\ref{prop:invers_fourier} such that
	\begin{equation*}
	\widehat{g_{\rho}}
	=\bla{\lambda_\rho}^{-1}.
	\end{equation*} 
	Then
	\begin{equation*}\label{eq:char_pot_thm}
	H^{\rho,p}(\R^n)
	=\{g_\rho*f:\, f\in L^p(\R^n)\}
	\end{equation*}
	and, for all $f\in L^p(\R^n)$,
	\begin{equation*}
	\norm{\mathcal{F}^{-1}(\bla{\lambda_\rho}(\widehat{g_\rho*f}))}{L^{p}(\R^n)}
	=\norm{f}{L^p(\R^n)}.
	\end{equation*}
\end{theorem}
\begin{proof}
We establish the statement by showing that the map
\[
 L^p(\R^n) \to H^{\rho,p}(\R^n) , \qquad f \mapsto g_\rho*f
\]
is a bounded linear isomorphism with inverse
\[
	H^{\rho,p}(\R^n) \to L^p(\R^n), 
	\quad u\mapsto \mathcal{F}^{-1}(\bla{\lambda_\rho}\hat{u}),
\]
where the last expression is understood in the sense of \(\mathcal{S}'(\R^n)\) and lies in \(L^p(\R^n)\) thanks to Theorem~\ref{thm:alt_char}. 
We first prove that $g_{\rho} * f \in H^{\rho,p}(\R^n)$ for any $f \in L^p(\R^n)$.
Firstly, by Young's inequality, $g_{\rho} * f \in L^p(\R^n) $. Secondly, for any \(\varphi\in \mathcal{S}(\R^n)\),
\begin{align*}
    \langle \mathcal{F}^{-1}(\bla{\lambda_\rho}(\widehat{g_\rho*f})), \varphi\rangle_{\mathcal{S}'\times \mathcal{S}}
    &= \langle g_\rho*f , \mathcal{F}(\bla{\lambda_\rho}\mathcal{F}^{-1}\varphi)\rangle\\
    &= \langle g_\rho*f , \mathcal{F}^{-1}(\bla{\lambda_\rho}\hat{\varphi})\rangle\\
    &= \langle f , g_\rho*\mathcal{F}^{-1}(\bla{\lambda_\rho}\hat{\varphi})\rangle\\
    &= \langle f, \mathcal{F}^{-1}(\widehat{g_\rho}\, \bla{\lambda_\rho}\hat{\varphi})\rangle =  \langle f , \varphi\rangle,
\end{align*}
where we used additionally the first identity of Lemma~\ref{lem:fourier_identities}~\textit{\ref{item:fi_symm}} in the second equality, and for the last equality that \(\mathcal{F}^{-1}(\bla{\lambda_\rho}\hat{\varphi})\in \mathcal{S}(\R^n)\).  
Thus, by Theorem \ref{thm:alt_char}, $g_{\rho} * f \in H^{\rho,p}(\R^n)$.  

To see the injectivity of \(f \mapsto g_\rho*f\), we assume that \(g_\rho*f=0\) for some \(f\in L^p(\R^n)\). 
It follows then directly from Theorem~\ref{thm:alt_char} that 
\begin{equation*}
    0=\norm{g_\rho * f}{H^{\rho,p}(\R^n)}
    \ge c_p\norm{\Fcal^{-1}(\bla{\lambda_\rho}\, \widehat{g_\rho}\widehat{f})}{L^p(\R^n)}
    = c_p \norm{f}{L^p(\R^n)},
\end{equation*}
implying that \(f\equiv 0\). Consequently, \(f \mapsto g_\rho*f\) is injective.
%On the one hand, for any sequence \((f_j)_j\subset \mathcal{S}(\R^n)\) such that \(f_j\to f\) in \(L^p(\R^n)\) it holds (\(f\mapsto g_\rho * f\) is continuous by Young's inequality)
%\begin{equation*}
%    \lim_{j\to \infty} g_\rho * f_j = 0 \qquad\text{in \color{red} $L^p(\R^n)$}.
%\end{equation*}
%On the other hand, by Theorem~\ref{thm:alt_char},
%\begin{equation*}
%    \norm{g_\rho * f_j}{H^{\rho,p}(\R^n)}
%    \ge c_p\norm{\Fcal^{-1}(\bla{\lambda_\rho}\, \widehat{g_\rho}\widehat{f_j})}{L^p(\R^n)}
%    = c_p \norm{f_j}{L^p(\R^n)}.
%\end{equation*}
%\commentCK{To check this conclusion!! And spaces of convergence?!}
%It follows that \(f\equiv 0\), i.e., \(f \mapsto g_\rho*f\) is injective.

Turning our attention to surjectivity, we note that by the density of \(C_c^\infty (\R^n)\) according to \cite[Theorem~3.9]{BMCS25}, there exists for every \(u\in H^{\rho,p}(\R^n)\) a sequence \((u_j)_j\subset C_c^\infty(\R^n)\) such that \(u_j \to u\) in \(H^{\rho, p}(\R^n)\) and, with Theorem~\ref{thm:alt_char},
\begin{equation*}
    \mathcal{F}^{-1}(\bla{\lambda_\rho}\widehat{u_j})\to \mathcal{F}^{-1}(\bla{\lambda_\rho}\hat{u})
    \quad \text{in } L^p(\R^n). 
\end{equation*}
Using that convolution with \( g_\rho\) is continuous by Young's inequality, we infer that
\begin{equation*}
    g_\rho *\mathcal{F}^{-1}(\bla{\lambda_\rho}\widehat{u_j})
    \to g_\rho *\mathcal{F}^{-1}(\bla{\lambda_\rho}\hat{u})
    \quad \text{in } L^p(\R^n),
\end{equation*}
while 
\begin{equation*}
    g_\rho *\mathcal{F}^{-1}(\bla{\lambda_\rho}\widehat{u_j})
    = u_j \to u 
    \quad \text{in } L^p(\R^n). 
\end{equation*}
We conclude that \(g_\rho *\mathcal{F}^{-1}(\bla{\lambda_\rho}\hat{u})= u\). The proof is finished.
\end{proof}
%\begin{remark}
%	\comment{Does this characterization provide quick noteworthy proofs of embeddings or inequalities?}
%   \commentC{Shieh and Spector only use it to record Theorem 2.2, which is Adams Chapter 7, p. 221. Warning: this is Adams first edition. In the second (and more famous) one (Adams Fournier) this doesn't appear. We have a copy in the UAM library. I never found it on the web (although the second edition is easily found)}
%\end{remark}

Having the sequence of functions \( ( g_{\rho,\ell})_{\ell}\) uniformly bounded in \(L^1(\R^n)\) at our disposal, an alternative proof of Theorem~\ref{thm:pot_char} can be obtained following the arguments in \cite[Chapter~V.3]{St70} (see also \cite[Chapter~7]{Mi96} and \cite[Theorem~1.7]{ShSp15}), where the coincidence of the nonlocal Sobolev spaces \(H^{s,p}(\R^n)\) with the Bessel potential spaces is shown. The only missing ingredient is the subsequent result regarding the fractions of the generalized Riesz and Bessel potentials. 
%\commentCK{insert little reference for Fourier(-Stieltjes) inverses or say in which of these references this appears. Maybe it is ok, though....} 
The proof is analogous to the one of Stein in \cite[Chapter~V.3, Lemma~2]{St70} in the case of the classical Riesz and Bessel potentials (see also \cite[Section~7, Lemma~2.1 and Lemma~2.2]{Mi96}).
        
        \begin{lemma}\label{lem:gen_pot_frac_inverse}
            Let \(p\in (1,\infty)\) and assume that \(\rho\) has compact support and satisfy assumptions (H0)--(H3). 
		  Then the fractions
         \begin{equation}\label{fractions98}
			\frac{\abs{\lambda_\rho}}{\bla{\lambda_\rho}}
			\quad\text{ and }\quad 
			\frac{\bla{\lambda_\rho}}{1+\abs{\lambda_\rho}}
			\end{equation}
			admit Fourier
            %(-Stieltjes) 
            inverses in the sense of finite measures.
        \end{lemma}
        \begin{proof}
            The starting point of the proof is the binomial series
            \begin{equation}\label{binomial}
                (1-t)^{1/2}
                = 1+\sum_{j=1}^\infty a_jt^j ,
                \qquad \abs{t} < 1 ,
            \end{equation}
for suitable coefficients \((a_j)_j\subset \R\) with \(\sum_{j=1}^\infty|a_j|<\infty\). This enables us to write
            \begin{equation*}
                \frac{\abs{\lambda_\rho(\xi)}}{\bla{\lambda_\rho(\xi)}}
                = 1+\sum_{j=1}^\infty a_j(1+|\lambda_\rho(\xi)|^2)^{-j}, 
                \qquad \xi \in \R^n\setminus \{0\},
            \end{equation*}
            where we chose \(t=1/\bla{\lambda_\rho(\xi)}^2\) for $\xi \in \R^n \setminus \{0\}$ in~\eqref{binomial}.
Let $(g_{\rho,\ell})_{\ell}$ be the family of  Proposition~\ref{prop:invers_fourier}.
Then $(g_{\rho,\ell})_{\ell}$ is bounded in $L^1 (\R^n)$ and \( \widehat{g_{\rho,2j}} = (1+|\lambda_\rho|^2)^{-j}\) for all \(j\in \N\). Consequently, defining 
                \begin{equation*}
                    \mu_1
                    \coloneqq \delta_0+\sum_{j=1}^\infty a_jg_{\rho,2j}\dl \mathcal{L}^n,
                \end{equation*}
                where \(\delta_0\) is the Dirac measure at the origin and \(\mathcal{L}^n\) denotes the Lebesgue measure on \(\R^n\), it holds
                \begin{equation*}
                    \mu_1(\R^n)<\infty
                    \quad\text{and}\quad
                    \widehat{\mu_1}
                    = \frac{\abs{\lambda_\rho}}{\bla{\lambda_\rho}}.
                \end{equation*}

As for the second fraction in~\eqref{fractions98},  we invoke Wiener's theorem (see e.g., \cite[Theorem~2.3]{Mi96}), according to which for any \(\Phi_1\in L^1(\R^n)\) such that \(1+\widehat{\Phi_1}(\xi)\ne 0\) for all \(\xi\in \R^n\), there exists \(\Phi_2\in L^1(\R^n)\) such that  
                \begin{equation*}
                    (1+\widehat{\Phi_1})^{-1}
                    = 1+\widehat{\Phi_2}.
                \end{equation*}
                In particular, for 
                \begin{equation*}
                    \Phi_1
                    \coloneqq g_{\rho,1}+\sum_{j=1}^\infty a_jg_{\rho,2j},
                \end{equation*}
                it holds
                \begin{equation*}
                    1+\widehat{\Phi_1}
                    = \frac{1+\abs{\lambda_\rho}}{\bla{\lambda_\rho}}>0.
                \end{equation*}
            Inverting this expression with the help of Wiener's theorem yields a function \(\Phi_2\in L^1(\R^n)\) such that, with 
                \begin{equation*}
                    \mu_2
                    \coloneqq \delta_0+\Phi_2\dl\mathcal{L}^n,
                \end{equation*}
            we deduce 
                \begin{equation*}
                    \widehat{\mu_2}
                    = \frac{\bla{\lambda_\rho}}{1+\abs{\lambda_\rho}},
                \end{equation*}
                which concludes the proof.
        \end{proof}
  
          We close with a remark regarding an alternative proof strategy  (under slightly stronger assumptions) for the  Fourier multiplier result Proposition~\ref{prop:fourier_multipliers}. 
        \begin{remark} 
        In the case that \(\rho\) has compact support and fulfills (H0)--(H3), in order to prove Proposition~\ref{prop:fourier_multipliers}, one can also argue more directly by means of Proposition~\ref{prop:invers_fourier} and Lemma~\ref{lem:gen_pot_frac_inverse}. Indeed, for Proposition~\ref{prop:fourier_multipliers}~\textit{(iii)} one can use Proposition~\ref{prop:invers_fourier} to write for any \(u\in \mathcal{S}(\R^n)\) (with \({g_{\rho}}\coloneqq {g_{\rho,1}}\))
		\begin{equation*}
			\mathcal{F}^{-1}(\bla{\lambda_\rho}^{-1} \hat{u})
			= g_{\rho}*u,
		\end{equation*}
		and by Young's convolution inequality it holds immediately that
		\begin{equation*}
			\norm{{g_{\rho}}*u}{L^p(\R^n)}
			\le \norm{{g_{\rho}}}{L^1(\R^n)}\norm{u}{L^p(\R^n)}.
		\end{equation*}
		%As usual, one concludes that the multiplier \((1+4\pi^2 \abs{\lambda_{\rho} }^2)^{-{1}/{2}}\) can indeed be extended to a bounded operator in \(L^p(\R^n)\).

        A similar approach proves Proposition~\ref{prop:fourier_multipliers}~\textit{(i)} and \textit{(ii)}, since  the respective fractions actually admit Fourier inverses in the sense of finite measures by Lemma~\ref{lem:gen_pot_frac_inverse}. 
    \end{remark}

	%%%%%%%%%%%%%%%%%%%%%%%%%%%%%%%%%%%%%%%%%%%%%%%%%%%%%%%%%%%%%%%%%
    \appendix

\addtocontents{toc}{\protect\setcounter{tocdepth}{1}}

\setcounter{equation}{0}
\renewcommand{\thesubsection}{\Alph{subsection}}

\section{Proofs of Lemmas~\ref{lem:estimates} and \ref{lem:estimates_exp}}\label{sec:appendix_ests}

Recall that \(\lambda_\rho\), the Fourier symbol of the nonlocal gradient \(D^\rho_0\), is given by the expression (see Section~\ref{sec:fourier})
     \begin{align*}
     \lambda_\rho(\xi) = 2\pi i \xi \widehat{Q_{\rho}}(\xi), \qquad \xi\in \R^n. 
     \end{align*}  
     
In the following, we will denote for some multi-index \(\alp\in \N_0^n\) by \(C_\alp>0\) a positive constant that may change from line to line and only depends on \(\alp\) (and implicitly on the dimension \(n\in \N\)). %\commentCK{Careful about notation of multiindices and definition of $\N$. With the current convention (part of the in notations section??) would need to write $\alpha\in \N^n_0$. }
	
	\begin{proof}[Proof of  Lemma~\ref{lem:estimates}]
        For \(\alp=0\) it is readily seen that all inequalities are fulfilled with suitable constants \(C_0>0\). We point out that regarding~\eqref{eq:fourth_estimate}, one needs to use that \(|\xi|\widehat{Q_\rho}(\xi)\ge C\) for $\xi\in B_1^c$ due to the assumptions and \eqref{eq:lower_prelim}, and thus, \(1+\abs{\lambda_{\rho}}^2\le C_0 \abs{\lambda_{\rho}}^2\) in $B_1^c$. 
        
		Pick now some \(\alpha\in \N^n_0\setminus\{0\}\).  We repeat that, since by assumption \(\rho\) has compact support, the same is true for \(Q_\rho\), and \(\widehat{Q_\rho}\) is smooth and, in particular, all subsequently computed derivatives will be well defined. 
		
		We start by recalling a version of Faà di Bruno's formula for the derivatives of compositions, which we will repeatedly rely on it throughout the proof. Precisely, one has for any \(f\in C^\infty(\R^n)\) and \(g:\R\to \R\) such that \(g\circ f \in C^\infty(\R)\), and any multi-index \(\alp\in \N^n_0\setminus\{0\}\) the estimate 
		\begin{equation}\label{eq:FdB}
		|\p^\alp (g\circ f)(\xi)|
		\le C_\alp \sum_{k=1}^{|\alp|} |((\p^k g)\circ f)(\xi)|
		\sum _{\{\alp'_1,\dots,\alp'_k\}}  \prod_{l=1}^k |\p^{\alp'_l} f(\xi)|, \qquad \xi\in \R^n,
        %B_{\alp^'_1\cdots\alp^'_k},
		\end{equation}
		where  
		%\begin{equation*}
	%	B_{\alp^'_1\cdots\alp^'_k}
	%	= \prod_{l=1}^k |\p^{\alp^'_l} f(\xi)|,
	%	\end{equation*}
		and the second sum runs over all partitions \(\{\alp'_1,\dots,\alp'_k\}\) of \(\alp\) into \(k\) multi-indices. By a partition of \(\alp\), we mean a set \(\{\alp'_1, \dots,\alp'_N\}\) of multi-indices \(\alp'_i\in \N^n_0\setminus\{0\}\), $N\in \N$ with \(N\le |\alp|\), such that \(\alp'_1+\ldots+\alp'_N=\alp\), where the summation is understood component-wise in the sense of multi-indices. In particular, the partitions \(\{\alp'_1,\dots,\alp'_k\}\) in the inner sum in \eqref{eq:FdB} are subject to the constraint
	\begin{equation}\label{eq:FdB_sum_constraint}
		|\alp'_1|+\ldots +|\alp'_k|=|\alp|
		\end{equation}
		for each $k\in\N$ with $k\leq |\alpha|$. Note that a partition \(\{\alp'_1,\dots,\alp'_k\}\) might contain the same multi-index multiple times. 
        
        As a particular case of the general formula, we will use in the following that for positive \(f\in C^\infty(\R^n)\) and \(g=(\cdot)^\gamma\) with some \(\gamma \in\R\) it holds
		\begin{equation}\label{eq:FdB_exp}
		|\p^\alp f^\gamma(\xi)|
		\le C_\alp \sum_{k=1}^{|\alp|} |f^{\gamma -k}(\xi)|
		\sum _{\{\alp'_1,\dots,\alp'_k\}} \prod_{l=1}^k |\p^{\alp'_l} f(\xi)|, \qquad \xi\in \R^n.
		\end{equation}

        Note that we will additionally use these elementary estimates, valid for every $\beta \in \mathbb{N}_0^n$,
        %we will additionally use the following elementary estimates. Let \(\beta\in \N^n_0\), then it holds 
		\begin{equation}\label{eq:first_aux_est}
        |\p^\beta \xi|\le |\xi|^{1-|\beta|},
		|\p^\beta |\xi|^2|\le 2|\xi|^{2-|\beta|},
		\qquad \xi\in \R^n\setminus\{0\}.
		\end{equation}
%		And then, second, for any \(\gamma<1\) and any \(k\ge 1\) we have the immediate inequality 
%		\begin{equation}\label{eq:second_aux_est}
%		|(1+4\pi^2|\xi|^2)^{\gamma-k}|
%		\le |\xi|^{2(\gamma-k)}.
%		\end{equation}

Now we are in the position to derive the four stated estimates~\eqref{eq:second_estimate}-\eqref{eq:fifth_estimate}, which are proven separately in the following steps.
		
		\textit{Step~1: estimate~\eqref{eq:second_estimate}.}
        With the aid of Leibniz' formula we compute 
		\begin{equation}\label{eq:leibniz_first}
		|\p^\alp (\xi\widehat{Q_\rho}(\xi))|
		\le C_\alp \sum_{\beta\le \alp}|\p^\beta \xi||\p^{\alp-\beta}\widehat{Q_\rho}(\xi)|, \qquad \xi\in \R^n.
		\end{equation}
		In combination with \eqref{eq:first_aux_est} and \eqref{eq:uest_qrho_alp_pre}, i.e. \(|\p^{\alp-\beta}\widehat{Q_\rho}(\xi)|\le C_\alp |\xi|^{-|\alp|+|\beta|}\widehat{Q_\rho}(\xi)\) for all \(  \xi\in B_1^c\), it follows that 
		\begin{equation*}
		|\p^\alp (\xi\widehat{Q_\rho}(\xi))|
		\le C_\alp \sum_{\beta\le \alp} |\xi|^{1-|{\beta}|}|\xi|^{-|\alp|+|\beta|}\widehat{Q_\rho}(\xi)
		\le C_\alp|\xi|^{1-|\alp|}\widehat{Q_\rho}(\xi),
		\qquad \xi\in B_1^c.
		\end{equation*}

		As a variant of \eqref{eq:leibniz_first}, we also have
		\begin{equation*}
		|\p^\alp (|\xi|\widehat{Q_\rho}(\xi))|
		\le C_\alp \sum_{\beta\le \alp}|\p^\beta |\xi|||\p^{\alp-\beta}\widehat{Q_\rho}(\xi)|, \qquad \xi\in \R^n.
		\end{equation*}
		With Faà di Bruno's formula in \eqref{eq:FdB_exp} for \(\gamma=1/2\) and \eqref{eq:first_aux_est}, one readily obtains that 
		\begin{align*}
		\begin{split}
		    |\p^{\tilde{\alp}} |\xi||
		= |\p^{\tilde{\alp}} (|\xi|^2)^{1/2}|
		&\le C_{{\alp}} \sum_{k=1}^{|\tilde{\alp}|}(|\xi|^2)^{(1/2-k)}\sum _{\{\tilde{\alp}'_1,\dots,\tilde{\alp}'_k\}} \prod_{l=1}^k |\p^{\tilde{\alp}'_l} |\xi|^2|\\
		&\le  C_{{\alp}} \sum_{k=1}^{|\tilde{\alp}|}(|\xi|^2)^{(1/2-k)}\sum _{\{\tilde{\alp}'_1,\dots,\tilde{\alp}'_k\}}\prod_{l=1}^k  |\xi|^{2-|\tilde{\alp}'_l|}\\
		&\le C_{{\alp}} |\xi|^{1-|\tilde{\alp}|},
        \qquad \xi\in \R^n
		\end{split}
		\end{align*}
		for any $\tilde \alpha\in \N_0^n$ with partition $\{\tilde{\alpha}'_1,\ldots , \tilde{\alpha}'_k\}$ and \(\tilde{\alp}\le\alp\).
		Thus, another application of \eqref{eq:uest_qrho_alp_pre} leads to
		\begin{equation*}
		|\p^\alp (|\xi|\widehat{Q_\rho}(\xi))|
		\le C_\alp \sum_{\beta\le \alp} |\xi|^{1-|{\beta}|}|\xi|^{-|\alp|+|\beta|}\widehat{Q_\rho}(\xi)
		\le C_\alp|\xi|^{1-|\alp|}\widehat{Q_\rho}(\xi),
		\qquad \xi\in B_1^c.
		\end{equation*}
		
		\textit{Step~2: estimate~\eqref{eq:third_estimate}.} %\comment{add general \(\gamma\) case}
		%Now we invoke Faà di Bruno's formula in \eqref{eq:FdB_exp} for \(\gamma=-1\). Thus, 
		Now we invoke Faà di Bruno's formula in \eqref{eq:FdB_exp} for general  \(\gamma\in \R\). Thus, 
        \begin{align*}
		|\p^\alp (|\xi|\widehat{Q_\rho}(\xi))^{\gamma}|
		&\le C_{{\alp}} \sum_{k=1}^{|{\alp}|}(|\xi|\widehat{Q_\rho}(\xi))^{\gamma-k}\sum _{\{\alp'_1,\dots,\alp'_k\}} \prod_{l=1}^k |\p^{\alp'_l} (|\xi|\widehat{Q_\rho}(\xi)|)\\
		&\le C_{{\alp}} \sum_{k=1}^{|{\alp}|}(|\xi|\widehat{Q_\rho}(\xi))^{\gamma-k}\sum _{\{\alp'_1,\dots,\alp'_k\}} \prod_{l=1}^k |\xi|^{1-|\alp'_l|}\widehat{Q_\rho}(\xi)\\
		&\le C_\alp |\xi|^{\gamma-|\alp|}\widehat{Q_\rho}^{\gamma}(\xi),
		\qquad \xi\in B_1^c,
		\end{align*}
		validating the first estimate of ~\eqref{eq:third_estimate}. The second estimate of ~\eqref{eq:third_estimate}, concerning bounds for \(\p^\alp (1+|\xi|\widehat{Q_\rho}(\xi))^{\gamma}\), is obtained similarly. It suffices to observe that \((1+|\xi|\widehat{Q_\rho}(\xi))^{\gamma-k}\le C_\alp (|\xi|\widehat{Q_\rho}(\xi))^{\gamma-k}\) for all \(\gamma\in \R\) and all \(k\in\N\) with \(k\le |\alp|\), where in the case that \(\gamma-k>0\) we use that \(1+\abs{\lambda_{\rho}}^2\le C_\alp \abs{\lambda_{\rho}}^2\) in $B_1^c$ according \eqref{eq:lower_prelim}. 
		
		\textit{Step~3: estimate~\eqref{eq:fourth_estimate}.}
		Once more appealing to Leibniz' formula, entails
		\begin{equation}\label{eq:sec_est_initial}
		|\p^\alp (1+\abs{\lambda_{\rho} (\xi)}^2)|
		\le C_\alp \sum_{\beta\le \alp}|\p^\beta |\xi|^2||\p^{\alp-\beta}\widehat{Q_\rho}^2(\xi)|,
        \qquad \xi\in \R^n.
		\end{equation}
		To further estimate the right-hand side, we apply Faà di Bruno's formula in \eqref{eq:FdB_exp} to the function \(\widehat{Q_\rho}^2\). Together with the estimate \eqref{eq:uest_qrho_alp_pre} for the derivatives of \(\widehat{Q_\rho}\), we infer that 
		\begin{align*}
		|\p^{\tilde{\alp}}\widehat{Q_\rho}^2(\xi)|
		&\le C_{\tilde{\alp}} \sum_{k=1}^{|\tilde{\alp}|} \widehat{Q_\rho}^{2 -k}(\xi)
		\sum _{\{\tilde{\alp}'_1,\dots,\tilde{\alp}'_k\}} \prod_{l=1}^k |\p^{\tilde{\alp}'_l} \widehat{Q_\rho}(\xi)|\\
		&\le  C_{\tilde{\alp}} \sum_{k=1}^{|\tilde{\alp}|} \widehat{Q_\rho}^{2 -k}(\xi)
		\sum _{\{\tilde{\alp}'_1,\dots,\tilde{\alp}'_k\}} \prod_{l=1}^k |\xi|^{-|\tilde{\alp}'_l|}\widehat{Q_\rho}(\xi)\\
		&= C_{\tilde{\alp}}  \widehat{Q_\rho}^{2}(\xi)|\xi|^{-|\tilde{\alp}|},
		\qquad \xi\in B_1^c,
		\end{align*}
	where the last equality invoked again the sum constraint from \eqref{eq:FdB_sum_constraint}. Coming back to~\eqref{eq:sec_est_initial} and using also the second part in \eqref{eq:first_aux_est}, we conclude with  \(|\alp-\beta|=|\alp|-|\beta|\) for \(\beta\le \alp\) that
		\begin{align*}
		|\p^\alp (1+\abs{\lambda_{\rho} (\xi)}^2)|
		\le C_\alp \sum_{\beta\le \alp} |\xi|^{2-|\beta|}\widehat{Q_\rho}^{2}(\xi)|\xi|^{-(|\alp|-|\beta|)}
		= C_\alp|\xi|^{2-|\alp|}\widehat{Q_\rho}^{2}(\xi),
        \qquad \xi\in B_1^c.
		\end{align*}
		This finishes the proof of~\eqref{eq:fourth_estimate}.
		
		\textit{Step~4: estimate~\eqref{eq:fifth_estimate}.}
		We apply once again Faà di Bruno's formula in the version of \eqref{eq:FdB_exp}. In that regard, we note that for \(\gamma\in \R\) and \(k\ge 1\), 
		\begin{equation}\label{eq:third_aux_est}
		(1+\abs{\lambda_{\rho} (\xi)}^2)^{\gamma -k}
		\le |\xi|^{2\gamma-2k}\widehat{Q_\rho}^{2\gamma-2k}(\xi), \qquad \xi\in B_1^c,
		\end{equation}
        where in the case that \(\gamma-k>0\), we used again that, by \eqref{eq:lower_prelim}, it holds \(1+\abs{\lambda_{\rho} (\xi)}^2\le C \abs{\lambda_{\rho} (\xi)}^2\).
		Combining this with the just proved third estimate \eqref{eq:fourth_estimate}, we obtain from \eqref{eq:FdB_exp} and \eqref{eq:FdB_sum_constraint} that
		\begin{align*}
		&|\p^\alp (1+\abs{\lambda_{\rho} (\xi)}^2)^{\gamma}|\\
		&\quad \le  C_\alp \sum_{k=1}^{|\alp|} (1+\abs{\lambda_{\rho} (\xi)}^2)^{\gamma-k}(\xi)
		\sum _{\{\alp'_1,\dots,\alp'_k\}} \prod_{l=1}^k |\p^{\alp'_l} (1+\abs{\lambda_{\rho} (\xi)}^2)|\\
		&\quad \le  C_\alp \sum_{k=1}^{|\alp|} |\xi|^{2\gamma-2k}\widehat{Q_\rho}^{2\gamma-2k}(\xi)
		\sum _{\{\alp'_1,\dots,\alp'_k\}} \prod_{l=1}^k|\xi|^{2-|\alp'_l|}\widehat{Q_\rho}^{2}(\xi)\\
		&\quad \le  C_\alp  |\xi|^{2\gamma-|\alp|}\widehat{Q_\rho}^{2\gamma}(\xi),
        \qquad \xi\in B_1^c,
		\end{align*}
        concluding the proof.
		%		With the help of \eqref{eq:lest_qrho}, which gives the upper bound \(|\widehat{Q_\rho}^{-1}(\xi)|\le C|\xi|^{1-\sigma}\), we finally conclude
		%		\begin{align*}
		%			|\p^\alp (1+4\pi^2\abs{\lambda_{\rho} (\xi)}^2)^{-{1}/{2}}|
		%			\le C_\alp  |\xi|^{-1-|\alp|}\widehat{Q_\rho}^{-1}(\xi)
		%			\le C_\alp |\xi|^{-\sigma-|\alp|}.
		%		\end{align*}
	\end{proof}

  %\commentCK{ Give the reduced version taylored to what is used in the proofs, that is,~\eqref{eq:first_der:est}, \eqref{eq:new_est}, \eqref{eq:der_small_xi}.}
  %\commentCK{Check again range of validity $B_1^c$, what is the role of $\delta_0$?}
   % \commentCK{Adjust introtext text and merge with Step 3...} 
    Next, we collect here some estimates of similar kind  for the function (see~\eqref{Theta_rhodelta})
    	\begin{equation*}
		\redconvker = e^{-\abs{\lambda_{\rho} }^2\delta / (4 \pi^2)}
		\end{equation*} 
   % \begin{equation*}
   %     \redconvker(\xi)
	%	\coloneqq e^{-\abs{\lambda_{\rho} (\xi)}^2\delta / (4 \pi^2)}, \qquad \xi\in \R^n
   % \end{equation*}
    for any fixed \(\delta>0\), which were exploited in the proof of Proposition~\ref{prop:invers_fourier}.
    %The proof can be found again in Appendix~\ref{sec:appendix_ests}.

%\color{blue}    \begin{lemma}\label{lem:estimates_exp}
%        Let \(p\in (1,\infty)\). Assume that \(\rho\) has compact support and satisfies assumptions (H0)--(H2) and (H3'). Then, 
%\begin{equation}\label{eq:lem_gen_der_est}
%	 |\p^\alp \redconvker(\xi)|
%		\le C_\alp |\xi|^{-|\alp|}
%		\sum_{k=1}^{|\alp|}\delta^k \redconvker (\xi)|\lambda_\rho(\xi)|^{2k},
%        \qquad \xi\in B_1^c,
%\end{equation}
%and, for any $\sigma>0$,        
%        \begin{equation}\label{eq:lem_gen_der_est_cov}
%        \Big|\p^\alp \redconvker\left(\frac{\xi}{\delta^{1/(2\sigma)}}\right)\Big|
%		\le C_\alp |\xi|^{-|\alp|}
%		\sum_{k=1}^{|\alp|} \delta^k\redconvker\left(\frac{\xi}{\delta^{1/(2\sigma)}}\right) \Big|\lambda_\rho \left(\frac{\xi}{\delta^{1/(2\sigma)}}\right)\Big|^{2k},
%        \qquad \xi\in B_1^c,
%        \end{equation}
%    for all multi-indices \(\alp\in \N_0^n\), and for some constant \(C_\alp>0\).
%    \end{lemma}

%    \color{black}
    \begin{lemma}\label{lem:estimates_exp}
        Let \(p\in (1,\infty)\) and \(\alp\in \N_0^n\). Assume that \(\rho\) has compact support and satisfies assumptions (H0)--(H3), 
        %\commentCK{We now assume (H3) instead of (H3'), ok?}
        and let
$ \gamma = (4 - 3 \sigma)^{-1}\in (1/4, 1)$, where $\sigma$ is the constant from (H3). Then, for all multi-indices \(\alp\in \N_0^n\),
\begin{alignat}{2}
&|\p^\alp \redconvker\left(\xi\right)|
\le C_\alp |\xi|^{-|\alp|},
\qquad &&\xi\in B_1^c,\, \delta\in (0,\infty),\label{eq:first_der:est}\\
&\Big|\p^\alp \redconvker(\delta^{-1/(2\sigma)}\xi)\Big|
\le C_{\alp} |\xi|^{-(|\alp|+\gamma \sigma)},
\qquad &&\xi\in B_{\delta^{1/(2\sigma)}}^c,%\in B_1^c
\, \delta\in (0,\infty), \label{eq:new_est}\\
&|\p^\alp \redconvker\left(\xi\right)|
\leq C_{\alp} , 
\qquad &&\xi\in B_1,\, \delta\in (0,1), \label{eq:der_small_xi_new1}\\
 &\Big|\p^\alp \redconvker(\delta^{-1/2}\xi)\Big|
\le C_{\alp}   (1+|\xi|)^{2|\alp|} e^{-c|\xi|^2},
\qquad &&\xi\in B_{\delta^{1/2}},\, \delta\in (1,\infty),\label{eq:large_delta1}\\
 &\Big|\p^\alp \redconvker(\delta^{-1/2}\xi)\Big|
\le C_{\alp}   e^{-c|\xi|^{2\sigma}},
\qquad  &&\xi\in B_{\delta^{1/2}}^c,\, \delta\in (1,\infty), \label{eq:large_delta2}
\end{alignat}
     % with $\gamma \in (1/4, 1)$  defined as $ \gamma \coloneq (4 - 3 \sigma)^{-1}$ with $\sigma$ from (H3'), 
     %and 
     %old-old
%\begin{equation}\label{eq:lem_gen_der_est}
%	 |\p^\alp \redconvker(\xi)|
%		\le C_\alp |\xi|^{-|\alp|}
%		\sum_{k=1}^{|\alp|}\delta^k \redconvker (\xi)|\lambda_\rho(\xi)|^{2k},
%        \qquad \xi\in B_1^c,
%\end{equation}
%and, for any $\sigma>0$,        
%        \begin{equation}\label{eq:lem_gen_der_est_cov}
%        \Big|\p^\alp \redconvker\left(\frac{\xi}{\delta^{1/(2\sigma)}}\right)\Big|
%		\le C_\alp |\xi|^{-|\alp|}
%		\sum_{k=1}^{|\alp|} \delta^k\redconvker\left(\frac{\xi}{\delta^{1/(2\sigma)}}\right) \Big|\lambda_\rho \left(\frac{\xi}{\delta^{1/(2\sigma)}}\right)\Big|^{2k},
 %       \qquad \xi\in B_1^c,
 %       \end{equation}
 %old
%\begin{equation}\label{eq:der_small_xi}
%    	|\p^\alp \redconvker(\xi)|
%	\leq C_{\alp} \delta, \qquad B_1\setminus \{0\}, 
%		\end{equation}
       % \comment{what I think we finally need:}
 with constants \(C_\alp,c>0\).
    \end{lemma}
    \begin{proof}[Proof of Lemma~\ref{lem:estimates_exp}]
    We start by deriving some auxiliary estimates for \(\alp\in \N_0^n\setminus\{0\}\). By Faà di Bruno's formula in the version of \eqref{eq:FdB}, one has
		\begin{equation}\label{eq:theta_ders}
		|\p^\alp \redconvker(\xi)|
		\le C_\alp \sum_{k=1}^{|\alp|}  \redconvker(\xi)
		\sum _{\{\alp'_1,\dots,\alp'_k\}} \prod_{l=1}^k |\delta\p^{\alp'_l} (\abs{\lambda_{\rho} (\xi)}^2)|, \qquad \xi\in \R^n,
		\end{equation}
        with the sum taken over the different partitions $\{{\alpha}'_1,\ldots , {\alpha}'_k\}$ of \(\alp\) for \(k\in \N\) such that \(k\le |\alp|\).
       % \commentCK{Check if $\xi$ are ok in this proof.}
        From~\eqref{eq:fourth_estimate} of Lemma~\ref{lem:estimates}, we infer that 
        \begin{align*}
		|\p^{\beta} \abs{\lambda_{\rho} (\xi)}^2|
		  \le  C_{\beta}|\xi|^{2-|{\beta}|}\widehat{Q_\rho}^{2}(\xi), \qquad \xi\in B_1^c
        \end{align*}
        for all \(\beta\in \N_0^n\),
        and it follows with \eqref{eq:theta_ders} that 
        \begin{align*}
            |\p^\alp \redconvker(\xi)|
		      &\le C_\alp \sum_{k=1}^{|\alp|} \delta^k\redconvker(\xi)
		\sum _{\{\alp'_1,\dots,\alp'_k\}} \prod_{l=1}^k |\xi|^{2-|\alp'_l|}\widehat{Q_\rho}^{2}(\xi)\\
            &\le C_\alp |\xi|^{-|\alp|}
		\sum_{k=1}^{|\alp|}\delta^k \redconvker(\xi)|\xi|^{2k}\widehat{Q_\rho}^{2k}(\xi),
        \qquad \xi\in B_1^c,
        \end{align*}
        where we used \eqref{eq:FdB_exp}  in the last inequality. This proves 
        \begin{equation}\label{eq:lem_gen_der_est}
	 |\p^\alp \redconvker(\xi)|
		\le C_\alp |\xi|^{-|\alp|}
		\sum_{k=1}^{|\alp|}\delta^k \redconvker (\xi)|\lambda_\rho(\xi)|^{2k},
        \qquad \xi\in B_1^c,
\end{equation}
which serves as the basis for verifying the estimates which are stated on the complement of a ball around the origin. 

Furthermore, since \(\xi\mapsto e^{-\xi}\)
%$\redconvker$ \commentCK{FIXME! isn't it rather the exponential function with $-$? Phrase properly... } 
decays faster than any polynomial 
%(because of \eqref{eq:lowerboundlambda}/\eqref{eq:lest_qrho_pre} \commentCK{Check again this argument, also above!}\comment{I actually think, this is not needed}), 
we have that, for any $l \in \N_0$, 
\begin{align}
        &\redconvker(\delta^{-1/(2\sigma)}\xi)
			\le C_l %\left(\Big|\lambda_\rho \left(\frac{\xi}{\delta^{1/(2\sigma)}}\right)\Big|^{2}\delta\right)^{-\gamma/2-l}
            \frac{1}{(|\lambda_\rho (\delta^{-1/(2\sigma)}\xi)|^{2}\delta)^{\gamma/2+l}},\label{eq:red1estimate}
            \quad \\
            &\redconvker (\delta^{-1/2}\xi)
		      \le C_l \frac{e^{-|\xi|^2\widehat{Q_\rho}^2(\delta^{-1/2}\xi)/2}}{(\abs{\lambda_{\rho} (\delta^{-1/2}\xi)}^2\delta)^l},
            \qquad \xi\in \R^n.\label{eq:red1estimate2}
\end{align}
		%		and 
		%		\begin{equation*}
		%		g\left(\frac{\xi}{\delta^{1/2}}\right)
		%		\le C_\alp \frac{1}{(|\xi|^2\widehat{Q_\rho}^2(\xi/\delta^{1/2}))^k},
		%		\end{equation*}

%\comment{and also:
%\begin{equation}\label{eq:red1estimate_new}
%    \redconvker\left(\frac{\xi}{\delta}\right)
%			\le C_l 
%            \frac{1}{\left(\Big|\lambda_\rho \left(\frac{\xi}{\delta}\right)\Big|^{2}\delta\right)^{l}},
%            \qquad \xi\in \R^n,
%\end{equation}
%which behaves like \(\delta^l\).
%}
We now turn to the verification of the stated inequalities \eqref{eq:first_der:est}--\eqref{eq:large_delta2}. 
Estimate \eqref{eq:first_der:est} is a direct consequence of inserting the first inequality in \eqref{eq:red1estimate} into \eqref{eq:lem_gen_der_est}.
%by Lemma~\ref{lem:estimates_exp}, for any $\alp \in \N^{n}_0$,~\eqref{eq:first_der:est}
%\commentCK{Check definition of multiindices.}

To show~\eqref{eq:new_est}  for \(\alp=0\) it suffices to combine \eqref{eq:red1estimate} for \(l=0\) with condition (H3) and the resulting inequality \eqref{eq:lest_qrho_pre}, which yields exactly
\begin{align*}
			\redconvker\left(\delta^{-1/(2\sigma)}\xi\right)
			%\le C_0 \frac{1}{\left(\Big|\lambda_\rho \left(\frac{\xi}{\delta^{1/(2\sigma)}}\right)\Big|^{2}\delta\right)^{\gamma/2}}
			\le C_0 |\xi|^{-\gamma\sigma},
            \qquad \xi\in B_{\delta^{1/(2\sigma)}}^c.
		\end{align*}
 For \(\alp\in \N_0^n\setminus\{0\}\), we argue that applying the chain rule
      % The remaining inequality \eqref{eq:lem_gen_der_est_cov} follows now from 
     to \eqref{eq:lem_gen_der_est} gives
     %by applying the chain rule,    i.e.,
        \begin{align}\label{aux865}
        \begin{split}
		\Big|\p^\alp \redconvker\left(\delta^{-1/(2\sigma)}\xi\right)\Big|
		&= \Big|(\p^\alp \redconvker)\left(\delta^{-1/(2\sigma)}\xi\right)\Big|
		\delta^{-|\alp|/(2\sigma)}\\
		&\le C_\alp |\xi|^{-|\alp|}
		\sum_{k=1}^{|\alp|}\delta^k \redconvker\left(\delta^{-1/(2\sigma)}\xi\right) \Big|\lambda_\rho \left(\delta^{-1/(2\sigma)}\xi\right)\Big|^{2k},
        \qquad \xi\in  B_{\delta^{1/(2\sigma)}}^c.
        \end{split}
		\end{align}
        %Furthermore, using condition (H3) and the resulting inequality \eqref{eq:lest_qrho_pre}, 
        Furthermore, using \eqref{eq:red1estimate} and then \eqref{eq:lest_qrho_pre}, there also holds, for any $k \in \N$,
		\begin{align*}
			\redconvker\left(\delta^{-1/(2\sigma)}\xi\right)
			%\le C_\alp \frac{1}{\left(\Big|\lambda_\rho \left(\frac{\xi}{\delta^{1/(2\sigma)}}\right)\Big|^{2}\delta\right)^{\gamma/2+k}}
			\le C_\alp \frac{1}{\Big|\lambda_\rho \left(\delta^{-1/(2\sigma)}\xi\right)\Big|^{2k}\delta^{k}}
			|\xi|^{-\gamma\sigma},
            \qquad \xi\in  B_{\delta^{1/(2\sigma)}}^c.
		\end{align*}
 Combining this with~\eqref{aux865}, 
 %the second estimate proved in Lemma~\ref{lem:estimates_exp}, which is
	%	\begin{align*}
	%	\Big|\p^\alp \redconvker\left(\frac{\xi}{\delta^{1/(2\sigma)}}\right)\Big|
	%	\le C_\alp |\xi|^{-|\alp|}
	%	\sum_{k=1}^{|\alp|} \delta^k\redconvker\left(\frac{\xi}{\delta^{1/(2\sigma)}}\right) \Big|\lambda_\rho \left(\frac{\xi}{\delta^{1/(2\sigma)}}\right)\Big|^{2k},
   %     \qquad \xi\in B_1^c,
	%	\end{align*}
        we get the refined estimate~\eqref{eq:new_est}.  
        
               %\comment{New:}
        %Recalling \eqref{eq:der_abs},
         %\begin{align*}
		%|\p^{\tilde{\beta}} |\xi||
		%\le C_{{\tilde{\beta}_1}} |\xi|^{1-|\tilde{\beta}_1|}
        %\le C_{\tilde{\beta}_1}|\xi|,
        %\qquad \xi\in \R^n,
		%\end{align*}
        %\comment{New:}
        Next we prove \eqref{eq:der_small_xi_new1} with \(\delta\in (0,1)\). We use that
        \begin{equation}\label{eq:xisqaured_new}
            |\p^{\tilde{\beta}} |\xi|^2|\le 2|\xi|^{\min\{2-|\tilde{\beta}|,0\}}\le 2,
		\qquad \xi\in B_1%\setminus\{0\}
        \end{equation} 
        for any \(\tilde{\beta}\in \N_0^n\). Then, \eqref{eq:der_small_xi_new1} with \(\alp=0\) follows immediately from \(\redconvker(\xi)\le 1\) for all \(\xi\in \R^n\).
        For \(\alp\in \N_0^n\setminus\{0\}\), we first conclude by inserting \eqref{eq:xisqaured_new} into \eqref{eq:sec_est_initial}, and with the smoothness of \(\widehat{Q_\rho}\), that 
        \begin{equation}\label{eq:ingredient_cr}
            |\p^\beta \abs{\lambda_{\rho} (\xi)}^2|
		\le C_\beta \sum_{\tilde{\beta}\le \beta}|\p^{\tilde{\beta}} |\xi|^2||\p^{\beta-\tilde{\beta}}\widehat{Q_\rho}^2(\xi)|
        \le C_\beta,
        \qquad \xi\in B_1%\setminus\{0\}
        \end{equation}
        for any \(\beta\in \N_0^n\).
        Together with \eqref{eq:theta_ders}, one has
        \begin{equation*}
            |\p^\alp \redconvker(\xi)|
		\le C_\alp \sum_{k=1}^{|\alp|} \delta^k 
        \le C_\alp , \qquad \xi\in B_1, \, \delta\in (0,1),
        \end{equation*}
        validating \eqref{eq:der_small_xi_new1}.
        
        In view of \eqref{eq:large_delta1}, we note that it holds
         \begin{equation*}
            |\p^{\tilde{\beta}} |\xi|^2|\le 2(1+|\xi|)^{2},
		\qquad \xi\in \R^n 
        \end{equation*} 
        for every \(\tilde{\beta}\in \N_0^n\).
        Hence, one deduces with the Leibniz rule from the smoothness of \(\widehat{Q_\rho}\) that
        \begin{align*}
		|\p^{\beta} \abs{\lambda_{\rho} (\delta^{-1/2}\xi)}^2|
		  \le  C_{\beta}\delta^{-1}|(1+\xi|)^{2}, \qquad\xi\in B_{\delta^{1/2}},\, \delta\in (1,\infty)
        \end{align*}
        for all \(\beta\in \N_0^n\). 
        It follows similar to \eqref{eq:theta_ders} that        
        \begin{equation*}
            \Big|\p^\alp \redconvker(\delta^{-1/2}\xi)\Big|
		\le C_{\alp}   \redconvker(\delta^{-1/2}\xi)(1+|\xi|)^{2|\alp|},
        \qquad \delta\in (1,\infty),\, \xi\in B_{\delta^{1/2}}.
        \end{equation*}
        Inserting the estimate \(\widehat{Q_\rho}\ge c>0\) on \(B_1\) and the resulting inequality \(\redconvker(\delta^{-1/2}\xi)\le e^{-c|\xi|^2}\) for \(\xi\in B_{\delta^{1/2}}\), finishes the proof of \eqref{eq:large_delta1}. 

        Finally, turning to \eqref{eq:large_delta2}, \eqref{eq:lem_gen_der_est} implies that 
        \begin{align*}
            |\p^\alp \redconvker(\delta^{-1/2}\xi)|
            &= \Big|(\p^\alp \redconvker)(\delta^{-1/2}\xi)\Big|
		\delta^{-|\alp|/2}\\
    		&\le C_\alp |\xi|^{-|\alp|}
    		\sum_{k=1}^{|\alp|}\delta^k \redconvker (\delta^{-1/2}\xi)|\lambda_\rho(\delta^{-1/2}\xi)|^{2k}\\
            &\le C_\alp 
    		 e^{-|\xi|^2\widehat{Q_\rho}^2(\delta^{-1/2}\xi)/2},
            \qquad \xi\in B_{\delta^{1/2}}^c,
        \end{align*}
        where we used for the last inequality that \(|\xi|^{-|\alp|}\le 1\) for \(\xi \in B_{\delta^{1/2}}^c\), together with 
        %\(t^ke^{-t}\le C_k e^{-t/2}\) for all \(t\ge 0\). 
        \eqref{eq:red1estimate2}. Moreover, it holds with \eqref{eq:lest_qrho_pre} that
        \begin{equation*}
            |\xi|^2\widehat{Q_\rho}^2(\delta^{-1/2}\xi)
            \ge \delta^{1-\sigma} |\xi|^{2\sigma}
            \ge c|\xi|^{2\sigma},
        \end{equation*}
        and so we obtain 
        \begin{equation*}
            |\p^\alp \redconvker(\delta^{-1/2}\xi)|
            \le C_\alp e^{-c|\xi|^{2\sigma}},
            \qquad \xi\in B_{\delta^{1/2}}^c
        \end{equation*}
        as claimed in \eqref{eq:large_delta2}.
        This finishes the proof.

    \end{proof}

    \section{Proof of Lemma~\ref{lem:lap_est}}\label{sec:appendix_aux}
    
    In this section, we prove Lemma~\ref{lem:lap_est}, stating a second-order version of the nonlocal Leibniz rule in \cite[Lemma~3.8]{BMCS25}.

    \begin{proof}[Proof of Lemma~\ref{lem:lap_est}]
    Throughout the proof,  \(C\) will denote a constant possibly changing from line to line, and possibly depending on the kernel \(\rho\).

    %The fact that \(g\in H^{\rho,p}(\R^n)\) is a direct consequence of Lemma~\ref{lem:lap_to_grad}.
    We first turn to the proof of estimate \eqref{eq:est_lap}. It was established in \cite[Lemma~3.8]{BMCS25} that 
    \begin{equation}\label{Leibniz33}
        D^\rho (fg)
        = f D^\rho g
            + K_f g    \end{equation}
    for 
    \begin{equation*}
			K_fg(x)
			\coloneqq \int_{\R^n} \frac{f(x)-f(y)}{|x-y|}\frac{x-y}{|x-y|}\rho(x-y)g(y)\dl y , \qquad x\in \R^n .
	\end{equation*} With the change of variables \(y\mapsto x+y\) and the radial symmetry of \(\rho\) according to (H0), \(K_f g\) can be equivalently expressed as
    \begin{equation}\label{eq:kfg_cov}
			K_fg(x)
			= \int_{\R^n} \frac{f(x+y)-f(x)}{|y|}\frac{y}{|y|}\rho(y)g(x+y)\dl y, \qquad x\in \R^n.
		\end{equation}
    Moreover, as shown in the proof of \cite[Lemma~3.8]{BMCS25}, for the nonlocal divergence operator, one finds similarly
    \begin{equation}\label{Leibniz_div}
         \diver^\rho (fG)
        = f\diver^\rho G
            + \tilde{K}_f G,
            \qquad  G\in C_c^\infty(\R^n;\R^n)
    \end{equation}
    for 
    \begin{equation*}
			\tilde{K}_fG(x)
			\coloneqq \int_{\R^n} \frac{f(x)-f(y)}{|x-y|}\frac{x-y}{|x-y|}\rho(x-y)\cdot G(y)\dl y, \qquad x\in \R^n.
	\end{equation*}
    Actually, the compact support of $G$ is not relevant for the validity of the equality, and so the result extends directly to \(G\in C^\infty(\R^n;\R^n)\). 

    With~\eqref{Leibniz33} and~\eqref{Leibniz_div}, we get the equality 
    \begin{equation}\label{eq:lap_decompose}
        \Delta^\rho_0 (fg)-f\Delta^\rho g
        = \tilde{K}_f D^\rho g 
        + \diver^\rho K_f g,
    \end{equation}
    and it remains to estimate the expressions on the right-hand side. Using that \(\rho\in L^1(\R^n)\), the term \(\tilde{K}_f D^\rho g \) can be estimated mimicking the steps in the proof of \cite[Eq. (3.4), Lemma~3.8]{BMCS25}, yielding 
    \begin{equation}\label{eq:lap_first:est}
        \norm{\tilde{K}_f D^\rho g }{L^p(\R^n)}
        \le C \norm{\nabla f}{L^{\infty}(\R^n;\R^n)}\norm{D^\rho g}{L^{p}(\R^n)}.
    \end{equation}

    Next, we use that with \eqref{eq:kfg_cov} we have
    \begin{equation*}
    (\diver^\rho  K_f g) (x)
        = \int_{\R^n} D^\rho\left(\frac{f(x+y)-f(x)}{|y|}g(x+y)\right)\cdot\frac{y}{|y|}\rho(y)\dl y.
    \end{equation*}
    Another application of \cite[Lemma~3.8]{BMCS25} gives for $x, y\in \R^n$, 
    \begin{align*}
        D^\rho\left(\frac{f(x+y)-f(x)}{|y|}g(x+y)\right)
        = \frac{f(x+y)-f(x)}{|y|}D^\rho g(x +y)
        + K_{\frac{f(x+y)-f(x)}{|y|}} g(x +y), 
    \end{align*}
    which can be estimated with the help of
    \begin{equation*}
        \left| \frac{f(x+y)-f(x)}{|y|}D^\rho g(x +y)\right|
        \le \norm{\nabla f}{L^\infty(\R^n;\R^n)}|D^\rho g(x +y)|
    \end{equation*}
     and 
    \begin{equation*}
        |K_{\frac{f(x+y)-f(x)}{|y|}} g(x +y)|
        \le C\norm{\nabla f}{W^{1,\infty}(\R^n;\R^n)}(|g|*\rho)(x +y);
    \end{equation*}
    the latter follows as in the proof of \cite[Eq. (3.4), Lemma~3.8]{BMCS25}.
    Altogether, taking into account again the spherical symmetry of \(\rho\), it follows
    \begin{equation*}
        \left|\diver^\rho  K_f g\right|
        \le C\norm{\nabla f}{W^{1,\infty}(\R^n;\R^n)} (|D^\rho g|*\rho + (|g|*\rho)*\rho).
    \end{equation*}
    Finally, recalling that \(\rho\in L^1(\R^n)\), a repeated application of Young's convolution inequality leads to the estimate 
    \begin{equation}\label{eq:lap_third:est}
        \norm{\diver^\rho  K_f g}{L^{p}(\R^n)}
        \le C \norm{\nabla f}{W^{1,\infty}(\R^n;\R^n)}\norm{g}{H^{\rho,p}(\R^n)}.
    \end{equation}
    The sought inequality \eqref{eq:est_lap} follows now from \eqref{eq:lap_decompose} together with the estimates \eqref{eq:lap_first:est} and \eqref{eq:lap_third:est}.
    	%It remains to prove the convergence in \eqref{eq:simul_approx_lem}, where $\{\sigma_j \}_{j>0}$ is a smooth sequence of cut-off functions obtained by setting \(\sigma_j(\cdot)\coloneqq \sigma(\cdot/j)\) for some initial \(\sigma \in C_c^\infty(\R^n)\) with \(\sigma\equiv 1\) in \(B_1\). We note that \(\sigma_j\Delta^{\rho} g\to \Delta^{\rho} g\) by dominated convergence, and that \(\norm{\nabla \sigma_j}{W^{1,\infty}}\le C/j\). The convergence follows then immediately from the above-proved estimate~\eqref{eq:est_lap} for the choice \(f= \sigma_j\), since
%\begin{equation*}
%		\norm{\Delta^{\rho}_0 (\sigma_j g)-\sigma_j\Delta^{\rho} g}{L^p(\R^n)}
%		\le C\norm{\nabla \sigma_j}{W^{1,\infty}(\R^n;\R^n)}\norm{g}{H^{\rho,p}(\R^n)}.
%\end{equation*} 
%    The proof is concluded.
\end{proof}

%Finally, we state and prove here a variant of \cite[Proposition~4.4]{BMCS25}, generalizing it from smooth to non-smooth test functions of the type \(|\lambda_\rho|\Fcal^{-1}\varphi\), \(\varphi\in C_c^\infty(\R^n)\). To that end, we introduce the operator
%\begin{equation*}
%    \mathcal{R}: L^p(\R^n;\R^n)\to L^p(\R^n), 
%    \quad U\mapsto \sum_{j=1}^n  \mathcal R_j U_j,
%\end{equation*}
%where the \(\mathcal{R}_j\) are the \(j\)-th Riesz transforms defined through multiplication with the Fourier symbol \(\xi\mapsto -i\xi_j/|\xi|\). %Note that the Riesz transforms are bounded operators on \(L^p(\R^n)\) for \(p\in(1,\infty)\) (see \cite[Corollary~5.2.8]{Gr14}). 

%\begin{lemma}\label{lem:aux_ftoc}  Assume that \(\rho\) has compact support and satisfies assumptions (H0)--(H2) and (H3').
%	For all $u, \varphi \in C^{\infty}_c (\R^n)$ it holds
    %\[
    %\int_{\R^n} \hat{u} |\lambda_\rho| \hat \varphi \dl \xi 
    %		= \int_{\R^n} \frac{-i\xi }{2\pi|\xi|^2\hat Q_\rho} \cdot \widehat{D^\rho_0 u}|\lambda_\rho|\hat\varphi\dl \xi .
    %\]
%\begin{equation*}
%    \langle \hat{u},|\lambda_\rho|\Fcal^{-1}\varphi\rangle
%	= \langle \Fcal (\mathcal{R} D^\rho_0 u),\Fcal^{-1}\varphi\rangle.
%\end{equation*}
%
%	\end{lemma}

    %%%%%%%%%%%%%%%%%%%%%%%%%%%%%%%%%%%%%%%%%%%%%%%%%%%%%%%%%%%%%%%%%
    \addtocontents{toc}{\protect\setcounter{tocdepth}{-1}}
    \section*{Declarations}
    \textit{Declaration of Generative AI and AI-Assisted Technologies:}
	During the preparation of this manuscript, the authors used ChatGPT, Claude, Mistral and Copilot to assist with literature search and spell-checking, and to improve the fluency of selected portions of the text. 
    %\comment{...} was used to generate or refine LaTeX code for improving the document layout, typography, and overall presentation. 
	No generative AI tools were used in the development of the mathematical ideas, the derivation of the theoretical results, or the proofs. The authors are fully responsible for the content of the manuscript. % and for the accuracy, validity, and integrity of all results presented. \commentCK{Shorten second part about the results? }
     \addtocontents{toc}{\protect\setcounter{tocdepth}{1}}

     %%%%%%%%%%%%%%%%%%%%%%%%%%%%%%%%%%%%%%%%%%%%%%%%%%%%%%%%%%%%%%%%%
    \addtocontents{toc}{\protect\setcounter{tocdepth}{-1}}
    \section*{Acknowledgments}
    The authors are grateful to Jos\'e Carlos Bellido and Guillermo Gac\'ia-S\'aez for helpful discussions on the topic. This research was part of the Joint Spanish-German Research Project ``Nonlocal Gradients in Variational Analysis and Materials Modelling: Limits, Kernels, Boundaries" funded by the Agencia Estatal de Investigaci\'on (AEI) and the Deutsche Forschungsgemeinschaft (DFG). Specifically,  \textsc{L.H.} is employed under Project PCI2024-155023-2 funded by MICIU/AEI/10.13039/501100011033/UE and the European Union. \textsc{C.K.} acknowledges financial support by DFG via Project 541520348. %\commentCK{FIXME!} \color{black}
    The work of \textsc{C. M.-C.} was supported through Project PID2024-158664NB-C2 funded by MICIU/AEI/10.13039/501100011033/FEDER/UE.
    \addtocontents{toc}{\protect\setcounter{tocdepth}{1}}

%project number 541520348.
%project number PCI2024-155023-2 
  %  \textsc{L.H.} acknowledges support from Project PCI2024-155023-2 funded by MICIU/AEI/10.13039/501100011033/UE and the European Union. 
    %Grant PCI2024-155023-2 funded by the European Union and MCIN/AEI/10.13039/501100011033. 

 %   \color{red}
%All three authors acknowledge support by the AEI-DFG project ``Nonlocal Gradients in Variational Analysis and Materials Modelling: Limits, Kernels, Boundaries"through the  Agencia Estatal de Investigaci\'on (AEI) project number PCI2024-155023-2 and the Deutsche Forschungsgemeinschaft (DFG)  project number 541520348. \color{black}

	%%%%%%%%%%%%%%%%%%%%%%%%%%%%%%%%%%%%%%%%%%%%%%%%%%%%%%%%%%%%%%%%%
	%%%%%%%%%%%%%%%%%%%%%%%%%%%%%%%%%%%%%%%%%%%%%%%%%%%%%%%%%%%%%%%%%
	
	\bibliography{Bibliography.bib}
	\bibliographystyle{abbrvurl}
\end{document}